\documentclass[11pt,psamsfonts]{amsart}

\newtheorem{theorem}{Theorem}[subsection]

\usepackage[all]{xy}
\usepackage[german,english]{babel}
\usepackage{amssymb}
\usepackage{latexsym}
\usepackage{xspace}                         
\usepackage{stmaryrd}                       
\usepackage{calc}
\usepackage{ifthen}
\usepackage[margin=1.3in]{geometry}

\newcommand{\numiii}{\renewcommand{\labelenumi}{(\roman{enumi})}}
\numberwithin{equation}{section}

\theoremstyle{plain}

\newtheorem*{theorem*}{Theorem}
\newtheorem{proposition}[theorem]{Proposition}
\newtheorem{corollary}[theorem]{Corollary}
\newtheorem{lemma}[theorem]{Lemma}

\newtheorem{booklemma}[theorem]{Booklemma}
\newtheorem{bookproposition}[theorem]{Bookproposition}

\theoremstyle{definition}
\newtheorem{definition}[theorem]{Definition}
\newtheorem{example}[theorem]{Example}

\newtheorem{bookexample}[theorem]{Bookexample}
\newtheorem{bookremark}[theorem]{Bookremark}

\newtheorem{construction}[theorem]{Construction}
\newtheorem{notation}[theorem]{Notation}

\newtheorem{remark}[theorem]{Remark}

\newtheorem{situation}[theorem]{Situation}

\newtheorem{convention}[theorem]{Convention}

\theoremstyle{remark}

\newboolean{book} 
\setboolean{book}{false}

\newcommand{\new}{\newcommand}
\newcommand{\ren}{\renewcommand}

\new{\extraskip}{\bigskip} 
\new{\extranewpage}{\newpage}
\new{\extravspace}{\vspace}

\newcommand{\raa}{b} 
\newcommand{\ral}{a} 
\newcommand{\ram}{b}
\newcommand{\rar}{c}
\newcommand{\ras}{a} 
\newcommand{\rat}{b} 

\newcommand{\rak}{{k}}

\newcommand{\ran}{{n}}

\newcommand{\ramat}{r}

\newcommand{\annmath}{Ann. Math.}

\renewcommand{\AA}{\mathbb{A}}

\newcommand{\CC}{\mathbb{C}}

\newcommand{\FF}{\mathbb{F}}

\newcommand{\II}{\mathbb{I}}

\newcommand{\LL}{\mathbb{L}}
\newcommand{\MM}{\mathbb{M}}

\newcommand{\NN}{\mathbb{N}}

\newcommand{\RR}{\mathbb{R}}

\newcommand{\UU}{\mathbb{U}}

\newcommand{\VV}{\mathbb{V}}

\newcommand{\XX}{\mathbb{X}}

\newcommand{\YY}{\mathbb{Y}}

\newcommand{\ZZ}{\mathbb{Z}}

\newcommand  {\shB}     {\mathcal{B}}

\newcommand  {\shE}     {\mathcal{E}}
\newcommand  {\shF}     {\mathcal{F}}
\newcommand  {\shG}     {\mathcal{G}}
\newcommand  {\shH}     {\mathcal{H}}

\newcommand  {\shI}     {\mathcal{I}}

\newcommand  {\shK}     {\mathcal{K}}
\newcommand  {\shM}     {\mathcal{M}}

\newcommand  {\shL}     {\mathcal{L}}
\newcommand  {\shP}     {\mathcal{P}}

\newcommand  {\shS}     {\mathcal{S}}
\newcommand  {\shT}     {\mathcal{T}}

\newcommand  {\foa}     {\mathfrak{a}}

\newcommand  {\fob}     {\mathfrak{b}}

\newcommand  {\fom}     {\mathfrak{m}}

\newcommand  {\fon}     {\mathfrak{n}}

\newcommand  {\fop}     {\mathfrak{p}}

\newcommand  {\foq}     {\mathfrak{q}}

\newcommand  {\an}      {{\text{an}}}
\newcommand  {\Ann}     {\operatorname{Ann}}

\newcommand  {\ara}     {\operatorname{ara}}

\newcommand  {\Ass}     {\operatorname{Ass}}

\newcommand  {\Cl}      {\operatorname{Cl}}

\newcommand  {\dual}    {\vee}

\newcommand  {\ext}     {\operatorname{ext}}

\newcommand  {\Hom}     {\operatorname{Hom}}

\newcommand  {\id}      {\operatorname{id}}
\newcommand  {\im}      {\operatorname{im}}

\renewcommand  {\ker }  {\operatorname{kern}}

\newcommand  {\liminv}  {\varprojlim}

\newcommand  {\modu}     {\operatorname{mod}}

\newcommand  {\Mor}     {{\operatorname{Mor}}}

\newcommand  {\nor}     {{\operatorname{nor}}}

\renewcommand{\O}       {\mathcal{O}}

\newcommand  {\Proj}    {\operatorname{Proj}}

\newcommand  {\Ra}      {\Rightarrow}

\newcommand  {\rad}      {\operatorname{rad}}
\newcommand  {\rk}    {\operatorname{rk}}
\newcommand  {\red}     {{\operatorname{red}}}
\newcommand  {\reg}     {\operatorname{reg}}

\newcommand{\glpoint}{P}

\newcommand  {\sat}     {{\operatorname{sat}}}

\newcommand  {\Spec}    {\operatorname{Spec}}

\newcommand  {\support} {\operatorname{supp}}

\newcommand  {\Syz}     {\operatorname{Syz}}
\newcommand  {\Tor}     {\operatorname{Tor}}

\newcommand{\zar}{{\operatorname{Zar}}}

\newcommand{\comdots}{ , \ldots , }
\newcommand{\komdots}{ , \ldots , }
\newcommand{\plusdots}{ + \ldots + }
\newcommand{\uplusdots}{ \uplus \ldots \uplus }

\newcommand{\timesdots}{ \times \ldots \times }

\newcommand{\capdots}{ \cap \ldots \cap }

\newcommand{\ltodots}{ \lto \ldots \lto }

\newcommand{\todots}{ \to \ldots \to }

\newcommand{\op}[1]{\operatorname{#1}}
\newcommand{\cal}[1]{\mathcal{#1}}

\renewcommand{\bar}[1]{\overline{#1}}

\newcommand{\longto}{\longrightarrow}

\renewcommand{\phi}{\varphi}        
\renewcommand{\epsilon}{\varepsilon}

\newcommand{\height}{\op{ht}}
\newcommand{\depth}{\op{depth}}

\newcommand{\dirlim}{\varinjlim}

\newcommand{\tensor}{\otimes}         

\renewcommand{\to}[1][]{\xrightarrow{\ #1\ }}

\newcommand{\usc}[1][m]{\underline{\phantom{#1}}}

\usepackage{amscd}
\usepackage{amssymb}

\newboolean{extra} 

\setboolean{extra}{false}

\ren{\sec}[1]{{ {#1}'}}

\newcommand{\setj}{{J}}

\newcommand{\elea}{{a}} 
\newcommand{\eleb}{{b}}
\newcommand{\elec}{{c}}

\newcommand{\elef}{{f}}
\newcommand{\eleg}{{g}}

\newcommand{\elel}{{l}}
\newcommand{\elem}{{m}}
\newcommand{\elen}{{n}}

\newcommand{\eleq}{{q}}
\newcommand{\eler}{{r}}
\newcommand{\eles}{{s}}
\newcommand{\elet}{{t}}

\newcommand{\elev}{{v}}

\newcommand{\eley}{{y}}

\newcommand{\eletest}{{z}} 

\ren{\elef}{t} 

\ren{\eleg}{s}

\ren{\elel}{t} 

\ren{\elem}{s} 

\ren{\elen}{t}

\ren{\elet}{y}

\newcommand{\elerp}{{h}}

\newcommand{\elelp}{{t'\!\!}}

\newcommand{\elere}{\tilde{r}}

\newcommand{\elema}{{u}}

\newcommand{\elemp}{{\tilde{s}}}

\newcommand{\elemaa}{{v}}

\newcommand{\numk}{{k}}

\newcommand{\numm}{{m}}
\newcommand{\numn}{{n}}

\newcommand{\numsec}{{\ell}}

\newcommand{\expok}{{k}}
\newcommand{\expon}{{e}}
\newcommand{\expom}{{m}}
\newcommand{\expot}{\tau}
\newcommand{\expor}{{r}}

\newcommand{\ind}{{\kappa}}

\newcommand{\indsec}{{j}}

\newcommand{\indi}{i}

\newcommand{\indj}{j}

\newcommand{\indk}{k}

\newcommand{\inda}{{i}}

\newcommand{\indo}{{j}}

\newcommand{\indu}{{k}}

\newcommand{\indao}{{\inda \indo}}

\newcommand{\indou}{{\indo \indu}}

\newcommand{\indau}{{\inda \indu}}

\newcommand{\indaou}{{\inda \indo \indu}}

\newcommand{\runn}{{n}}

\newcommand{\setind}{{I}}
\newcommand{\indinset}{{\ind \in \setind}}
\newcommand{\indseize}{{i}}  
\newcommand{\indseizesec}{{j}}

\newcommand{\indobj}{{i}}  

\newcommand{\indob}{{\indobj}}

\newcommand{\catab}{{\cal A}}

\newcommand{\catsch}{Sch} 

\newcommand{\obset}{{A}}

\newcommand{\obab}{{A}}

\newcommand{\obabsec}{{ \tilde{A}}}

\newcommand{\colimab}{{A}}

\newcommand{\abso}{{\rm abs}}

\newcommand{\catt}{\shT}

\newcommand{\catc}{\cal C}

\newcommand{\cattt}{\shT}

\newcommand{\arphi}{{\varphi}} 
\newcommand{\arpsi}{{\psi}} 

\newcommand{\armr}{{\psi_1}}

\newcommand{\armrsec}{{\psi_2}}

\newcommand{\mor}{{\varphi}}

\newcommand{\morleft}{{\gamma}}

\newcommand{\morj}{{j}}

\newcommand{\posp}{{p}} 

\newcommand{\potq}{{q}}

\newcommand{\frobpotq}{{[\potq]}}

\newcommand{\Frob}{F} 

\newcommand{\frobf}{{F}} 

\newcommand{\frobe}{{{}^{\expoe}\!}}

\newcommand{\frobesec}{{{}^{\expoesec}\!}}

\newcommand{\frobeplussec}{{{}^{\expoe + \expoesec}\!}}

\newcommand{\expoe}{{e}}

\newcommand{\expoesec}{{\tilde{e}}}

\newcommand{\perf}{{\infty}} 

\newcommand{\shef}{\theta} 

\newcommand{\grocont}{G} 
\newcommand{\catopen}{\shT}

\newcommand{\catopenspax}{\catopen_\spax}

\newcommand{\modcat}{\shM}

\newcommand{\comp}{{\operatorname{comp}}}

\newcommand{\compic}{{\operatorname{comp}}}

\newcommand{\differ}{\delta} 

\newcommand{\catringtest}{\cal D} 

\newcommand{\subcatring}{\catringtest}

\newcommand{\ringtest}{{T}}

\newcommand{\ringtestincat}{\ringtest \in \catringtest}

\newcommand{\testring}{{\ringtest}}

\newcommand{\subcatringcl}{{\subcatring-{\rm cl}}}

\newcommand{\patha}{\gamma} 

\newcommand{\toponw}{} 

\newcommand{\sect}{\psi} 

\newcommand{\colim}{\dirlim}

\newcommand{\xexact}{$\spax$-exact}

\newcommand{\pre}{{\operatorname{pre}}} 

\newcommand{\heighd}{{d}}

\newcommand{\automor}{{\alpha}}

\newcommand{\diag}{\triangle}

\newcommand{\submor}{\varphi} 

\newcommand{\Ab}{\mathsf{Ab}}
\newcommand{\PSh}{\mathsf{PSh}}
\newcommand{\Sh}{\mathsf{Sh}}

\newcommand{\obu}{{U}} 

\newcommand{\oburight}{{U'}}

\newcommand{\obuleft}{{\tilde{U}}}

\newcommand{\obusec}{{U'}}

\newcommand{\obur}{{U'}} 

\newcommand{\obtop}{{U}} 

\newcommand{\obtopa}{{U}} 

\newcommand{\obtopo}{{V}}

\newcommand{\obtopsup}{{\spax - \obtop}} 

\newcommand{\obcov}{{U}} 
\newcommand{\obcovp}{U} 

\newcommand{\obcovindic}{{\obcov_\indcov}}

\newcommand{\obcovindicsec}{{\obcov_\indcovsec}}

\newcommand{\namecov}{{\cal U}}

\newcommand{\namecovfine}{{\cal V}} 

\newcommand{\morfine}{{\theta}}

\newcommand{\obcovsec}{{V}}

\newcommand{\morcov}{{\varphi}} 

\newcommand{\morcovsec}{{\psi}}

\newcommand{\obv}{{V}}

\newcommand{\obx}{{X}}
\newcommand{\obw}{{W}}

\newcommand{\obtest}{{T}}

\newcommand{\obl}{U}
\newcommand{\obm}{V}
\newcommand{\obr}{W}
\newcommand{\obsec}{W}

\newcommand{\obsecsec}{U}

\newcommand{\openzar}{{U}}

\newcommand{\closedzar}{{Y}}

\newcommand{\opzar}{{\openzar}} 

\newcommand{\obpart}{{\spax}} 

\newcommand{\indpart}{{i}} %

\newcommand{\clopart}{{W}} 

\newcommand{\numpart}{{k}} %

\newcommand{\morpart}{{\struto}} %

\newcommand{\morpartindic}{{\morpart_\indpart}}

\newcommand{\morpartindictospax}{{ \morpartindic: \obpartindic \to \spax}}

\newcommand{\obpartsec}{{\spav}}

\newcommand{\setindpart}{{I}} %

\newcommand{\indpartinset}{{\indpart \in \setindpart}} %

\newcommand{\obpartindic}{{\obpart_\indpart}}

\newcommand{\partitionofspax}{{\obpartindic \to \spax}} %

\newcommand{\indpartsec}{{j}}

\newcommand{\setindpartsec}{{J}}

\newcommand{\indpartsecinset}{{\indpartsec \in \setindpartsec}}

\newcommand{\algforc}{B}

\newcommand{\mata}{D}
\newcommand{\entry}{{r}}
\newcommand{\indcri}{{i}}
\newcommand{\indcrj}{{j}}

\newcommand{\dvd}{V} 

\newcommand{\modul}{{M}}

\newcommand{\modulsec}{{N}}

\newcommand{\submod}{{N}}
\newcommand{\submodul}{{\submod \subseteq \modul}}
\newcommand{\quotmod}{{\modul/\submod}}

\newcommand{\submodzweit}{U}
\newcommand{\modulzweit}{V}
\newcommand{\submodulzweit}{\submodzweit \subseteq \modulzweit}

\newcommand{\modl}{N} 
\newcommand{\modm}{M} 
\newcommand{\modr}{L} 

\newcommand{\finsubmod}{\submod '}
\newcommand{\betweenmod}{W}

\newcommand{\indmod}{{j}}

\newcommand{\setindmod}{{J}}

\newcommand{\indmodinset}{{\indmod \in \setindmod}}

\newcommand{\modfree}{{F}}

\newcommand{\modra}{\mu}
\newcommand{\submodra}{\nu}

\newcommand{\coa}{\alpha} 
\newcommand{\cob}{\beta}

\newcommand{\freisurn}{\nu}

\newcommand{\lmr}{\modl \stackrel{\coa}{\to} \modm \stackrel{\cob}{\to} \modr}
\newcommand{\lmrr}[1]{\modl {#1} \stackrel{\coa}{\to} \modm {#1} \stackrel{\cob}{\to} \modr {#1}}
\newcommand{\lmrtensor}[1]{\modl \tensor_\ring {#1}
 \stackrel{\coa}{\to} \modm \tensor_\ring {#1} \stackrel{\cob}{\to} \modr\tensor_\ring {#1}}

\newcommand{\chch}{\check{H}} 

\newcommand{\algcat}{B}
\newcommand{\algcata}{\bar{B}}
\newcommand{\algcatb}{B}

\newcommand{\stalkabs}{H}

\renewcommand{\freisurn}{\sigma}
\newcommand{\fgen}{f_1 \komdots f_n}

\newcommand{\tensordvd}{\tensor_{\dvd}}

\newcommand{\lto}{\longrightarrow}

\newcommand{\lmrglobal}[2]{\Gamma(#1,\modl {#2}) \stackrel{\coa}{\to}
 \Gamma(#1,\modm{#2}) \stackrel{\cob}{\to} \Gamma(#1,\modr{#2})}
\newcommand{\shekok}{\shK}
\newcommand{\shker}{\shT}

\newcommand{\elemzweit}{z}

\newcommand{\covadma}{A}

\newcommand{\fixata}{H}

\newcommand{\gammashf}[1]{\Gamma(#1,\shF)}
\newcommand{\timesx}{\times_{\spax}}
\newcommand{\spax}{X}

\newcommand{\covmor}{\varphi}

\newcommand{\sheafify}{{\#}}

\newcommand{\sheafplus}{{+}} 

\newcommand{\sheaf}{{\shF}}

\newcommand{\specring}{\Spec \ring}
\newcommand{\spaxeqspecring}{\spax = \Spec \ring }

\newcommand{\alg}{{A}}
\newcommand{\algb}{{B}}

\newcommand{\nzdring}{{\ring ^*}}

\newcommand{\rings}{S}

\newcommand{\iindex}{{ i \in I}}

\newcommand{\indcov}{{i}} 

\newcommand{\indcovsec}{{j}}

\newcommand{\indcovtri}{{k}}

\newcommand{\indcovdob}{{\indcov \indcovsec}}

\newcommand{\setindcovsec}{{J}}

\newcommand{\indcovsecinset}{{\indcovsec \in \setindcovsec}}

\newcommand{\indcovinsetsec}{{ \indcovsec \in \setindcovsec}}

\newcommand{\setindcov}{{I}}

\newcommand{\indcovinset}{{ \indcov \in \setindcov}}  

\newcommand{\numcov}{{n}}

\newcommand{\indcovfine}{{j}}

\newcommand{\obcovfine}{{V}}

\newcommand{\setindcovfine}{{J}}

\newcommand{\indcovfineinset}{{ \indcovfine \in \setindcovfine}}  

\newcommand{\covering}{{ \obcov_\indcov \to \obcovp}}

\newcommand{\coveringring}{{ \ring \to \ring_\indcov }}

\newcommand{\shel}{\shG}
\newcommand{\shem}{\shF}
\newcommand{\sher}{\shE}

\newcommand{\sheafrightcomplexx}[1]{\shel{#1} \to \shem{#1} \to \sher{#1} \to 0}
\newcommand{\olmrosheaf}{0 \to \shel \to \shem \to \sher \to 0}
\newcommand{\olmrosheafstack}{0 \to \shel \stackrel{\coa}{\to} \shem \stackrel{\cob}{\to}
\sher \to 0}

\newcommand{\rest}{\rho} 

\newcommand{\tensorr}{{\tensor_\ring}}

\newcommand{\lmroo}[1]{ \modl{#1} \to \modm {#1} \to \modr {#1} \to 0}
\newcommand{\olmro}{0 \to \modl \to \modm \to \modr \to 0}
\newcommand{\olmroo}[1]{0 \to \modl {#1} \to \modm {#1} \to \modr {#1} \to 0}

\newcommand{\lmrsheaf}{\shel \to \shem \to \sher}
\newcommand{\lmrsheafstack}{\shel \stackrel{\coa}{\to} \shem \stackrel{\cob}{\to} \sher}

\newcommand{\lmrsheafstackk}[1]{\shel {#1} \stackrel{\coa}{\to} \shem {#1} \stackrel{\cob}{\to} \sher {#1}}

\newcommand{\lmrositemorsheafstack}{ \sitemorsheaf (\shel) \stackrel{\coa}{\to} \sitemorsheaf (\shem )
\stackrel{\cob}{\to} \sitemorsheaf (\sher) \to 0 }

\newcommand{\shinv}{\shL}

\newcommand{\popi}{n}
\newcommand{\pipo}{\pi^\popi}

\newcommand{\gammaf}[1]{\Gamma(#1 , \shF)}

\newcommand{\contr}{{\rm cont}} 

\newcommand{\contin}{{\rm cont}} 

\newcommand{\kersyz}{\rm Syz}
\newcommand{\syzgroup}{{\rm SYZ}}

\newcommand{\prop}{{\rm prop}}

\newcommand{\dwd}{W}

\newcommand{\fuglob}{{f}} 

\newcommand{\seminorm}{{\rm semi}}
\newcommand{\tensortwist}[1]{{#1 \tensor 1 - 1 \tensor #1}}

\newcommand{\tensortwistfuglob}{{\tensortwist{\fuglob}}}

\newcommand{\fin}{{\rm fin}}
\renewcommand{\finsubmod}{L}

\newcommand{\tensorre}{\tensorr \frobe R}

\newcommand{\morsite}{{\sitemor}}

\newcommand{\sitemor}{ \varphi}

\newcommand{\sitemortop}{ \sitemor^{*}} 
\newcommand{\sitemorinv}{ \sitemor^{-1}} 
\newcommand{\sitemorsheaf}{ \sitemor^{*}} 
\newcommand{\sitemorpre}{ \sitemor^{\pre}} 
\newcommand{\sitemormodule}{ \sitemor^{*}} 
\newcommand{\sitemorc}{ {c}} 

\newcommand{\mortop}{ \varphi}

\newcommand{\mortoptop}{ \mortop^{*}} 
\newcommand{\mortopmodule}{ \mortop^{*}} 
\newcommand{\mortopmod}{ \mortopmodule} 

\newcommand{\struto}{{j}} 

\newcommand{\strutoinv}{ \struto^{-1}} 
\newcommand{\strutopre}{ \struto^{\pre}} 
\newcommand{\strutosheaf}{ \struto^{*}} 

\newcommand{\strutosec}{{\rho}} 

\newcommand{\presh}   {\shF} 

\newcommand{\preshb}{\shG}

\newcommand{\phili}{ \rho}

\newcommand{\philileft}{\tilde{ \rho}}

\newcommand{\xili}{ \sigma}

\newcommand{\topc}{ {\rm tc }}

\newcommand{\presheafmor}{ {f }}

\newcommand{\torclass}{{\cal T}}
\newcommand{\freeclass}{{\cal F}}

\newcommand{\modtor}{{T}}

\newcommand{\modtorsec}{{S}}

\newcommand{\modfreesec}{{G}}

\newcommand{\sitemorpf}{{\sitemor_*}}
\newcommand{\frobmo}{{\Phi}}
\newcommand{\mofi}{{\psi}}

\newcommand{\shehom}{{\mu }}

\newcommand{\shell}{\shH}
\newcommand{\triv}{{\rm triv}}
\newcommand{\inc}{\alpha}

\newcommand{\covmap}{\varphi}
\newcommand{\cmpbi}[1]{ \covmap_i^*(#1)}
\newcommand{\covmapcov}{\covmap_i: X_i \to X}
\newcommand{\iini}{i \in I}
\newcommand{\cmpbimodul}{{ \cmpbi {\modul }}}
\newcommand{\gammaxmodul}{\Gamma(\spax,\modul)}
\newcommand{\secmap}{\psi}
\newcommand{\secpb}{\secmap^*}

\newcommand  {\topo}     {{\operatorname{top}}} 

\newcommand  {\frob}     {{\operatorname{frob}}}

\newcommand{\surj}{{\operatorname{surj}}}

\newcommand{\tors}{{\rm tors}}

\newcommand{\nzd}{{\rm nzd}}

\newcommand{\zftop}{{\rm zftop}}

\newcommand{\absol}{{\rm abs}}

\newcommand{\idealsup}{{J}}

\newcommand{\idealsupsec}{{I}}

\newcommand{\ctopc}{{\rm co-tc}}

\newcommand{\dense}{{\rm dense}}

\newcommand{\divtop}{{\rm div}}

\newcommand{\jac}{{\rm Jac}}

\newcommand{\arb}{{\rm arb}}

\newcommand{\topdim}{{\rm top}} 

\newcommand{\clop}{c} 

\newcommand{\cont}{{\rm cont}}

\newcommand{\topocl}     {{\topo-{\rm cl}}}

\newcommand{\tocl}{{\topocl}}

\newcommand{\soc}{{\rm soc}}

\newcommand{\hommod}{{\mu}} 

\newcommand{\homproj}{{p}} 

\newcommand{\modprim}{{P}}

\newcommand{\varrad}{{\widetilde{\rad}}}

\newcommand{\cons}{{\rm cons}}

\newcommand{\lmrfrobe}{ \frobe \modl \stackrel{\coa}{\to}   \frobe \modm \stackrel{\cob}{\to } \frobe \modr }

\newcommand{\algfilt}{{H}}

\newcommand{\inem}{\emph}

\newcommand{\eptau}{\tau}

\newcommand{\synchronize}[2]{\renewcommand{#1}{{#2}}}

\newcommand{\closed}{{A}}

\newcommand{\prim}{{p}} 

\newcommand{\ulm}{{u}} 

\newcommand{\num}{{m}}

\newcommand{\numcoho}{{k}} 
\newcommand{\mapind}{{\psi}} 
\newcommand{\chcyc}{{s}} 

\newcommand{\choiceindcovnumcoho} {{\indcov_0 \ldots
\indcov_\numcoho}} 

\newcommand{\choicecommaindcovnumcoho} {{\indcov_0 \comdots
\indcov_\numcoho}}

\newcommand{\choiceindcovzero} {{\indcov_0 }}

\newcommand{\choicecommaindcovzero} {{\indcov_0}}

\newcommand{\choiceindcovzeroone} {{\indcov_0 \indcov_1 }}

\newcommand{\choicecommaindcovzeroone} {{\indcov_0 ,\indcov_1
}}

\newcommand{\choiceindcovzerotwo} {{\indcov_0 \indcov_1 \indcov_2
}}

\newcommand{\choicecommaindcovzerotwo} {{\indcov_0 ,\indcov_1,
\indcov_2}}

\newcommand{\choiceindcovzerothree} {{\indcov_0 \indcov_1 \indcov_2
\indcov_3}}

\newcommand{\choicecommaindcovzerothree} {{\indcov_0 ,\indcov_1,
\indcov_2,\indcov_3}}

\newcommand{\choicebreak}{{t}}

\newcommand{\catind}{{\Lambda}} 

\newcommand{\indexcat}{{\catind}}

\newcommand{\catideal}{{\triangle}} 

\newcommand{\indec}{{\lambda}}
\newcommand{\inded}{{\lambda'}}

\newcommand{\indee}{{\lambda_0}}

\newcommand{\indfin}{{\lambda_0}}

\newcommand{\indm}{{\lambda}} 
\newcommand{\indmo}{{\lambda}}
\newcommand{\indmu }{{\lambda'}}
\newcommand{\indr}{{\mu}}

\newcommand{\inds}{{s}}

\newcommand{\indl}{{\kappa}}

\newcommand{\indcomp}{{i}}

\new{\idbla}{{H}}
\new{\idblasec}{{G}}
\new{\idblatri}{{E}}

\newcommand{\compindic}{{ \spax_\indcomp}}

\newcommand{\setindcomp}{{I}}

\newcommand{\indcompinset}{{\indcomp \in \setindcomp}}

\newcommand{\indeeincat}{{\indee \in \catind}}
\newcommand{\indmincat}{{\indm \in \catind}}

\newcommand{\indlincat}{{\indl \in \catind}}

\newcommand{\elmo}{{a}}

\newcommand{\elmu}{{a'}}

\newcommand{\spafix}{{Z}} 

\newcommand{\filt}{{F}}

\newcommand{\openfil}{{U}} 

\newcommand{\obfil}{{V}}

\newcommand{\obfilsec}{{U}}

\newcommand{\obfiltri}{{U'}}

\newcommand{\obfilu}{{U}}

\newcommand{\obfilw}{{W}}

\newcommand{\obfilinfil}{{\obfil \in \filtzar}}

\newcommand{\filtzar}{{F}}

\newcommand{\filtgab}{{GF}} 

\newcommand{\imu}{{i}} 

\newcommand{\setfin}{{I}}

\newcommand{\setfinl}{{\tilde{I}}}
\newcommand{\setfinsec}{{I'}}

\newcommand{\setfintri}{{\tilde{I}}}

\newcommand{\arlsec}{{\psi}}
\newcommand{\arl}{{\varphi}}

\newcommand{\tupf}{{f}}
\newcommand{\tupfsec}{{g}}

\newcommand{\vecunit}{{e}} 

\newcommand{\spay}{Y}
\newcommand{\spaz}{Z}

\newcommand{\spav}{{V}}

\newcommand{\spaa}{{A}}

\newcommand{\spau}{{Z}}

\newcommand{\spaw}{{W}}

\newcommand{\incl}{{i}}

\renewcommand{\curve}{{C}}

\newcommand{\sing}{{S}}

\newcommand{\morbc}{{\psi}} 

\newcommand{\moryx}{{\morbc : \spay \to \spax}}

\newcommand{\numnor}{{n}} 

\newcommand{\point}{{P}}

\newcommand{\pointy}{{y}}

\newcommand{\embpoint}{{i}} 

\newcommand{\project}{{p}}

\newcommand{\projcech}{{q}} 

\newcommand{\inject}{{i}}

\new{\normcomp}{{|}}

\newcommand{\filtcatindtocatopenspax}{{ \filt:\indexcat \to
\catopenspax}}

\new{\gen}{{g}}  

\new{\classres}{{u}} 

\new{\classcoho}{{c}}

\new{\homcon}{{\delta}} 

\new{\ringnorm}{{\ring^\nor}}

\new{\spaxnorm}{{\spay}}

\new{\linel}{{L}}

\newcommand{\group}{{G}}

\newcommand{\autg}{{G}}

\newcommand{\field}{{K}}

\newcommand{\quotfield}{{Q}}

\new{\cloalg}{\overline}

\newcommand{\fieldl}{{L}}

\newcommand{\fieldres}{{\kappa}}

\newcommand{\fieldmax}{{k}} 

\newcommand{\ring}{{R}}

\newcommand{\ringl}{{L}}

\newcommand{\ringtop}{{R'}} 

\newcommand{\ringsec}{{S}}

\newcommand{\ringsecsec}{{R''}}

\newcommand{\ringcomp}{{\rings_0}}  

\newcommand{\ringtri}{{R''}}

\newcommand{\opentop}{{U}} 

\newcommand{\ringt}{{T}}

\newcommand{\homtest}{{\varphi}}

\newcommand{\homlift}{{\psi}}

\newcommand{\homring}{{\varphi}}

\newcommand{\homsurj}{{\varphi}}

\newcommand{\ringcov}{{S}}

\newcommand{\ringcovindic}{{\ringcov_\indcov}}

\newcommand{\fua}{{a}}
\newcommand{\fub}{{b}}

\newcommand{\fuf}{{f}}
\newcommand{\fug}{{g}}
\newcommand{\fuh}{{h}}

\newcommand{\fuq}{{q}}

\newcommand{\funum}{{a}}
\newcommand{\fudeno}{{f}} 

\newcommand{\fucov}{{\fudeno}} 

\newcommand{\fudenosec}{{g}} 
\newcommand{\fucovsec}{{\fudenosec}}

\newcommand{\runfuf}{{\fuf_1 \comdots \fuf_\numgen}}

\newcommand{\indfu}{{i}} 

\newcommand{\numvar}{m}

\newcommand{\var}{z}

\newcommand{\varx}{z}
\newcommand{\vary}{w}
\newcommand{\varz}{v}

\newcommand{\varu}{u}

\newcommand{\vara}{{u}}
\newcommand{\varasec}{{s}}
\newcommand{\varusec}{{t}}

\new{\expoideal}{{d}}

\newcommand{\expoid}{{d}}

\new{\varcoef}{{t}}  

\new{\fucoef}{{g}}

\new{\indcoef}{{\nu}}

\new{\vargen}{{z}}

\newcommand{\obcomp}{{Z}} 

\newcommand{\indcompsec}{{j'}}

\newcommand{\obcompindic}{{\obcomp_\indcomp}}

\newcommand{\numcomp}{{m}}

\newcommand{\numcompbreak}{{s}}

\newcommand{\ideal}{{I}}

\newcommand{\idealsec}{{J}}

\newcommand{\idealsheaf}{{\shI}} 

\newcommand{\idealsubring}{{\ideal \subseteq \ring}}

\newcommand{\numgen}{{n}} 

\newcommand{\numfuf}{{n}}

\newcommand{\numfufbreak}{{k}}

\newcommand{\indgen}{{i}}

\newcommand{\indgenset}{{I}}

\newcommand{\indgeninset}{{\indgen \in \indgenset}}

\newcommand{\idealminor}{{J}}

\newcommand{\intclo}{\overline }

\new{\fuarb}{{g}}

\newcommand{\vart}{{T}} 

\begin{document}

\numberwithin{equation}{section}


{
\renewcommand{\topo}{\heartsuit}
\renewcommand{\clop}{\clubsuit}%

\ren{\field}{{\FF}}

\ren{\ring}{{\RR}}

\ren{\alg}{{\AA}}

\ren{\rings}{{\mathbb{S}}}

\ren{\ideal}{{\II}}

\ren{\modul}{{\MM}}

\ren{\submod}{{\NN}}

\ren{\spax}{{\XX}}

\ren{\spay}{{\YY}}

\ren{\spaz}{{\ZZ}}

\ren{\modl}{{\NN}}

\ren{\modm}{{\MM}}

\ren{\modr}{{\LL}}

\ren{\obcov}{{\UU}}

\ren{\openfil}{{\UU}}

\ren{\opentop}{{\UU}}

\ren{\obtopa}{{\UU}}

\ren{\obtopo}{{\VV}}}

\setboolean{book}{true}

\setboolean{book}{false}

\title[Grothendieck topologies and ideal closure operations]
{Grothendieck topologies and ideal closure operations}

\author[Holger Brenner ]{Holger Brenner}
\address{ \hspace{-0.9mm}Department \hspace{-0.9mm} of \hspace{-0.9mm} Pure \hspace{-0.9mm} Mathematics,\hspace{-0.9mm}
University \hspace{-0.9mm} of \hspace{-0.9mm} Sheffield,\hspace{-0.9mm}
Hicks \hspace{-0.9mm} Building,\hspace{-0.9mm} Houns\-field Road, Sheffield S3 7RH, United Kingdom}
\email{H.Brenner@sheffield.ac.uk}





\begin{abstract}
We relate closure operations for ideals and for submodules to
non-flat Gro\-then\-dieck topologies. We show how a Gro\-then\-dieck
topology on an affine scheme induces a closure operation in a
natural way, and how to construct for a given closure operation
fulfilling certain properties a Grothendieck topology which induces
this operation. In this way we relate the radical to the surjective
topology and the constructible topology, the integral closure to the
submersive topology, to the proper topology and to Voevodsky's
$h$-topology, the Frobenius closure to the Frobenius topology and
the plus closure to the finite topology. The topologies which are
induced by a Zariski filter yield the closure operations which are
studied under the name of hereditary torsion theories. The
Grothendieck topologies enrich the corresponding closure operation
by providing cohomology theories, rings of global sections, concepts
of exactness and of stalks.
\end{abstract}

\maketitle



\section*{introduction}

\ren{\var}{{X}}

Hilbert's Nullstellensatz states that for a polynomial
$\fuf \in \CC[\var_1 \komdots \var_\numvar]$
which vanishes on the common zero set
$V(\runfuf) \subseteq \CC^\numvar$
of a set of polynomials there exists
$\expok \in \NN$ such that
$\fuf^\expok \in (\runfuf) =\ideal$,
or, in other words, that $\fuf$ belongs to the
\inem{radical} $\rad (\ideal)$ of $\ideal$.
This theorem is probably the first instance in history where an ideal is replaced by a bigger
ideal which reflects better certain properties of the ideal in the
given context.
The \inem{integral closure} $\intclo{\ideal}$ of an
ideal in a commutative ring is the biggest ideal with the property
that it has the same Hilbert-Samuel multiplicity as the given ideal.
It was first studied by Zariski in the context of equisingularity problems.
Northcott and Rees \cite{northcottreesreduction} enriched
this notion introducing the theory of reductions and analytic
spread.

In the last 20 years, the \inem{tight closure} $\ideal^*$, as
defined by Hochster and Huneke in positive characteristic
(\cite{hochsterhunekebriancon}, \cite{hunekeapplication}),
has proved useful in homological algebra, commutative
algebra, invariant theory, singularity theory and algebraic
geometry. Tight closure comes along with \inem{plus closure}
(conjecturally the same),
\inem{Frobenius closure}, \inem{regular closure}, \inem{diamond closure} and \inem{solid
closure}.
It can be extended to rings containing a field of
characteristic $0$ by reduction to positive characteristic and there
have been several attempts to obtain closure operations with similar
features also in mixed characteristic
(\cite{hochstersolid}, \cite{brennerparasolid}, \cite{heitmannplus},
\cite{heitmannplusextension}, \cite{heitmanndirectsummand}). These
constructions are related to the search for \inem{big Cohen-Macaulay
modules} and \inem{algebras}.
Other closure operations include the \inem{weak subintegral closure}, the
\inem{blow-up closure}, Ratliff's $\catideal$-\inem{closure}, the $2,3$-\inem{closure}
and the \inem{continuous closure}
(\cite{leahyvitulliweaksubintegralclosure},
\cite{faridiblowup}, \cite{ratliffdeltaclosure}, \cite{odayoshida23}, \cite{brennercontinuous})

These closure operations have more or less natural extensions to
submodules. Also, a coherent subsheaf of a locally free sheaf is
often replaced by its \inem{reflexive hull} or by its
\inem{saturation} in the sense that the quotient is torsion-free.
Moreover, some constructions in commutative algebra can be thought
of as a closure of $0$ inside an $\ring$-module $\modul$, like its
\inem{torsion} or the \inem{submodule of sections with support}
inside a given closed subset of $\specring$, as it is used in the
definition of local cohomology. The idea of torsion has been
systematically studied under the name of \inem{hereditary torsion
theories} in the categorial context, but in particular for modules
over a ring (\cite{blandtorsiontheory}, \cite{golantorsiontheory}).

On the other hand, \inem{Grothendieck topologies}, in particular the
\inem{\'{e}tale topology} for schemes, were introduced to give for
varieties over a finite field a suitable replacement for the complex
topology, in order to define algebraically a cohomology theory with
similar properties as the singular cohomology in algebraic topology
(\cite{artingrothendieck}, \cite{SGA4}, \cite{milne}).
This development culminated in Deligne's proof of the Weil
conjectures.

The concept of a Grothendieck topology axiomatizes the idea of an
\inem{open set} and of \inem{coverings}. However, the open sets need
neither be subsets nor is a covering something which can be tested
by containment of points. Instead, coverings are given as a
collection of morphisms $\obcov_\indcov \to \obcovp$,
$\indcovinset$, fulfilling certain structural conditions, and points
need not to exist at all. With this notion of coverings at hand, one
can talk about \inem{local} and \inem{global} properties and one can
define \inem{presheaves}, \inem{sheaves}, \inem{sheafification},
\inem{exactness} and \inem{cohomology}.

\medskip
The basic idea of this (and subsequent) paper is to describe closure
operations for ideals and submodules over a commutative ring $\ring$
as a \inem{sheafification procedure} in a suitable Grothendieck
topology over $\spaxeqspecring$. More generally, for a morphism
between \inem{ringed sites} (i.e. categories with a Grothendieck
topology and a sheaf of rings) we can talk about \inem{extension}
and \inem{contraction} of modules and hence of a closure operation.
For $\ring$-modules $\submodul$, a Grothendieck topology on $\spax\!$
\--denoted by $\!\spax_\topo$\--
yields sheafified $\O_\topo$-modules
$\submod_\topo \to \modul_\topo$ on $\spax_\topo$ and an image sheaf
$\submod^\topo \subseteq \modul_\topo$. Then we declare for
$\elem \in \modul$ that
$\elem \in \submod^\clop$ if and only if
$\shef(\elem) \in \Gamma(\spax_\topo, \submod^\topo)$, where $\shef:
\modul \to \Gamma(\spax_\topo, \modul_\topo)$ is the sheafification
homomorphism. This means, in a certain sense, that
$\elem \in \submod^\clop$ holds if ``$\elem \in \submod$'' holds \inem{locally}
in the Grothendieck topology, and this means that it holds for a
covering. In particular, $\elem \in 0^\clop$, if $\elem =0$ holds
locally in the topology.

The classical Grothendieck topologies, like the Zariski, the flat and the
\'{e}tale topology
(this is also true for the Nisnevich \emph{topology}
\cite{nisnevich}, which is between the
Zariski and the \'{e}tale topology, and for the \emph{primitive
topology} of M. Walker \cite{walkerprimitive}) do not yield in this
way anything new, i.e.
$\submod^\clop = \submod$.
This is due to the fact that the assignment
$\obv \mapsto \Gamma(\obv, \struto^* \submod)$ for
$\struto:\obv \to \spaxeqspecring$ flat
(denoted by $W(\modul)$ in \cite[Example VII.2c)]{SGA4} and in \cite{milne}),
is already a sheaf in the flat topology \cite[Proposition I.2.18 and
Corollary II.1.6]{milne}, and since a faithfully flat homomorphism
$\ring \to \rings$ is \inem{pure},
which means that
$\modul \to \modul \tensorr\rings$ is injective for
every $\ring$-module $\modul$. If
$\homring:\ring \to \rings$ is however not pure,
then the inclusion
$\ideal \subset \homring^{-1}(\ideal \rings)$ might be strict, and if
we declare $\homring$ to be a covering, then
$\homring(\fuf)\in \ideal \rings$ says that $\fuf$ belongs locally to the ideal $\ideal$,
hence to $\ideal^\clop$.

Lets explain more concretely how a Grothendieck topology induces a
closure operation and what properties it has. Lets restrict for
simplicity in this introduction to the \inem{affine}
\inem{single-handed} case. This means that we have a category of
$\ring$-algebras, closed under composition and tensor product, and
certain $\ring$-algebra homomorphisms are declared to be covers (so
all coverings consist of just one morphism). In this situation the
closure operation given by sheafification in the topology boils down
to saying that for $\fuf
\in \ring$ and an ideal $\idealsubring$ the containment $\fuf \in
\ideal^\clop$ holds if and only if there exists a cover
$\homring:\ring \to \rings$ such that
$\homring(\fuf) \in \ideal\rings$.
More generally, for $\ring$-submodules $\submodul$ and
$\elem \in \modul$, the containment $\elem \in \submod^\clop$ holds
if and only there exists such a cover such that $\elem \tensor 1$ is
in the image of
$\submod \tensorr \rings \to \modul\tensorr \rings$.
It follows immediately that
$\submod \subseteq \submod^\clop =(\submod^\clop)^\clop$
holds and that the closure operation is
\inem{order preserving}, that is for
$\submod \subseteq \submod'$ in
$\modul$ one has
$\submod^\clop \subseteq (\submod')^\clop$.

A great deal of this project is in finding the right concept of
covers for a given closure operation and to study properties of
these morphisms instead (or rather as a way) of studying the closure
operation. For some closure operations it is immediately clear which
ring homomorphisms should be considered to be a cover, but not for
all. We go through the closure operations which we will deal with in
this and in subsequent papers
(\cite{blicklebrennersubmersion}, \cite{blicklebrennertighttight},
\cite{brennercontinuous}, \cite{brennersemiintegraltest})
and we describe the natural covers for them.

Let's look first at the \inem{Frobenius closure} in positive
characteristic $\posp$, which is defined by $\fuf \in \ideal^\frobf$
if and only if there exists a power of the Frobenius
$\frobmo^\expoe :\ring \to \ring$,
$\fuf \mapsto \fuf^{\posp^\expoe}$,
such that
$\fuf^{\posp^\expoe} \in \frobmo^\expoe (\ideal)
= :\ideal^{[\posp^\expoe]}$.
In this case one wants to allow exactly
the $\frobmo^\expoe$ to be covers to define the \inem{Frobenius
topology}. However, in order to get a theory with good properties
regarding base change and also to have only covers of finite type,
it is better to define a homomorphism $\ring \to \rings$ of finite
type to be a cover in the Frobenius topology if there exists a
factorization
$\frobmo^\expoe :\ring \to \rings \to \ring$ for some
$\expoe$ such that $\Spec \ring \to \Spec \rings$ is surjective.

The \inem{plus closure} of an ideal $\ideal$ in a domain $\ring$
consists of all elements
$\fuf \in \ring$ such that there exists a
finite extension of domains
$\ring \subseteq \rings $ with
$\fuf \in\ideal \rings$.
Here the appropriate concept of a cover, which works
also for rings with zero divisors, is to impose that $\ring \to
\rings$ is finite and with surjective morphism
$\Spec \rings \to\specring$. This gives the \inem{finite topology}.

What is the appropriate concept of a cover for the \inem{radical} of
an ideal (and what is the radical of a submodule)? If we take
Hilbert's Nullstellensatz as our guide, then $\fuf$ belongs to $\rad
(\ideal)$ if and only if $V(\fuf) \supseteq V(\ideal)$ in
$\specring$. This means for all prime ideals
$\ideal \subseteq \fop$
that also $\fuf \in \fop$ holds,
or that all homomorphisms $\ring \to
\fieldres(\fop)$ which map $\ideal$ to $0$ also map $\fuf$ to $0$.
This property is then true for all homomorphisms from $\ring$ to a
field. Which ring homomorphisms $\ring \to \rings$ behave nicely
with respect to homomorphism to fields? There are several
possibilities, but we will use the fact that
$\Spec \rings \to\specring$ is surjective if and only if for all
$\ring \to \field $
we have
$\rings \tensorr \field \neq 0$
(or $\Spec (\rings \tensorr \field) \neq \emptyset$).

Declaring \inem{spec-surjective morphisms} to be covers is in fact a
way to obtain the radical via a Grothendieck topology, the
\inem{surjective topology}.
For, if $\ring \to \rings   $ is
spec-surjective and
$\fuf \in \ideal \rings$, then for all
homomorphisms to a field
$\mofi: \ring \to \field$ it follows by
the faithful flatness of
$\field \to \rings \tensorr \field $
($\neq 0$) that
$\mofi (\fuf) \in \ideal \field$.
Hence if $\ideal\field=0$, then also $\mofi(\fuf) = 0$.

For the other direction we have to construct for $\fuf \in \rad
(\ideal)$ a spec-surjective homomorphism
$\homring:\ring \to \rings $
with
$\homring(\fuf) \in \ideal \rings$. For this we use
\inem{forcing algebras}, as introduced by Hochster
\cite{hochstersolid} in connection with solid closure in an (at the
end not successful) attempt to construct a closure operation in all
characteristics with tight closure properties \cite{hochstersolid}.
Fix ideal generators $\ideal = (\runfuf)$. The forcing algebra for
these data is
$\algforc
= \ring[\vart_1 \komdots \vart_\numfuf]/
(\fuf_1\vart_1 \plusdots \fuf_\numfuf\vart_\numfuf+\fuf)$.
In $\algforc$, $\fuf \in \ideal \algforc$ holds, and every
$\ring$-algebra $\rings$ with $\fuf \in \ideal \rings$ factors (non
uniquely) through $\algforc$. It is clear by this property that
forcing algebras should play an important role in the study of
closure operations. In the context of the radical, $\fuf \in \rad
(\ideal)$ holds if and only if the corresponding forcing algebra
induces a spec-surjective morphism.

There are several equivalent characterizations of the \inem{integral
closure} of an ideal. If $\ring$ is noetherian, then $\fuf \in
\intclo{\ideal}$ if and only if $\fuf \in \ideal  \dvd$ for all ring
homomorphisms $\ring \to \dvd$ to discrete valuation domains. Hence
a forcing algebra
$\Spec \algforc \to \specring$ should declared to
be a cover if and only if for all $\Spec \dvd \to \specring$ there
exists a lifting
$\Spec \dvd \to \Spec \algforc$. For forcing
algebras, but not in general, this path-lifting property is
equivalent to say that for every $\ring \to \dvd$ the homomorphism
$\dvd \to \algforc \tensorr \dvd$ has a $\dvd$-module section or
that it is pure. There exists even a topological characterization
for an $\dvd$-algebra $\ringt$ to be pure, namely that over the
prime ideals in $\Spec \dvd$, i.e. $0 \subset \fom_\dvd$, there
exists a pair of prime ideals
$\foq \subset \fop$ in $\Spec \ringt$
mapping to $0$ and $\fom_\dvd$ respectively.

This topological observation ties in purity after base change to
discrete valuation domains with submersions. A \inem{submersion}
$\mor : \Spec \rings \to \specring$ is a surjective mapping such
that the Zariski topology of $\specring$ is the image (or quotient)
topology under the morphism, that is, a subset
$\obw \subseteq \specring$ is open if and only if its preimage
$\mor^{-1}(\obw)$
is open. This definition can be found in SGA 1
(see \cite[D\'{e}finition IX.2.1]{SGA1} and \cite[D\'{e}finition 3.10.1]{EGAI}) and also in
\cite[Definition 3.1.1]{voevodskyhomology}
under the name of \inem{topological epimorphism}.
A morphism
$\mor :\Spec \rings \to \specring$ ($\ring$ noetherian) is
\inem{universally submersive} if and only if for all $\ring \to \dvd$ to
discrete valuation domains
$\dvd \to \rings \tensorr \dvd$ is pure
(\cite[Remarque IX.2.6]{SGA1}, \cite[Th\'{e}or\`{e}me 37]{picavet}).
Hence the \inem{submersive topology}, in which the universal
submersions are covers, gives a Grothendieck topology which induces
the integral closure. Surprisingly, this topology is very closely
related to Voevodsky's $h$-\inem{topology}, which arose as a new
approach to understand the algebraic topology of varieties and which
was the first systematic study of a non-flat Grothendieck topology
(\cite{voevodskyhomology}, \cite{suslinvoevodskysingularhomology}).
We will study submersions and integral closure in a separate paper
\cite{blicklebrennersubmersion}.

Recall that for a ring of positive characteristic $\posp$ the
\inem{tight closure} of an ideal $\ideal$ is by definition
$\ideal^*
= \{ \fuf \in \ring:\, \mbox{there exists } \eletest \in \ring^o,
\eletest \fuf^\potq \in \ideal^{[\potq]} \mbox{ for all } \potq =
\posp^\expoe\}$, where $\ring^o$ denotes the complement of all
minimal primes. Tight closure agrees in positive characteristic (if
$\ring$ has a completely stable test element) with \inem{solid
closure} \cite[Theorem 8.6]{hochstersolid}.
This means that the containment $\fuf \in \ideal^*$ for an $\fom$-primary ideal
$\ideal$, where $\fom$ is a maximal ideal of height $\heighd$, is
characterized by the property that the local cohomology module
$H^\heighd_\fom(\algforc)$ is not zero, where $\algforc$ is the
forcing algebra \cite[Corollary 2.4]{hochstersolid}. This property
makes sense in every (even mixed) characteristic, leading to
\inem{solid closure} in general and to the \inem{solid topology}. It
is however known that for rings containing a field of characteristic
$0$ solid closure does not exhibit all features one expects from a
tight closure type theory, in particular, ideals in a regular ring
are not always solidly closed \cite{robertscomputation}. A variant
of solid closure, called \inem{parasolid closure}, introduced in
\cite{brennerparasolid}, is defined characteristic free by the
property that all cohomology classes
$1/\fug_1 \cdots \fug_\heighd \in H^\heighd_\fom (\ring)$, where
$\fug_1 \komdots \fug_ \heighd$
is a system of parameters in $\ring$, do not vanish in
$H^\heighd_\fom (\algforc)$. We will continue the search for the
\inem{tight topology} along these lines in
\cite{blicklebrennertighttight}.

Another source for closure operations is given by \inem{torsion
theories}. Torsion theories have been studied in various categorial
contexts; for the category of $\ring$-modules over a commutative
ring it means that we have two classes $\torclass$ and $\freeclass$
of modules (`\inem{torsion modules}' and `\inem{torsion-free
modules}' in an abstract sense) such that $\torclass$ is closed
under quotients and $\freeclass$ is closed under submodules and such
that for every $\ring$-module $\modul$ there exists a short exact
sequence
$0 \to \modtor \to \modul \to \modfree \to 0$
($\modtor \in \torclass$, $\modfree \in \freeclass$).
In particular then, every module has a \inem{biggest torsion submodule},
and this can be considered as the closure of $0$. A torsion theory
is called \inem{hereditary} if the torsion class is also closed
under submodules.

A hereditary torsion theory is basically given by its so-called
\inem{Gabriel filter}, which consists of all ideals
$\idealsec \subseteq \ring$ such that $\ring/\idealsec$ is torsion.
In the noetherian case, this Gabriel filter is stable under radical,
so it is essentially a (Zariski-)topological filter.
Such a \inem{Zariski filter} gives rise to a (non-affine) Grothendieck topology,
where the Zaiski open subset
$\obfil \subseteq \specring$ is said to be a cover if
$\obfil$ belongs to this filter.
So here the covers are flat, but they are not surjective.
This is quite a natural construction, since the relevant
information lies often already in certain open subsets, like e.g.
the divisor class group depends only on points of codimension one.
The induced closure of a submodule $\submod$ consists of all
$\elem \in\modul$ such that the restriction
$\elem|\obfil \in \submod$ for some open subset in the filter.
For $\submod=0$ and one fixed $\obu = D(\foa)$ this also equals
$\Gamma_{V(\foa)} (\modul) = \Gamma_{\foa}(\modul)$,
the \inem{module of sections with support} in $\obv
(\foa)$, as used in the definition of local cohomology.
In general, a typical feature for hereditary torsion theories is that the
closure operation $0 \mapsto 0^\clop$ is left exact.

{\ren{\hommod}{{\varphi}}

So for a lot of closure operations we can describe a Grothendieck
topology which induces it, but one can also give an exact
characterization of the closures for which one can construct a
Grothendieck topology.
This construction uses again forcing algebras.
Assume that
$\fuf \in (\fuf_1 \komdots \fuf_\numfuf)^\clop$ in a given closure operation $\clop$. Then it
is natural to declare the corresponding forcing algebra $\algforc$
to be a cover.
It may then happen that $\fug \in \idealsec
\algforc$, where $\idealsec$ is another ideal, so $\fug$ will belong
to the induced closure of $\idealsec$ in every topology where $\ring
\to \algforc$ is a cover.
So in order that this construction does
not produce to many elements in the closure we need that $\fug \in
\idealsec^\clop$.
We express this requirement by saying that $\clop$-\inem{admissible forcing
algebras} must be $\clop$-\inem{pure}.
This property characterizes
exactly the closure operations which come from a Grothendieck
topology.
It implies that the closure is \inem{module persistent},
i.e. for a module homomorphism
$\hommod: \modul \to \modul'$ one has
$\hommod(\submod^\clop) \subseteq \hommod(\submod)^\clop$ and that
it is \inem{independent of the presentation},
i.e. if
$\hommod :\modul \to \modul/\submod$,
then $\submod^\clop =\hommod^{-1}(0^\clop)$.}

An important feature of closure operations are \inem{persistence
properties} with respect to ring changes, which should be reflected
by the Grothendieck topologies. A closure operation is said to be
persistent with respect to a class of ring homomorphisms if for
$\homring :\ring \to \ringsec $ and $\idealsubring$ one has
$\homring (\ideal^\clop) \subseteq \homring (\ideal)^\clop$.
Most closure operations mentioned above
(with the exception of torsion theories)
are persistent with respect to all ring homomorphisms; for tight
closure this is a non-trivial fact.
If we have a Grothendieck topology on a category of rings,
then a ring homomorphism $\ring \to \ringsec  $
might or might not respect the Grothendieck topology, i.e. induces a
\inem{site morphism}
$(\Spec \ringsec )_\topo \to (\specring)_\topo$
or not.
If it does, then the closure operation is persistent
with respect to this homomorphism.

\medskip
After having explained how a Grothendieck topology induces a closure
operation we shall look now at some other features which a
Grothendieck topology provides. Some of these features have occurred
in the study of closure operations in some disguise, and
Grothendieck topologies give a conceptual view on these.

The containment in the Frobenius closure can also be expressed by
saying that $\fuf \in \ideal  \ring^\perf$, where $\ring^\perf$ is
the \inem{perfect closure} of $\ring$. The containment in the plus
closure is equivalent with $\fuf \in \ideal \ring^+$, where
$\ring^+$ is the integral closure of $\ring$ in an algebraic closure
of its quotient field $\quotfield( \ring)$ ($\ring$ a domain),
called the \inem{absolute integral closure}. The study of tight
closure is closely related to the search for \inem{big
Cohen-Macaulay modules} or \inem{algebras}. The common feature here
is that there exists one ``big'' algebra which encompasses in a
certain sense all relevant covers, like in topology the universal
covering space encompasses all unramified covers. We systemize this
observation with the concept of an \inem{absolute filter} in a
Grothendieck topology, at which we can evaluate a (pre-)sheaf to get
an \inem{absolute stalk}. This is similar to the \'{e}tale topology
of a local ring, where a suitable colimit over all \'{e}tale covers
yields the \inem{strict Henselization} of the ring. An absolute
stalk should carry essentially the whole information of the
Grothendieck topology and of the closure operation, but its
construction and its explicit description might be difficult, and
its functorial properties are weak in general.

A Grothendieck topology provides a concept of \inem{exactness} for a
complex of $\ring$-modules $\modl \to \modm \to \modr$, namely the
exactness of the complex of sheaves $\modl_\topo \to \modm_\topo \to
\modr_\topo$ on the underlying category of the Grothendieck
topology.
This notion depends not only on the concept of a covering,
but also on the seize of the category.
Quite often we will work with
a \inem{small site}, where (roughly speaken) only objects occur
which are relevant for some covering of the final object
$\specring$.
This exactness on the small site is still a strong
condition and we will also consider weaker notions which are closely
related to the property that the kernel in the complex of modules is
contained in the closure of the image.
This property has been studied in tight
closure theory under the name of \inem{phantom homology}
\cite[Chapter 10]{hunekeapplication} and for the integral closure
(for the Koszul complex) by D. Katz \cite{katzdcomplex}.
In the surjective topology, the sheaf exactness is equivalent to the
property that $\modl \tensorr \field \to \modm \tensorr \field \to
\modr \tensorr \field $ is exact for every ring homomorphism $\ring
\to \field $ to a field.

A Grothendieck topology on $\spaxeqspecring$ gives also for an
$\ring$-module $\modul$ via sheafification a \inem{global module of
sections}, namely $\Gamma( \spax_\topo, \modul_\topo)$
(with a homomorphism from $\modul$ to it),
and in particular a \inem{global ring of sections}
$\Gamma(\spax_\topo, \O_\topo)$.
An element in this ring is represented by a
compatible element $\fuglob \in \rings$, where $\ring \to \rings $
is a cover and where compatible means that
$\fuglob \tensor 1 - 1 \tensor \fuglob \in \rings\tensorr\rings$
is $0$ in the topology, i.e. in some cover
$\rings \tensorr\rings \to \ringt$.
The assignment
$\ring \mapsto\Gamma(\spax_\topo, \O_\topo)$
is functorial for ring homomorphisms
$\rings \to \ringsec$ which respect the Grothendieck topology.

For the Frobenius topology this construction yields the perfect
closure, since for $\fuglob \in \frobe \ring$ the element $\fuglob
\tensor 1 - 1 \tensor \fuglob $ is nilpotent. In this case the ring
of global sections equals the absolute stalk. In the surjective
topology (in characteristic $0$), a global function is essentially
given by an affine \inem{constructible partition} $\spax= \obv_1
\uplusdots \obv_\numpart$, where $\obv_\indpart$ is an intersection
of a closed and an open subset, and functions $\fuglob_\indpart \in
\Gamma(\obv_\indpart, \O_{\obv_\indpart})$ (an absolute stalk is
given as the product of the algebraic closures of all residue class
fields).

For other topologies the ring of global sections is closer to the
ring itself, e.g. for the plus closure and in the submersive
topology (in characteristic $0$), it is the so-called
\inem{seminormalization}. In particular, $\fuglob \in \ideal^\clop$
does not mean that
$\fuglob \in \ideal \Gamma(\spax_\topo,\O_\topo)$,
which is a global and much stronger condition.

As the module of global sections is the $0$th cohomology,
we arrive at \inem{sheaf cohomology} in general.
Cohomology can be defined quite generally via injective resolutions,
and at least the first cohomology might be computed by
\inem{\v{C}ech cohomology}.
The cohomology of a coherent $\O_\spax$-module over an
affine scheme might be non-trivial, contrary to the flat,
\'{e}tale or Zariski topology.
A typical short exact sequence arises from
ideal generators $\ideal =(\runfuf)$ and gives
$0 \to \Syz \to \O_\topo^\numfuf \stackrel{\runfuf}{\to} \ideal^\topo \to 0$,
where $\ideal^\topo \subseteq \O_\topo$ is the image sheaf of $\ideal_\topo$ inside
$\O_\topo$ and where $\Syz$ is defined as the kernel of the
surjection (not as $\Syz_\topo$).
The corresponding long cohomology sequence is
$0 \to \Gamma(\spax_\topo, \Syz) \to \oplus_\numfuf
\Gamma(\spax_\topo, \O_\topo) \stackrel{\runfuf}{\to}
\Gamma(\spax_\topo, \ideal^\topo)
\to H^1(\spax_\topo, \Syz) \to \oplus_\numfuf H^1(\spax_\topo, \O_\topo)$.
In the case of the Frobenius topology and of the surjective topology we can prove that
$H^1(\spax_\topo, \O_\topo)=0$ and also that
$H^1(\spax_\topo, \Syz)=0$,
and so all $\fuglob \in \ideal^\clop$ come from a `global solution'
$\fuglob \in \ideal \Gamma(\spax_\topo, \O_\topo)$.

In cases where
$\Gamma(\spax_\topo, \O_\topo)=\ring$ we have
$\Gamma(\spax_\topo,\ideal^\topo)= \ideal^\clop \subseteq \ring$ and
the image from the left is
$\ideal \Gamma(\spax_\topo, \O_\topo) =\ideal $,
so $\ideal^\clop/ \ideal \subseteq H^1(\spax_\topo,\Syz)$,
which is a bijection provided that $H^1(\spax_\topo,\O_\topo)=0$.
We expect this to hold in the submersive topology for a normal domain
in characteristic $0$. The relation between the first cohomology of
the sheaf of syzygies and the closure quotient
$\ideal^\clop/\ideal$ extends to a relation between the naive
sequence
$0 \to \ideal^\clop/\ideal \to \ring/\ideal
\to \ring/\ideal^\clop \to 0$ and the long exact sequence
$\,\,\, \ldots \to H^1(\spax_\topo, \Syz) \to H^1(\obu_\topo, \Syz)
\to H^2_{\spax_\topo- \obu_\topo}(\Syz) \to \ldots \,\,\,$
relating local and global cohomology, where $\obu=D(\ideal)$.

If a closure operation or a Grothendieck topology is defined for a
class of rings
(usually the class of all noetherian rings),
then it is natural to ask for which rings the Grothendieck topology is
\inem{pure} and the closure is trivial.
Purity means here that the sheafification does not annihilate anything,
and for a single-handed Grothendieck topology it just means that
every cover is pure.
A stronger condition is that every cover allows a module section or
even that every cover allows a ring section, like in the \'{e}tale
topology the algebraically simply connected varieties are
characterized by the property that every \'{e}tale cover has a
section.
The class of rings (or singularities) which are pure with
respect to a certain Grothendieck topology might be quite
significant, which is in particular true for tight closure, where
the so-called $\frobf$-\inem{regular} rings are closely related to
rational singularities.

\medskip
We believe that Grothendieck topologies provide a natural framework
to study closure operations. It gives a conceptual view on certain
features which have been loosely accompanying closure operations
ever since (like exactness or absolute covers). Once the `right'
Grothendieck topology for a given closure operation for ideals is
found we can immediately extend the closure to submodules.

The focus on morphisms coming along with studying the appropriate
coverings brings in a geometric flavor to questions on closure of
ideals. For example, the complex-topological characterization of
submersions yields at once the complex-analytic characterization of
the integral closure of submodules. Also, statements on closures
have often natural generalizations to morphisms, in such a way that
the original statement is regained by looking at the morphisms given
by forcing algebras.

Grothendick topologies provide also a new way to construct closure
operations. Here the most challenging problem is to give in mixed
characteristic an operation with similar properties as tight
closure. This would solve most remaining homological conjectures.
Less ambitious is the task to give a tight closure like operation in
equal characteristic zero without reduction to positive
characteristic. This was essentially accomplished in
\cite{brennerparasolid}, but a complex-analytic characterization
extending the relation between tight closure test ideals and
multiplier ideals remains to be done.

\medskip
We give a brief summary of the content and the structure of this
paper.

We start in part 1 with the concepts of Grothendieck topologies,
sheaves (\ref{presheafsubsection}), sheafification
(\ref{sheafifysubsection}) and of morphisms of (ringed) sites
(\ref{sitesubsection}).
A morphism of (ringed) sites is the natural
setting in which it makes sense to talk about extensions and
contractions of (sub-)modules, which gives a closure operation
(\ref{extensioncontractionsubsection}).
This approach has the advantage of being a relative setting, so that we can sheafifiy and
get the closures in several steps. Later on our site morphisms will
be between a scheme with its Zariski topology
(or even its trivial topology, where only isomorphisms are coverings)
and another larger (or finer) topology. In
(\ref{grothendieckschemesubsection})
we discuss such Grothendieck topologies on
a scheme and how a module defines a (pre-)sheaf. In
(\ref{globalsectionsubsection}) we discuss modules of global
sections and provide some lemmata which will later be of use in
computations for several Grothendieck topologies.

Starting with (\ref{filterstalksubsection}) we look at filters in a
Grothendieck topology and the stalks they lead to. Basically, a
filter is a cofiltered system (or diagram) in the topology. An easy
way to produce filters is by a fixation
(Construction \ref{fixationconstruction}). In some sense, a filter, and in
particular the irreducible ones, are a replacement for points
(\ref{irreduciblesubsection}). A filter which consists only of
covers and has the property that for every cover it has a
factorization is called \inem{absolute}
(\ref{absolutestalksubsection}). They carry in some sense the whole
information of the Grothendieck topology. In
(\ref{exactnesssubsection}) we discuss several notions of exactness
and in (\ref{cohomologysubsection}) we discuss cohomology, \v{C}ech
cohomology and local cohomology in a site. It might be a good idea
to read initially only (\ref{extensioncontractionsubsection}) of the
first part and go back if needed.

The second part relates Grothendieck topologies on a scheme, in
particular on an affine scheme, with closure operations. In
(\ref{forcingsubsection}) we start with the basics about forcing
algebras. In (\ref{admissiblesubsection}) we introduce
\inem{admissible closure operations}
which will turn out to be the ones for which one can construct a
Grothendieck topology.
We also introduce certain finiteness conditions for a closure operation.
In (\ref{exampleadmissiblesubsection}) we give examples of admissible
closure operations and also of non-admissible closures like the
Ratliff-Rush closure. In (\ref{closuregrothendiecksubsection}) we
describe, on the basis of (\ref{extensioncontractionsubsection}),
how a Grothendieck topology on  a scheme induces an admissible
closure operation. In (\ref{puregrothendiecksubsection}) we
characterize the pure topologies and in
(\ref{grothendieckconstructionsubsection}) we construct for a given
admissible closure a Grothendieck topology which induces it. In
(\ref{exactnessclosuresubsection}) we relate the notion of exactness
with properties of the closure, and in
(\ref{syzygycohomologysubsection}) we deal with the sheaf of
syzygies and how the closure operation is reflected in the
cohomology of this sheaf.


The following parts are dedicated in detail to specific closure
operations and their ``natural'' Grothendieck topologies and can be
read independently from each other.

Part \ref{surjectivesection} deals with the radical and the
surjective topology. In (\ref{radicalsubsection}) we recall the
radical of an ideal and introduce the \inem{radical of a submodule}.
The surjective topology is introduced in
(\ref{surjectivetopologysubsection}) and it is shown that it induces
the radical (Proposition \ref{surjectiveinduceradical}). In
(\ref{jacobsonsubsection}) we discuss the \inem{Jacobson radical} and the
\inem{Jacobson topology}, where a cover need only be surjective on the
closed points. It follows from a version of the Nakayama Lemma that
for a finitely generated module the Jacobson radical and the radical
is not the whole module (Corollary \ref{radicalnotall}). In
(\ref{surjectivefieldsubsection}) we deal with the surjective
topology over a field and show that the ring of global sections is
the perfect closure of the field. Note that though over a field a
ring is faithfully flat if and only if it is spec-surjective, the
surjective topology over a field is not the same as the faithfully
flat topology, since also the covers of $\alg \otimes_\field \alg$
are crucial for the topology. In
(\ref{surjectiveabsolutesubsection}) we show that the direct product
over the algebraic closures of all residue class fields of
$\specring$ provides an absolute stalk (Proposition
\ref{surjectiveabsolutestalk}).

In (\ref{constructiblesubsection}) we discuss
\inem{constructible partitions} of a scheme, which are the easiest non-trivial
surjective mappings, and which are rich enough to induce also the
radical.
Constructible partitions are easier to handle than arbitrary
surjections, yet they contain often the full information of the
surjective topology and many questions can be reduced to the
\inem{constructible topology}. Section
(\ref{surjectiveconstructiblesubsection}) provides some lemmata for
this. In
(\ref{surjectiveglobalsubsection}) we show that the ring of global
sections in the surjective topology for a ring of finite
(Krull) dimension over a field of characteristic $0$ is the same as in
the constructible topology (Theorem \ref{surjectiveglobalsection}).
So a global function is given by a collection of algebraic functions
on a constructible partition of the scheme. The exactness in the
surjective topology is characterized in
(\ref{surjectiveexactsubsection}) and is equivalent to the exactness
after tensoration with an arbitrary field. It is also equivalent to
global exactness, which is rather untypical for a Grothendieck
topology. Finally, in (\ref{surjectivecohomologysubsection}) we
prove that the first cohomology of coherent sheaves vanishes in the
surjective topology (Theorem \ref{surjectivecohomologytheorem} and
Corollary \ref{surjectivecohomologycoherent}).

Part 4 deals with the \inem{Frobenius topology} and the \inem{Frobenius
closure}
in positive characteristic. In (\ref{frobeniussubsection}) we recall
the definition of Frobenius closure of modules and in
(\ref{frobeniustopologysubsection}) we introduce the Frobenius
topology which induces the Frobenius closure (Proposition
\ref{frobeniusfrobenius}) and we characterize the purity of the
Frobenius topology (Proposition \ref{frobeniuspureequivalent}). In
(\ref{frobeniusglobalsubsection}) we compute the ring of global
sections in the Frobenius topology to be the perfect closure
(Theorem \ref{frobeniusglobal}). In (\ref{frobeniusexactsubsection})
we characterize exactness in the Frobenius topology and in
(\ref{frobeniuscohomologysubsection}) we show that the first
cohomology of coherent modules vanishes (Proposition
\ref{frobeniuscohomology} and Corollary
\ref{frobeniuscohomologycorollary}).

Part 5 deals with the \inem{plus closure} and the \inem{finite topology}.
In (\ref{plusclosuresubsection}) we recall the plus closure and
introduce the finite topology which induces the plus closure.
In (\ref{finiteglobalsubsection}) we show that the ring of global
sections for a domain of characteristic zero in the finite topology
is the seminormalization, i.e. the maximal subring of the
normalization which does not separate points
(Proposition \ref{finiteglobalsection}).
The absolute stalks in the finite
topology, which are just products of the absolute integral closures,
are described in (\ref{finitestalksubsection}).
Our results on the computation of cohomology in the finite topology are rather modest
(\ref{finitecohomologysubsection}).

Part 6 is concerned with \inem{torsion theories} and topologies given by
\inem{Zariski filters}. These topologies are flat but not surjective (and
not affine).
A topological filter in the Zariski topology defines a
Grothendieck topology by declaring the open sets inside the filter
to be covers (\ref{zariskifiltersubsection}). For a filter given by
just one open set the closure of $0$ is nothing but the module of
sections with support in the closed complement, as used in the
definition of local cohomology (\ref{submodulesupportsubsection}).
Section (\ref{torsiontheorysubsection}) provides the link between
Grothendieck topologies given by a Zariski filter and hereditary
torsion theories (Theorem
\ref{grothendieckleftexacttorsiontheory}). In
(\ref{divisorialsubsection}) we look at the \inem{divisorial topology},
where the Zariski filter consists of all open subsets which contain
all points of codimension one.
The module of global sections is under certain condition
the \inem{reflexive hull} of the module, i.e. its bidual
(Lemma \ref{divisorialglobalreflexive}), and the
closure of a submodule inside a reflexive module over a normal ring
is also the reflexive hull
(Proposition \ref{divisorialnormalclosure}).
We also describe the relation to \inem{symbolic powers} and we show that the
divisor class group is the Picard group in the divisorial topology
(Proposition \ref{divisorialdivisor}).

In the final part 7 we treat some further Grothendieck topologies.
In (\ref{completionsubsection}) we show how to define, starting from
an ideal $\foa$ in a ring $\ring$, a Grothendieck topology such that
the ring of global sections is the $\foa$-\inem{adic completion} of
$\ring$. Here a covering needs infinitely many objects. In
(\ref{propersubsection}) we treat the \inem{proper topology} (first
introduced by S.-I. Kimura in \cite{kimuraalexander}), where a cover
is given by a proper morphism. We show in Proposition
\ref{properintegral} that the induced closure is the \inem{integral
closure}. We will come back to this in the study of submersive
topologies \cite{blicklebrennersubmersion}. Section
(\ref{ratliffsubsection}) finally studies the
so-called $\catideal$-\inem{closure} of Ratliff, which is given by a
multiplicatively closed set of ideals. This closure is induced by
taking the blow-ups of these ideals as covers (Proposition
\ref{ratliffblowup}).

We do not touch in this paper the following topics: the topos of
modules defined by a Grothendieck topology, e.g. structure theory
for the $\O_\topo$-modules, the question when scheme morphisms
induce isomorphisms of sites, free resolutions in the topology,
minimal number of generators, reductions, core, projective
dimension, Picard group in the Grothendieck topology (the
multiplicative theory), other non-quasicoherent sheaves like
constant sheaves, the algebro-topological or homotopical content of
the Grothendieck topologies, extension problems for covers.

I thank R. Buchweitz, A. Conca, G. Dietz, D. Eisenbud, N. Epstein, T. Gaffney, S. Goto,
R. Hart\-shorne, M. Hashimoto, R. Heitmann, A. Kaid,  M. Katzman, E.
Kunz, M. Morales, P. Roberts, E. Sbarra, R. Sharp, V. Srinivas, M. Vitulli,
K.I. Watanabe for their interest and useful remarks. Most of all I
thank Manuel Blickle for critical remarks and many helpful
discussions during the whole project.

\section{Grothendieck topologies}

\subsection{Presheaves, Grothendieck topologies and sheaves}
\label{presheafsubsection}

\

\medskip
An (abelian) \emph{presheaf} on a category $\catt$ is a
contravariant functor $  \presh: \catt \to \catab$, where $\catab$
is the category of abelian groups
(or another abelian category).
For a morphism $\morj: \obu \to \obv$ we will often denote
$\presh(\morj): \presh (\obv) \to \presh(\obu)$ by $\rest^\obv_\obu(\morj)$
or simply by $\rest^\obv_\obu$ and call this the \emph{restriction}
from
$\obv$ to $\obu$.
\ifthenelse{\boolean{book}}{We denote by $\PSh
(\catt)$ the category of abelian presheaves on $\catt$, where a
morphism of presheaves $\presh \to \preshb$ is a natural
transformation of functors. The usual constructions in $\catab$ can
be carried out in the category $\PSh (\catt)$ of presheaves on
$\catt$ by doing them for each $\obu \in \catt$. Thus for example
the direct sum of presheaves $\presh$ and $\shG$ is defined by
$(\presh \oplus \preshb)(\obu) = \presh (\obu) \oplus \shG( \obu)$,
or the kernel of a morphism of presheaves
$\presheafmor: \presh \to \preshb$ is given by
$(\ker \presheafmor)(\obu) = \ker \presheafmor (\obu)$ and so forth.
This endows the category $\PSh(\catt)$ of
presheaves on $\catt$ with the structure of an abelian category.}{}
For a topological space, a sheaf is a presheaf which satisfies
certain properties with respect to coverings.
To talk about sheaves in a given category one has to define a notion of coverings with
certain structural properties
(see  \cite{artingrothendieck}, \cite{SGA4}, \cite{milne}, \cite{vistoligrothendieck}).
This is the basic idea of a Grothendieck topology.

Recall that a \emph{product} in a category for $\obu \to \obw$, $\obv \to \obw$
is an object $\obu \times_\obw \obv$ with projections to $\obu$ and
$\obv$ compatible over $\obw$ such that for every $\obtest \to \obu$, $\obtest \to \obv$
compatible over $\obw$ there exists a unique morphism
$\obtest \to \obu \times_\obw \obv$.

\begin{definition}
\label{Grothendieckdef}
Let $\catt$ be a category with products. A
\emph{Grothendieck topology}
(sometimes called a pretopology,
see \cite[D\'{e}finition II.1.1]{SGA4} or \cite[Definition 2.23 and
Remark 2.24]{vistoligrothendieck}) on
$\catt$ consists of a collection of coverings, that is, families
$\morcov_\indcov: \obcovindic \to \obu$, $\indcovinset$, where in each family the
target $\obu \in \catt$ is fixed. These coverings should satisfy.

\numiii
\begin{enumerate}

\item
Isomorphisms are coverings.

\item
If $\morcov_\indcov: \obcov_\indcov \to \obu$, $\indcovinset$,
and
$\morcov_{\indcovdob}:\obcov_\indcovdob \to \obcov_\indcov$,
$\indcovsec \in \setindcovsec_\indcov$,
are coverings for every $\indcov$,
then so is the composition
$ \morcov_\indcov \circ \morcov_\indcovdob : \obcov_\indcovdob \to \obu$,
$(\indcov,\indcovsec) \in \biguplus_\indcovinset
\setindcovsec_\indcov$.

\item
If $\obcov_\indcov \to \obu$ is a covering and
$\obv \to \obu$ is any morphism in
$\catt$, then $\obv \times_\obu \obcovindic \to \obv$ is a covering.
\end{enumerate}
\end{definition}

\ifthenelse{\boolean{book}}{
\begin{bookremark}
A Grothendieck topology does not necessarily have the property that
if $\obcov_\indcov \to \obu$ is a covering, then also the enrichment by
further open sets is again a covering.
\end{bookremark}}{}

\ifthenelse{\boolean{book}}{
\begin{bookexample}
The first example for a Grothendieck topology is given by the
category $\cattt$ of open subsets in a topological space $\spax$
such that the coverings are the families $\obcov_\indcov \subseteq
\obu$ of open subsets such that $\bigcup_{\indcovinset} \obcov_\indcov = \obu$.
Since in this context the product $\obu_1 \times_\obu \obu_2$ of two
open subsets is just their intersection $\obu_1 \cap \obu_2$, the
above properties are easily checked.
\end{bookexample}}

\ifthenelse{\boolean{book}}{
\begin{bookremark}
The definition of a Grothendieck topology above is sometimes called
a pretopology. The corresponding topology is then enriched by
declaring a family to be also a cover whenever there exists a
refinement which is a cover. E.g. starting with the Zariski
topology, one declares also a disjoint union of open sets which
cover a covering morphism (provided it belongs to the category)
\end{bookremark}}{}

\new{\doublelongrightarrow}{{\rightrightarrows}}

Each covering
$\obcovindic \to \obu$ in a Grothendieck topology induces
the diagram
(which is the first part of the \v{C}ech diagram)

\vspace{3.mm}
\begin{picture}(20,10)

\unitlength 1cm

\put(4.5,0){${\displaystyle \biguplus_{\indcov, \indcovsec}
\obcov_\indcov \times_\obu \obcov_\indcovsec} $}
\put(6.64,0.07){$\longto$}
\put(6.64,-0.07){$\longto$}
\put(7.45,0){${\displaystyle \biguplus_\indcov \obcov_\indcov} $}
\put(8.45,0.){$\longto$}
\put(9.25,0){$\obu \,$}
\end{picture}

\vspace{5mm}
\noindent
(the disjoint unions do not have to exist in $\catt$),
where the arrows
$\biguplus_\indcovdob \obcov_\indcov \times_\obu \obcov_\indcovsec
\to \biguplus_\indcov \obcov_\indcov $ are on each
open component either the first or the second projection. For a
presheaf $\presh$ we get the \v{C}ech diagram of the covering with
values in $\presh$,
$\presh(\obu) \to \prod_\indcov \presh ( \obcov_\indcov)\rightrightarrows
\prod_{\indcov,\indcovsec} \presh(\obcov_\indcov \times_\obu \obcov_\indcovsec)$.

\ifthenelse{\boolean{book}}{
\begin{bookexample}
A sheaf on a topological space is a presheaf $\sheaf $ such that for
every covering $ U= \bigcup_{\indcovinset} \obcov_\indcov$ the
sequence
$$ \sheaf (U) \to \prod \sheaf (\obcov_\indcov) \rightrightarrows \prod\sheaf (\obcov_\indcov \cap
U_j) $$ is exact. By which is meant that the left is the equalizer
of the maps on the right (and the left map is injective). With the
generalized definition of topology this also leads to a notion of
sheaves.
\end{bookexample}}{}

\begin{definition}
\label{sheafdefinition}
Let $\catt$ be a Grothendieck topology. Then
a presheaf $\sheaf  $ is called a \emph{sheaf},
if for all coverings
$\obcov_\indcov \to \obu$
in the topology the sequence

\vspace{3.mm}
\begin{picture}(20,10)

\unitlength 1cm

\put(3.05,0){$\displaystyle{
 \sheaf (\obu)}$}

\put(4.02, 0.){ $ \longto$}

\put(4.95,0){${\displaystyle
 \prod_{\iindex} \sheaf (\obcov_\indcov)}$}

 \put(6.63,0.07){$ \longto$}

 \put(6.63,-0.07){$ \longto$}

\put(7.45,0) {${\displaystyle \prod_{(\indcov,\indcovsec)
\in \setindcov \times \setindcov} \sheaf (\obcov_\indcov \times_\obu
\obcov_\indcovsec)}$}

\end{picture}

\vspace{.6cm}


\noindent
is exact (which also means that the first map is injective).
\end{definition}

\begin{remark}
Explicitly, a sheaf must satisfy the following two properties
.

\numiii
\begin{enumerate}

\item
If $\elem \in \sheaf (\obu)$ is such that
$\rest^\obu_{\obcovindic}(\elem)=0$ for all $\indcovinset$, then
$\elem=0$.

\item
If $\elem_\indcov \in \sheaf (\obcovindic)$ is such that
$(\rest^{\obcovindic}_{\obcovindic \times_\obu
\obcov_\indcovsec}(\project_1)) (\elem_\indcov)
= (\rest^{\obcov_\indcovsec}_{\obcovindic \times_\obu
\obcov_\indcovsec}(\project_2))(\elem_\indcovsec)$
for all $\indcov,\indcovsec$,
then there exists an $\elem \in \sheaf (\obu)$ such that
$\rest^\obu_{\obcovindic}(\elem)=\elem_\indcov$.
\end{enumerate}
A presheaf which satisfies only the first condition is called a
\emph{separated pre\-sheaf}. Note that the sheaf condition is in
this general setting, contrary to the case of a topological space,
also non-trivial for a single cover $\obv \to \obu$
(we will usually prefer the term cover for a covering consisting
of only one object), because the two projections
$\project_1, \project_2: \obv \times_\obu \obv \to \obv$ are different.
\end{remark}

\ifthenelse{\boolean{book}}{We make the category of sheafs
$\Sh(\catt)$ on $\catt$ into a full subcategory of $\PSh(\catt)$.
Note that, whereas the presheafs only depend on the underlying
category of the topology $\catt$, the notion of a sheaf does depend
on the notion of coverings given by the topology $\catt$.
Nevertheless it can be shown that $\Sh(\catt)$ is an abelian
category as well, though it is \emph{not} an abelian subcategory of
$\PSh(\catt)$ in general, as the cokernel in $\PSh(\catt)$ of a map
of sheaves is generally not a sheaf. To construct cokernels in
$\Sh(\catt)$ one needs to be able to sheafify.}{}

\subsection{Sheafification}
\label{sheafifysubsection}

\

\medskip
One needs to construct a \emph{sheafification} of a presheaf
$\sheaf $, that is, a sheaf $\sheaf ^\sheafify$ together with a presheaf
homomorphism
$\shef: \sheaf  \to \sheaf ^\sheafify$ which is
universal for homomorphisms from $\sheaf $ to sheaves.
That is, whenever the solid arrow to the sheaf $\shG$ is given,
$$ \xymatrix{  &  \shG  \\
\sheaf  \ar[ur] \ar[r]^\shef & \sheaf ^\sheafify \ar@{.>}[u]}$$ the
dotted arrow exists uniquely.
As we will use sheafification explicitly we recall its construction.
The sheafification is a two step process,
see \cite[Theorem II 2.11]{milne}, \cite[Expos\'{e} II.3]{SGA4} or
\cite[Section 2.3.7]{vistoligrothendieck}.
In the first step one defines a subpresheaf
$\sheaf _0 \subseteq \sheaf $ by
$$\sheaf _0(\obu)\!
:= \! \{\elem \in \sheaf (\obu) : \, \exists \mbox{ a covering } \obcovindic \to
\obu, \indcovinset, \mbox{ such that } \rest^\obu_{\obcovindic}(\elem) = 0 \,\,
\forall \indcov \} \,.$$
Then one sets
$\sheaf _1(\obu) := \sheaf (\obu)/\sheaf _0(\obu)$.
This presheaf $\sheaf _1$ is separated.

To describe the second step we pause to recall some categorial
notions, in particular colimits over a (quasi)filtered system.
Beside the sheafification procedure we need this to define the
inverse presheaf under a site morphism
(Section \ref{sitesubsection})
and to talk about (absolute) filters and stalks in a Grothendieck
topology
(Section \ref{filterstalksubsection}).
We are explicit about this point also since not everything
written down in the literature is correct
(see Example \ref{stalkquasifilternotfilter}).

\synchronize{\indec}{\indmo}
\synchronize{\inded}{\indmu}

\begin{definition}
\label{cofilteredcatdef}
Let $\indexcat$ ba a category. We call
$\indexcat$ \emph{quasicofiltered}, if the following property holds.

(i)
For every pair $\indmo \to \indr  $ and
$\indmu \to \indr$ in $ \indexcat$
there exists $\indl \in \indexcat$ and morphisms $\indl
\to \indmo  $ and $\indl \to \indmu $ such that the compositions
are the same; and for every $\indmo, \indmu \in \indexcat$ there
exists $\indl$ with arrows $\indl \to \indmo, \indl \to \indmu$
(this follows from the first condition if $\indexcat$ has a final
object).

We call it \emph{cofiltered} if it has also the property that

(ii) For every pair $\armr ,\armrsec    : \indm \rightrightarrows
\indr  $ in $\indexcat$ there exists $\arphi: \indl \to \indmo  $
such that $ \armr \circ \arphi = \armrsec \circ \varphi $ (so
$\arphi$ `unifies' the two morphisms).
\end{definition}


\begin{remark}
If the arrows are reversed the corresponding notions are called
\emph{quasifiltered} and \emph{filtered}. A
$\indexcat$-\emph{diagram} is a functor defined on a category
$\indexcat$. We will deal with covariant topological diagrams which
lead by evaluating a (pre)sheaf to contravariant diagrams in an
abelian category. The \emph{colimit} of such a contravariant diagram
$\indexcat \to \catab$ is the object $\colimab$ in $\catab$
characterized by the property that there exist homomorphisms
$\rest_\indm  : \obab_\indm \to \colimab$, which are compatible with
respect to the homomorphisms $\obab_\indr \to \obab_\indm$ indexed
by $\indm   \to \indr  $, and that for every object $\obabsec$ in
$\catab$ with such compatible homomorphisms there exists a unique
factorization $\colimab \to \obabsec$.

The existence of colimits and their description depend on properties
of the indexing category and of the category $\catab$. Colimits
exist for $\ring$-modules over a small category \cite[Proposition A
6.2]{eisenbud}. Over a small quasifiltered category the colimit of
sets is the disjoint union of all sets $\obset_\indm $ modulo the
equivalence relation which identifies two elements $\elmo \in
\obset_{\indmo }$ and $\elmu \in \obset_{\indmu}$ if they map to the
same element in some $\obset_\indl$ under
$\arpsi, \arpsi': \obset_{\indmo }, \obset_{\indmu} \to \obset_\indl$.
For a quasifiltered diagram of abelian groups (or $\ring$-modules) the
colimit in the category of groups is different from the colimit in
the category of sets (\cite[Proposition A 6.3]{eisenbud} is wrong,
see Example \ref{stalkquasifilternotfilter}). For a filtered diagram
however the two notions coincide, and the colimit in the category of
groups has the easy set theoretical description.
\end{remark}

\ifthenelse{\boolean{book}}{

Hence, if $\presh$ is a presheaf of $\ring$-modules and $\filt$ is a
quasifilter, then the colimit (the stalk) exists as an
$\ring$-module, every element is represented by an element $s \in
\Gamma(U_\indmo  , \presh)$, and two elements $s \in \Gamma(U_\indmo
, \presh)$ and $s' \in \Gamma(U_\indmu  , \presh)$ represent the
same element in the stalk if and only if there exist $\psi:\delta
\to \indmo   $ and $\psi':\delta \to \indmu  $ in $\indexcat$ such
that $\psi(s)=\psi'(s')$ in $\Gamma(U_\delta, \presh)$. In
particular, two elements $s,s' \in \Gamma(U_\indmo  ,\presh)$ define
the same element in the stalk if and only if there exists $\psi_1,
\psi_2 : \delta \to \indmo  $ such that $\psi_1(s) = \psi_2(s')$ in
$\Gamma(U_\delta,\presh)$. If moreover $F$ is a filter, then these
two morphisms $\psi_1$ and $\psi_2$ can be united to one morphism
$\varphi: k \to \delta $ such that $(\varphi \circ \psi_1)(s)=
(\varphi \circ \psi_2)(s')$. So in the case of a filter two elements
in $\Gamma(U_\delta,\presh)$ with the same image in the stalk are
identified by one homomorphism.

\begin{bookremark}
Another difference between quasifilters and filters is with respect
to a sheaf of rings. In both cases the colimit exists, however for
quasifilters it needs not be the same as the corresponding colimit
of sets or $\ring$-modules. For if only one object, say a ring
$\alg$ is given, and two morphisms $\psi,\id: A \to A$, then the
module-colimit is the $\ring$-module $\alg$ modulo the subgroup (or
$\ring$-submodule) $(a-\psi(a), a \in A)$, but the colimit as a ring
is the residue class ring of $\alg$ modulo the ideal generated by
these elements, see \cite[Proposition A6.7]{eisenbud}. For a filter
however the set-theoretic colimit is also the ring-theoretic
colimit, see \cite[Appendix A, Corollary 2]{milne}. The point is
that for $\psi_1,\psi_2: B \to B'$ (rings), the two images
$\psi_1(b)$ and $\psi_2(b)$ come together again under $\varphi_1,
\varphi_2: B' \to C$, so they define the same element in the stalk;
however an element in the ideal generated by $\psi_1(b)-\psi_2(b)$
does not necessarily map to $0$ somewhere and need not be $0$ in the
colimit of $\ring$-modules, but it must be $0$ in the colimit of
rings.
\end{bookremark}}{}

\ren{\presh}{{\shP}}

We return now to the construction of the sheafification. In the
second step one sets
$\presh^\sheafplus (\obu)= \colim
\chch(\namecov, \presh)$
(applied then to the separated presheaf
$\presh= \sheaf _1$),
where $\chch( \namecov, \presh)$ consists of the elements
$\elem_\indcov \in \Gamma(\obcov_\indcov ,\presh)$
which are compatible in the sense that
$\rest^{\obcov_\indcov}_{\obcov_\indcov \times_\obu
\obcov_\indcovsec}(\elem_\indcov)
 = \rest^{\obcov_\indcovsec}_{\obcov_\indcov \times_\obu
\obcov_\indcovsec} (\elem_\indcovsec)$ and where the colimit is
taken over all coverings of $\obu$ ordered by refinements. So the
objects in the indexing category are coverings
$\namecov =(\obcov_\indcov \to \obu, \indcovinset)$,
and the morphisms are refinements of coverings.
A covering
$\obcovsec_\indcovsec \to \obu$, $\indcovinsetsec$,
is a \emph{refinement} of
$\obcov_\indcov \to \obu$, $\indcovinset$, if there exists a mapping
$\mapind: \setindcovsec \to \setindcov$ and morphisms
(in $\catopenspax $)
$\morfine_\indcovsec:\obcovsec_\indcovsec \to \obcov_{\mapind
(\indcovsec)}$
(\cite[III \S2]{milne} or \cite[Definition 2.43]{vistoligrothendieck}).

A refinement $\namecovfine \to \namecov$ induces a homomorphism
$\chch( \namecov, \presh) \to \chch (\namecovfine, \presh)$ by
sending $(\elem_\indcov)$, $\indcovinset$, to
$(\morfine_\indcovfine(\elem_{\mapind (\indcovfine)}))$,
$\indcovfineinset$,
$\morfine_\indcovfine=\rest(\morfine_\indcovfine):
\Gamma(\obcov_{\mapind(\indcovfine)} , \presh)
\to \Gamma(\obcovfine_\indcovfine, \presh)$.
Because of the compatibility
property this homomorphism is independent of $\mapind$ and of
$\morfine_\indcovfine$ \cite[Corollaire V.2.3.5]{SGA4}.
Therefore this is a filtered diagram in the
category of abelian groups (though on the topological level it is
not even quasicofiltered) and its colimit equals the colimit in the
category of sets. In particular, the homomorphisms $\presh (\obu)
\to \presh ^\sheafplus (\obu)$ are injective for a separated
presheaf.

Combining the two steps gives us the sheafification $ \sheaf ^\sheafify
= ( \sheaf _1)^\sheafplus$ of a given presheaf $\sheaf $. It is also true
that $\sheaf ^\sheafify =( \sheaf ^\sheafplus)^\sheafplus$, but we
will use the first description.

\ifthenelse{\boolean{book}}{ If we want that a refinement has the
proprety that $\obcovsec_\indcovsec \to \obcov_\indcov$, $\mapind(
\indcovsec)= \indcov$, is a covering, then this is a quasicofiltered
category on the topological level, nu not a cofiltered one. For the
abelian diagram this is however automatically true, since there
exists at most one homomorphism between the coverings.

The construction of $\sheaf ^+$ is as follows:
For $\obu \in \catt$ to
define $\sheaf ^+(\obu)$ let $J_\obu$ be the category of coverings
$\obcovindic \to \obu$ where the maps are \emph{refinements} of
coverings.
A refinement of a cover $\obcovindic \to \obu$, $\indcovinset$, is a
cover $\obv_j \to \obu$, $j \in J$, and maps
$\mapind: J \to \setindcov$ and morphisms (in $\catopenspax    $)
$\obv _j \to \obu_{\mapind (j)}$ (over $\obu$).
Clearly, any presheaf $\sheaf $ on $\catt$ induces
a functor $\sheaf _\obu: J_\obu \to \Ab$ by defining
$$\sheaf _\obu(\obcovindic \to \obu)
:= \ker(\prod \sheaf (\obcovindic) \rightrightarrows \prod
\sheaf (\obcovindic \times_\obu \obu_j))
= \chch^0(\obcovindic \to \obu,\sheaf ).$$
Now just set
$$ \sheaf ^+(\obu) := \colim_{J_\obu} \sheaf _\obu$$
which can be thought of the $0$-th \v{C}ech cohomology
$\chch^0(\obu;\sheaf )$ of $\sheaf $ on $\obu$, after all it is the
colimit over all covers of $\obu$.

To take the colimit over all of $J_\obu$ is somewhat overkill since
for two different refinements $f,g : ( \obcov_\indcov \to \obu ) \to
(V_i
\to \obu)$ the resulting maps on the \v{C}ech complex of the
coverings are homotopic, thus on cohomology, in particular on
$H^0(\usc,\sheaf )$, they are equal. Thus $J_\obu$ can be replaced by
the partially ordered set $\bar{J}_\obu$ given $( \obv_j \to \obu)
\leq ( \obcov_\indcov \to \obu)$ if $ \obv_j \to \obu$ is a refinement of $(
\obcov_\indcov \to \obu)$, the specific maps giving the refinement are not
important.

Using that each two elements of $\bar{J}_\obu$ have a minimum, it
follows that $+:\PSh \to \PSh$ is a left exact functor. Furthermore,
for every covering $\obcov_\indcov \to \obu$ the natural maps
$\sheaf (\obu)
\to \prod \sheaf (\obcov_\indcov)$ induce a map $\sheaf  \to \sheaf ^+$. If $\sheaf $
was already a sheaf this map is an isomorphism. This shows that any
map $\sheaf  \to \shG$ to a sheaf $\shG$ must factor through
$\sheaf ^+$.
Now one only has to show that this construction has the claimed
properties of above. We note some properties of the sheafification
process.}{}

\ifthenelse{\boolean{book}}{
\begin{bookproposition}
\label{sheafificationproperties} Let $\catt$ denote a Grothendieck
topology. Then the following holds.

\numiii
\begin{enumerate}

\item The inclusion $i_\catt: \Sh(\catt) \to \PSh( \catt)$ is left exact.

\item The functor $\#_\catt: \PSh(\catt) \to \Sh(\catt)$ is exact and left adjoint to
   $i_\catt$.

\item The category $\Sh(\catt)$ has enough injectives. An injective sheaf is
injective as a presheaf.

\item
The section functor $\Gamma(U, -) : \Sh(\catt) \to \Ab$ defined by
$\Gamma(U,\sheaf )=\sheaf (U)$ is left exact.
\end{enumerate}
\end{bookproposition}
\begin{proof}
See \cite[Theorem 2.15]{milne}.
\end{proof}}{}


\medskip
The kernel presheaf of a sheaf morphism $\shG \to \sheaf $ is a sheaf,
but the presheaf image is not. The \emph{image sheaf} of a sheaf
morphism is by definition the sheafification of the presheaf image.
A complex of sheaves $\lmrsheafstack$ is called \emph{exact} if the
identity of sheaves
$\ker (\cob)= \im (\coa)$ holds;
this means that for every $\obu \in \catopenspax    $ and every
element
$\eles \in\Gamma(\obu, \shem)$ mapping to $0$ in
$\Gamma(\obu,\sher)$ there exists a covering $\obcov_\indcov \to
\obu$, $\indcovinset$, and sections
$\elel_\indcov \in \Gamma(\obcov_\indcov,\shel)$,
$\indcovinset$,
(without compatibility condition)
mapping to the restrictions
$\rest^\obu_{\obcov_\indcov}(\eles) \in \Gamma(\obcov_\indcov,\shem )$.
With these definitions the category of sheaves
on $\catt$ is an abelian category.

A complex of presheaves is called \emph{exact} if for every $\obu$ in
$\catopen$ the evaluated complex of abelian groups is exact. This
notion is often too strong, and it is independent of the
Grothendieck topology. We will use also the following definition.

\begin{definition}
\label{catexactdef}
Let $\catt$ be a Grothendieck topology and let $\opentop \in \catt$
be an object. Let $\lmrsheafstack$ be a complex of presheaves of abelian
groups. Then we say that the complex is
$\catopen$-\emph{exact} on $\opentop$ or $\opentop$-\emph{exact} if for
every $\elem \in \Gamma(\opentop, \shem)$ mapping to $0$ in
$\Gamma(\opentop,\sher)$ there exists a covering
$\obcovindic \to \opentop$, $\indcovinset$,
and sections
$\elel_\indcov \in \Gamma(\obcovindic, \shel)$ mapping to the
restrictions $\rest^\opentop_{\obcovindic}(\elem)$ in $\Gamma(\obcovindic, \shem)$.
We say that the complex is $\catt$-\emph{exact} if it is
$\opentop$-exact
for every $\opentop \in \catt$.
\end{definition}

\begin{proposition}
\label{exactexact} Let $\lmrsheafstack$ be a complex of presheaves
on a Grothendieck topology $\catt$. Then the following holds.

\numiii

\begin{enumerate}

\item
If the complex is exact as a complex of presheaves, then it is
$\catt$-exact.

\item
If the sheafified complex is $\obu$-exact, then the complex of
presheaves is $\obu$-exact.

\item
The complex is $\catt$-exact if and only if the associated complex
of sheaves is $\catt$-exact.

\item
A complex of sheaves is exact if and only if it is $\catt$-exact.
\end{enumerate}
\end{proposition}
\begin{proof}
(i) is clear. (ii). Suppose that the sheafified complex is
$\obu$-exact, and let $\elem \in \Gamma(\obu, \shem)$ mapping to $0$. Then
there exists a covering $\obcov_\indcov \to \obu$, $\indcovinset$,
and elements
$\elel_\indcov \in \Gamma(\obcov_\indcov, \shel^\sheafify)$ mapping to
$\elem_\indcov= \rest^\obu_{\obcov_\indcov}(\shef (\elem))
\in \Gamma(\obcov_\indcov, \shem^\sheafify)$.
The elements $\elel_\indcov$ are represented by sections of $\shel$ in a
finer covering and also the identities hold already in $\shem$ in a
finer covering.

(iii).
One direction follows from (ii).
For the other direction let
$\elem \in \Gamma(\obu,\shem^\sheafify)$ mapping to $0$ in
$\Gamma(\obu,\sher^\sheafify)$.
Then there exists a covering
$\obcov_\indcov \to \obu$, $\indcovinset$, and
$\elem_\indcov \in \Gamma(\obcov_\indcov ,\shem)$ representing $\elem$ and
coverings $\obcovsec_\indcovsec \to \obcov_\indcov$,
$\indcovsecinset_\indcov$,
such that the restrictions of
$\cob(\elem_\indcov)$ in
$\Gamma(\obcovsec_\indcovsec, \sher)$ are $0$.
Since the statement is local we may assume that
$\obcov_\indcov = \obcovsec_\indcovsec=\obu$.
But then again the assumptions on the presheaves yield immediately a
covering
$\obcov_\indcov \to \obu$ and sections
$\elel_\indcov \in \Gamma(\obcov_\indcov, \shel) \to \Gamma(\obcov_\indcov,
\shel^\sheafify)$ which map to the restrictions of $\elem$.
(iv) follows from the definitions.
%
\end{proof}

\ifthenelse{\boolean{book}}{
\begin{bookremark}
In the Zariski topology on $\spax$, both an exact complex of sheaves
as well as a globally exact complex of sheaves are $\spax$-exact.
Examples of exact, but not globally exact, and examples of globally
exact, but not exact complexes show that the converse is not true.
\end{bookremark}}{}

\ren{\catind}{{I}}

\subsection{Site morphisms}
\label{sitesubsection}

\

\medskip
We will call a category $\catopen$ with a final object $\spax$, with
arbitrary products and endowed with a Grothendieck topology a
\emph{site}, and we will denote it by $(\spax,\catopenspax    )$ or
sometimes just by $\spax$. Later on $\spax$ will be a scheme, the
products will be products of schemes and $\catopenspax    $ will be
a subcategory of schemes over $\spax$.

If $\catopenspax    $ and $\catopen_\spay$ are two sites with final
objects $\spax$ and $\spay$, then a \emph{site morphism} $
\sitemortop :\catopenspax     \to \catopen_\spay$, for which we will
briefly write $\sitemor: \spay \to \spax$, is a functor
$\catopenspax \to \catopen_\spay$ with
$\sitemortop (\spax) =\spay$ and
which respects products and coverings, so that if
$\obcovindic \to \obu$, $\indcovinset$,
is a covering in $\catopenspax $, then
$\sitemortop (\obcovindic) \to \sitemortop (\obu)$, $\indcovinset$,
is a covering in $\catopen_\spay$.

\begin{construction}
\label{presheafinvers}
A site morphism $\sitemor : \spay \to \spax$ associates to a
presheaf $\presh$ of abelian groups on $\spax$ a presheaf
$\sitemorinv (\presh)$ on $\spay$
(called its \emph{inverse image})
in the following way
(cf. the construction of $u_!$ in \cite[Proposition I.5.1]{SGA4} in connection
with \cite[Proposition III.1.6]{SGA4}). For
$\obv\in\catopen_\spay$, let
$\catind^\sitemor_\obv$ be the category of all pairs $(\obu,\phili)$
with
$\obu \in \catopenspax $ and
$\phili \in \Mor_{\catopen_\spay}(\obv, \sitemortop(\obu))$.
The morphisms
$\Mor_{\catind^\sitemor_\obv}((\obu_1,\phili_1),(\obu_2,\phili_2))$
in this category are the maps
$\xili \in \Mor_{\catopenspax}(\obu_1,\obu_2)$
such that
$\sitemortop (\xili) \circ \phili_1 =\phili_2$.
This category is quasicofiltered, since for
$\phili_1: \obv \to \sitemortop (\obu_1)$
and
$\phili_2: \obv \to \sitemortop (\obu_2)$
mapping to
$\phili_0: \obv \to \sitemortop (\obu_0)$
via
$ \xili_1 :\obu_1 \to \obu_0$
and
$\xili_2 : \obu_2 \to \obu_0$
there exists
$\phili: \obv \to \sitemortop(\obu_1) \times_{\sitemortop(\obu_0)}
\sitemortop (\obu_2)
= \sitemortop (\obu_1 \times_{\obu_0} \obu_2)$
over $\phili_1$ and $\phili_2$.
We take the colimit
$$ \sitemorinv (\presh)(\obv)
 : = \colim_{(\obu,\phili) \in \catind^\sitemor_\obv} \Gamma(\obu,  \presh ) \, .$$
in the category of abelian groups.
This indeed defines a presheaf on
$\catopen_\spay $,
since a map
$\eptau \in \Mor_{\catopen_\spay }
(\obv,\obw)$
induces a functor
$\catind^\sitemor_\obw \to \catind^\sitemor_\obv$ by sending
$(\obu,\phili)$ to
$(\obu, \phili\circ \eptau)$
and thus induces a map of the colimits above.

The categories $\catind^\sitemor_\obv$ are in general not cofiltered, but they are under the following
\emph{cofiltered condition}. For every pair
$\arpsi_1,\arpsi_2 : \obu \rightrightarrows \oburight$ in
$\catopen_\spax$ and $\phili: \obv \to \sitemortop(\obu)$ in  $\catopen_\spay$
there exists $\morleft :\obuleft \to \obu$ which unifies $\arpsi_1$
and $\arpsi_2$ and
$\tilde{\phili}: \obv \to \sitemortop (\obuleft)$ such that
$ \sitemortop(\morleft)\circ \philileft =\phili$.
Under this condition $\sitemorinv(\presh)$ has the easy
set-theoretical description.
\end{construction}

The presheaf $\sitemorinv (\presh)$ constructed in
\ref{presheafinvers}
can be sheafified, but for us it will be more important to sheafify
with respect to a structure sheaf of rings on the sites (we need the
`image r\'{e}ciproque au sens modules' instead of the `image
r\'{e}ciproque ensembliste', cf.
\cite[Expos\'{e}s IV.11-13]{SGA4}).

\begin{definition}
A site $\catopen$ is called a \emph{{\rm(}pre-{\rm)}ringed site} if
there exists a fixed {\rm(}pre{\rm)}sheaf
$\obu \mapsto \Gamma(\obu,\O) $, $\obu \in \catopen$, of commutative rings, called its
\emph{structure {\rm(}pre{\rm)}sheaf}.

A presheaf $\presh$ of abelian groups on such a
{\rm(}pre-{\rm)}ringed site is called an $\O$-\emph{premodule}, if
it has additionally an $\O$-multiplication such that $\Gamma(\obu,
\presh)$ is an $\Gamma(\obu, \O)$-module for every $\obu \in
\catopen$. If moreover $\presh$ is a sheaf, then we call it an
$\O$-module.
\end{definition}

The sheafification of the structure presheaf $\O$ on a preringed
site makes it into a ringed site. The sheafification of an
$\O$-premodule $\sheaf $ is an $\O$-module and an
$\O^\sheafify$-module, since the multiplication $\O \times \sheaf  \to
\sheaf $ sheafifies to $\O^\sheafify \times \sheaf ^\sheafify \to
\sheaf ^\sheafify$.

\begin{definition}
If $(\catopenspax, \O_\spax)$ and
$(\catopen_\spay, \O_\spay )$
are two {\rm(}pre{\rm)}-ringed sites, then a site morphism $
\sitemortop :\catopenspax \to \catopen_\spay $, for which we
will write again
$\sitemor: \spay \to\spax$, is called a
\emph{morphism of {\rm(}pre{\rm)}ringed sites}, if there is given a
family of ring homomorphisms $\Gamma(\obu, \O_\spax) \to \Gamma(
\sitemortop (\obu), \O_\spay )$ compatible with the restrictions in
the topologies.
\end{definition}

\begin{remark}
A morphism of (pre-)ringed sites $\sitemor:\spay \to\spax$ induces
on $ \catopen_\spay $ a morphism of presheaves of rings
$\sitemorinv(\O_\spax) \to \O_\spay $ by
$$ \Gamma(\obv,  \sitemorinv(\O_\spax))
= \colim_{(\obu,\phili) \in \catind^\sitemor_\obv} \Gamma(\obu, \O_\spax) \lto
\colim_{(\obu,\phili) \in \catind^\sitemor_\obv} \Gamma(\sitemortop
(\obu), \O_\spay ) \stackrel{}{\lto} \Gamma(\obv, \O_\spay ) \, .
$$
(Note that $\sitemorinv(\O_\spax)$ is here the inverse sheaf in the
category of rings, so the colimit we have to take in its
construction is the colimit inside the category of rings.)
On the other hand, if such a presheaf morphism is given, then the
composition
$\Gamma(\obu, \O_\spax)
\to \Gamma(\sitemortop(\obu),\sitemorinv (\O_\spax))
\to \Gamma(\sitemortop(\obu), \O_\spay )$
yields the ring homomorphisms between the sites.
\end{remark}

In the case of a morphism between (pre-)ringed sites
$\sitemor:\spay \to \spax$ the inverse image
$\sitemorinv(\sheaf )$ of an
$\O_\spax$-(pre)module $\sheaf $ is in general not an $\O_\spay
$-(pre)module.

\begin{definition}
\label{sitemorpredef}
Let $\sitemor: \spay \to\spax$ be a morphism between (pre-)ringed
sites.
We define $\sitemorpre (\sheaf )$ to be the
$\O_\spay $-premodule given by
$$ \obv \longmapsto \Gamma(\obv, \sitemorinv (\sheaf ))
\tensor_{\Gamma(\obv, \sitemorinv(\O_\spax))}  \Gamma(\obv,
\O_\spay ) \, .$$ We denote its sheafification by $\sitemormodule
(\sheaf )$ (cf. \cite[ IV.13.2.2]{SGA4}).
\end{definition}

$\sitemorpre (\sheaf )$ is an $\O_\spay $-premodule and
$\sitemor^*(\sheaf )$ is an $\O_\spay $-module.

\begin{lemma}
\label{rightexactpullback}
Let $\sitemor: \spay \to\spax$ be a morphism between ringed sites.
Then $\sitemormodule$ is a right exact functor for
$\O_\spax$-modules.
\end{lemma}
\begin{proof}
See \cite[Corollaire IV.13.6]{SGA4}.
\end{proof}

\ifthenelse{\boolean{book}}{

Der erste Teil gilt nur unter der cofiltered Bedingung.

\begin{lemma}
\label{rightexactpullback}
Let $\sitemor: \spay \to\spax$ be a site
morphism and let $\lmrsheafstack $ be a $\catopenspax    $-exact
sequence of abelian presheaves on $\catopenspax    $. Then
$\sitemorinv (\shel) \to \sitemorinv (\shem) \to \sitemorinv
(\shel)$ is a $\catopen_\spay $-exact sequence of abelian
presheaves. If
$\sitemor: \spay \to\spax$ is a morphism between
ringed sites, then $\sitemormodule$ is a right exact functor for
$\O_\spax$-modules {\rm(}cf. \cite[Corollaire IV.13.6]{SGA4}{\rm)}.
\end{lemma}
\begin{proof}
Let
$\elem \in \Gamma(\obv, \sitemorinv (\shem))$ mapping to $0$ in
$\Gamma(\obv, \sitemorinv (\sher))$
and let $\elem$ be represented by $\elem \in \Gamma(\obu, \shem)$,
where
$\obv \to \sitemortop(\obu)$
belongs to $\catopen_\spay $.
Since $\cob (\elem)=0$ there
exists another morphism
$\obv \to \sitemortop (\obu')$ such that
$\cob(\elem)=0 $ in
$\Gamma(\obu \timesx \obu', \sher)$. We may
write $\obu$ instead of $\obu \timesx \obu'$. By assumption there
exists a covering
$\obcovindic \to \obu$, $\indcovinset$, in
$\catopenspax $ and
$\elel_\indcov \in \Gamma(\obcovindic, \shel)$
mapping to
$\rest^\obu_{\obcovindic}(\elem)$.
Hence
$\sitemortop (\obcovindic) \to\sitemortop (\obu)$
is a covering in $\catopen_\spay $ and therefore
$\obv_\indcov
=\obv \times_{\sitemortop (\obu)} \sitemortop(\obcovindic) \to\obv$
is a covering.
The $\elel_\indcov$ give elements
$\elel_\indcov \in
\Gamma(\obv \times_{\sitemortop (\obu)}
\sitemortop (\obcovindic),\sitemorinv (\shel))$
mapping to $\rest^\obv_{\obv_\indcov} (\elem)$.

Let now $\olmrosheafstack $ be an exact sequence of
$\O_\spax$-modules which gives a $\catopen_\spay $-exact sequence
$0\to \sitemorinv (\shel)
\to \sitemorinv (\shem)
\to \sitemorinv (\sher) \to 0$ of presheaves,
which are also $\sitemorinv(\O_\spax)$-premodules.

Now in general, if $\spaz$ is a preringed site and if
$0 \to \shH_1\to\shH_2 \to \shH_3 \to 0$
is a
$\catopen_\spaz$-exact sequence of
$\O_\spaz$-premodules, then
$\shH_1 \tensor_{\O_\spaz} \shP
\to\shH_2\tensor_{\O_\spaz} \shP
\to \shH_3\tensor_{\O_\spaz} \shP \to 0$ is
$\catopen_\spaz$-exact for every $\O_\spaz$-premodule $\shP$. For this we
may replace by Proposition \ref{exactexact}(iv) the presheaves by
their sheafification. The right exactness of tensorization with a
sheaf (in the sheaf category) follows from the universal properties
of tensor products as in the case of modules
(see \cite[Satz 82.9]{schejastorch2},
cf. \cite[Propositions IV.12.7 and IV.12.10 and
Corollaire IV.12.12]{SGA4}
for general properties of the tensor product for sheaves).
Taking
$\O_\spaz = \sitemorinv (\O_\spax)$ and
$\shP= \O_\spay $
gives the result.
\end{proof}}{}



\subsection{Extension and contraction under a site morphism}
\label{extensioncontractionsubsection}

\

\medskip
A ring homomorphism $\homring: \ring \to \rings$ defines for
an
ideal $\idealsubring$ the \emph{extended ideal} $\ideal \rings$ in
$\rings$ and then the \emph{contraction}
$\homring^{-1}(\ideal \rings)$, which
contains $\ideal$.
Basically the same construction applies to an
arbitrary site morphism.
Later on we will specialize to the situation
where the site morphism
$(\specring)_\topo \to \specring$ is given by a Grothendieck
topology on an affine scheme $\specring$.
This gives then a closure operation for ideals in $\ring$.

\begin{definition}
\label{extendedsubmoduledef}
Let $\spax$ and $\spay$ stand for two
ringed sites $(\catopenspax,\O_\spax)$ and
$(\catopen_\spay ,\O_\spay)$ and let
$ \sitemor:\spay \to \spax$ denote a morphism of ringed
sites. Let $\shG \subseteq \sheaf $ denote
$\O_{\spax}$-modules. This induces an $\O_{\spay}$-module
homomorphisms
$\sitemormodule(\shG) \to \sitemormodule(\sheaf )$.
The image sheaf of $\sitemormodule(\shG)$ inside
$\sitemormodule(\sheaf )$ is called the \emph{extended submodule} of
$\shG$ and is denoted by
$\shG^{\ext}$. In particular, for an ideal sheaf
$\idealsheaf \subseteq \O_\spax$ we call the image sheaf of
$\sitemor^*(\idealsheaf)$ inside
$\O_\spay = \sitemormodule(\O_\spax)$ the \emph{extended ideal sheaf}, denoted
by $\idealsheaf^{\ext}$.
\end{definition}

Recall that a site morphism
$\sitemor: \spay \to \spax$ defines a
functor
$\shH \mapsto \sitemorpf (\shH)$ for presheaves by setting
$\Gamma(\obu,\sitemorpf (\shH)):= \Gamma(\sitemortop(\obu), \shH)$
for $\obu \in \catopenspax    $.
This \emph{push forward} functor
sends sheaves to sheaves and is left exact. It also sends $\O_\spay
$-modules to $\O_\spax$-modules. In particular, for an
$\O_\spay $-submodule $\shH \subseteq \shH'$
also
$\sitemorpf (\shH) \subseteq\sitemorpf(\shH')$
is an $\O_\spax$-submodule. For an
$\O_\spax$-module $\sheaf $ we have a natural sheaf homomorphism
$\shef: \sheaf  \to \sitemorpf (\sitemorsheaf (\sheaf ))$ given by
$$\shef: \Gamma(\obu, \sheaf ) \to \Gamma(\obu, \sitemorpf (\sitemormodule (\sheaf )))=
\Gamma( \sitemortop (\obu), \sitemormodule (\sheaf )) \, .$$


\begin{definition}
Let $\sitemor:\spay \to \spax$ be a site morphism between two ringed
sites, let $\sheaf $ be an $\O_\spax$-module. An
$\O_{\spay}$-submodule
$\shH \subseteq \sitemormodule(\sheaf )$ yields an $\O_\spax$-submodule
$\sitemorpf (\shH) \subseteq  \sitemorpf (\sitemormodule (\sheaf ))$.
We call the $\O_\spax$-submodule of $\sheaf $ given as the preimage
sheaf of $\sitemorpf (\shH)$ under the natural sheaf homomorphism
$\shef: \sheaf  \to \sitemorpf (\sitemormodule (\sheaf ))$ the
\emph{contracted module}, denoted by $\shH^\contr$.
Explicitly, we have
$\Gamma(\obu, \shH^\contr)= \{ \elem \in \Gamma(\obu, \sheaf ):\,
\shef(\elem) \in \Gamma(\obu, \sitemorpf (\shH))
= \Gamma(\sitemortop(\obu), \shH) \}$.
\end{definition}

\begin{definition}
Let $\sitemor: \spay \to \spax$ be a site morphism between two
ringed sites, let $\shG \subseteq \sheaf $ be $\O_\spax$-modules. We
call the contraction of the extension the \emph{closure submodule}
of $\shG$ inside $\sheaf $ via $\sitemor$, and denote it by
$\shG^{\sitemor -\sitemorc} =\shG^\sitemorc= (\shG ^{\ext})^\contr$.
\end{definition}

\begin{proposition}
\label{sitemorcprop}
Let $\sitemor: \spay \to \spax$ be a site
morphism between two ringed sites, let $\shG \subseteq \sheaf $ be
$\O_\spax$-modules. Then the following holds.

\numiii
\begin{enumerate}

\item
The closure submodule $\shG^\sitemorc$ is an $\O_\spax$-module and
in particular a sheaf.

\item
We have $\shel \subseteq \shel^\sitemorc =
(\shel^\sitemorc)^\sitemorc$.

\item
Let $\shel' \subseteq \shem '$ be another pair of $\O_\spax$-modules
and let
$\shehom: \shem \to \shem '$
be an $\O_\spax$-module homomorphism with
$\shehom (\shel) \subseteq \shel '$. Then
$\shehom(\shel^\sitemorc)\subseteq (\shel')^\sitemorc$.

\item
If $\shel \subseteq \shell \subseteq \shem$ is another
$\O_\spax$-submodule, then
$\shel ^\sitemorc \subseteq \shell^\sitemorc$.

\item
Consider the short exact sequence $\olmrosheafstack $ of
$\O_\spax$-modules, where
$\sher = \shem/\shel$. Then
$\shG^\sitemorc = \cob^{-1} (0^\sitemorc)$.
\end{enumerate}
\end{proposition}
\begin{proof}
(i)
is clear from the corresponding properties of extension and
contraction.

(ii).
The inclusion $\shG \subseteq \shG^\sitemorc$ and so also
$\shG^\sitemorc \subseteq (\shG^\sitemorc)^\sitemorc$ is clear.
For the other direction note that in general
$(\shH^\contr)^{\ext} \subseteq \shH$ for a submodule
$\shH \subseteq \sitemor^*(\sheaf )$.
This is because
$\sitemorsheaf (\shH^\contr) \to \sitemorsheaf (\sheaf )$
factors through $\shH$, since we have
$\shH^\contr \to \sitemorpf (\shH)$ by definition and so
$\sitemorsheaf (\shH^\contr) \to \sitemorsheaf \sitemorpf (\shH)$,
and by adjunction we get
$\sitemorsheaf \sitemorpf (\shH) \to \shH \to \sitemorsheaf (\sheaf )$.

(iii)
follows from applying $\sitemorsheaf$ to the commutative
diagram
$$ \xymatrix{ \shel \ar[r] \ar[d] &  \shem \ar[d]^\shehom \\
\shel' \ar[r] & \shem' \, .  }$$
(iv) is a special case of (iii).

(v).
The inclusion $\subseteq $ follows from (iii). For the other
inclusion note that we have by Lemma \ref{rightexactpullback} an
exact sequence of sheaves $ \lmrositemorsheafstack $ in
$\catopen_\spay $, where $\shG^{\ext} = \im (\coa)= \ker (\cob) $
holds. If $\elem \in \Gamma(\obu, \shem)$ is such that
$\cob (\elem)\in 0^\sitemorc$
(inside $\sher = \shem/\shel$), then
$\shef (\cob(\elem)) = 0$ in
$\Gamma( \sitemortop (\obu), \sitemorsheaf (\sher))$,
since $0^{\ext}=0$. But then
$\shef(\elem)\in \Gamma(\sitemortop(\obu), \sitemorsheaf(\shem))$ maps to $0$ in
$\Gamma(\sitemortop (\obu), \sitemorsheaf (\sher))$ and belongs to
$\Gamma(\sitemortop (\obu), \shel^{\ext})$.
\end{proof}

\ifthenelse{\boolean{book}}{

\begin{remark}
Property (ii) in Proposition \ref{sitemorcprop} says that the
closure submodule gives indeed a closure operation. The property
expressed in (iii) implies for $\shG'= \shehom (\shG)$ that this
closure is \emph{module-persistent}. Property (iv) is called
\emph{order preserving} or \emph{monotony}. Property (v) is called
\emph{independence of presentation}. This independence allows us to
restrict to the $0$-submodule, where the sheafification is $0$ as
well and an element $\elem$ lies in the closure module of $0$ if and
only if $\shef(\elem) =0$. This is the case if and only if there
exists a covering such that $\elem$ is mapped to $0$ on the presheaf
level. On the other hand, if $\O_\spax^\numn \to \shE \to 0$ is a
presentation, then the computation of $0^\sitemorc$ in
$\shE$ can be reduced to the computation of a closure submodule
inside a free module.
%
\end{remark}}{}


We describe the containment in the closure submodule more
explicitly.

\renewcommand{\covmor}{\tau}

\begin{lemma}
\label{explicitsite}
Let $\sitemor: \spay \to \spax$ be a morphism
between two {\rm(}pre-{\rm)}ringed sites and suppose that it
fulfills the cofiltered condition of \ref{presheafinvers}.
Let
$\shG \subseteq \sheaf$ be
$\O_\spax$-{\rm(}pre-{\rm)} modules.
Let
$\elem \in\Gamma(\obu,\sheaf )$, $\obu \in \catopenspax $.
Set
$\obv=\sitemortop(\obu)$.
Then
$\elem \in \Gamma(\obu,\shG^\sitemorc)$
if and only if there exists a covering
$\obv_\indcov \stackrel{\covmor_\indcov}{\to} \obv$, $\indcovinset$,
and morphisms
$\rho_\indcov: \obv_\indcov \to\sitemortop(\obcov_\indcov) $
in $\catopen_\spay $ and
$\sigma_\indcov: \obcov_\indcov \to \obu$ in $\catopenspax $
such that
$\covmor_\indcov = \sitemortop (\sigma_\indcov) \circ \rho_\indcov $
and sections
$\elel_\indcov \in \Gamma(\obcov_\indcov, \shG)
\tensor_{\Gamma(\obcov_\indcov, \O_\spax)} \Gamma(\obv_\indcov,\O_\spay )$
such that
$\rest^\obu_{\obcov_\indcov}(\elem) \tensor 1 = \inc(\elel_\indcov)$
in
$\Gamma(\obcov_\indcov, \sheaf )
\tensor_{\Gamma(\obcov_\indcov,\O_\spax)} \Gamma(\obv_\indcov,\O_\spay )$,
where $\inc$ denotes the inclusion
$\inc: \shG \to \sheaf $.
\end{lemma}
\begin{proof}
The element
$\elem \in \Gamma(\obu, \sheaf )$ defines the element
$\shef(\elem) \in \Gamma(\obv, \sitemormodule(\sheaf ))$
as the image of $\elem \tensor 1$ under
$$\Gamma(\obu, \sheaf ) \tensor_{\Gamma(\obu,\O_\spax)} \Gamma(\obv, \O_\spay )
\lto  \Gamma(\obv, \sitemorinv(\sheaf ))
\tensor_{\Gamma(\obv, \sitemorinv (\O_\spax))} \Gamma(\obv,\O_\spay )
\lto \Gamma(\obv, \sitemormodule(\sheaf )) \, .$$
This element is in the image sheaf of
$\sitemormodule(\shG) \to \sitemormodule(\sheaf )$
if and only if there exists a covering
$\obv_\indcov \to \obv$, $\indcovinset$,
and sections
$\eler_\indcov \in \Gamma(\obv_\indcov, \sitemormodule (\shG))$
mapping to
$\rest^\obv_{\obv_\indcov}(\shef(\elem))$.
Each $\eler_\indcov$ is represented by a covering
$\obw_\indcovsec \to \obv_\indcov$,
$\indcovsec \in \setindcovsec_\indcov$,
and elements
$\elere_\indcovsec \in \Gamma(\obw_\indcovsec ,\sitemorpre (\shG))
= \Gamma(\obw_\indcovsec, \sitemorinv (\shG))
\tensor_{\Gamma(\obw_\indcovsec, \sitemorinv (\O_\spax))}
\Gamma(\obw_\indcovsec, \O_\spay ) $
(fulfilling some compatible condition, which we do not need).
We denote the composed covering of $\obv$ again by $\obv_\indcov$.
By the definition of inverse presheaf,
the sections
$\eler_\indcov \in \Gamma(\obv_\indcov, \sitemorinv (\shG))
\tensor_{ \Gamma(\obv_\indcov, \sitemorinv (\O_\spax))}
\Gamma(\obv_\indcov, \O_\spay ) $
are represented by
$\elel_\indcov \in \Gamma(\obcov_\indcov, \shG)
\tensor_{\Gamma(\obcov_\indcov,\O_\spax)}
\Gamma( \obv_\indcov, \O_\spay  )$,
where
$\rho_\indcov : \obv_\indcov \to \sitemortop (\obcov_\indcov)$
are morphisms in
$\catopen_\spay $ and
$\obcov_\indcov \in \catopenspax $.
These elements define sections
$\inc (\elel_\indcov) \in \Gamma(\obcov_\indcov, \sheaf )
\tensor_{\Gamma(\obcov_\indcov, \O_\spax)}\Gamma( \obv_\indcov, \O_\spay )$,
which map by assumption to
$\rest^\obv_{\obv_\indcov}(\shef(\elem))
\in \Gamma(\obv_\indcov, \sitemorsheaf (\sheaf ))$.

We may replace $\obcov_\indcov$ by
$\obcov_\indcov \timesx \obu$
(denoted again by $\obcov_\indcov$)
so that we have the projections to
$\obu$, the morphisms commute in the described way and
$\elem_\indcov =\rest^\obu_{\obcov_\indcov}(\elem)
\in \Gamma(\obcov_\indcov, \sheaf )$
makes sense. So we have elements
$\inc(\elel_\indcov)$ and $\elem_\indcov \tensor 1$ in
$\Gamma(\obcov_\indcov, \sheaf )
\tensor_{\Gamma(\obcov_\indcov, \O_\spax)}
\Gamma(\obv_\indcov, \O_\spay)$
and we know that they define the same element in the
sheafification
$\Gamma(\obv_\indcov,\sitemor^*(\sheaf ))$.
This means that there exist coverings
$\obw_\indcovsec \to \obv_\indcov$, $\indcovsec \in \setindcovsec_\indcov$, such that
their difference is $0$ in
$\Gamma(\obw_\indcovsec, \sitemorpre(\sheaf ))$.
This means by the cofiltered condition that there exist morphisms
$\obw_\indcovsec \to \sitemortop (\obx_\indcovsec)$,
$\obx_\indcovsec \in \catopenspax $,
such that the restriction of
$\elem_\indcov \tensor 1$ in
$\Gamma(\obcov_\indcov \timesx \obx_\indcovsec, \sheaf )
\tensor_{\Gamma(\obcov_\indcov \timesx \obx_\indcovsec, \O_\spax) }
\Gamma(\obw_\indcovsec, \O_\spay )$
equals the restriction of
$\inc(\elel_\indcov)$.
So the covering
$\obw_\indcovsec \to \obv$, $\indcovsec \in \bigcup_\indcov \setindcovsec_\indcov$,
has all the desired properties.
\end{proof}

\ifthenelse{\boolean{book}}{
\begin{example}
Without the cofiltered condition the preceeding Lemma does not hold.
Let $\spax$ be a point and let $\catopen_\spax$ be the category of
finite sets, where we allow only surjective morphisms, which are
also coverings. Let $\spay$ be also a point and let $\catopen_\spay$
consist only in $\spay$ with one identical morphism.
Declare $\sitemor:\spay \to \spax$ by $\sitemortop(\obu)=\spay$
for all $\obu \in \catopen_\spax$.
Make both sites to a ringed site with the constant sheaf $\ZZ$
(so that we have a morphism of ringed sites;
tensoring with the structure sheaf has no effect).
Consider on $\spax$ the $\ZZ$-module $\sheaf$ given by
$\obu \mapsto \ZZ^\obu$ (locally constant) and the submodule
$\shG=0$. Then $\sitemorinv(\shF)=0$, since
$\colim_{(\obu , \phili) \in \catind^\sitemor_\obv} \Gamma(\obu,\sheaf)
= \colim_{\obu \in \catopen_\spax} \Gamma(\obu,\sheaf) =0$. For this
we may assume that $\obu$ has at least three elements and we look at
two different morphisms to another finite set with one element less.
Hence $\elem \in \Gamma(\obu,\sheaf)$ belongs to the closure of $0$,
but the explicit condition in the Lemma is not fulfilled.
\end{example}
}{}

\begin{example}
Let $\spax$ denote an algebraic variety over the complex numbers
$\CC$ and let $\spax(\CC)$ denote the corresponding complex space
with the sheaf of holomorphic functions. Then we have a morphism of
ringed spaces
$\spax(\CC) \to \spax$, which is in particular a site morphism.
Extension and contraction yields the same ideal back, since
the algebraic stalks and the analytic stalks at a closed point have
the same completion, which is faithfully flat. So if
$\fuf \in \ideal^\an $,
then this holds in every analytic stalk, in their
completions and therefore already in the algebraic setting.

If we consider however $\spax(\CC)$ with the sheaf of continuous
functions, denoted by $\spax_\contin$, then extension and
contraction may give back a larger ideal $\ideal^\contin$, which we
call the \emph{continuous closure}.
For example, in $\CC[\varx,\vary]$ we have
$\varx^2\vary^2 \not \in (\varx^3, \vary^3)$, but $\varx^2\vary^2 \in (\varx^3, \vary^3)^\cont$.
We deal with this continuous closure in
\cite{brennercontinuous},
where we give a precise description of it for monomial ideals.
\end{example}

\subsection{Grothendieck topologies on a scheme}
\label{grothendieckschemesubsection}

\

\medskip
Let $\spax$ denote a scheme, and let $\catsch_\spax$ denote the
category of schemes over $\spax$, where the morphisms are supposed
to be compatible with the structure morphisms. A Grothendieck
topology on a subcategory
$\catopenspax \subseteq \catsch_\spax$
containing $\spax$ and closed under products is called a
Grothendieck topology on $\spax$. We will denote the scheme $\spax$
endowed with such a topology by $\spax_\topo$, where $\topo$ is a
name for the corresponding topology. Note that the products
appearing in Definition \ref{Grothendieckdef} are now products of
schemes. The site
$(\spax_\topo, \catopenspax)$
is immediately a preringed site
$(\spax_\topo, \catopenspax,\O^\pre)$ by setting
$\Gamma(\obu,\O^\pre) := \Gamma(\obu, \O_\obu)$
for $ \obu \in \catopenspax $, and sheafification in the topology
gives a ringed site
$(\spax_\topo, \catopenspax, \O_\topo)$.

A site morphism between schemes $\spay_\topo$ and $\spax_\topo$ is a
scheme morphism
$\sitemor: \spay \to \spax$ such that
$\obu \in \catopenspax $ implies that
$\sitemortop (\obu)= \spay \timesx \obu \in \catopen_\spay $
and such that a covering
$\obcovindic \to \obu$, $\indcovinset$,
yields a covering
$ \spay \timesx \obcovindic \to \spay \timesx \obu$, $\indcovinset$.
It is then also a morphism of (pre)ringed sites by
taking
$\Gamma(\obu,\O^\pre)
= \Gamma (\obu,\O_\obu)
\to \Gamma(\sitemortop (\obu), \O_{\sitemortop(\obu)} )
= \Gamma(\sitemortop (\obu), \O^\pre)$
and its sheafification. A scheme morphism may or may not induce a
site morphism. For our treatment we will also adopt the following
convention.

\begin{convention}
\label{basistopology}
In saying that a scheme $\spax$ has a Grothendieck topology
$\catopen_\spax$ we mean that there exists a basic (or standard)
topology $\shB_\spax$ on $\spax$ which is always inside the Zariski
topology
(it will be usually either the trivial topology (Example
\ref{trivialsite}) or the Zariski topology)
and which we will not mention in general, of which the given
Grothendieck topology is a refinement.
In particular we have a site morphism
$\mortop: \spax_\topo \to \spax$.
We assume also that
if $\struto: \obv \to \spax$ belongs to
$\catopen_\spax$ and factors through an open subscheme
$\obv \to\openzar \subseteq \spax$,
then also $\obv \to \openzar \in \catopen_\spax$.
We denote the $\O_\topo$-sheafification of an $\O_X$-module $\modul$ by
$\modul_\topo =\sitemorsheaf (\modul) $.
\end{convention}

\ifthenelse{\boolean{book}}{

We are interested in what happens to  $\O_\spax$-modules $\sheaf$
under the refinement morphism
$\mortop: \spax_\topo \to \spax$ given
by a Grothendieck topology on $\spax$.

\begin{lemma}
\label{sheafifyexactproperty}
Let $\spax$ be a scheme and let
$\catopenspax    $ be a Grothendieck topology on $\spax$. Then the
module-sheafification in  the Grothendieck topology is right-exact
on $\O_\spax$-modules, that is, an exact complex of
$\O_\spax$-modules $\olmrosheaf $ yields an exact complex
$\sheafrightcomplexx {_\topo} $ of $\O_\topo$-modules sheaves on
$\spax_\topo$.
\end{lemma}
\begin{proof}
This follows from Lemma \ref{rightexactpullback}.
\end{proof}}{}


\begin{example}
\label{trivialsite}
Let $\ring$ be a commutative ring and consider the category
consisting of one object $\spax$ (`$\specring$') and one identical
morphism. We define a sheaf of rings by
$\Gamma(\spax,\O_\spax) = \ring$
and call this the \emph{trivial site} for $\ring$. An $\ring$-module
$\modul$ yields immediately an
$\O_\spax$-module, which we also denote by $\modul$.

Let $\catopen_\ring $ be a subcategory of $\ring$-algebras
(or rather its opposite category, so we will write $\Spec \rings $ instead of
$\rings$) containing $\ring$, closed under tensor products and
endowed with a certain Grothendieck topology. Then
$\catopen_\ring$ is immediately a preringed site by setting
$\Gamma(\Spec \rings,\O^\pre) = \rings$,
and we get a morphism of preringed sites by the inclusion
$\sitemortop : \{ \spax \} \to \catopen_\ring$
(with $\sitemortop (\spax) = \spax$).
The ring homomorphism is given by the identity
$\Gamma( \spax,\O_\spax)
= \ring \to \Gamma( \sitemortop ( \spax ),\O^\pre) = \ring$.
The inverse image presheaf of an $\O_\spax$-module $\modul$ is given
by the constant presheaf since
$$\Gamma(\Spec \rings, \sitemorinv (\modul))
= \colim_{ (\obu, \phili) \in \ideal_\rings^\sitemor } \Gamma (\obu, \modul)
= \colim_{ \phili:\ring \to \rings   } \Gamma( \spax, \modul)
= \modul \, .$$
Here the (only) indexing morphism $\ring \to \rings   $ is the
structure homomorphism of the $\ring$-algebra $\rings$. The
$\O_\spax$-module $\modul$ yields the corresponding
$\O^\pre$-premodule $\sitemorpre (\modul)$ by sending
$\Spec \rings \mapsto \modul \tensorr \rings$. The sheafification
$\modul_\topo$
(and also $\O_\topo$) might be difficult to compute,
depending on the topology.
\end{example}

\begin{example}
Let $\spaxeqspecring$ be an affine scheme with the
trivial topology and let $\spax_\zar$ be the affine scheme with the
Zariski topology and with the structure sheaf $\O_\zar$, so that
$\morsite: \spax_\zar \to \spax_\triv$ is a morphism of ringed sites.
If $\modul$ is an $\ring$-module, then the inverse image presheaf
$\sitemorinv(\modul)$ is the constant presheaf
$\obv \mapsto \modul$. The associated $\O_\zar$-premodule is given by
$\obv \mapsto \modul \tensorr  \Gamma(\obv, \O_\zar)$.
The pull-back of $\modul$ is the sheaf associated to this presheaf
and gives the quasicoherent module $\tilde{\modul}$
\cite[Proposition II.5.1]{haralg}.
\end{example}

\ifthenelse{\boolean{book}}{
\begin{example}
\label{zariskisiteexample}
Assume that a scheme $\spax$ is endowed with a Grothendieck topology
which is a refinement of the Zariski topology, and let
$\sitemor:
\spax_\topo \to \spax_\zar$ be the corresponding site morphism. An
$\O_\zar$-module $\modul$ on
$\spax_\zar$ yields the inverse image presheaf $\sitemorinv
(\modul)$ and the
$\O^\pre$-premodule $\sitemorpre( \modul)$ is then
given by
$$\obv \mapsto \Gamma(\obv , \sitemorinv (\modul))
\tensor_{\Gamma(\obv , \sitemorinv ( \O_\zar ))}
\Gamma(\obv, \O_\obv) \, $$
for $\obv \in \catopenspax $
(Definition \ref{sitemorpredef}).
We denote its sheafification $\sitemorsheaf (\modul) $, which is an
$\O_\topo$-module, by $\modul_\topo$. Via the sheafification we have
a natural homomorphism
$ \shef: \Gamma(\spax, \modul) \lto \Gamma(\spax_{\topo}, \modul_{\topo})$.
An $\O_\zar$-module homomorphism
$\modulsec \to \modul$ on $\spax$
yields an
$\O_\topo$-module homomorphism
$\modulsec_{\topo} \to \modul_{\topo}$ on $\spax_\topo$.
\end{example}}{}

\begin{definition}
\label{singledef}
Let $\catopenspax$ denote a Grothendieck topology on a scheme
$\spax$.


We call $\catopenspax$ or $\spax_\topo$ \emph{affine}, if all
structural morphisms
$\obu \to \spax$ in $\catopenspax$ are
affine
(then all morphisms are affine).

We call it \emph{of finite type} if all morphisms $\obu \to \obv$ in
$\catopenspax    $ are of finite type.

We call it \emph{quasicompact} if for every covering
$\obcov_\indcov \to\obu$, $\indcovinset$,
there exists a finite subset
$\setindcovsec \subseteq \setindcov$
such that
$\obcov_\indcov \to \obu$, $\indcov \in \setindcovsec$, is also a covering.

We call it \emph{single-handed} if it is of finite type and if every
covering is given by a single morphism $\obu \to \obv$.

We call it \emph{covering} if every morphism in $\catopen_\spax$ is
a cover.
\end{definition}

{\ren{\catind}{{I}}

\begin{proposition}
\label{sheafzariskicompare}
Let $\spax$ be a scheme
and let $\catopenspax$ be a Grothendieck topology on $\spax$
{\rm(}recall Convention \ref{basistopology}{\rm)}.
Let $\sitemor: \spax_\topo \to \spax$ be the corresponding site
morphism.
Let $\struto:\obv \to \spax$ be in $\catopenspax$, considered as a
Zariski morphism, and let $\modul$ be an $\O_\spax$-module. Then the following hold.

\numiii

\begin{enumerate}

\item
We have
$\Gamma(\obv, \sitemorinv (\modul))
= \Gamma(\obv, \strutoinv (\modul))$.

\item
We have
$\Gamma(\obv, \sitemorpre (\modul))
= \Gamma(\obv,\strutopre (\modul))$.

\item
The sheafification
$\sitemormodule (\modul)= \modul_\topo$
equals the
$\catopenspax $-sheafification of the $\O^\pre$-premodule
$\obv \mapsto \Gamma(\obv, \strutopre (\modul))$.

\item
If $\catopen_\spax$ is a refinement of the Zariski topology or an
affine topology on an affine scheme, then $\modul_\topo$ equals also
the sheafification of
$\obv \mapsto \Gamma(\obv, \strutosheaf (\modul))$,
where $\strutosheaf$ is the pull-back of modules in the
Zariski topology.
\end{enumerate}
\end{proposition}
\begin{proof}
(i).
We have
$$\Gamma(\obv, \sitemorinv (\modul))
=\colim_{(\openzar, \phili) \in \catind^\sitemor_\obv}
\Gamma(\openzar , \modul)
= \colim _{\struto (\obv) \subseteq\openzar} \Gamma(\openzar , \modul)
= \Gamma(\obv, \strutoinv(\modul))\, ,$$
because there exists by our convention a morphism
$\obv\to \openzar = \sitemortop (\openzar)$
in $\catopenspax $ if and only if
$\struto (\obv) \subseteq \openzar$.

(ii). We have by Definition \ref{sitemorpredef} and part (i)
\begin{eqnarray*}
\Gamma(\obv, \sitemorpre (\modul))
&=&\Gamma(\obv, \sitemorinv(\modul)) \tensor_{\Gamma(\obv,\sitemorinv
(\O_\zar))} \Gamma(\obv, \O^\pre) \cr
&=& \Gamma(\obv, \strutoinv (\modul))
\tensor_{\Gamma(\obv,\strutoinv (\O_\zar))} \Gamma(\obv,\O_\obv)
= \Gamma(\obv, \strutopre (\modul)) \, .
\end{eqnarray*}

(iii) follows from (ii).
(iv) follows also from (ii), since in the first case
$\obv \mapsto \Gamma(\obv, \strutosheaf (\modul))$ is the
Zariski sheafification of the presheaf
$ \obv \mapsto \Gamma(\obv,\sitemorinv (\modul))
= \Gamma(\obv, \strutoinv (\modul))$,
and the topology is a refinement of the Zariski topology.
In the second case
$\Gamma(\obv, \strutosheaf (\modul))
= \Gamma(\obv, \strutopre (\modul))$
anyway.
\end{proof}

\ifthenelse{\boolean{book}}{ \begin{bookexample} The condition holds
if for example $\catopenspax    $ is a full subcategory. It does not
hold in general. Let $\spax= Spec R \supset D(f)=U$, and let $V= U
\times \AA$ be an object in $\catopenspax    $ and assume that $j :
U \times \AA \to \spax$ belongs to the category, but the projection
to $U$ not. Then $\Gamma(U \times \AA, \sitemorinv
(\O_\spax))=\ring$, but $\Gamma(U \times \AA, j^{-i} (\O_\spax)) =
\ring_f$.
\end{bookexample}}{}

\begin{corollary}
\label{pullbacksheafify}
Let $\spax$ denote a scheme endowed with a
Grothendieck topology $\catopen_\spax$, and suppose that it is
either a refinement of the Zariski topology or affine over an affine
scheme. Let
$\struto:\spay \to \spax$ be in $\catopenspax $ and let
$\modul$ be an $\O_\spax$-module on $\spax$. Let
$\modul_\topo$ be its sheafification in $\spax_\topo$ and let
$\struto^*(\modul)$ be the Zariski pull-back to $\spay$. Let
$\catopen_\spay$ be the category consisting of all morphisms
$\spaz \to \spay$ in $\catopen_\spax$
as objects and the morphisms in
$\catopen_\spax$ over $\spay$ as morphisms {\rm(}and the induced
coverings{\rm)}.
Then the restriction sheaf
$(\modul_\topo)|\catopen_\spay$ given by
$(\spaz \to \spay) \mapsto \Gamma(\spaz,\modul_\topo)$ equals the sheafification
$(\struto^*(\modul))_\topo$ in $\catopen_\spay$.
\end{corollary}
\begin{proof}
We have to show for
$\strutosec: \spaz \to \spay$
(in $\catopen_\spay $) that
$\Gamma(\spaz, (\struto^*(\modul))_\topo)
=\Gamma(\spaz, \modul_\topo)$.
By Proposition \ref{sheafzariskicompare}(iv)
both sheaves come from the presheaf
$( \strutosec : \spaz \to \spay) \mapsto \strutosec^*(\struto^* (\modul))$,
one time sheafified in $\catopen_\spax$ and one time sheafified in
$\catopen_\spay$. The objects $\spaz \to \spay$ have the same
coverings in $\catopen_\spax$ and in $\catopen_\spay$, because a
covering $\spaz_\indcov \to \spaz$ over $\spax$ can also be
considered as a covering with basis $\spay$ via $\spaz \to
\spay$. Hence the associated separated presheaves are the same and
then also their sheafifications.
\end{proof}

\begin{example}
Proposition \ref{sheafzariskicompare} and Corollary
\ref{pullbacksheafify} are not true without Convention
\ref{basistopology} on Zariski factorizations. To give an example,
let $\ring$ be a normal domain of dimension $\geq 2$. Let $\spax=
\Spec \ring$ have the trivial topology and let $\catopen_\spax$ be
given by the Zariski open subsets which contain all points of
codimension one (and Zariski coverings). Then
$\Gamma(\obv,\O_\spax)=\ring$ for all
$\obv \in \catopenspax$, and $\spax_\topo$
is a ringed site with constant structure sheaf. For an
$\ring$-module $\modul$ the inverse image presheaf under
$\sitemor:\spax_\topo \to \spax$ is constant, and so also
$\sitemorpre(\modul)$ is constant $=\modul$. Let now
$\spay \in \catopen_\spax$
and $\modul$ be an $\ring$-module such that
$\Gamma(\spay,\tilde{\modul}) \neq \modul$, e.g. take
$\struto: \spay = \spax -\{\point\} \to \spax$ and
$\modul= \fieldres (\point)$, where
$\point \in \spax$ is a closed point. Then
$\struto^*(\modul)=0$ and so
$(\struto^*(\modul))_\topo=0$, but
$\Gamma(\spay ,\modul_\topo)= \modul$.
\ifthenelse{\boolean{book}}{
A similar example shows that we need the assumptions in Proposition \ref{sheafzariskicompare}(iv).
Let $\spax$ be the affine plane with the trivial topology and
$\catopen_\spax$ consisting of $\spax$ and $\spau=\spax -
\{\point\}$ ($\point$ a closed point), again with the trivial
topology.
Then for $\modul=\fieldres(\point)$ the inverse presheaf is just
$\fieldres(\point)$, and in the given topology this is also the
sheaf pull-back. In the Zariski topology however we have
$\^struto^*(\modul)=0$.
}{}
\end{example}
}

\ifthenelse{\boolean{book}}{

\begin{example}
\label{recallflattopology}
We recall the definition of the
\emph{flat} and of the \emph{\'{e}tale} Grothen\-dieck topology on a
scheme $\spax$
(we describe the so-called small site).
The category
$\catopenspax$ is given by the full subcategory of
$\catsch_\spax$ consisting of flat (\'{e}tale) morphisms $\obu \to
\spax$ of finite type
(there are variants like the fppf and the fpqc topology,
see \cite[Expos\'{e} IV.6.3]{SGA3} and \cite[Example 2.30 and
Definition 2.32]{vistoligrothendieck}) and such that a family
$ \obcov_\indcov \to\obu$, $\indcovinset$,
is a covering if the (open) images of the $\obcovindic$ do
cover $\spax$ in the set-theoretical sense.

A common feature of the classical Grothendieck topologies
(the flat, the \'{e}tale and the Zariski topology;
this feature also holds for the Nisnevich \emph{topology}
\cite{nisnevich}, which is between the
Zariski and the \'{e}tale topology, and for the \emph{primitive
topology} of M. Walker \cite{walkerprimitive})
is that the presheaf
$\obv \mapsto \Gamma(\obv, \strutosheaf (\modul))$
(denoted by $W(\modul)$ in \cite[Example VII.2c)]{SGA4} and in \cite{milne})
associated to a quasicoherent $\O_\spax$-module $\modul$ is already
a sheaf \cite[Proposition I.2.18 and Corollary II.1.6]{milne}.
Related to this are the facts that a
Zariski-exact complex is also exact in the flat (\'{e}tale) topology
(in the small site) and that quasicoherent sheaves on an affine
scheme have trivial flat (\'{e}tale) cohomology
(\cite[Corollaire VII.4.4]{SGA4} or \cite[Proposition III.3.7]{milne}).
These properties will not hold
for the topologies which give interesting closure operations.}


\ifthenelse{\boolean{book}}
{Another common property of these
classical topologies is that they are what is called
\emph{subcanonical}. This means that every presheaf (of sets,
groups) on the flat (\'{e}tale) site which is representable by a
scheme is already a sheaf, i.e. that the assignment $U \mapsto
\Mor_\spax(U,Z)$ is a sheaf for a given scheme $Z$ over $\spax$.
This relies on the fact that a faithfully flat morphism $\spay \to
\spax$ of finite type is a \emph{strict epimorphism}, meaning that
the complex of sets
$$\Mor(\spax,Z) \lto \Mor(\spay,Z) \rightrightarrows \Mor(\spay \times_\spax \spay,Z) $$
is exact for every scheme $Z$ \cite[Theorem I.2.17]{milne}.
\end{example}}{}

\ifthenelse{\boolean{book}}{

\begin{remark}
Non flat Grothendieck topologies seem to have first occurred in the
work of V. Voevodsky in an attempt to understand the algebraic
topology of varieties (\cite{voevodskyhomology}, \cite{mazzaweibel}
\cite{levineintroduction} \cite{goodwillielichtenbaum}
\cite{suslinvoevodskysingularhomology}), namely the so-called
$h$-\emph{topology} and the $qfh$-\emph{topology} \cite[Section
6.1]{levineintroduction}. In the $h$-topology a cover is given by a
universal (Zariski-)topological epimorphism (a universal
submersion). These include in particular proper surjective maps. The
\emph{proper topology} was studied by S.-I. Kimura (see
\cite[Definition 3.1]{kimuraalexander}) in the context of Alexander
schemes. We will deal with such submersive topologies and its
relation to the integral closure in Section \ref{propersubsection}
and in \cite{blicklebrennersubmersion}.

One advantage of this topology is that
every excellent scheme has a smooth cover (by resolution of
singularities in characteristic zero and by de Jongs alterations in
positive characteristic). Suitable comparison results \cite[Theorems
8.3, 9.1 and 10.1]{levineintroduction} allow in the end to show that
Suslins algebraic singular (co)homology for finite (constant)
coefficients coincides with \'{e}tale cohomology \cite[Corollary
10.3]{levineintroduction} for a variety over an algebraically closed
field.

We will encounter among other non-pure topologies also the
$h$-topology, but our viewpoint is what happens to quasi-coherent
modules and their submodules in such topologies.
\end{remark}}{}

%



To a large extent we will restrict to affine, single-handed
Grothendieck topologies on an affine scheme $\spaxeqspecring$, so
that
$\catopen_\ring = \catopenspax$ is given by a category of
$\ring$-algebras and the coverings are given by a certain class of
$\ring$-algebra homomorphisms.

\ifthenelse{\boolean{book}}{
\begin{bookremark}
Any property of a ring homomorphism which holds for isomorphisms,
which is stable under composition and under base change defines an
affine single-handed Grothendieck topology by declaring the
homomorphisms $\ring \to \rings   $ with this property to be a
cover. There are several possibilities for $\catopen_\ring$; one can
take the full category $\catsch_\ring$, or the category of
$\ring$-algebras such that the structural homomorphism has the
property, or allow only homomorphisms with the property.
\end{bookremark}}{}

{\ren{\group}{{H}}

\begin{remark}
\label{zariskisinglehanded}
We will often consider the Zariski topology on a
scheme $\spax$ as an (affine) single-handed Gro\-then\-dieck
topology, called the (\emph{affine}) \emph{single-handed
Zariski topology},
where
$\catopenspax$ consists of all finite disjoint unions of (affine) open
immersions
$\obu \to \spax$ and where the only coverings are given by
$\biguplus \obcovindic \to \obusec$, where the $\obcovindic$ do cover $\obusec$ as a set.
If
$\spay= \biguplus_{\iindex} \obcovindic \to \spax$
is a Zariski covering and $\sheaf $ is a quasicoherent sheaf on
$\spax$,
then
\begin{eqnarray*}
\gammashf {\spay\times_\spax \spay }
&=& \gammashf {( \biguplus_\indcovinset \obcovindic) \timesx
(\biguplus_\indcovinset \obcovindic)} \cr
&=& \gammashf {\biguplus_{(\indcov,\indcovsec) \in \setindcov \times \setindcov}
\obcovindic \timesx \obcov_\indcovsec}\cr
&=& \bigoplus_{(\indcov, \indcovsec) \in \setindcov \times \setindcov}
\gammashf {\obcovindic \timesx \obcov_\indcovsec}
 = \bigoplus_{(\indcov,\indcovsec) \in \setindcov \times \setindcov}
\gammashf {\obcovindic \cap \obcov_\indcovsec} \, .
\end{eqnarray*}
This holds more generally for every sheaf which respects disjoint
Zariski unions (as in Lemma \ref{pullbacksheafdisjoint}) and it
holds also for products with more than two factors. From this it
follows that the \v{C}ech cohomology for such sheaves is the same in
the Zariski site and in the single-handed Zariski site.

\end{remark}
}

\ifthenelse{\boolean{book}}{
In one respect one has to be careful with the Zariski site and the
single Zariski site.
The constant presheaf
$\openzar \mapsto \group$ for all $\openzar \in \catopenspax    $
($\emptyset \mapsto 0$) given by a group $\group$
is already a sheaf in a single-handed Grothendieck topology, whereas
this is not true in the Zariski (or flat, \'{e}tale) site, where its
sheafification is given by sending
$\openzar \mapsto \group^\setj$, $\setj$
being the set of connected components of $\openzar$.
}{}

Although we will mainly deal with single-handed Grothendieck
topologies it is good to know when they respect Zariski covers.

\begin{definition}
\label{respectzariskidef}
Let $\spax$ be a scheme with a
Grothendieck topology $\catopen _\spax$.
We say that the Grothendieck topology \emph{respects} (\emph{disjoint})
\emph{Zariski coverings} if a family of morphisms $\obcovindic \to \obu$,
$\indcovinset$,
(inside $\catopenspax $) is a covering if and only if for an
arbitrary (disjoint) Zariski cover
$ \obu = \bigcup_\indcovsecinset \obcovsec_\indcovsec$
($\obu = \biguplus_\indcovsecinset \obcovsec_\indcovsec$)
the family
$\obcovindic \times_\obcovsec \obcovsec_\indcovsec
\to \obcovsec_\indcovsec$,
$\indcovinset$, is a covering of $\obcovsec_\indcovsec$ for every $\indcovsec$.
\end{definition}

{
\ren{\mor}{{\psi}}

\ren{\elet}{{t}}

\begin{lemma}
\label{pullbacksheafdisjoint}
Let
$\spax = \biguplus_\indcovinset \obcovindic$ be a finite disjoint
union of open subschemes of a scheme $\spax$.
Let an affine single-handed Grothendieck topology
$\spax_\topo$ be given which is a refinement of the single-handed Zariski topology
and which respects disjoint Zariski coverings
.
Let $\modul$ be an $\O_\spax$-module
on $\spax$. Then
$\Gamma(\spax_\topo, \modul_\topo)
=\bigoplus_{\indcovinset} \Gamma(\obcov_\indcov,\modul_\topo)
=\bigoplus_{\indcovinset} \Gamma(\obcov_\indcov,(\morcov_\indcov^*\modul)_\topo)  $
($\morcov_\indcov:\obcov_\indcov \to \spax$).
\end{lemma}
\begin{proof}
The second equation follows from Corollary \ref{pullbacksheafify}.
The morphisms
$\morcov_\indcov: \obcov_\indcov \to \spax$ induce the restrictions
$\Gamma(\spax_\topo, \modul_\topo)
\to \Gamma(\obcov_{\indcov \,\toponw},\modul_\topo)$.
An element
$\eles \in \Gamma(\spax_\topo, \modul_\topo)$
is represented by
$\elet \in \Gamma(\spay, \mor ^* (\modul))$
(Proposition \ref{sheafzariskicompare}(iv)),
where
$\mor: \spay \to \spax$
is a cover in the topology. Suppose that $\eles$ restricts to $0$ in
$\Gamma(\obcov_\indcov, \modul_\topo)$ for all open subsets
$\obcov_\indcov$.
This means that
$\shef(\elet) \in \Gamma(\spay, \modul_\topo)$
restricts to $0$ in
$\Gamma( \obcov_\indcov \times_\spax \spay, \modul_\topo)$
and so there exist covers
$\spaz_\indcov \to \obcov_\indcov \timesx \spay$
such that the pull-backs of $\elet$ are
$0$.
These covers glue together to the cover
$\biguplus_\indcovinset \spaz_\indcov \to \spay$
(since being a cover can be tested on the disjoint subsets of the union),
hence also
$\shef(\elet) = 0$ in $\Gamma(\spay, \modul_\topo)$
and so $\eles=0$ in
$\Gamma(\spax, \modul_\topo)$.
On the other hand, let
$\eles_\indcov \in \Gamma(\obcov_\indcov, \modul_\topo)$ be given.
Say $\eles_\indcov$ is represented by compatible elements
$\elet_\indcov \in \Gamma(\spaz_\indcov,\mor_\indcov^*(\modul))$,
where
$\mor_\indcov:\spaz_\indcov \to \obcov_\indcov$
is a cover for each $\indcovinset$.
Then again these covers glue together to a cover of
$\spax$, and the tuple
$(\elet_\indcov)$
is also compatible over
$\spax$, so it defines a global section.
\end{proof}
}

\ifthenelse{\boolean{book}}{

\begin{remark}
\label{zariskidisjointremark}
If a Grothendieck topology respects disjoint Zariski covers,
then disjoint covers glue together:
if $ \obcovindic \to \obu$, $\indcovinset$,
and $\obcovsec_\indcovsec \to \obv$, $\indcovsecinset$, are coverings, then
also
$\obcovindic, \obcovsec_\indcovsec \to \obu \uplus \obv$,
$\indcov, \indcovsec \in \setindcov \uplus \setindcovsec$,
is a covering, because the restrictions to $\obu$ and $\obv$ give
the coverings back
(at least if adding empty sets does not destroy
the covering property).
\end{remark}

\begin{corollary}
\label{pullbacksheafdisjoint}
Let
$\spax = \biguplus_\indcovinset \obcovindic$ be a finite disjoint
union of open subschemes. Let a Grothendieck topology $\spax_\topo$
be given which respects disjoint covers and such that
$\obcov_\indcov \in \catopenspax$.
Let $\modul$ be an $\O_\spax$-module
on $\spax$. Then
$\Gamma(\spax_\topo, \modul_\topo)
=\bigoplus_{\indcovinset} \Gamma(\obcov_\indcov,\modul_\topo)$.
\end{corollary}
\begin{proof}
The statement is true for arbitrary sheaves if
$\obcov_\indcov \to \spax$, $\indcovinset$,
is a covering in the topology, but this is here not part of
the assumptions. The morphisms
$\morcov_\indcov: \obcov_\indcov \to \spax$ induce the restriction homomorphisms
$\Gamma(\spax_\topo, \modul_\topo)
\to \Gamma(\obcov_{\indcov \,\toponw},\modul_\topo)
= \Gamma(\obcov_{\indcov \,\toponw},(\morcov_\indcov^*(\modul))_\topo)$
(by Corollary \ref{pullbacksheafify}).
An element
$\eles \in \Gamma(\spax_\topo, \modul_\topo)$
is represented by
$\eles_\indcovsec \in \Gamma(\spax_\indcovsec, \morcovsec_\indcovsec^* (\modul))$,
where
$\morcovsec_\indcovsec: \spax_\indcovsec \to \spax$,
$\indcovsecinset$, is a covering in the topology. Suppose that $\eles$
restricts to $0$ in
$\Gamma(\obcov_\indcov, \modul_\topo)$ for all open subsets
$\obcov_\indcov$.
This means that $\shef(\eles_\indcovsec)$ restricts to $0$ in
$\Gamma( \obcov_\indcov \times_\spax \spax_\indcovsec, \modul_\topo)$
and so there exist coverings
$\obv_{\indcovtri_{(\indcov,\indcovsec)}} \to \obcov_\indcov \timesx \spax_\indcovsec$
such that the pull-backs of $\eles_\indcovsec$ are
$0$.
These coverings glue together to a covering of
$\spax_\indcovsec$
(Remark \ref{zariskidisjointremark}),
hence also
$\shef(\eles_\indcovsec) = 0$ in $\Gamma(\spax_\indcovsec, \modul_\topo)$
and so $\eles=0$ in
$\Gamma(\spax, \modul_\topo)$.
On the other hand, let
$\eles_\indcov \in \Gamma(\obcov_\indcov, \modul_\topo)$ be given.
Say $\eles_\indcov$ is represented by compatible elements
$\eles_{\indcov, \indcovsec_\indcov}
\in \Gamma(\obv_{\indcovsec_\indcov},\morcov_{\indcovsec_\indcov}^*(\modul))$,
where
$\morcov_{\indcovsec_\indcov}: \obv_{\indcovsec_\indcov}
\to \obcov_\indcov$,
$\indcovsec_\indcov \in \setindcovsec_\indcov$,
is a covering.
Then again these coverings glue together to a covering
$\obv_\indcovsec \to \spax$, $\indcovsec \in \setindcovsec
= \biguplus_{\indcovinset} \setindcovsec_\indcov$.
The elements
$\eles_{\indcov, \indcovsec_\indcov}$
are also compatible over
$\spax$, so they define a global section.
\end{proof}}

\subsection{Rings and modules of global sections}
\label{globalsectionsubsection}

\

\medskip
In the flat (or the \'{e}tale) Grothendieck topology on an affine scheme
$\specring$ we have
$\Gamma(\specring, \O_\topo)= \ring$.
This does not hold for non-pure topologies in general.
We have to compute in a given Grothendieck topology $\catopenspax$ the sheafification
$\obu \mapsto \Gamma(\obu, \O_\topo)$, $\obu \in \catopenspax $,
associated to the presheaf
$\obu \mapsto \Gamma(\obu, \O_\obu)$,
and this computation is a non-trivial matter.

\begin{proposition}
Suppose that we have a Grothendieck topology on a certain
subcategory $\grocont$ of schemes, meaning that for every $\spax \in
\grocont$ we have a Grothendieck topology $\catopenspax$ and
every scheme morphism $\spay \to \spax$ in $\grocont$ induces a site
morphism
$\spay_\topo \to \spax_\topo$.
Then we get a functor
$\spax \mapsto \Gamma(\spax,\O_\topo)$, $\spax \in \grocont$, together with
a natural transformation
$\shef: \Gamma(\spax, \O_\spax) \to \Gamma(\spax, \O_\topo)$.
\end{proposition}
\begin{proof}
This is clear.
\end{proof}

\begin{remark}
In this way we will encounter the perfect closure (Theorem
\ref{frobeniusglobal}), the semi-normalization (Proposition
\ref{finiteglobalsection}), the ring of constructible sections
(Theorem \ref{surjectiveglobalsection}), the completion of a local
ring (Proposition \ref{completeproperties}) and more.
In particular, we have
$\Gamma((\specring)_\topo, \O_\topo) \neq \ring$.
From this it
follows that the additive group scheme $\AA^1$ does not define a
sheaf in the topology, because
$\Mor(\spax, \AA^1) \cong
\Gamma(\spax,\O_\spax ) \neq \Gamma(\spax_\topo, \O_\topo)$
\cite[Exercise 2.4]{haralg}.
Hence such topologies are not \emph{subcanonical}, i.e. not every representable functor is a sheaf
(see \cite[Definition 2.54]{vistoligrothendieck} or \cite[D\'{e}finition II.2.5]{SGA4}).
\end{remark}

\ifthenelse{\boolean{book}}{
\begin{booklemma}
Suppose that $\catopenspax $ is an affine Grothendieck topology on
$\spaxeqspecring$ and that $\modul$ is an $\ring$-module. Then for
$\ring \to \ringsec  $ in $\catopen_\ring$ we have
$$(\modul \tensorr \ringsec  )_0
 \!=\! \{\elem \in \modul \tensorr \ringsec   :
\exists \mbox{ a covering } \varphi_i: \ringsec   \to \ring_i, i \in I,
\mbox{ such that } \varphi_i(\elem)
=0 \} \, $$
and the associated separated presheaf is
$$\Gamma(\ringsec  , \modul_1) := (\modul \tensorr \ringsec  )_1
 = \modul \tensorr \ringsec   /(\modul \tensorr \ringsec  )_0 \, .$$
Moreover, $\Gamma(\ring, \modul_\topo)$ is given as the limit of all
tuples $\elem_i \in \modul \tensorr \ring_i$, $\ring \to \ring_i$ a
covering, such that the restrictions of $\elem_i$ and $\elem_j$ in
$(\modul \tensorr \ring_i \tensorr \ring_j)_1$ are the same.
\end{booklemma}
\begin{proof}
This follows from general facts about sheafification, see Section
\ref{sheafifysubsection}.
\end{proof}}{}

A global section $\fuglob \! \in \! \Gamma(\spax, \O_\topo)$ is
represented by a compatible family
$\fuglob_\indcov \! \in \!
 \Gamma(\obcovindic, \O_{\obcovindic})$, $ \indcovinset$,
where
$\obcovindic \to \spax$, $ \indcovinset$,
is a covering in the Grothendieck topology, and where the
compatibility condition means that
$\project_1^*(\fuglob_\indcov) - \project_2^*(\fuglob_\indcovsec)
\in \Gamma(\obcovindic \timesx \obcovindicsec)$ has the property that it vanishes on a covering of
$\obcovindic \timesx \obcovindicsec $.
A similar characterization holds for an $\O_\spax$-module $\modul$.
Quite often we will deal with a single cover $\spay \to \spax$,
and a global section is then
given by $\fuglob \in \Gamma(\spay, \O_\spay )$ such that
$\project_1^*(\fuglob) - \project_2^*(\fuglob) \in \Gamma(\spay \timesx \spay, \O)$
is zero in a covering of $\spay \timesx \spay$.

\begin{lemma}
\label{sectionlemmaglobal}
Suppose that a Grothendieck topology on a
scheme $\spax$ has the property that every covering
$\morcov_\indcov : \obcovindic \to \spax$, $\indcovinset$,
has a section
$\sect_\indcov: \spax \to \obcovindic$ for some $\indcov$
{\rm(}M. Walker calls such $\spax$ \emph{acyclic};
cf. \cite[Lemma 3.7]{walkerprimitive}{\rm)}.
Then for an $\O_\spax$-module $\modul$
the natural homomorphism
$\shef: \Gamma(\spax,\modul) \to \Gamma(\spax, \modul_\topo)$
is an isomorphism.
\end{lemma}
\begin{proof}
Let $\elem \in \Gamma(\spax,\modul)$. If $\elem = 0$ in
$\Gamma(\spax, \modul_\topo)$, then there exists a covering
$\morcov_\indcov: \obcovindic \to \spax$, $\indcovinset$, such that
$\morcov_\indcov^*(\elem) =0$ in
$\Gamma(\obcovindic, \cmpbimodul)$
for all $\indcovinset$. The factorization
$\spax \to \obcovindic \to \spax$ shows that $\elem =0$.

On the other hand, let $\elem \in \gammaxmodul$ be represented by
$\elem_\indcov \in \Gamma(\obcovindic, \cmpbimodul) $ for a covering
$\obcov_\indcov \to \spax$, $\indcovinset$.
Then via
$\spax \to \obcovindic \to \spax$ the element
$\elem_\indcov$ gives an element
$\secpb_\indcov(\elem_\indcov) \in \gammaxmodul$,
and since
$\id: \spax \to \spax$ is a covering, we may represent
$\elem$ also by
$\secpb_\indcov(\elem_\indcov)$.
\end{proof}

\ifthenelse{\boolean{book}}{
\begin{bookexample}
\label{nisnevich} The Nisnevich topology (or completely decomposed
topology) is a variant of the \'{e}tale topology where a covering is
an \'{e}tale covering and fulfills additionally the condition that
every residue class field has a lifting, see \cite{nisnevich}. This
implies in particular that in the Nisnevich topology of a field
every covering has a section as in Lemma \ref{sectionlemmaglobal}.
The stalks in the Nisnevich topology are the Henselizations of local
rings \cite[D\'{e}finition 18.6.5]{EGAIV}.
\end{bookexample}}{}

\begin{lemma}
\label{autoglobal}
Let $\spay \to \spax$ be a cover in a
Grothendieck topology on $\spax$, and let
$\fuglob \in \Gamma(\spay,\O_\spay )$ be a compatible element,
defining a global element in
$\Gamma(\spax,\O_\topo)$.
If $\automor :\spay \to \spay$ is an
$\spax$-automorphism, then $\automor^*(\fuglob)$ is also compatible and
it defines the same global element.
\end{lemma}
\begin{proof}
It is clear that $\automor^* (\fuglob)$ is also compatible, since an
automorphism preserves the compatibility conditions. Suppose that
the pull-backs of $\project_1^*(\fuglob)$ and
$\project_2^*(\fuglob)$ are the same under the covering
$\spaz_\indcov \to \spay \timesx \spay$,
$\iini$.
The automorphism $\id \times \automor: \spay \timesx \spay \to \spay
\timesx \spay$ shows that also
$(\id \times \automor)^*
(\project_1^*(\fuglob))
 = \project_1^*(\fuglob)$ and
 $(\id \times \automor)^*
(\project_2^*(\fuglob))
 = \project_2^*(\automor^*(\fuglob))$ are the same in some
covering.
Therefore also
$ \project_2^*(\fuglob)$ and
$\project_2^*(\automor^*(\fuglob))$ are the same in some covering. Hence
$\fuglob$ and $\automor^*(\fuglob)$ are the same in some covering and so
they define the same element in $\Gamma(\spax, \O_\topo)$.
\end{proof}

If a Grothendieck topology $\catopenspax$ on $\spaxeqspecring$ has
the property that the reduction defines a cover, then it is
important to know for an element
$\fuglob \in \rings$, $\ring \to \rings $
a cover, whether
$\fuglob \tensor 1 -1 \tensor \fuglob \in \rings \tensorr \rings$
is nilpotent. In this case $\fuglob$ is
compatible and defines a global element in
$\Gamma(\ring_\topo,\O_\topo)$.

\begin{lemma}
\label{nilpotenttest}
Let $\rings$ denote an $\ring$-algebra,
$\fuglob \in \rings$. Then
$\fuglob \tensor 1 - 1 \tensor \fuglob \in\rings \tensorr \rings$
is nilpotent if and only if the following hold: For all
homomorphisms
$\homtest: \ring \to \field   $, $\field$
a field, and any two extensions
$\homlift_1, \homlift_2: \rings\to \field
$ of $\homtest$, we have
$\homlift_1(\fuglob)= \homlift_2(\fuglob)$.
\end{lemma}
\begin{proof}
The element
$\differ (\fuglob)
=\fuglob \tensor 1 - 1 \tensor \fuglob
\in \rings \tensorr \rings$
is nilpotent if and only if for all ring homomorphisms
$ \rings \tensorr \rings \to \field$, $\field$ a field, the image of $\differ
(\fuglob)$ is $0$.
These homomorphisms are in one to one
correspondence with the data described in the lemma, so the result
follows.
\end{proof}


\ifthenelse{\boolean{book}}{
\begin{lemma}
Let $\spax_\topo$ denote a scheme endowed with a Grothendieck
topology. Let $\spay \to \spax$ be a cover and let $\fuglob \in
\Gamma(\spay, \O_\spay )$ denote a compatible element. Then the
corresponding element $\fuglob \in \Gamma(\spay, \O_\topo)$ maps to
the corresponding element $\fuglob' \in \Gamma(\spax, \O_\topo)$.
\end{lemma}
\begin{proof}
\end{proof}
}{}

{\ren{\point}{{x}}

\begin{lemma}
\label{compatiblevalue}
Let a Grothendieck topology be given on a
subcategory $\grocont$ of schemes such that every morphism in
$\grocont$ induces a site morphism. Suppose that for $\spax \in
\grocont$ all $\obu \to \spax$ {\rm (}in $\catopenspax    ${\rm)}
and all morphisms
$\Spec \fieldres(\point) = \point \to \spax$,
$\point \in \spax$ a point, belong to $\grocont$.
Suppose also that
$\Gamma(\point, \O_\topo) = \fieldres(\point)$
{\rm (}or $\fieldres(\point)^\perf${\rm)}
for every point $\point$ in $\grocont$.
Then a compatible element
$\fuglob \in \Gamma(\spay,\O_\spay )$,
$\spay \to \spax$ a cover, has the property that for
two points
$\pointy_1,\pointy_2 \in \spay$ over $\point \in \spax$
the values of $\fuglob$ in
$\fieldres(\pointy_1),\fieldres(\pointy_2) \supseteq \fieldres(\point)$
lie in
$\fieldres(\point)$
{\rm(}or in their perfect closures{\rm )}
and are identical.
\end{lemma}
\begin{proof}
We have the commutative diagram
$$ \xymatrix{ \Gamma(\spax, \O_\topo) \ar[r] \ar[d] &  \Gamma(\spay, \O_\topo) \ar[d]   \\
\Gamma(\point, \O_\topo)=\fieldres(\point)  \ar[r] &
\Gamma(\pointy_1, \O_\topo)=\fieldres(\pointy_1),\,
\Gamma(\pointy_2, \O_\topo)=\fieldres(\pointy_2) \, .}$$
In this diagram
the element $\fuglob \in \Gamma(\spay, \O_\spay ) \to
\Gamma(\spay,\O_\topo)$ comes from an element $\tilde{\fuglob} \in
\Gamma(\spax, \O_\topo)$, since
$\fuglob \in \Gamma(\spay,\O_\spay )$
is compatible over $\spax$ (both rings of global sections are
colimits of compatible elements over coverings). Hence the image of
$\fuglob$ in $\fieldres(\pointy_1)$ and $\fieldres(\pointy_2)$ comes
from
$\tilde{\fuglob}$ via $\fieldres(\point)$. The same argument
works if the fields are replaced by their perfect closures.
\end{proof}
}

The next lemma describes a situation where a lot of ``new'' global
elements arise in the Grothendieck topology. We will use this
argument in Corollary \ref{constructibleglobalsectioncorollary},
Theorem \ref{frobeniusglobal}, Proposition \ref{zariskifilterproperties}
and Proposition \ref{ratliffglobalstalk}.

\begin{lemma}
\label{antistrictlemma}
Suppose that $\spay \to \spax$ is a cover in a Grothendieck topology
on $\spax$ such that the diagonal
$\diag: \spay \to \spay \times_\spax \spay$
is a cover.
Let $\sheaf $ denote a sheaf of abelian groups on
$\spax_\topo$.
Then the restriction homomorphism
$\Gamma(\spax,\sheaf ) \to \Gamma(\spay,\sheaf )$
is an isomorphism.
This holds in particular if one {\rm(}then both{\rm)} projection
$\spay \timesx \spay \to \spay$ is an isomorphism.
\end{lemma}
\begin{proof}
We look at the exact complex
$$\Gamma(\spax, \sheaf ) \lto \Gamma(\spay,\sheaf )
\stackrel{\rest_1 -\rest_2}{\lto} \Gamma(\spay \timesx \spay,\sheaf )\, $$
and have to show that every element
$\elem \in \Gamma(\spay,\sheaf)$ maps to $0$ (the first mapping is injective anyway).
The restriction homomorphism
$\rest(\diag) :\Gamma(\spay \timesx \spay,\sheaf ) \to \Gamma(\spay, \sheaf)$
comes from a cover and
$\rest(\diag) (\rest_1 (\elem) - \rest_2(\elem))=0$, so
$\rest_1 (\elem) -\rest_2(\elem)=0$
as well.

If one projection is an isomorphism, then also the diagonal
$\diag :\spay \to \spay \timesx \spay$ is an isomorphism, hence a cover. Then both
projections equal
$ \diag^{-1}$.
\end{proof}

\subsection{Filters and stalks}

\label{filterstalksubsection}

\ren{\point}{{x}}

\ren{\catind}{{\Lambda}}

\

\medskip
In the \'{e}tale topology, the \emph{Henselization} of the local
ring $\O_\point$ of a point $\point \in \spax$ with separably closed
residue class field can be recovered by taken a suitable colimit
over all \emph{\'{e}tale neighborhoods} of that point
(\cite[D\'{e}finition VIII.4.3]{SGA4}, \cite[D\'{e}finition 18.6.5]{EGAIV}, \cite[Remark I.4.11]{milne}).
We will describe here how to get suitable colimits
over reasonably filtered systems in $\catopenspax $ and in which
sense it is possible to bring all coverings together into a single
absolute cover.
Motivating examples for this are that all finite
extensions of a domain are inside its absolute integral closure
\cite{hochsterhunekeinfinitebig},
or all Frobenius powers are inside the perfect closure.

\synchronize{\indec}{\indmo}
\synchronize{\inded}{\indmu}

\begin{definition}
\label{filterdefinition}
Let $\indexcat$ denote a small quasicofiltered category
(Definition \ref{cofilteredcatdef})
with a final object and let $\catopen_\spax$ be a Grothendieck
topology on a scheme $\spax$.
Then we call a covariant functor
$\filt:\indexcat \to \catc$
which respects the final object a \emph{quasifilter} in
$\catopen_\spax$.
If moreover $\indexcat$ is cofiltered, then we
call $\filt$ a \emph{filter}.
We write often $\openfil_\indmo  $ or
$\openfil_\indmo \to \spax$ instead of $\filt(\indmo )$.
\end{definition}


\ren{\presh}{{\shP}}

Suppose that $\presh$ is a presheaf on a Grothendieck topology
$\catopen_\spax$ with values in a category $\catab$ and that $\filt:
\indexcat \to \catopen_\spax$ is a quasifilter. Then we get also a
contravariant (quasifiltered) diagram in $\catab$ by looking at
$$\presh \circ \filt : \indexcat \lto \catab, \,\, \indm
\longmapsto \Gamma(\openfil_\indm  , \presh) \, .$$
If $\rest: \indm  \to \indr $ belongs to $\indexcat$, then we will
denote the induced restriction morphism
$\Gamma(\obu_{\indr}, \presh) \to \Gamma(\obu_{\indm },\presh)$
also by $\rest$ or
$\rest^{\indr }_\indm  $. If $\catab$ is such that
colimits over (certain) diagrams exist, then we denote the colimit
by
$$ \presh_\filt
:= \colim_{\indm   \in \indexcat} \Gamma(\openfil_\indm ,\presh)\,$$
and call it the \emph{stalk} of $\presh$ in $\filt$.

\ifthenelse{\boolean{book}}{

Hence, if $\presh$ is a presheaf of $\ring$-modules and $F$ is a
quasifilter, then the colimit (the stalk) exists as an
$\ring$-module, every element is represented by an element $s \in
\Gamma(U_\indmo  , \presh)$, and two elements $s \in \Gamma(U_\indmo
, \presh)$ and $s' \in \Gamma(U_\indmu  , \presh)$ represent the
same element in the stalk if and only if there exist $\psi:\delta
\to \indmo   $ and $\psi':\delta \to \indmu  $ in $\indexcat$ such
that $\psi(s)=\psi'(s')$ in $\Gamma(U_\delta, \presh)$. In
particular, two elements $s,s' \in \Gamma(U_\indmo  ,\presh)$ define
the same element in the stalk if and only if there exists $\psi_1,
\psi_2 : \delta \to \indmo  $ such that $\psi_1(s) = \psi_2(s')$ in
$\Gamma(U_\delta,\presh)$. If moreover $F$ is a filter, then these
two morphisms $\psi_1$ and $\psi_2$ can be united to one morphism
$\varphi: k \to \delta $ such that $(\varphi \circ \psi_1)(s)=
(\varphi \circ \psi_2)(s')$. So in the case of a filter two elements
in $\Gamma(U_\delta,\presh)$ with the same image in the stalk are
identified by one homomorphism.

\begin{bookremark}
Another difference between quasifilters and filters is with respect
to a sheaf of rings. In both cases the colimit exists, however for
quasifilters it needs not be the same as the corresponding colimit
of sets or $\ring$-modules. For if only one object, say a ring
$\alg$ is given, and two morphisms $\psi,\id: A \to A$, then the
module-colimit is the $\ring$-module $\alg$ modulo the subgroup (or
$\ring$-submodule) $(a-\psi(a), a \in A)$, but the colimit as a ring
is the residue class ring of $\alg$ modulo the ideal generated by
these elements, see \cite[Proposition A6.7]{eisenbud}. For a filter
however the set-theoretic colimit is also the ring-theoretic
colimit, see \cite[Appendix A, Corollary 2]{milne}. The point is
that for $\psi_1,\psi_2: B \to B'$ (rings), the two images
$\psi_1(b)$ and $\psi_2(b)$ come together again under $\varphi_1,
\varphi_2: B' \to C$, so they define the same element in the stalk;
however an element in the ideal generated by $\psi_1(b)-\psi_2(b)$
does not necessarily map to $0$ somewhere and need not be $0$ in the
colimit of $\ring$-modules, but it must be $0$ in the colimit of
rings.
\end{bookremark}}{}

\ifthenelse{\boolean{book}}{

\subsection{Quasifilter given by a fixation}
\label{fixationsubsection}

\

\medskip}{}

\ifthenelse{\boolean{book}}{
 If $U_{\indmo  ,\topo} \to
U_{\indec ,\zar}$ is a site morphism for every $\indmo   \in
\indexcat$, then we have morphisms
$$\Gamma(U_\indmo  ,\O_{U_\indmo  }) \lto \Gamma(U_\indmo  ,\O_\topo)\lto S \, ,$$
which are compatible for every morphism $U_\indmo   \to U_\indmu  $
indexed by $F$. This gives for $U_{\indmo  }$ affine compatible
scheme morphisms $Z=\Spec \rings  \to U_\indmo  $ for every $\indmo
\in \indexcat$.

On the other hand, if we fix a scheme $Z \to \spax$ and consider all
scheme morphisms $Z \to U$, $U \in \catopenspax    $, then we get a
quasifilter, as the following construction shows. }{}

\begin{construction}
\label{fixationconstruction}
Let $\catopen_\spax$ be a Grothendieck
topology on $\spax$ and fix a scheme morphism $\spafix \to \spax$
(we call this a \emph{fixation}).
This defines a quasifilter in the following way:
Let the indexing category
$\indexcat=\indexcat_\spafix=\indexcat_{\spafix \to \spax}$
be given by (isomorphism types of) all commutative diagrams
$$ \xymatrix{  &   \obu \ar[d] \\
 \spafix \ar[ur] \ar[r] & \spax}$$
as objects, where $\obu \to \spax$ is in $\catopenspax$
(so the objects are the $\spafix$-pointed open sets in $\catopenspax$), and
by all the morphisms
$\obu \to \obusec$ in $\catopenspax    $ which
commute with the fixations $\spafix \to \obu$ and
$\spafix \to \obusec$.
The final object is given by
$\spafix \to \spax \stackrel{\id}{\to}\spax$.
If two diagrams $\spafix \to \obu \to \spax$ and
$\spafix \to \obusec \to\spax$ are given over
$\spafix \to \obu_0 \to \spax$
(that is, $\obu, \obusec \to \obu_0$ fixing $\spafix$),
then they come together in
$\spafix \to \obu \times_{\obu_0} \obusec \to \spax$
(arbitrary products exist in $\catopenspax$).
Hence $\indexcat$ is a quasicofiltered category.
The corresponding quasifilter is given by sending
$$ \indexcat \lto \catopenspax ,\, \, (\spafix \to \obu \to \spax) \longmapsto
(\obu \to \spax) \, ,$$
i.e. by forgetting $\spafix \to \obu$.

If $\catopen_\spax$ is the full category of finite type over $\spax$ (a
big site) and the topology is affine, then an affine fixation
$\spafix =\Spec \fixata$ yields even a filter. If $\armr ,\armrsec :
\obu \rightrightarrows \obur$ are given corresponding to $\armr
,\armrsec : \algb' \rightrightarrows \algb$, then
$\algb \tensor_{\algb'}\algb/
(\armr (\elec) \tensor 1 - 1 \tensor \armrsec (\elec), \elec \in \algb' )$
gives the equalizer which belongs then to the category.
\end{construction}

\ifthenelse{\boolean{book}}{
\begin{bookremark}
Note that for every $U_\indmo  $, $\indmo   \in \indexcat$, we have
a fixed morphism $Z \to U_\indmo  $, but an object $U \in
\catopenspax    $ may occur for several indices with different $Z
\to U$. The morphisms $Z \to U_\indmo  $ do all commute with respect
to the morphisms $U_\indmo   \to U_\indmu  $ which are indexed by
$\indmo \to \indmu  $ in $\indexcat$, but of course not with respect
to all morphisms $U_\indmo   \to U_\indmu  $ in $\catopenspax    $.
\end{bookremark}}{}

\ifthenelse{\boolean{book}}{
\begin{bookremark}
The construction above has several variants. For example one may
impose restrictions on the indexing objects $Z \to U \to \spax$ and
also on the morphisms $U \to U'$. For example for $Z= \Spec \field
\to \spaxeqspecring$ corresponding to an inclusion $R \subseteq
\field $ in an algebraically closed field one may restrict to
inclusions $\ring \to \algcat\to K$, $U= \Spec \algcat$. In this
case it is convenient to replace the tensor product $\algcat
\tensorr A'$ by the $\ring$-algebra $R[\algcat,\algcat'] \subseteq
\field $ and to impose that this algebra belongs to the Grothendieck
category.
\end{bookremark}}{}

\ifthenelse{\boolean{book}}{
\begin{bookremark}
Let $\catopen_\ring$ be given by affine single coverings and
consider the quasifilter in $\catopen_\ring$ given by a fixation
$\ring \to A$. Then we have a ring homomorphism $\colim _{\indmo
\in \indexcat_A} B_\indmo   \to A$, but in general there will be no
ring homomorphism $\colim_{\indmo   \in \indexcat_A} \Gamma(B_\indmo
, \O) \to A$.
\end{bookremark}}{}

\begin{example}
\label{geometricpoint}
An important case is when $\spafix= \Spec
\field$ is the affine spectrum of an algebraically closed field
$\field$, in which case $\spafix \to \spax$ is called a
\emph{geometric point}. Goodwillie and Lichtenbaum call in
\cite[Section 2]{goodwillielichtenbaum} a scheme $\glpoint$ a
\emph{point} in a Grothendieck topology if for every covering
$\obcov_\indcov \to \obcovp$, $\indcovinset$, and every morphism
$\glpoint \to \obcovp$ there exists a lifting $\glpoint \to
\obcov_\indcov$ for some $\indcov$. If the topology is such that all
coverings are surjective and of finite type, then every geometric
point is also a point in this sense.
\end{example}

\ifthenelse{\boolean{book}}{ Hier geht covering und nicht covering
durcheinander:
\begin{example}
\label{zariskifilterquasifilter} Let $\spaxeqspecring$ be an affine
scheme with the Zariski topology. Then the fixation given by
$\specring_\fop \to \specring$ does in general yield only a
quasifilter, not a filter. Let $\fop \in U \subset \spax$, $U$ open.
We consider the covers $\spay= U \uplus \spax$ and $\spay'=U \uplus
\spax \uplus \spax$ with the fixation given by $\fop \in U$. Then
the two morphisms $\armr   ,\armrsec     :\spay \to \spay'$ given by
the identity on $U$ and the two possible embeddings of $\spax$ into
$\spax \uplus \spax$ are not unified by any $Z \to \spay$. For a
point $\point \in \spax$, $\point \not\in U$, there exists a point
$\point' \in Z$ (over $\point$), which must be sent to $\point \in
\spax \subset \spay$, and $\armr $ and
$\armrsec    $ have different values there.
\end{example}

auch hier:

\begin{bookexample}
\label{quasinotfilterexample}
Let $\spax= \AA^1$ be the affine line
endowed with a Grothen\-dieck topology where every object $\spay \to
\spax$ is a surjection (or flat) and where the coverings are given
by single surjections (or faithfully flat morphisms). Consider the
fixation given by the zero point $0 \to \AA^1$ and the corresponding
quasifilter (Construction \ref{fixationconstruction}). Let $\indmo $
be the identity $0 \to \AA^1 \to \AA^1$ and let $\indmu  $ be $0 \to
(0,0) \in \AA^2 \stackrel{p_1}{\to} \AA^1$, and consider the two
morphism $\armr   , \armrsec    :\indmo   \to \indmu  $ given by the
diagonal and the antidiagonal.

\vspace{.6cm}

\unitlength1cm
\begin{picture}(12,5)


\put(8,3){\line(1,0){2}} \put(8,3){\line(0,1){2}}
\put(8,5){\line(1,0){2}} \put(10.,3){\line(0,1){2}}

\put(9,4){\circle*{0.1}}

\put(3,4){\line(1,0){2}}

\put(4,4){\circle*{0.1}}

\put(4.8,3.6){$\spax$}

\put(5.7,4.1){\vector(3,1){1.6}}

\put(5.7,3.9){\vector(3,-1){1.6}}

\put(6.3,4.55){$\armrsec    $}

\put(6.3,3.3){$\armr   $}

\put(8,3){\line(1,1){2}}

\put(8,5){\line(1,-1){2}}

\put(9,2.7){\vector(0,-1){.8}}

\put(9.1,2.2){$\indmu  $}

\put(9,1.6){\circle*{0.1}}

\put(8,1.6){\line(1,0){2}}

\put(9.8,1.2){$\spax$}

\end{picture}

\vspace{-0.8cm}

Then there does not exist an object (the ``$k$'' in Definition
\ref{cofilteredcatdef}) $0 \to \spay \to \spax$ which unifies $\armr
$ and $\armrsec    $, because the only possible morphism
$\varphi:\spay \to\AA^1$ is the surjective structure morphism, and
so $ \armr \circ \varphi \neq \armrsec     \circ \varphi$. Hence
this quasifilter is not a filter.
\end{bookexample}

Diese ganzen Beispiele machen nur Sinn, wenn nur Produkte ueber der
Basis wieder dazu gehoeren!

\begin{bookexample}
Let $\ring \to \field    $ a be homomorphism to an algebraically
closed field, giving rise to an absolute quasifilter in the
primitive or in the faithfully flat or in the surjective topology.
This is in general not a filter, the unifying property does not hold
for the polynomial ring (in one variable). Consider in this
quasifilter the two (indexing) objects given by the polynomial ring
$R[T]$ and the fixations $T \mapsto a ,a' \in K$. A homomorphism
$R[T] \to R[T]$, $T \mapsto F$ belongs to the indexing category if
and only if $F(a)=a'$. Let two such homomorphisms $ \armr   ,
\armrsec    : R[T] \to R[T]$ be given by $T \mapsto F_1, F_2$. A
morphism $\varphi :R[T] \to S$, $ T \mapsto s$ ($S$ an
$\ring$-algebra with a $\field$-fixation), is unifying if and only
if $F_1(s)=F_2(s) \in S$. Hence $\varphi$ must factor through
$R[T]/(F_1 - F_2)$. If e.g. $a=a' =0$ and $F_1=T$, $F_2=T+g$, where
$g \in \ring$ mapping to $0 \in K$, but not nilpotent, then this is
$\ring/(g)[T]$ and cannot define a surjective morphism.
\end{bookexample}}{}

\ifthenelse{\boolean{book}}{
\begin{bookremark}
The easiest way to get via the construction a filter is to make $Z$
big enough such that there exists at most one morphism between $Z\to
U \to \spax$ and $Z \to U' \to \spax$ fixing the ends. Finding such
a $Z$ is however in some sense the same as looking for a reasonable
stalk bring all the different coverings together.

A reasonable way for $\ring$ a domain is probably to look at $\ring
\to L$, $L$ a big field, and to consider only $\ring$-subalgebras $R
\subseteq \algcat \subseteq L$ such that $\Spec \algcat \to
\specring$ is a covering. Then there exists at most one homomorphism
$\algcat \to \algcat'$ (an inclusion) fixing $L$ (and belonging to
$\catopen$). Two such subalgebras $\algcat,\algcat' \subseteq L$
come then together in $R[\algcat,\algcat'] \subseteq L$ (make sure
that this is a covering). Such a filter is absolute (see below) if
for every covering $\ring \to \algcatb$ there exists $\ring \to
\algcatb \to \algcata \subseteq L$, where $\ring \to \algcata$ is a
covering and $\algcatb \to \algcata$ is in $\catopen$.
\end{bookremark}}{}

\begin{example}
\label{stalkquasifilternotfilter}
In a Grothendieck topology the
quasifilter given by a fixation is in general not a filter. Consider
the category of finite non-empty sets where only surjective
morphisms are allowed and are covers. This is the same as the affine
single-handed covering Zariski topology on $\spax=\Spec \field$,
$\field$ a field. The identity $\spax \to \spax$ gives a fixation
and an (absolute) quasifilter, the objects in this quasifilter are
the pointed sets. Let $\obr$ consist of the two points $P_0$ (fixed)
and $P_1$ and let $\obm$ consist of three points $Q_0$ (fixed),
$Q_1$ and $Q_2$. Let $\armr $ and $\armrsec $ be given by $\armr
(Q_0)= \armrsec (Q_0) = P_0$, $\armr (Q_1)= \armrsec (Q_1) = P_1$,
and $\armr (Q_2)= P_0$, $\armrsec (Q_2)= P_1$. Then there does not
exist a surjective map $\obl \to \obm$ which unifies the two
morphisms.

\ren{\setfin}{\spay} \ren{\setfinsec}{{\bar{\spay}}}
\ren{\setfinl}{{\tilde{\spay}}}
\ren{\elef}{{f}} \synchronize{\setfintri}{\setfinl}

Consider the sheaf
$\obm \to \field^\obm= \Gamma(\obm, \O_\topo)$
and the diagram induced by the given quasifilter. 
The colimit of this system as a diagram of sets is
$(\biguplus_{\obm} \field^\obm)/ \sim$, where $\sim$ is the
equivalence relation  in which $ \tupf \in \field^\obm$ and
$\tupfsec \in \field^\obsec $ are equivalent if and only if there
exists a set $\obl$ and surjective mappings
$\arl:\obl \to \obm $,
$\arlsec:\obl \to \obsec $ with
$\arl^*(\tupf) = \arlsec^*(\tupfsec)$.
This is true if and only if $\tupf$ and
$\tupfsec$ have the same values.
Therefore the equivalence classes are given by
finite unordered tupels of different elements in $\field$.

This is however not a group, since it is not possible to define a
group operation. Two elements $\tupf \in \field^\obm$ and $\tupfsec
\in \field^\obsec$ come together in some $\field^\obsecsec$, but
there are several choices, so their sum is not well defined. For
each summand it is possible to bring the different choices together,
but not for both summands simultaneously. This example shows that
Proposition \cite[Proposition A6.3]{eisenbud} is wrong. The colimit
in the category of groups is just $\field$ given by the evaluation
at the fixed point.
\end{example}

\ifthenelse{\boolean{book}}{
\begin{bookexample}
Another example for a quasifiltered but not a filtered category is
given by the category of one object $\point$ together with a group
of automorphisms $\autg$. An operation of $\autg$ as group
automorphisms on an abelian group $H$ gives a covariant functor. The
colimit as sets is the space of orbits of $H$ under the action,
which is not a group.
\end{bookexample}}{}

\subsection{Irreducible filter}
\label{irreduciblesubsection}

\

\begin{definition}
\label{irreduciblefilterdef}
We call a (quasi)filter
$\filt:\indexcat \to \catopenspax $
in a Grothendieck topology
$\catopenspax$ on a scheme $\spax$ \emph{irreducible} if for every
$\indmo \in \indexcat$ and every covering
$\obcov_\indcov \to \obcovp_\indm $, $\indcovinset$,
there exists $\indcovinset$ and
$(\indl \to \indm )\in \indexcat$ such that
$\filt(\indl \to \indm) = \obu_\indl \to \obu_\indm $
factors through
$\obcov_\indcov \to \obcovp_\indm  $ via some
$(\obu_\indl \to \obcov_\indcov) \in \catopenspax$.
\end{definition}

\begin{example}
In an usual topological space, the filter of open neighborhoods of a
point is irreducible, because if $\bigcup_{\indcovinset}
\obcov_\indcov$ is a neighborhood of a point $\point$, then
$\obcov_\indcov$ must be a neighborhood of $\point$ for at least one
$\indcov$.
\end{example}

\begin{example}
\label{geometricpointirreducible}
If a Grothendieck topology over a
scheme $\spax$ has the property that every covering is also
set-theoretically a covering, then the quasifilter corresponding to
a geometric point
$\spafix = \Spec \field \to \spax$ is irreducible.
For if $\Spec \field \to \obu \to \spax$ is a given neighborhood of
$\point$ and $\covering$, $\indcovinset$, is a covering, then there
exists $\indcov$ and
$\Spec \field \to \obcov_\indcov$ compatible
with the given neighborhood.
\end{example}

The following lemma is a variant of the fact that for a point in a
topological space the stalk of a presheaf and of its sheafification
are the same.

\begin{lemma}
\label{stalkpresheafsheaf}
Let $\catopenspax $ denote a Grothendieck
topology on $\spax$, let
$\filt: \indexcat \to \catopenspax $ denote
an irreducible filter. Let $\presh$ denote a presheaf of abelian
groups with sheafification
${\presh}^{\sheafify}$.
Then
$\presh_\filt=({\presh}^{\sheafify})_\filt$.
\end{lemma}
\begin{proof}
We have a natural homomorphism
$$ \presh_\filt
= \colim_{\indm \in \indexcat} \Gamma(\openfil_\indm , \presh)
\lto \colim_{\indm \in \indexcat}\Gamma(\openfil_\indm,{\presh}^\sheafify)
= ({\presh}^{\sheafify})_\filt \, .$$
Suppose
that $\elem \in \presh_\filt$, represented by $\elem \in
\Gamma(\openfil_\indm , \presh)$, maps to $0$ on the right. We may then
assume that $\elem=0 $ in
$\Gamma(\openfil_\indm  ,\presh^\sheafify)$
(because the colimit over a filter is the set theoretical colimit).
This means that there exists a covering
$\obcov_\indcov \to \obcovp_\indm $,
$\indcovinset$, such that the restrictions
$\rest^{\openfil_\indm}_{\obcov_\indcov}(\elem )
\in \Gamma(\obcov_\indcov,\presh)$ are zero.
Since $(\openfil_\indl \to
\obcov_\indcov \to \openfil_\indm ) \in \filt$
for some $\indcovinset$, it follows that $\elem=0$ in
$\Gamma(\openfil_\indl, \presh)$ and hence in
$\presh_\filt$.

Now suppose that
$\elem \in ({\presh}^{\sheafify})_\filt$ is
represented by
$\elem \in \Gamma(\openfil_\indm, \presh^\sheafify)$,
$\indm \in \indexcat$.
This means that there exists a covering
$\obcov_\indcov \to \obcovp_\indm  $, $\indcovinset$, and elements
$\elem_\indcov \in \Gamma(\obcov_\indcov,\presh)$ which are
compatible in the sense that
$\elem_\indcov = \elem_\indcovsec$ in
$\Gamma(\obcov_\indcov \times_{\obu_{\indm}} \obcov_\indcovsec,\presh_1)$.
Let
$(\openfil_\indl \to \obcov_\indcov \to \openfil_\indm  )\in \filt$.
Then $\elem = \elem_\indcov$ on $\obcov_\indcov$ and
also on $\openfil_\indl $,
hence $\elem$ comes from the left.
\end{proof}

\begin{corollary}
\label{irreduciblefiltermodulestalk}
Let $\spaxeqspecring$ be
endowed with an affine Grothendieck topology, let
$\filt:\indexcat \to \catopenspax$ denote an irreducible filter with
structure stalk $\O_\filt$.
Then for an $\ring$-module $\modul$ we have
$(\modul_\topo)_\filt = \modul \tensorr \O_\filt$.
\end{corollary}
\begin{proof}
We have by Lemma \ref{stalkpresheafsheaf} and since tensor products
commute with colimits the identities
(setting $\openfil_\indm = \Spec \ring_\indm  $)
\begin{eqnarray*}
(\modul_\topo)_\filt
=\colim_{\indm \in \indexcat}
\Gamma(\openfil_\indm , \modul_\topo)
&=& \colim_{\indm   \in\indexcat}  \modul \tensorr \ring_\indm
= \modul \tensorr (\colim_{\indm \in \indexcat}\ring_\indm ) \cr
&=& \modul \tensorr (\colim_{\indm \in \indexcat}
\Gamma(\openfil_\indm ,\O_\topo))
= \modul \tensorr \O_\filt \,.
\end{eqnarray*}
\end{proof}

\ifthenelse{\boolean{extra}}{
\begin{lemma}
Tensor products and colimits commute, that is, $M \tensorr
(\colim_\Lambda N_\lambda) = \colim_\Lambda ( M \tensorr
N_\lambda)$.
\end{lemma}
\begin{proof}
The homomorphisms $M \tensorr N_\lambda \to M \tensorr
\colim_\Lambda N_\lambda$ induce a homomorphism $$ \colim_\Lambda (
M \tensorr N_\lambda) \lto M \tensorr \colim_\Lambda N_\lambda \,
.$$ The surjectivity is clear. For injectivity we go back to the
definition of tensor product as generated by symbols module
relations as is \ref{tensorhilfslemma}. In fact this Lemma there
should be formulated more generally that the tensor product commute
with limits in every place. Then that an element $t \in M \tensorr
N_\lambda$ which becomes $0$ on the right means that a certain
linear combination of symbols belong to a certain submodul, and this
is then also true in some $M \tensorr N_{\lambda'}$.
\end{proof}

We take the following definition of
\cite[Section 2]{goodwillielichtenbaum}

\begin{definition}
\label{pointdef}
Let a Grothendieck topology be given on a category
of schemes. We say that a scheme $\glpoint$ is a \emph{point} in the
topology if for every covering $\obcovindic \to \obu$,
$\indcovinset$, and every morphism $\glpoint \to \obu$ there exists a
lifting $\glpoint \to \obcovindic$ for some $\indcov$.
\end{definition}

\begin{example}
If the topology is such that all coverings are surjective and of
finite type, then $\Spec \field$ is a point for every algebraically
closed field.
\end{example}

\begin{remark}
A scheme $\glpoint$ is a point if and only if for all fixations
$\glpoint \to \obu$ the corresponding quasifilter given by
$\indexcat_{\glpoint \to \obu}$ is irreducible.
\end{remark}

Suppose now that $Z$ is also endowed with a Grothendieck topology
$\catopen_Z$ and that every scheme morphism $Z \to U$ induces a site
morphism. Then we get in particular ring homomorphisms
$$\Gamma( U_\indmo   , \O_\topo) \lto \Gamma(Z_\topo , \O_\topo ) \, $$
which are compatible with $U_\indmo   \to U_\indmu  $, and therefore
we get a ring homomorphism
$$  \colim_{\indmo   \in \indexcat} \Gamma(U_\indmo  , \O_\topo)
\lto  \Gamma(Z_\topo , \O_\topo ) \, .$$

An important case is when $\spafix = \Spec \field$ is the affine
spectrum of an algebraically closed field $\field$, in which case
$\spafix \to \spax$ is called a \emph{geometric point}. Suppose that
we have such a geometric point and that $\Gamma((\Spec
\field)_\topo, \O_\topo) = \field$. In this case we get a ring
homomorphism
$$ \colim_{\indmo   \in \indexcat} \Gamma(U_\indmo  , \O_\topo) \lto \field \, .$$
}{}

\begin{proposition}
\label{filtertest}
Let $\catopen_\ring$ denote a Grothendieck
topology on $\specring$, let $\ring \to \fixata$ be a ring
homomorphism defining a quasifilter $\filt$ by the small indexing
category
$$\indexcat_{\ring \to \fixata}
= \{ [\ring \to \algcat \to \fixata], (\ring \to\algcat ) \in \catopen_\ring \} \,$$
{\rm(}here $[\, ]$ denotes the isomorphy class{\rm)}.
Suppose that $\filt$ is irreducible and let $\O_\filt$ denote the stalk of
the structure sheaf in $\filt$. Suppose in the first two statements
that $\filt$ is a filter. Then the following holds.

\numiii
\begin{enumerate}

\item
There exists a natural ring homomorphism $\O_\filt \to \fixata$.

\item
Suppose that $\ring \to \ring[\vart]$ is a cover in the topology.
Then $\O_\filt \to \fixata$ is a surjection.

\item
Suppose now that for each
$\ring \to \algcat \to \fixata$ in
$\indexcat_{\ring \to \fixata}$ the image
$\bar{\algcat}$ of
$\algcat$ in $\fixata$ has the property that the induced
homomorphism
$\algcat \to \bar{\algcat}$ belongs also to
$\catopen_\ring$.
Then $\filt$ is a filter and the stalk $\O_\filt$
is a subring of $\fixata$.
\end{enumerate}
\end{proposition}
\begin{proof}
(i).
We always have a ring homomorphism
$\colim_{\indm \in \indexcat} \algcat_\indm \to \fixata$,
and
$\O_\filt = \colim_{\indm \in \indexcat} \algcat_\indm$
by Lemma \ref{stalkpresheafsheaf} if
$\filt$ is an irreducible filter.
(ii).
For $\elea \in \fixata$ consider the index $\indm$ given by
$$\ring \lto \ring[\vart]
\lto \fixata, \, \vart \longmapsto \elea \, .$$
Then the image of $\vart \in \algcat_\indm$ in
the colimit
$\colim_{\indm \in \indexcat} \algcat_\indm$
(which is $\O_\filt$) maps to $\elea$.

We have to verify property (ii) in Definition \ref{cofilteredcatdef}
for the indexing category $\indexcat_\fixata$.
So assume that the first two rows of the following diagram commute.

\vspace{8.31cm}

\begin{picture}(20,10)

\unitlength 1cm

\put(4.28,8){$\ring$}
\put(4.72,8.15){\vector(1,0){.83}}
\put(5.7,8){$\algcat'$}
\put(6.23,8.15){\vector(1,0){.83}}
\put(7.2,8){$\fixata$}
\put(9.82,8.05){$ \indr$}

\put(4.42,7.85){\vector(0,-1){.83}}
\put(4.10,7.25){$^{=}$}
\put(5.41,7.20){$^{\armr}$}
\put(5.92,7.85){\vector(0,-1){.83}}
\put(5.78,7.85){\vector(0,-1){.83}}
\put(6.03,7.20){$^{\armrsec}$}
\put(7.36,7.85){\vector(0,-1){.83}}
\put(7.5,7.25){$^{=}$}
\put(9.85,7.02){\vector(0,1){.83}}
\put(9.99,7.02){\vector(0,1){.83}}

\put(4.28,6.59){$\ring$}
\put(4.72,6.74){\vector(1,0){.83}}
\put(5.7,6.59){$\algcat$}
\put(6.23,6.74){\vector(1,0){.83}}
\put(7.2,6.59){$\fixata$}
\put(9.77,6.615){$ \indm$}

\put(4.42,6.44){\vector(0,-1){.83}}
\put(4.10,5.84){$^{=}$}
\put(5.88,6.44){\vector(0,-1){.83}}
\put(6.02,5.79){$^{\homring}$}
\put(7.36,6.44){\vector(0,-1){.83}}
\put(7.5,5.84){$^{=}$}
\put(9.93,5.61){\vector(0,1){.83}}

\put(4.28,5.18){$\ring$}
\put(4.72,5.33){\vector(1,0){.83}}
\put(5.7,5.16){$\overline{\algcat}$}
\put(6.23,5.33){\vector(1,0){.83}}
\put(7.2,5.18){$\fixata$}
\put(9.82,5.23){$ \indl$}

\end{picture}

\vspace{-4.93cm}

%

\noindent
This can be commutatively extended to the third row, and the
homomorphisms
$\homring \circ \armr , \homring \circ \armrsec:
\algcat' \to \bar{\algcat}$
are identical, since they are identical as homomorphisms to
$\fixata$.

Suppose that $\fuf \in \O_\filt$ maps to $0$ in $\fixata$,
represented by $\fuf \in \algcat_\indm  $.
Then the image of $\fuf$ in
$\algcat_\indm \to \overline{\algcat_\indm  } \to \fixata$
is $0$ and so $\fuf=0$ in $\O_\filt$.
\end{proof}

\ifthenelse{\boolean{extra}}{
\begin{example}
Let $\ring$ be a ring and consider the Grothendieck topology
$\catopen_\ring$ generated by $\ring \to \rings   =
R[\epsilon]/(\epsilon^2)$. So the object are $$\rings_k=S \tensorr
\ldots \tensorr S \cong R[\epsilon_1 \komdots
\epsilon_k]/(\epsilon_1^2 \komdots \epsilon_k^2) \, .$$ (but various
morphisms). Consider the absolute stalk given by the identity $\ring
\to \ring$ (filter?). Suppose that $\ring$ is reduced, so that there
exists only one $\ring$-algebra homomorphism $\rings_k \to \ring$.
All these homomorphisms (including the identity) factor through $S$,
and this is the stalk.
\end{example}}{}

\subsection{Absolute filters and absolute stalks}
\label{absolutestalksubsection}

\

\medskip
With absolute filters and absolute stalks we want to introduce a
concept which is supposed to control all covers in a single-handed
Grothendieck topology and contains in some sense all the information
of the topology. The absolute stalk contains in particular the
information about the induced closure operation
(Corollary \ref{topclosureabsolute}).
A leading and motivating example is the
role of the absolute integral closure $\ring^+$ of a domain which
controls the plus closure, see Part \ref{plussection}.

\begin{definition}
\label{absolutefilterdef}
We call a (quasi)filter
$\filt: \indexcat\to \catopenspax $ in a Grothendieck topology $\catopenspax$
over $\spax$ \emph{absolute} if
$\openfil_\indmo   \to \spax$ is a
cover for every $\indmo \in \indexcat$ and if for every cover
$\spay \to \openfil_\indmo$ in $\catopenspax $
($\indmo \in \indexcat$)
there exists
$\indl \in \indexcat$ and a factorization
$\openfil_\indl \to \spay \to \openfil_\indmo$
(so that the composition comes from $\indl \to \indmo$).
We call the stalk of a presheaf $\shS$ in an
absolute filter an \emph{absolute stalk}.
\end{definition}

\begin{remark}
\label{filterabsoluteremark}
An absolute quasifilter has in particular
the property that for every cover $\spay \to \spax$ there exists
$\openfil_\indm \to \spay \to \spax$, as
$\spax = \openfil_\indfin$ for the
final element $\indfin \in \indexcat$. We will also encounter the
condition that every morphism
$\openfil_\indl \to \openfil_\indm$
indicated by $\filt$ is a cover, and the stronger property that
every morphism
$\spay \to \openfil_\indm$ in $\catopen_\spax$ is a cover for every
$\indm \in \indexcat$.
See Proposition \ref{exactnessabsolute} and
Examples \ref{finabsoluteintegral} and \ref{finabsoluteextra} in the
finite topology, where these conditions occur naturally.
If $\filt$ is given in an affine single-handed Grothendieck topology
by a fixation $\ring \to \fixata$, then this quasifilter is absolute
if
$\ring \to \rings \to \fixata $ implies that $\ring \to \rings$
is a cover and that for given $\rings \to \fixata$ and a cover
$\rings \to \ringt$ there exists an $\rings$-algebra homomorphism
$\ringt \to \fixata$.
\end{remark}

\begin{definition}
\label{absolutecoverdef}
A scheme morphism $\spafix \to \spax$
(not necessarily in $\catopenspax $)
is called an \emph{absolute cover}
if a morphism $\obu \to \spax$ in $\catopenspax$ is a cover if
and only if there exists a factorization $\spafix \to \obu \to
\spax$.
\end{definition}

\begin{lemma}
\label{absoluteirreducible}
Let $\catopen_\spax$ be a single-handed
Grothendieck topology on a scheme $\spax$.
Then the following hold.

\numiii

\begin{enumerate}
\item
An absolute quasifilter $\filtcatindtocatopenspax$ is irreducible.

\item
An irreducible quasifilter $\filtcatindtocatopenspax$ with the
property that every $\openfil_\indm$, $\indmincat$, is a cover of
$\spax$,
is absolute.
\end{enumerate}
\end{lemma}
\begin{proof}
(i).
If $ \obu \to \openfil_\indm$ is a cover (there are only single
covers), $\indmincat$, then there exists $\indl \in \indexcat$ and
morphisms
$\openfil_\indl \to \obu \to \openfil_\indm$,
where
$(\obu_{\indl } \to \openfil_\indmo  )= \filt(\indl \to \indm)$.

(ii).
We have to check the second condition of an absolute filter,
so let $\spay \to \openfil_\indm$ be a cover. Then by irreducibility
there exists a $\indlincat$ and a factorization $\openfil_\indl \to
\spay \to \openfil_\indm$.
\end{proof}

Lemma \ref{absoluteirreducible}(i) is not true for non
single-handed Grothendieck topologies, see Proposition
\ref{completeproperties}(iv) for an example.

\ifthenelse{\boolean{book}}{
\begin{bookremark}
If $\filt: \indexcat \to \catopenspax    $ is an absolute filter, say in
an affine single-handed Grothendieck topology, then for every
$\ring$-algebra $\alg$ in the category there exists a homomorphism
$A \to A_{\indmo  } $ and hence also from $\alg$ to the absolute
stalk $H= \colim_{\indmo   \in \indexcat} A_\indmo   $, but these
homomorphisms are not uniquely determined. E.g. every
$\ring$-automorphism yields another map. We can always enrich the
given filter with stalk $H$ by the (quasi)filter given by the
fixation of $H$. This separates the different homomorphisms to $H$.
\end{bookremark}}{}

\begin{lemma}
\label{globalstalkinjective}
Let $\spax$ denote a scheme with a
Grothendieck topology $\catopenspax $ and an absolute filter
$\filt:\indexcat \to \catopenspax $. Let $\shS$ denote a sheaf of
abelian groups on $\spax_\topo$.
Then the natural homomorphism
$\Gamma(\spax_\topo, \shS) \to \shS_\filt$ is injective.
\end{lemma}
\begin{proof}
If $\elem \in \Gamma(\spax_\topo, \shS)$ becomes $0$ in the absolute
stalk, then it must be $0$ in
$\Gamma(\openfil_\indm   , \shS)$ for some
$\openfil_\indm   \to \spax$, $\indmo   \in \indexcat$.
Because
$\openfil_\indm \to \spax$ is a cover, $\elem$ must itself be $0$ by
the local nature of a sheaf.
\end{proof}

\begin{corollary}
\label{stalkpresheaf}
Let $\catopenspax\! $ denote a single-handed
Grothendieck topology on a scheme $\spax$. Let $\modul$ denote an
$\O_\spax $-module and let $\modul_{\topo}$ denote the corresponding
sheaf {\rm(}of abelian groups{\rm)} in the topology. Let
$\indexcat \to \catopenspax $ denote an absolute filter. Then the absolute
stalk is
$ \colim_{\indmo \in \indexcat}
\Gamma(\openfil_\indmo,\struto_\indmo^* (\modul))
= \colim_{\indmo \in \indexcat}
\Gamma(\openfil_\indmo,\modul_{\topo})$.
In the affine case this means that
$ \colim_{\indmo   \in \indexcat} \modul \tensorr \ring_\indmo
= \colim_{\indmo   \in \indexcat} \Gamma(\ring_\indmo,\modul_{\topo})$.
\end{corollary}
\begin{proof}
This follows from Lemma \ref{stalkpresheafsheaf},
Lemma \ref{absoluteirreducible}(i) and Proposition \ref{sheafzariskicompare}(iv).
\end{proof}

\ifthenelse{\boolean{book}}{
\begin{bookremark}
\label{filterfunctorial} A site morphism $\varphi: \spax' \to \spax$
and a quasifilter $\filt: \indexcat \to \catopenspax    $ induce by
pull-back a quasifilter $\filt': \indexcat \to \catopen_{\spax'}$,
$\indmo \mapsto \varphi^{-1} (U_\indmo  )$. If $\indexcat$ and
$\filt$ are given by a fixation $Z \to \spax$, then $Z'=Z
\times_\spax\spax' \to \spax'$ gives a fixation of $\spax'$, and a
corresponding indexing category $\indexcat_{Z'}$, which are related
by
$$\indexcat_Z \lto \indexcat_{Z'}, \,\,\,
\indmo  =(Z \to U \to \spax) \longmapsto (Z' \to U
\times_\spax\spax' \to \spax')=\indmu   \, .$$ If $\shS$ is a sheaf
on $\spax$ with pull-back $\shS'$ on $\spax'$, then the
homomorphisms $\Gamma(U_\indmo , \shS) \to \Gamma(U_\indmu  ,\shS')$
induce a stalk homomorphism
$$ \colim_{\indmo   \in \indexcat_Z} \Gamma(U_\indmo  , \shS)
\lto \colim_{\indmu   \in \indexcat_{Z'}} \Gamma(U_{\indmu  }, \shS') \, .$$
This holds in particular for the structure sheaf. Note however that
the pull-back of an absolute filter is in generally not absolute
anymore.
\end{bookremark}}{}

\begin{example}
\label{zariskiabsolutefilter}
Let $\spaxeqspecring$ be an affine
scheme with the single-handed Zariski topology (Example \ref{zariskisinglehanded}).
Then every global function
$\fug=(\fug_\indcov)$,
$\fug_\indcov \in \Gamma(\openzar_\indcov, \O)$,
over a cover
$\biguplus_\indcovinset \openzar_\indcov \to \spax$
will yield an element in an absolute stalk, so this cannot be an easy object. We look at
the fixation given by
$\spafix= \Spec (\prod_{\fop \in \spax}\ring_\fop)$.
An object in the corresponding quasifilter is given (after
refinement) by a covering
$\biguplus_\indcov \openzar_\indcov \to \spax$,
$\openzar_\indcov =\specring_{\fucov_\indcov}$, together
with a morphism from $\spafix$ to the covering.
Such a morphism is
given by a family of ring homomorphisms
$\prod_\indcovinset \ring_{\fucov_\indcov} \to \ring_\fop$ for every $\fop$, and they
are given by an $\indcov$th projection
(where $\indcov=\indcov(\fop)$ is such that $\fop \in D(\fucov_\indcov)$)
and followed by
$\ring_{\fucov_\indcov} \to \ring_\fop$. This
quasifilter is absolute, since for a homomorphism
$\prod_\indcovinset \ring_{\fucov_\indcov} \to \ring_\fop$
and a cover
$ \prod_\indcovinset \ring_{\fucov_\indcov} \to
 \prod_\indcovinsetsec \ring_{\fucovsec_\indcovsec}$
there exists
$\prod_\indcovinsetsec \ring_{\fucovsec_\indcovsec} \to \ring_\fop$.

Assume now that $\ring$ is a domain.
The image of
$\prod_\indcovinset \ring_{\fucov_\indcov}$ inside
$\prod \ring_\fop$ is
$\prod_{\indcov \in \setindcovsec} \ring_{\fucov_\indcov}$ for a subset
$\setindcovsec \subseteq\setindcov $
so that
$\biguplus _{\indcov \in \setindcovsec}\openzar_\indcov \to \spax$
is also cover.
So by Proposition \ref{filtertest}(iii) this quasifilter is a
filter. Hence the absolute stalk is given by
$$\{(\fuh_\fop),\fuh_\fop
\in \ring_\fop, \exists (\fucov_1 \comdots\fucov_\numcov)=(1)
\mbox{ and } \funum_\indcov, \expok_\indcov \mbox{ such that }
\fuh_\fop =\funum_\indcov/\fucov_\indcov^{\expok_\indcov}
\mbox{ for some } \indcov 
\} \, .$$
In this construction one can also restrict to the product over all local rings at maximal
ideals, but the resulting absolute stalk is different. For a local
ring $\ring$ the fixation by $\ring$ is an absolute filter with
absolute stalk $\ring$. If we fix by
$\ring \to \ring/ \fom_\ring$
then this defines the same filter with again $\ring$ as stalk (but
the condition in Proposition \ref{filtertest}(iii) is not
fulfilled). If say $\ring$ is a discrete valuation domain and we fix
by
$\ring \times \quotfield(\ring)$
(as described in the first paragraph), then the resulting absolute
stalk is $\ring \times \quotfield(\ring)$. 
\end{example}

\ifthenelse{\boolean{book}}{
\begin{bookexample}
\label{finitesetsexample} We illustrate how the quasifilters and
stalks depend on the morphisms allowed in $\catopen$. Let $\field$
denote a field and consider the topology where the coverings are
given by all finite non-empty disjoint unions of $P=\Spec \field $
(this is for $\field$ algebraically closed the finite topology, see
Section \ref{plussection}). So the covering objects are just $\Spec
(K^I)$, where $\ideal$ is a finite, non-empty set (this Grothendieck
topology may be identified with the category of non-empty finite
sets).

The morphisms in $\catopen_P$ will be either all morphisms or only
the spec-surjective morphisms. We consider first the quasifilter
which is given by the identical fixation $\Spec \field \to \Spec
\field $ (thus the indexing category consists of pointed sets
$(I,i)$, where $i \in I$). This is an absolute quasifilter. Whether
it is a filter depends on the morphisms which we allow. If we allow
all morphisms and if $\indmo  =(P \to I \to P)$ and $\indmu   = (P
\to I' \to P)$ and $\armr   , \armrsec    :I \to I'$ (commuting with
the fixations), then we can take $k=(P \to P \to P)$ and the natural
mappings $k \to \indmo  $ unifies $\armr   $ and $\psi_2$. If
however only surjective mappings are allowed, then for
$\ideal=I'=\{P,Q\}$ and $\armr = \id$, $\armrsec    =P$ constant,
there does not exist such a unification, because for $\varphi:
k=(P\to J \to P) \to (P \to I \to P)$ the mapping $J \to I$ need be
surjective, hence $\armr    \circ \varphi \neq \armrsec     \circ
\varphi$.

We claim that the stalk in both cases (the colimit over this
quasifilter) is just $\field$: We have $(a,0,0 \komdots
0;c)=(a,a,0,0 \komdots  0 ;c)$ by doubling the $a$-value, but also
$(a,0,0 \komdots ,0;c)=(a,0,0,0 \komdots 0 ;c)$ by doubling a
$0$-value. So the difference, which is $(0,a,0 \komdots 0 ;0)$, must
be $0$ in the colimit.

Consider now the natural numbers $\NN$ and the quasifilter indexed
by the surjective morphisms $\indmo  : \NN \to I$ which have only
one infinite fiber. With this fixation there is at most one morphism
$\indmo   \to \indmu  $. Note that such a morphism indexes a
surjective mapping $\ideal \to I'$ (independent of the morphisms
allowed in $\catopen_K$). On the other hand every covering $\ideal$
occurs in this indexing. So we have an absolute filter.

These indexing mappings $\index: \NN \to I$ induce injections $K^I
\to A \subset K^{\NN}$, where $\alg$ consists of all sequences which
are almost constant. On the other hand every element in $\alg$ is
the pull-back of an element in $K^I$ under a suitable mapping $\NN
\to I$. Hence $\alg$ is the absolute stalk of this filter. Of course
an element in $K^I$ has a lot of different images in $\alg$
depending on the indexing $\NN \to I$.
\end{bookexample}
}{}

\subsection{Exactness}

\label{exactnesssubsection}

\

\medskip
An
exact sequence of $\O_\spax $-modules yields an exact sequence of
sheaves in the flat (or the \'{e}tale) topology (in the small site),
but not in an arbitrary Grothendieck topology. The exactness of a
complex of sheaves in a non-flat topology is a phenomenon in its own
right, which we study here and in later sections for several
topologies.

The exactness in the topology is in general quite a strong
condition, in particular if $\catopen_\spax$ is big (this is already
true for the big Zariski site).
For example, if $\fuf \in \ring$ is a
non-zero divisor, then
$0 \to \ring \stackrel{\fuf}{\to} \ring$ is
exact, but as soon as there is one object
$\rings \in \catopen_\ring$ where $\fuf$ is a zero divisor, the corresponding
complex of sheaves is not exact anymore.
Therefore we will also deal with weaker exactness properties.

{
\ren{\ringsecsec}{{\ringsec'}}

\begin{corollary}
\label{exactnessintop}
Let $\spaxeqspecring$ denote an affine scheme
endowed with an affine single-handed Grothendieck topology
$\catopenspax $. Let $\lmr$ denote a complex of $\ring$-modules, and
let $\lmrr{_\topo}$ denote the corresponding complex of sheaves in
the topology. Then the following are equivalent.

\numiii
\begin{enumerate}

\item
The sheaf complex $\lmrr{_\topo}$ is exact.

\item
For every ring homomorphism $\ring \to \ringsec$ in $\catopen_\ring$
and every
$\elem \in \Gamma(\ringsec,\modm_\topo)$ such that
$\cob(\elem)= 0$ in
$\Gamma(\ringsec,\modr_\topo)$ there exists a cover
$\ringsec \to \ringsecsec$ and $\elel \in
\Gamma(\ringsecsec,\modl_\topo)$ such that $\coa(\elel)$ equals the
restriction of $\elem$ to
$\Gamma(\ringsecsec,\modm_{\topo})$.

\item
For every ring homomorphism
$\ring \to \ringsec$ in $\catopen_\ring$
and every
$\elem \in \modm \tensorr \ringsec$ such that $\cob(\elem)= 0$ in
$\modr \tensorr \ringsec$
there exists a cover
$\ringsec\to\ringsecsec$ and
$\elem \in \modl \tensorr \ringsecsec$
such that
$\elem=\coa(\elel)$ in $\modm \tensorr \ringsecsec$.
\end{enumerate}
\end{corollary}
\begin{proof}
(ii)
is just the formulation of the exactness of a complex of sheaves
(Section \ref{sheafifysubsection}) in the given affine single-handed
Grothendieck topology.
(iii)
is the formulation for
$\catopen_\ring$-exactness (Definition \ref{catexactdef}) for the complex of presheaves given by
tensoration with $\ring \to \ringsec  $.
Hence the equivalence (ii)$\Leftrightarrow$ (iii)
follows from Proposition \ref{exactexact}(iii)
and Proposition \ref{sheafzariskicompare}(iii).
\end{proof}
}


A complex $\lmrsheaf$ of sheaves on a Grothendieck topology
$\catopenspax$ on $\spax$ yields also a complex in a
(quasi-)filter
$\filtcatindtocatopenspax$. Its exactness can be
characterized as follows.

{ \ren{\indee}{{\indl}}
\begin{lemma}
\label{stalkexactexplicit}
Let $\spax$ denote a scheme together with
a Grothendieck topology $\catopenspax$, let
$\filtcatindtocatopenspax$ denote a filter and let $\lmrsheafstack $
denote a complex of abelian sheaves. Then the following are
equivalent.

\numiii
\begin{enumerate}

\item
The complex is exact in $\filt$, that is, the stalk complex
$\lmrsheafstackk {_\filt}$
is
exact.

\item
For all
$\elem \in \Gamma(\openfil_\indm  , \shem)$ such that
$\cob(\elem)=0$ in
$\Gamma(\openfil_\indm  , \sher)$
{\rm(}$\indm \in \indexcat${\rm)}
there exists
$\openfil_\indl \to \openfil_\indm  $
{\rm(}indexed by $\indl \to \indm  $ in $\indexcat${\rm)}
and
$\elel\in \Gamma(\openfil_\indl,\shel)$
such that
$\coa(\elel)= \rest^\indm _\indl(\elem)$ in
$\Gamma(\openfil_\indl,\shem)$.
\end{enumerate}
\end{lemma}
\begin{proof}
(i) $\Rightarrow$ (ii).
Since $\cob(\elem)=0$ in
$\Gamma(\openfil_\indm , \sher)$, this holds also in $\sher_\filt$.
Hence there exists
$\elel \in \shel_\filt$ such that $\coa(\elel)=\elem$ in
$\shem_\filt$. So there exists
$\obu_\indee$, $\indeeincat$, such
that $\elel$ is represented by
$\elel \in \Gamma(\obu_\indee,\shel)$
and such that
$\coa(\elel)=\rest^\indm _\indee(\elem)$ holds in
$\Gamma(\obu_\indee,\shem)$.

(ii) $\Rightarrow$ (i).
Let $\elem \in \shem_\filt$ map to $0$ in
$\sher_\filt$. We may assume that
$\elem \in \Gamma(\openfil_\indm ,\shem)$ and that
$\cob(\elem)=0$ in
$\Gamma(\openfil_\indm, \sher)$.
So by (ii) there exists
$\obu_{\indl} \to\openfil_\indm$ and
$\elel\in\Gamma(\openfil_\indl, \shel)$ mapping to
$\rest^{\indm}_\indl (\elem)$,
and this is then also true in the stalk.
\end{proof}
}

\begin{lemma}
\label{exactirreducible}
Let $\catopenspax $ denote a Grothendieck
topology on $\spax$ and let $\lmrsheafstack $ denote an exact
complex of sheaves of abelian groups on $\spax_\topo$. Then the
complex is exact in every irreducible filter. In particular it is
exact in an absolute filter, if $\catopenspax$ is single-handed.
\end{lemma}
\begin{proof}
Let
$\filtcatindtocatopenspax$ denote an irreducible filter. Let
$\elem \in \shem_\filt$ be represented by
$\elem \in \Gamma(\openfil_\indm,\shem)$, $\indmincat$, and suppose that
$\cob(\elem)=0 $ in $\sher_\filt$.
We may assume that $\cob(\elem)=0$ in
$\Gamma(\openfil_\indm  ,\sher)$.
Due to exactness there exists a covering
$\obcov_\indcov \to \openfil_\indm$, $\indcovinset$,
and sections
$\elel_\indcov \in \Gamma(\obcov_\indcov, \shel)$, $ \indcovinset$, such that
$\coa(\elel_\indcov)
= \rest^{\openfil_\indm}_{\obcov_\indcov}(\eles)$
in
$\Gamma(\obcov_\indcov,\shem)$. Since $\filt$ is irreducible there
exists $\indcovinset$ and $\indlincat$ such that
$( \openfil_\indl \to \obcov_\indcov \to \openfil_\indm ) \in \filt$.
Then
$\rest^{\obcov_\indcov}_{\openfil_\indl} (\elel_\indcov)$ maps to
$\elem$ in $\shem_\filt$.
\end{proof}

\ren{\ringsec}{{\ring'}}

\ren{\ringtri}{{\ring''}}

We characterize the exactness in an absolute filter.
This is in general a weaker property than the exactness in the topology,
but under certain conditions they are equivalent.

\begin{proposition}
\label{exactnessabsolute}
Let $\spaxeqspecring$ denote an affine
scheme endowed with an affine single-handed Grothendieck topology
$\catopenspax $. Let $\lmr$ denote a complex of $\ring$-modules and
let
$\filt :\indexcat \to \catopen_\spax$ be an absolute filter with
absolute stalk $\stalkabs$. Then the following are equivalent.

\numiii
\begin{enumerate}

\item
The sheaf complex $\lmrr{_\topo}$ is exact in the absolute filter.

\item
The complex $\lmrr{\tensorr  \stalkabs}$ is exact.

\item
For every $\indm \in \indexcat$ and every $\elem \in
\Gamma(\openfil_\indm,\modm_\topo)$
mapping to $0$ in
$\Gamma(\openfil_\indm, \modr_\topo)$
there exists
$(\indl \to \indm) \in \indexcat$ and
$\elel \in\Gamma(\openfil_\indl, \modl_\topo)$ mapping to
$\rest^\indm_\indl(\elem) \in \Gamma(\openfil_\indl , \modm_\topo)$.

\item
For every $\indm \in \indexcat$ and every $\elem \in \modm \tensorr
\ring_\indm$ mapping to $0$ in
$ \modr \tensorr \ring_\indm$ there exists
$(\indl \to \indm) \in \indexcat$ and
$\elel \in \modl\tensorr \ring_\indl$ mapping to
$\elem \tensor 1 \in \modm \tensorr\ring_\indl$.

\smallskip\noindent\hspace{-1.5cm}
If the filter is covering
{\rm(}Remark \ref{filterabsoluteremark}{\rm)}
in the sense that for every
$\indl \to \indm $ in $\indexcat$ the morphism
$\openfil_\indl \to \openfil_\indm$ is a
cover, then this is also equivalent to the following.

\item
The sheaf complex $\lmrr{_\topo}$ is $\catopen$-exact for every
$\openfil_\indm$, $\indm \in \indexcat$.

\smallskip\noindent\hspace{-1.5cm}
If the filter is covering in the sense that for every $\indm \in
\indexcat$ every morphism
$\spay \to \openfil_\indm$ in $\catopen_\spax$ is a cover, then
this is also equivalent to the following.

\item
For every cover $\spay \to \spax$ and every
$\elem \in \Gamma(\spay,\modm_\topo)$
mapping to $0$ in
$\Gamma(\spay, \modr_\topo)$ there exists
$\indl\in \indexcat$, a morphism
$\obu_\ind \to \spay$ in $\catopen_\spax$
and
$\elel \in \Gamma(\openfil_\indl, \modl_\topo)$
mapping to
$\rest^\spay_{\openfil_\indl}(\elem)
\in \Gamma(\openfil_\indl,\modm_\topo)$.

\item
For every cover $\spay \to \spax$ and every
$\elem \in \Gamma(\spay,\modm_\topo)$
mapping to $0$ in
$ \Gamma(\spay, \modr_\topo)$
there exists a cover
$\spay' \to \spax$ and a morphism $\spay' \to \spay$
and
$\elel \in \Gamma(\spay', \modl_\topo)$
mapping to
$\rest^{\spay}_{\spay'}(\elem) \in \Gamma(\spay', \modm_\topo)$.

\smallskip\noindent\hspace{-1.5cm}
If the topology is covering {\rm(}Definition \ref{singledef}{\rm)}, then this is also equivalent to

\item
The sheaf complex $\lmrr{_\topo}$ is exact.
\end{enumerate}
\end{proposition}
\begin{proof}
(i) and (ii) are equivalent by
Corollary \ref{irreduciblefiltermodulestalk}.
(iii) is the explicit version of
(i) (Lemma \ref{stalkexactexplicit}) and (iv) is the explicit
version of (ii). 
The implications (viii) $\Rightarrow$ (v), (v) $\Rightarrow$
(i)-(iv), (i)-(iv)  $\Rightarrow$ (vi) and (vi) $\Leftrightarrow$
(vii) hold without further conditions, so we prove them first.
(viii) $\Rightarrow$ (v) is a restriction (Proposition
\ref{exactexact}).

(v) $\Rightarrow$ (iii). In the given situation of (iii) there
exists by (v) a cover $\spay \to \openfil_\indm$ such that the restriction of $\elem$ comes from the
left. By definition of an absolute filter there exists a morphism
$\openfil_\indl \to \spay$, and so the restriction of $\elem$ to
$\Gamma(\openfil_\indl, \modm_\topo)$ comes from an element in
$\Gamma(\openfil_\indl,\modl_\topo)$.

(iii) $\Rightarrow$ (vi).
In the given situation there exists $\indm
\in \indexcat$ and a morphism
$\openfil_\indm \to \spay$, since $\filt$
is an absolute filter. The condition holds also for
$\rest^\spay_{\openfil_\indm}(\elem)$ and by (iii) there exists
$\openfil_\indl \to \openfil_\indm$ such that $\elem$ comes from an element
in
$\Gamma( \openfil_\indl,\modl_\topo)$.
(vi) $\Rightarrow$ (vii)
is trivial ($\spay'=\openfil_\indl$).
For (vii) $\Rightarrow$ (vi) just take
$\openfil_\indl \to \spay' \to \spay$.

Now assume that the indicated morphisms in the filter are covers.
Then (iii) $\Rightarrow$ (v) is trivial.

(vii) $\Rightarrow$ (iii)
under the given assumption. Suppose the
situation of (iii).
There exists by (vii) a morphism
$\spay' \to \openfil_\indm$, $\spay' \to \spax$ a cover,
such that the restriction of $\elem$ to $\spay'$ comes from the
left.
By assumption on the filter,
$\spay' \to \openfil_\indm$ is a cover.
Hence by the definition of an absolute
filter there exists $\indl \to \indm$ and a factorization
$\openfil_\indl \to \spay' \to \openfil_\indm$.
If the topology is covering, then (vii) $\Rightarrow$ (viii) is
trivial.
\end{proof}

For an example where Proposition \ref{exactnessabsolute}(vi) and
(vii) hold, but not (i)-(v), see Example \ref{finabsoluteextra}.

\ifthenelse{\boolean{book}}{

\begin{remark}
\label{exactzariskismallbigwarning} The exactness of a complex
depends on the chosen category, not only on the concept of a
covering. This is already relevant for the Zariski topology. The
exactness of a complex of $\ring$-modules is equivalent to the
exactness in the small site (similar in the flat topology), where
only the open sets of $\specring$ occur, but not to the exactness in
the big site, where we allow all schemes over $\specring$.%
\end{remark}}{}

\ren{\presh}{{\shP}}

\subsection{Cohomology}
\label{cohomologysubsection}

\

\medskip
The category of sheaves of abelian groups in a Grothendieck topology
is an abelian category \cite[Theorem II.2.15]{milne} and every
abelian sheaf can be embedded into an injective sheaf
\cite[Proposition III.1.1]{milne}. In particular, a sheaf has an
injective resolution and so it is possible to define
\emph{cohomology} of a sheaf as the derived functor for the functor
of global sections
$\sheaf \mapsto \Gamma(\spax, \sheaf)$ \cite[Section III.1]{haralg}.

We will be interested in cohomology applied to a short exact
sequence
$0 \to \Syz \to \O_\topo^\numfuf \stackrel{\runfuf}{\to}\ideal^\topo \to 0$
given by ideal generators
$\ideal = (\runfuf) \subseteq \ring$,
as this gives under certain conditions a bijection
$\ideal^\topocl/\ideal = H^1((\Spec \ring)_\topo, \Syz)$
(Proposition \ref{syzygycohomology}).
So the first cohomology of the syzygy sheaf measures the difference
between an ideal and its closure induced by the Grothendieck
topology. We will also compute that in the surjective topology (in
characteristic zero) and in the Frobenius topology the first
cohomology of a quasicoherent module is zero
(Theorem \ref{surjectivecohomologytheorem}, Proposition
\ref{frobeniuscohomology}) using \v{C}ech cohomology.

To review \emph{\v{C}ech cohomology} let $\spax$ be a site and let
$\namecov= (\obcovindic \to \opentop$, $\indcovinset$) be a covering of
$\opentop \in \catopenspax $.
We write
$\obcov_\choiceindcovnumcoho
= \obcov_{\indcov_0} \times_\opentop \cdots \times_\opentop \obcov_{\indcov_\numcoho}$.
The \emph{\v{C}ech diagram} (of this covering) is

\vspace{3.mm}
\begin{picture}(20,10)

\unitlength 1cm

\put(1.17,0.0){$\ldots$}

\put(1.75,0){$ {\displaystyle\biguplus_\choicecommaindcovzerothree}
\obcov_\choiceindcovzerothree$}

\put(4.43,0.07){$\longto$}
\put(4.43,+0.21){$\longto$}
\put(4.43,-0.07){$\longto$}
\put(4.43,-0.21){$\longto$}

\put(5.24,0){${\displaystyle\biguplus_\choicecommaindcovzerotwo}
\obcov_\choiceindcovzerotwo$}

\put(7.39,+0.14){$\longto$}
\put(7.39,0.){$\longto$}
\put(7.39,-0.14){$\longto$}

\put(8.2,0){${{\displaystyle \biguplus_\choicecommaindcovzeroone}}
\obcov_\choiceindcovzeroone$}

\put(9.74,0.07){$\longto$}
\put(9.74,-0.07){$\longto$}

\put(10.58,0){${\displaystyle\biguplus_\choicecommaindcovzero }\,\obcov_\choiceindcovzero$}
\put(11.62,0.){$\longto$}
\put(12.46,0){$\opentop \, ,$}
\end{picture}

%
\vspace{.7cm}
\noindent
where the $\choicebreak$th arrow
$\biguplus_{\indcov_0 \comdots \indcov_{\numcoho +1}} \obcov_{\indcov_0
 \ldots \indcov_{\numcoho +1}} \to
 \biguplus_\choicecommaindcovnumcoho \obcov_\choiceindcovnumcoho$
($(\indcov_0 \komdots \indcov_{\numcoho+1}) \in
\setindcov^{\numcoho+1}$)
is on every open component just the projection
$\projcech_\choicebreak
= \project_{\indcov_0 \ldots \check{\indcov_\choicebreak} \ldots \indcov_{\numcoho +1}}:
\obcov_{\indcov_0 \ldots \indcov_{\numcoho +1}}
\to \obcov_{\indcov_0 \ldots \check{\indcov_\choicebreak}
\ldots i_{\numcoho +1}}$ (see \cite[V.1.10]{SGA4} or \cite[III \S2]{milne}).

Applying a presheaf $\presh$ of abelian groups to this diagram we
get the \v{C}ech diagram of the covering with values in $\presh$
\cite[V.2.3.3]{SGA4}, and taking the alternating sum over the arrows of this diagram
we get the \emph{\v{C}ech complex} of $\presh$ with respect to this
covering. Explicitly,
$C^\numcoho (\namecov, \presh) = \prod_\choicecommaindcovnumcoho
\Gamma(\obcov_\choiceindcovnumcoho, \presh)$
and the derivation is
$d^\numcoho : C^\numcoho (\namecov, \presh) \to C^{\numcoho +1} (\namecov, \presh)$
by sending
$ (d^\numcoho  \chcyc )_{\indcov_0 \comdots \indcov_{\numcoho +1}}
= \sum_{\choicebreak = 0}^{\numcoho + 1} (-1)^\choicebreak
\rest_\choicebreak (\chcyc_{\indcov_0 \cdots
\check{\indcov_\choicebreak} \cdots \indcov_{\numcoho +1}})$
, where
$\rest_\choicebreak
= \rest^{\indcov_0 \cdots \check{\indcov_\choicebreak}
 \cdots \indcov_{\numcoho +1}}_{\indcov_0 \comdots \indcov_{\numcoho +1}}$
stands for the restriction to the projection
$\projcech_\choicebreak$.
The homology groups of this complex are called the
\emph{\v{C}ech cohomology groups} and are denoted by
$\chch^\numcoho (\obcovindic \to \opentop, \presh)$.

If $\obcovsec_\indcovsec \to \opentop$, $\indcovsecinset$, is a second
covering which refines $\namecov$, then there exists a refinement
homomorphism
$\chch^\numcoho (\namecov , \presh)
\to \chch^\numcoho(\obcovsec_\indcovsec \to \opentop, \presh)$.
The \v{C}ech cohomology groups are defined as the colimit
$\chch^\numcoho(\opentop,\presh)
= \colim_\namecov  \chch^\numcoho (\namecov, \presh)$ over all
coverings $\namecov$ of $\opentop$ \cite[V.2.4.5]{SGA4}.
For a sheaf
$\presh$,
$\chch^1(\opentop,\presh)$
equals the derived functor cohomology $H^1(\opentop,\presh)$
\cite[Corollary III.2.10]{milne}.


The assumption in the next lemma holds in the constructible topology
(Lemma \ref{constructibleproductlemma}, Proposition
\ref{constructiblecohomology}) and in the Frobenius topology
(Proposition \ref{frobeniuscohomology}) and generalizes Lemma
\ref{antistrictlemma}.

{\ren{\presh}{{\sheaf}}

\begin{lemma}
\label{chchisoproj}
Suppose that a cover $\spay \to \spax$ in a Grothendieck topology on
a scheme $\spax$ has the property that also the diagonal
$\diag: \spay \to \spay \times_\spax \spay$ is a cover
{\rm(}this is in particular the case if a projection $\spay \timesx
\spay \rightarrow \spay$ is an isomorphism{\rm)}.
Then for a sheaf
$\presh$ of abelian groups on
$\spax_\topo$ we get in the \v{C}ech complex either the $0$-map or
the identity,
$$ \Gamma(\spax,\presh) \lto \Gamma(\spay,\presh)
\stackrel{0}{\lto} \Gamma(\spay\timesx \spay, \presh)
\stackrel{\id}{\lto} \Gamma(\spay \timesx \spay \timesx
\spay,\presh) \stackrel{0}{\lto} \, .$$
It follows that the \v{C}ech cohomology is trivial,
i.e. $\chch^0(\spay\to \spax, \presh)
=\Gamma(\spay,\presh)$ and
$\chch^\numcoho(\spay \to \spax,\presh)=0$
for $\numcoho \geq 1$.
\end{lemma}
\begin{proof}
We show that
$\rest_\choicebreak = \rest(\projcech_\choicebreak)$ is
the same mapping for every $\choicebreak$.
Since
$\projcech_\choicebreak \circ \diag_{\choicebreak, \choicebreak+1}
=\projcech_{\choicebreak+1} \circ \diag_{\choicebreak,
\choicebreak+1}$ we get
$\rest(\diag_{\choicebreak, \choicebreak+1})\circ \rest_\choicebreak
=\rest(\diag_{\choicebreak, \choicebreak+1})\circ \rest_{\choicebreak+1}$.
As $\diag$ is a cover, also
$\diag_{\choicebreak, \choicebreak+1}$
is a cover, hence it follows that
$\rest_\choicebreak =\rest_{\choicebreak+1}$.
All these mappings are identities, since for a cover
$\obu \to \obv$ with a covering section $\obv \to \obu$ the induced
maps on a sheaf are always isomorphisms.
In the \v{C}ech complex we have either
an even or an odd alternating sum of identities, so the sum is
either $0$ or the identity. The statements about the cohomology
follow.
\end{proof}

}

\synchronize{\obtop}{{\opentop}}

It is also possible to define \emph{local cohomology} for an open subset
$\obtop \to \spax$ in a Grothendieck topology on $\spax$
(see \cite[V.6.3]{SGA4}). For an
abelian sheaf $\sheaf$ define
$$\Gamma_\obtopsup (\sheaf)
= \{\elem \in \Gamma(\spax,\sheaf): \elem|_\obtop = 0 \}
\,= \ker( \rest^\spax_\obtop:
\Gamma(\spax, \sheaf)  \to \Gamma(\obtop,\sheaf))\, $$
($\obtopsup$ is just a symbol for us).
This is again a left exact functor
and gives rise to a derived functor, denoted by
$H^\numcoho_\obtopsup (-)$.
There exists a long exact sequence
(as in \cite[Corollaire I.2.9]{SGA2})
$0 \to \Gamma_\obtopsup(\sheaf) \to \Gamma (\spax,\sheaf) \to \Gamma
(\obtop,\sheaf) \to H^1_\obtopsup (\sheaf) \to H^1(\spax,\sheaf) \to
H^1(\obtop,\sheaf) $ \cite[Proposition V.6.5]{SGA4}. We will use this in Section
\ref{syzygycohomologysubsection}.

\ifthenelse{\boolean{book}}{
\begin{remark}
Note that we have a natural map $\Gamma(\spax, \presh) \to
\chch^0(\obcov_\indcov \to \obu, \presh)$, which is an isomorphism for a
sheaf. So we have in Lemma \ref{chchisoproj} also the identity
$\Gamma(\spax, \presh)= \Gamma(\spax', \presh)$ (for a sheaf!). This
lemma is relevant for several topologies: for the constructible
partitions in the surjective topology (Proposition
\ref{constructiblecohomology}), for the Frobenius topology (Part
\ref{frobeniussection}) and for the topologies induced by a Zariski
filter (Part \ref{zariskifiltersection}).
\end{remark}

\begin{remark}
Let $\filt: \indexcat \to \catopenspax    $, and let $\shem$ be a
(pre-)sheaf of abelian groups on $\spax_\topo$. Evaluation at the
filter gives a left exact functor $\Sh (\spax) \to \Ab$, $\shem
\mapsto \Gamma(F , \shem ) = \colim_{\indmo   \in \indexcat}
\Gamma(\openfil_\indmo  , \shem)$ from sheaves to groups. Its right
derived functor is the colimit of the cohomology groups, i.e.
$H^\numcoho (F, \shem)= \colim_{\indmo   \in \indexcat} H^\numcoho
(U_\indmo  , \shem)$. In particular we get for a short exact
sequence of sheaves of abelian groups $\olmrosheaf$ the long exact
sequence
$$0 \to \shel_F \to \shem_F \to \sher_F \to H^1(F,\shel) \to
H^1(F,\shem) \to \cdots \, .$$
\end{remark}

\begin{remark}
\label{cohomologyfunctoriality}
Let $\sitemor: \spay \to \spax$ be a
site morphism. There is a homomorphism $
\chch^{\numcoho}(\obcov_\indcov \to
\obu, \shem) \to \chch^{\numcoho}(\sitemortop (\obcov_\indcov) \to
\sitemortop (\obu), \sitemorsheaf (\shem))$, where $\obcov_\indcov \to \obu$
is a covering, and therefore there exists a homomorphism
$\chch^\numcoho (\obu, \shem) \to \chch^\numcoho (\sitemortop(\obu),
\sitemorsheaf (\shem))$. In particular, cohomology and \v{C}ech
cohomology yield presheaves on a site $\spax$. If
$\obu \mapsto \chch^\numcoho(\obu, \shem)$ is such a cohomological presheaf,
denoted by $\shH^\numcoho (\shem)$, and if $\sitemor: \spay \to
\spax$ is a morphism of pre-ringed sites, then we have two different
associated (pre)sheaves on $\spay$. On one hand we have the
pull-back $\O_\spay $-premodule $\sitemorpre (\shH^\numcoho(\shem))$
and the pull-back module $\sitemorsheaf (\shH^\numcoho(\shem))$,
where $\shH^\numcoho(\shem)$ is just considered as an arbitrary
$\O_\spax$-premodule. On the other hand we also have the assignment
$\obv \mapsto \chch^\numcoho (\obv, \sitemorsheaf (\shem))$. There
exists a homomorphism
$ \sitemorpre (\shH^\numcoho(\shem)) \to \chch^\numcoho (-, \sitemorsheaf (\shem))$.
In particular we have homomorphisms
$$ H^\numcoho (\spax, \shem) \longto \Gamma(\spay, \sitemorpre ( \shH^\numcoho( \shem)))
\longto H^\numcoho (\spay, \sitemorsheaf(\shem)) \, .$$
\end{remark}

\begin{remark}
If $\obu_\indec
 \to \obu$ is a cover, then there exist natural
homomorphisms
$$\chch^\numcoho( \openfil_\indmo
\to \obu, \presh) \longto \chch^\numcoho(\obu, \presh) \longto \chch^\numcoho(\openfil_\indmo  , \presh) \, ,$$
where the second homomorphism comes from the functoriality of
\v{C}ech cohomology. If the topology is single-handed and if $\filt:
\indexcat \to \catopenspax    $ is an absolute filter, then we also
have homomorphisms
$$  \colim_{\indmo   \in \indexcat} \chch^\numcoho( \openfil_\indmo   \to \spax, \presh)
= \chch^\numcoho(\spax, \presh) \longto
 \colim_{\indmo \in \indexcat} \chch^\numcoho (\openfil_\indmo, \presh)
 = \chch^\numcoho (F, \presh) \,  $$
(and the two sides should not be confused).
\end{remark}

\begin{remark}
\label{cechcohomologycautionremark}
Let $\sheaf $ be a sheaf on
$\spax=\spax_\topo$ and let a \v{C}ech cohomology class $c \in
\Gamma(\spax' \times_\spax\spax', \sheaf )$
be given (that is, a class mapping to $0$ in $\Gamma(\spax'
\times_\spax\spax'\times_\spax\spax', \sheaf )$).
This class is $0$ in the \v{C}ech
cohomology if there exists a refinement $V \to \spax'$ such that the
image of the cohomology class in $\Gamma(V\times_\spax V, \sheaf )$
comes from $\Gamma(V,\sheaf )$. It is possible however that there
exists a covering (!) $\spay \to \spax$ such that the pull-back of
$c$ gets $0$, but such that $c \neq 0$. The difference in these
conceptions should become clear from the following commutative
diagram ($\spay'=\spay \times_\spax\spax'$).
$$\xymatrix{\gammaf {\spax'} \ar[r]^{\delta_\spax} \ar[d]
& \gammaf {\spax' \timesx \spax'} \ar[d]\\
\gammaf {\spay'} \ar[r]^{\delta_\spax} \ar[d]_{=} & \gammaf {\spay' \timesx \spay'} \ar[d]\\
\gammaf {\spay'} \ar[r]^{\delta_\spay }  & \gammaf {\spay'
\times_\spay \spay'}\, .}$$ Here the image of $c \in \gammaf {\spax'
\timesx \spax'}$ may come on the bottom row from the left, but not
in the middle row.
\end{remark}}{}

\section{Closure operations and Grothendieck topologies}

We describe now how a Grothendieck topology on an affine scheme
$\spaxeqspecring$ defines a closure operation for $\ring$-submodules
$\submodul$ and we characterize which closure operations can be
obtained in this way. In general there will be several more or less
natural Grothendieck topologies which realize the same closure
operation. These Grothendieck topologies come along with a variety
of gadgets such as cohomology, ring and modules of global sections,
notions of stalks and exactness, which we try to understand in the
context of the closure operation.

\subsection{Forcing algebras}
\label{forcingsubsection}

\

\medskip
We recall the definition of (universal) forcing algebras
following Hochster \cite{hochstersolid}. Let $\ring$ be a
commutative ring and let $\submodul$ denote $\ring$-modules such
that $\quotmod$ is finitely presented, $\elem \in \modul$. This
means that we have a resolution
$$ \ring^\submodra \stackrel{\mata}{\lto} \ring^\modra \lto \quotmod  \lto 0 \, ,$$
where
$\mata=(\entry_{\indcri\indcrj})_{\indcri \indcrj}$ is a
$\modra \times \submodra $-matrix with entries in $\ring$. Let
$\tilde{\elem}= (\elem_1 \komdots \elem_\modra) \in \ring^\modra$
denote a lifting of $\elem$. In this situation we make the following
definition.

\begin{definition}
\label{forcingalgebradef}
The finitely generated $\ring$-algebra
$$ \algforc
 = \ring[\vart_1 \komdots \vart_\submodra ]/
 (\mata \vart -\tilde{\elem})
 = \ring[\vart_1\komdots \vart_\submodra ]/
 (\sum \entry_{1\indcrj}\vart_\indcrj -
\elem_1, \ldots, \sum \entry_{\modra \indcrj}\vart_\indcrj -
\elem_\modra)$$
is called a \emph{forcing algebra} for $(\modul,\submod,\elem)$.
Accordingly, $\Spec \algforc$ is called the \emph{forcing scheme}.
\end{definition}

\begin{remark}
\label{presentforcingremark}
Note that the three sets of data
$(\modul,\submod,\elem)$, $( \ring^\modra,\tilde{\submod},
\tilde{\elem})$, where $\tilde{\submod} $ is the kernel of
$\ring^\modra \to \quotmod$, and $(\bar{\modul},0,\bar{\elem})$,
where $\bar{\modul}= \modul/\submod$ and $\bar{\elem}$ denotes the
image of $\elem$, have the same forcing algebras. This means that we
may simplify the situation by either passing to a torsion-free
submodule inside a free module or to the $0$-submodule inside a
module. Or the other way round: forcing algebras can only be helpful
to understand such closure operations which are independent of the
presentation, see Proposition \ref{admissibleproperties}(iii).
\end{remark}

\begin{remark}
\label{forcingfunctorremark}
The construction of the forcing algebra
is functorial with respect to ring homomorphisms $\ring \to
\ringsec$. By right exactness we get immediately the resolution
$$\ringsec  ^\submodra  \stackrel{\mata'}{\lto} \ringsec^\modra
\lto (\quotmod) \tensorr \ringsec
 = \modul'/\submod'  \lto 0 \, ,$$
where $\submod'$ is the image of $\submod \tensorr \ringsec  $
inside $\modul' =\modul \tensorr \ringsec  $, and the forcing
algebra over $\ringsec  $ built with $\mata'$ is just $\algforc
\tensorr \ringsec  $.
\end{remark}

The forcing algebra has the following universal property.

\begin{proposition}
\label{universalforcing}
Let $\ring$ denote a commutative ring and
let $\algforc$ denote a forcing algebra for $\ring$-modules
$\submodul$ such that $\quotmod$ is finitely presented, $\elem \in
\modul$. Then the following holds.

\numiii
\begin{enumerate}

\item
The element $\elem \tensor 1 $ lies in the image of
$\submod \tensorr  \algforc \to \modul \tensorr \algforc$.

\item
If $\alg$ is an $\ring$-algebra such that
$\elem \tensor 1
\in \im(\submod \tensorr \alg \to \modul \tensorr \alg)$,
then there exists an $\ring$-algebra homomorphism
$\algforc \to \alg$.
\end{enumerate}
\end{proposition}
\begin{proof}
We use the notation of Remark \ref{presentforcingremark}. The
condition that $\elem \tensor 1 \! \in \im (\!\submod \tensorr \alg
\! \to \! \modul \tensorr \alg)$ is equivalent to the property that
$\bar{\elem} \otimes 1$ becomes $0$ in $(\quotmod) \tensorr \alg$,
and this is equivalent to the property that $\tilde{\elem} \otimes
1$ in $ \alg^\modra$ is in the image of $\alg^\submodra  \to
\alg^\modra$. This means that the inhomogeneous linear system of
equations $\sum \entry_{1\indcrj}\vart_\indcrj - \elem_1, \ldots,
\sum \entry_{\modra \indcrj}\vart_\indcrj - \elem_\modra $ has a
solution in $\alg$. For the forcing algebra $\algforc$ this is true
by the very definition, which gives (i). (ii). If a solution over
$\alg$ exists, then we just have to send the variables
$\vart_\indcrj$ to the corresponding coefficients in $\alg$, so
there exists an $\ring$-algebra homomorphism $\algforc \to \alg$.
\end{proof}

Recall that a ring homomorphism $\ring \to \alg$ is called
\emph{pure} if for every $\ring$-module $\modul$ the induced map
$\modul \to \modul \tensorr \alg$ is injective. For a pure
homomorphism the contraction of an extended ideal is the ideal
itself. Pure ring homomorphisms induce a spec-surjective morphism.
We note the following basic observation.

\begin{lemma}
\label{forcingequivalent} Let $\ring$ be a commutative ring. Let
$\submodul$ be $\ring$-modules such that $\quotmod$ is finitely
presented, let $\elem \in \modul$ be an element. Let
$$\algforc = R[\vart_1 \komdots \vart_\submodra ]/(\mata \vart-\tilde{\elem})$$
be the forcing algebra for $(\modul,\submod, \elem)$. Then the
following are equivalent.

\numiii
\begin{enumerate}

\item $\elem \in \submod$.

\item
There is a ring section $\algforc \to \ring$ for $\ring \to
\algforc$.

\item
$\ring \to \algforc$ splits as a homomorphism of $\ring$-modules
{\rm(}that is, $\ring$ is a direct summand in $\algforc${\rm)}.

\item $\ring \to \algforc$ is pure.
\end{enumerate}
\end{lemma}
\begin{proof}
(i) $\Rightarrow$ (ii) follows from the universal property of a
forcing algebra (Proposition \ref{universalforcing}(ii)). The
implications (ii) $\Rightarrow$ (iii) and (iii) $\Rightarrow$ (iv)
hold for arbitrary ring homomorphisms. So suppose that (iv) holds.
Apply the purity property to the $\ring$-module $\modul/\submod$.
Let $\bar{\elem}$ denote the image of
$\elem$ in $\modul/\submod$.
Then by Proposition \ref{universalforcing}(i) $\bar{\elem}$ becomes
$0$ in $(\modul/\submod) \tensorr  \algforc $, so it must be $0$ by
purety.
\end{proof}

\begin{remark}
Let $\ring$ denote a commutative ring, and let $\ideal = (\runfuf)$
denote a finitely generated ideal, let $\fuf \in \ring$. In this
case the resolution is just $\ring^\numgen
\stackrel{\fuf_\indfu}{\to} \ring \to \ring/ \ideal \to 0$ and the
forcing algebra
$$\algforc
= \ring[\vart_1 \komdots \vart_\numgen]/
(\fuf_1\vart_1 \plusdots\fuf_\numgen \vart_\numgen+ \fuf)$$
is given by one affine-linear
equation. In this case we can add to the equivalent conditions in
Lemma \ref{forcingequivalent} the property that $\ring \to \algforc$
is cyclic
pure
, which means by definition that
$\ideal \algforc \cap \ring = \ideal$ for every ideal $\idealsubring$.
\end{remark}

\ifthenelse{\boolean{book}}{Forcing algebras have the following
trivial gluing property.

\begin{corollary}
\label{forcingsectionlocal}
Let $\varphi: \spay \to\spaxeqspecring$
denote a forcing algebra and suppose that $\varphi$ has locally a
section, i.e. there exists a Zariski covering $\spax= \bigcup_{i \in
I} \obcov_\indcov$ and morphisms $s_i:\obcov_\indcov \to \varphi^{-1}(\obcov_\indcov)$ with
$\varphi_i \circ s_i = \id_{\obcov_\indcov}$. Then there exists also a global
section $s: \spax \to \spay$.
\end{corollary}
\begin{proof}
We may assume that the covering is affine, say $\obcov_\indcov=D(f_i)$, $i \in
I$. Let $\algforc$ be the forcing algebra for $(\modul,\submod,
\elem)$, $\spay = \Spec \algforc$. The existence of a section on
$\obcov_\indcov$ means by Lemma \ref{forcingequivalent} that $\elem \in
\submod_{f_i}$ in $\modul_{f_i}$. It follows that $\elem \in
\submod$ and again by Lemma \ref{forcingequivalent} we deduce the
existence of a global section.
\end{proof}

\begin{remark}
Suppose that $\spax$ is a scheme and $\submodul$ are
$\O_\spax$-modules with a presentation $\O_\spax^\submodra \to
\O_\spax^\modra \to \modul/\submod \to 0$. Let $\elem \in
\Gamma(\spax, \O_\spax^\modra)$. Then the affine forcing schemes
given by the forcing algebras on $U \subseteq \spax$ affine glue
together to form a \emph{forcing scheme} $\spay \to \spax$. A
forcing algebra $\algforc$ induces on every open subset $U \subseteq
\specring$ a forcing scheme $\Spec \algforc|U \to U$. Occasionally
we will use this more general viewpoint.
\end{remark}

\begin{remark}
\label{forcingglobalremark} Corollary \ref{forcingsectionlocal} is
not true on an arbitrary non-affine scheme $\spax$. If for example
$\spay=\Spec B \to \specring$ is a forcing algebra for $\fuf$ and an
ideal $\ideal$, and if $\psi:Z \to \specring$ is not affine (e.g. a
blow-up as in Section \ref{ratliffsubsection}), then $\spay_Z \to Z$
admits locally sections if and only if $\fuf$ is locally (and then
globally) in the extended ideal sheaf $\ideal \O_Z$. A global
section means however that $\fuf \in \ideal  \Gamma(Z,\O_Z)$.
\end{remark}}{}

\subsection{Admissible closure operations}
\label{admissiblesubsection}

\

\medskip
The closure operations we will consider will be defined for
submodules in a certain class of modules and will be persistent with
respect to certain ring homomorphisms. We fix the following
situation.

\begin{situation}
\label{closuresituation}
Let $\grocont$ denote a subcategory of
commutative rings and suppose that a subcategory
$\modcat=\modcat_\ring$ of $\ring$-modules is given for each $\ring
\in \grocont$, such that for every ring homomorphism
$\homring: \ring \to \rings$ in $\grocont$ the module
$\modul\tensorr \rings$
belongs to $\modcat_\rings$.
A \emph{closure operation} for
$\modcat$ is an assignment
$\submod \mapsto \submod^\clop$ from
$\modcat_\ring$ to $\modcat_\ring$ with
$\submod \subseteq \submod^\clop = (\submod^\clop)^\clop$.
The closure operation is $\grocont$-\emph{persistent} if for
$\homring : \ring \to \rings   $ in
$\grocont$ one has
$ \im ( \submod ^\clop) \subseteq (\submod')^\clop$
in $\modul \tensorr \rings$, where $\submod' = \im
(\submod \tensorr \rings)$.
\end{situation}

\begin{remark}
The category $\grocont$ will be the full category of all rings, or
all noetherian rings, all domains, all rings of characteristic
$\posp$ etc. We may also restrict the morphisms to homomorphisms of
finite type, inclusion of domains, etc.

The category of modules $\modcat $ will be either the category of
all ideals (then we do not insist that $\ideal \tensorr \rings$
belongs to $\modcat_\rings$), or the category of all finitely
generated $\ring$-module, or all torsion-free modules, or the
category of all modules.
\end{remark}

\ren{\covadma}{S}

\begin{definition}
Let an assignment $\clop$ on $\modcat$ as described in
\ref{closuresituation} be given. We call a forcing algebra
$\algforc$ for $(\modul, \submod, \elem)$ for $\ring$-modules
$\submodul$ such that $\quotmod$ is finitely presented
\emph{admissible} (with respect to the closure operation) or
$\clop$-\emph{admissible} if $\elem \in \submod^\clop$.

We call a ring homomorphism
$\homring: \ring \to \covadma$
\emph{$\clop$-pure} if for every $\submodulzweit$ in $\modcat_\ring$
we have that
$\elemzweit \tensor 1 \in
\im( \submodzweit \tensorr \covadma \to \modulzweit \tensorr \covadma)$
implies that
$\elemzweit \in \submodzweit^\clop$.
\end{definition}

\begin{proposition}
\label{admissiblecomposition}
Suppose the situation described in
\ref{closuresituation}. Then the following holds.

\numiii

\begin{enumerate}

\item
$\ring \to \covadma$ is $\clop$-pure if and only if for every
forcing algebra $\ring \to \algforc$ such that
$\covadma \to\algforc \tensorr \covadma =\algforc_\covadma $ is
$\clop$-admissible, then also $\ring \to \algforc$ is
$\clop$-admissible.

\item
The composition of two $\clop$-pure homomorphisms $\ring
\stackrel{\varphi}{\to} \rings \stackrel{\psi}{\to} \ringt$ is also
$\clop$-pure.
\end{enumerate}
\end{proposition}
\begin{proof}
(i). Let $\algforc$ be the forcing algebra for
$(\modul,\submod,\elem)$. Then $\elem \tensor 1$ is in the image of
$\submod \tensorr \covadma \to \modul \tensorr \covadma$ if and only
if $\elem \tensor 1$ is in the closure of the extended module
$\submod'$.  This is by definition equivalent to
$\covadma \to \algforc_\covadma$ being $\clop$-admissible. That $\ring \to
\algforc $ is $\clop$-admissible is equivalent to $\elem \in
\submod^\clop$.

(ii) follows direct or from (i), since if $\ring \to \algforc$ is a
forcing algebra such that $\ringt \to \algforc_\ringt$ is
admissible, then first $\rings \to \algforc_\rings$ is admissible
and then $\ring \to \algforc$ is admissible.
\end{proof}

\begin{definition}
\label{admissibleclosuredef} Suppose the situation described in
\ref{closuresituation}. We call $\clop$ $\grocont$-admissible, if
the following conditions hold.

\numiii

\begin{enumerate}

\item
$\clop$ is a closure operation, that is, $\submod \subseteq
\submod^\clop =(\submod^\clop)^\clop$.

\item
The admissible forcing algebras are $\clop$-pure.

\item
The closure operation is $\grocont$-persistent.
\end{enumerate}
\end{definition}

\ren{\betweenmod}{{L}}

\ren{\numcomp}{{k}}

\ren{\hommod}{{\varphi}}

\begin{proposition}
\label{admissibleproperties} Suppose that we have an admissible
closure operation on a subcategory of noetherian rings and for
finitely generated modules. Then the following hold.

\numiii

\begin{enumerate}

\item
{\rm(}Module-persistence{\rm)}
Let $\submodul $, $\elem \in \modul$
and let $\hommod: \modul \to \modul'$ be an $\ring$-module
homomorphism. Set $\submod' := \hommod(\submod)$ and suppose that
$\elem \in \submod^\clop$. Then also $\hommod(\elem) \in
(\submod')^\clop$.

\item
{\rm(}Independence of presentation{\rm)}
Let $\submodul $, $\elem
\in \modul$. Let
$\homproj:\tilde{\modul} \to \modul \to 0$ and
$\tilde{\submod} := \homproj^{-1}(\submod)$ and let $\tilde{\elem}
\in \tilde{\modul}$ be such that $\homproj(\tilde{\elem})=\elem$.
Then $\elem \in \submod^\clop$ if and only if $\tilde{\elem} \in
\tilde{\submod}^\clop$. In particular, $\elem\in \submod^\clop$ if
and only if $\bar{\elem} \in 0^\clop $ in $\modul/\submod$.

\item
The closure operation is order preserving,
that is, for submodules
$\submod \subseteq \betweenmod \subseteq \modul$ of a finitely
generated $\ring$-module $\modul$ we have $\submod^\clop \subseteq
\betweenmod^\clop$.

\item Suppose also that the closure is persistent with respect to homomorpisms of finite type.
Let $\submod_\indcomp \subseteq \modul_\indcomp, \indcomp = 1
\komdots \numcomp$, be submodules and set $\submod = \submod_1
\timesdots \submod_\numcomp \subseteq \modul_1 \timesdots
\modul_\numcomp = \modul$. Then $\submod^\clop = \submod_1^\clop
\timesdots \submod_\numcomp^\clop$, where $\submod_\indcomp^\clop$
is the closure
inside $\modul_i$.
\end{enumerate}
\end{proposition}
\begin{proof}
(i).
Let $\ring \to \algforc$ denote the (admissible) forcing
algebra for $(\modul,\submod,\elem)$. Then $\elem \tensor 1$ is in
the image of $\submod\tensorr  \algforc \to \modul \tensorr
\algforc$. But then $\hommod (\elem) \tensor 1$ is in the image of
$\submod'\tensorr \algforc \to \modul' \tensorr  \algforc$.

(ii).
One direction follows from (i). For the other direction
suppose that $\elem \in \submod^\clop$, let $\algforc$ denote a
forcing algebra for these data. Then $\algforc$ is also the forcing
algebra for $(\tilde{\modul}, \tilde{\submod},\tilde{\elem})$.

(iii) follows from (i) and (ii).


(iv).
By (ii) we may assume that $\submod =0$. Let $\elem =(\elem_1
\komdots \elem_\numcomp)$. If $\elem_\indcomp \in 0^\clop$ in
$\modul_\indcomp$ for every $\indcomp$, then the corresponding
forcing algebras $\algforc_\indcomp$ are $\clop$-admissible, hence
$\clop$-pure. Because of the persistence property this holds in the
universal sense. Therefore also their tensor product $\algforc$
(which is the forcing algebra of $\elem$) is $\clop$-pure and since
$\elem = 0$ in
$\modul \tensorr \algforc$ we get $\elem \in0^\clop$.
On the other hand, if $\elem \in 0^\clop$, then also
$\elem_\indcomp \in 0^\clop$ by applying (i) to the projections
$\modul \to \modul_\indcomp$.
\end{proof}


%
%
%
%
%

\begin{remark}
The module persistence of an admissible closure operation implies
that the assignment $\modul \mapsto 0^\clop$ is a functor from
$\ring$-Mod to itself (which is not a projection). An inclusion
$\submod \subseteq \modul$ yields to an inclusion $0^\clop_\submod
\subseteq 0^\clop_\modul$, but this functor is neither exact in the
middle nor does it respect surjections. The left exactness is
related to hereditary torsion theories; see Lemma \ref{leftexact}
and Section \ref{torsiontheorysubsection}.
\end{remark}

\begin{lemma}
\label{primarytriviallemma}
An order preserving closure operation
for ideals in a local noetherian ring $(\ring,\fom)$ is trivial if
it is trivial for all $\fom$-primary ideals.
\end{lemma}
\begin{proof}
Let $\idealsubring$ and suppose that $\fuf \in \ideal^\clop$, $\fuf
\not \in \ideal$. If $\fuf$ is a unit, then also $\fuf \in
\fom^\clop$. So we may assume that $\fuf$ is not a unit. Then there
exists $\expok \in \NN$ such that
$\fuf \not\in \ideal +\fom^\expok$
(by Krull's intersection theorem \cite[Corollary 5.4]{eisenbud}
applied to $\ring/ \ideal$).
Then $\fuf \in \ideal^\clop \subseteq
(\ideal + \fom^\expok)^\clop$, and the $\fom$-primary ideal $\ideal
+ \fom^\expok$ is not $\clop$-closed.
\end{proof}

\begin{remark}
If an order preserving closure operation is given on an
$\ring$-module $\modul$, then an arbitrary intersection of a family of closed submodules
$\submod_\indmod = \submod_\indmod^\clop$, $\indmodinset$,
is again closed. Therefore one can define a topology on $\modul$
(in the usual sense) by declaring finite unions of closed submodules
to be the closed subsets. As the union of two submodules is a
submodule only if one of them is contained in the other, the closed
submodules in this topology are exactly the closed submodules for
the operation.
\end{remark}

\ren{\numcomp}{{r}}

\ren{\numcompbreak}{{s}}

\begin{example}
\label{topdimensional}
There are several examples of closure
operations for ideals which are trivial for $\fom$-primary ideals,
but not in general, so these are not order preserving and hence not
admissible (and not induced by a Grothendieck topology). E.g. the
so-called \emph{top-dimensional closure} (or \emph{equidimensional
hull} \cite[Definition 5.7]{hunekeparameter}) of an ideal
$\idealsubring=\field[ \var_1 \komdots \var_\numvar]$ (or more
general) is given in the following way: let
$\ideal= \foq_1 \capdots \foq_\numcomp$ be a primary decomposition such that $\foq_\indcomp$
is $\fop_\indcomp$-primary and such that the $\fop_1 \komdots
\fop_\numcompbreak$
($\numcompbreak \leq \numcomp$) have minimal height
(or $\ring/\fop_\indcomp$ has maximal dimension) among the
$\fop_\indcomp$. Then $\ideal^\topdim= \foq_1 \capdots
\foq_\numcomp$. If $\ideal$ is $(\var_1 \komdots
\var_\numvar)$-primary, then of course $\ideal^\topdim = \ideal$.
\end{example}

\ren{\expon}{{n}}

\begin{example}
\label{symbolicremark}
The ``\emph{symbolic closure operation}'' $\ideal
\mapsto \bigcap_{\ideal \subseteq \fop \mbox{ minimal prime} } \ring
\cap \ideal \ring_{\fop}$ is not admissible, since it is trivial on
the primary ideals, but not in general. For a prime ideal $\fop$
this assignment maps the ordinary power $\fop^\expon$ to the
symbolic power $\fop^{(\expon)}$, \cite[Section 3.9]{eisenbud}.
See also Section \ref{divisorialsubsection}, where the divisorial
topology induces the symbolic powers for prime ideals of height one.
\end{example}

\begin{example}
\label{ratliffrush}
The Ratliff-Rush \emph{closure} of an ideal
$\idealsubring$ is defined as $\tilde{\ideal}= \bigcup_{\runn \geq
1} (\ideal^{\runn+1} : \ideal^{\runn})$, see
\cite{ratliffrushtwonotes}. This closure is not an admissible
closure operation, since it is not order preserving
\ifthenelse{\boolean{book}}{For example, the Ratliff-Rush closure of
the $0$-ideal is the whole ring, but there are also examples of
$\fom$-primary ideals where this closure does not respect inclusions
\cite{source}.}{}
\end{example}

\ren{\varx}{{x}}

\begin{example}
The socle $\soc (\modul) = \{\elem \in \modul:\, \elem \fom =0 \}$
for an $\ring$-module over a local ring $(\ring,\fom)$ is not a
closure operation. E.g. for $\ring = \field [\varx]/(\varx^3)$ we
get $\soc (\ring) =(\varx^2)$, but its ``socle-closure'' would be
$(\varx)$.
\end{example}

\begin{definition}
\label{finitepropdef}
Let $\clop$ denote a closure operation on
$\modcat$. Consider for a pair $\submodul$, $\elem \in \modul$,
$\elem \in \submod^\clop$ the properties

\numiii
\begin{enumerate}
\item
There exists an $\ring$-submodule $\submod \subseteq \modul'
\subseteq \modul$ such that $\elem \in \modul'$, $\elem \in
\submod^\clop$ inside $\modul'$ and $\modul'/ \submod$ is finitely
presented.

\item
There exists a finitely generated $\ring$-submodule $\submod'
\subseteq \submod$ such that $\elem \in (\submod')^\clop$.
\end{enumerate}
We say that the closure operation has the first finiteness property,
if (i) holds (for every $\submodul$) and the second finiteness
property, if (ii) holds.
\end{definition}

\ifthenelse{\boolean{book}}{
\begin{bookremark}
We will show in Lemma \ref{integralfinitelemma} below that the
integral closure has this property.
\end{bookremark}}{}

\ren{\covadma}{A}

\subsection{Examples of admissible closure operations}
\label{exampleadmissiblesubsection}

\

\medskip We give some examples of admissible closure operations.

\begin{example}
We consider the identical closure operation, where of course the
conditions (i) and (iii) (for the category of all rings) in
Definition \ref{admissibleclosuredef} are fulfilled. The ring
homomorphisms $\ring \to \covadma$ which are pure with respect to
the identical closure operation are exactly the pure ring
homomorphisms. An admissible forcing algebra for the trivial closure
operation is pure by Lemma \ref{forcingequivalent}, hence the
identical closure is admissible.
\end{example}

\begin{definition}
\label{zielringcategory}
Let $\subcatring$ be a subcategory of
rings. Such a category defines a closure operation for
$\ring$-submodules $\submodul $ by declaring for $\elem \in \modul$
that $\elem \in \submod^\clop $ if and only if
$$ \elem \tensor 1 \in \im(\submod \tensorr \testring  \to \modul \tensorr \testring)
\mbox{ for all ring homomorphisms } \ring \to \testring,\, \testring
\in \subcatring \, .$$ We call this the \emph{closure operation
associated} to $\subcatring$ (the rings in $\subcatring$ are the
testrings).
\end{definition}

\begin{lemma}
\label{testringtest} Let $\ring$ be a commutative ring, let
$\subcatring$ be a category of rings with the induced closure
operation $\clop$. Let $\submodul$ be $\ring$-modules such that
$\quotmod$ is finitely presented, $\elem \in \modul$, let $\algforc$
be the forcing algebra. Then the following are equivalent. \numiii
\begin{enumerate}

\item
$\elem \in \submod^\clop$.

\item
For every $\ring \to \testring$, $\testring \in \subcatring $, the
induced forcing algebra $\testring\to \algforc_\testring=\algforc
\tensorr \testring$ has a ring section.

\item
For every $\ring \to \testring$, $\testring \in \subcatring $, the
induced forcing algebra $\testring \to \algforc_\testring$ has an
$\testring$-module section.

\item
For every $\ring \to \testring$, $\testring \in \subcatring $, the
induced forcing algebra $\testring \to \algforc_\testring$ is pure.
\end{enumerate}
\end{lemma}
\begin{proof}
This follows from the functoriality of forcing algebras (Remark
\ref{forcingfunctorremark}) and Lemma \ref{forcingequivalent}.
\end{proof}

\begin{proposition}
\label{subcatringadmissible} Let $\submod \mapsto \submod^\clop$
denote a closure operation given by a subcategory $\subcatring$ of
rings. Then this closure operation is admissible with respect to all
ring homomorphisms.
\end{proposition}
\begin{proof} It is clear that $\submod \subseteq \submod^\clop=(\submod^\clop)^\clop$ and that the
closure is persistent for all ring homomorphisms $\ring \to \ringsec
$. By Lemma \ref{testringtest} a forcing algebra $\ring \to
\algforc$ is $\clop$-admissible if and only if for every $\ring \to
\testring$, $\testring \in \subcatring $, the induced mapping
$\testring\to \algforc_\testring = \algforc \tensorr \testring$ has
a ring-section. But an $\ring$-algebra $\algforc$ with this last
property is in general $\clop$-pure. For if $\elem \tensor 1 \in \im
(\submodzweit \tensorr \algforc \to \modulzweit \tensorr  \algforc)$
and if $\ring \to \testring$,
$ \testring \in \subcatring$ is given,
then we have a factorization $\ring \to \algforc \to \testring$
showing that
$\elem \tensor 1 \in \im (\submodzweit \tensorr
\testring \to \modulzweit \tensorr \testring)$.
\end{proof}

We describe some closure operations coming from such a subcategory
of rings.

\begin{example}
Take $\subcatring$ to be the category of fields. The corresponding
ideal closure operation is the radical $\rad(\ideal)$ of the ideal
$\ideal$, see Section \ref{radicalsubsection}. This closure
operation is related to the surjective topology, where the coverings
are given by spec-surjective morphisms.
\end{example}

\begin{example}
Suppose that $\ring$ is a (noetherian) ring and let $\subcatring$ be
the category of (discrete) valuation domains. Then the induced
closure operation is the integral closure, which we will study in
its relation to the submersive topology in
\cite{blicklebrennersubmersion}.
\end{example}

\begin{example}
The weak subintegral closure was introduced by Leahy and Vitulli
(\cite{leahyvitulliweaksubintegralclosure}, \cite{reidvitullimonomial}).
In \cite{brennersemiintegraltest} we show that \emph{crosses},
i.e. one dimensional schemes consisting of two normal components
meeting in one point transversally give a category of test rings for
this closure operation.
The \emph{scheme of axes},
i.e. one dimensional schemes consisting of normal components meeting
transversally in one point and having maximal embedding dimension,
yield also an interesting category of test rings, which induce the
so-called \emph{axes closure},
which we study in its relation to the
continuous closure in \cite{brennercontinuous}.
\end{example}

\begin{example}
Let $\subcatring$ be the category of regular rings. Then the induced
closure operation is the so-called \emph{regular closure}
\cite{hochsterhunekebriancon}.
The regular closure is clearly contained in the integral closure and in the radical.
It does contain the tight closure of an ideal, since tight closure is
persistent and every ideal in a regular ring is tightly closed
\cite[Theorem 2.3 and Theorem 1.3(e)]{hunekeapplication}.
Therefore
$\ideal^* \subseteq \ideal^{\reg} \subseteq \intclo{\ideal}
\subseteq \rad (\ideal)$.
\end{example}

\subsection{The closure operation induced by a Grothendieck topology}
\label{closuregrothendiecksubsection}

\

\medskip
Let $\spax$ denote a scheme endowed with a Grothendieck topology
which is a refinement of the trivial topology (Example
\ref{trivialsite}) (or of the single-handed Zariski topology (Remark
\ref{zariskisinglehanded}); which one should be clear from the
context). This defines a site morphism $\sitemor: \spax_\topo \to
\spax$, and we will specialize now the construction of extension and
contraction described in Section
\ref{extensioncontractionsubsection} to this situation. We could
essentially develop all our topologies as refinements of the Zariski
topology, however, this is sometimes unnecessary ballast, in
particular in the affine case.

\begin{notation}
Let $\spax$ be a scheme endowed with a Grothendieck topology such
that $\mortop : \spax_\topo \to \spax$ is a site morphism. For an
$\O_\spax$-module $\modul$ we denote the module pull-back (the
sheafification) by $\modul_\topo= \mortopmod$
Let $\submodul$ denote $\O_\spax$-modules.
Then the image sheaf of the induced $\O_{\spax_\topo}$-module
morphisms $\sitemortop (\submod)= \submod_\topo \to \sitemortop
(\modul)= \modul_\topo $ on $\spax_\topo$ is called the
\emph{extended submodule} of $\submodul$ in $\modul_\topo$, denoted
by $ \submod^\topo =\submod^{\ext}$. In particular, for an ideal
sheaf $\idealsheaf \subseteq \O_\spax $ we call the image sheaf of
$\idealsheaf_\topo$ inside $\O_\topo $ the \emph{extended ideal
sheaf}, denoted by $\idealsheaf^{\topo}$.

For an $\O_{\topo}$-submodule $\shH \subseteq \modul_\topo$ we call
the push forward sheaf on $\spax$ given by $\opzar \mapsto
\shef^{-1}(\Gamma(\mortoptop (\opzar), \shH))$, where $\shef: \Gamma
(\opzar, \modul ) \to \Gamma(\mortoptop(\opzar),\modul_\topo)$ is
the natural mapping, the \emph{contracted sheaf}, denoted by
$\shH^\contr$.

For $\submodul $ we call the contraction of the extension the
\emph{closure sheaf} or \emph{closure submodule} (induced by the
topology) of $\submod$ inside $\modul$ and denote it by
$\submod^\tocl= (\submod ^{\topo})^\contr$.
\end{notation}

\begin{proposition}
\label{toclprop} Let $\spax$ be a scheme endowed with a Grothendieck
topology $\spax_\topo \to \spax$, let $\submodul$ be $\O_\spax
$-modules. Then the following holds.

\numiii

\begin{enumerate}

\item
The closure submodule $\submod^\tocl$ is an $\O_\spax$-module
{\rm(}but in general not quasicoherent{\rm)}.

\item
We have $\submod \subseteq \submod^\tocl = (\submod^\tocl)^\tocl$.

\item
If $\hommod : \modul \to \modul'$ is an $\O_\spax$-module
homomorphism to another $\O_\spax$-module $\modul'$ and if $\submod'
\subseteq \modul'$ is a submodule, then $\hommod( \submod^\tocl)
\subseteq (\hommod(\submod'))^\tocl$.

\item
Let $\submod \subseteq \submod' \subseteq \modul$ be another
submodule. Then $\submod^\tocl \subseteq (\submod')^\tocl$.

\item
Let $\pi: \modul \to  \modul/\submod$. Then $\submod^\tocl =
\pi^{-1} (0 ^\tocl)$.

\end{enumerate}
\end{proposition}
\begin{proof}
This is a special case of Proposition \ref{sitemorcprop}.
\end{proof}

\begin{remark}
\label{topologyclosureoperation}
We will apply Proposition \ref{toclprop} mainly in the case of a
commutative ring (an affine scheme $\spaxeqspecring$ endowed with
the trivial topology) and an affine Grothendieck topology for
$\ring$ given by a certain class of $\ring$-algebras. For
$\ring$-submodules $\submodul $ we get
$$ \submod \longmapsto \submod^\pre \longmapsto \submod_{\topo} \longmapsto \submod^{\topo}
\longmapsto \Gamma(\spax_\topo, \submod^{\topo}) \longmapsto
\shef^{-1}(\Gamma(\spax_\topo, \submod^{\topo})) \, .$$ We denote
this $\ring$-submodule of $\modul$ by $\submod^\tocl$. For an ideal
$\idealsubring$ we get the assignment
$$ \ideal \longmapsto \ideal^\pre \longmapsto \ideal_{\topo} \longmapsto \ideal^{\topo}
\longmapsto \Gamma(\spax_\topo, \ideal^{\topo}) \longmapsto
\shef^{-1}(\Gamma(\spax_\topo, \ideal^{\topo}))\, $$
which is an ideal $\ideal^\tocl \subseteq \ring$.
In the case that
$\submod=0 \subseteq \modul$, we have
$\submod^\pre=\submod_\topo=\submod^\topo=0$ and the procedure
simplifies to $0^\tocl= \ker (\shef)$, i.e. an element $\elem \in
\modul$ is in $0^\tocl$ if and only if it is annihilated by the
global sheafification homomorphism $\shef: \modul \to
\Gamma(\spax_\topo, \modul_\topo)$.
\end{remark}

\ren{\elet}{{t}}


The following Lemma describes explicitly the closure operation given
by an affine Grothendieck topology.

\ren{\ringcov}{{\ring}}

\begin{lemma}
\label{explicit}
Let $\catopenspax    $ denote an affine
Grothen\-dieck topology on an affine scheme $\spaxeqspecring$. Let
$\submodul $ denote $\ring$-modules. Let $\elem \in \modul$. Then
$\elem \in \submod^\tocl$ if and only if there exists a covering
$\ring \to \ringcov_\indcov$, $\indcovinset$, and elements
$\elen_\indcov \in \submod \tensorr \ringcov_\indcov$ which map to
$\elem \tensor 1$ under the natural homomorphisms $\submod \tensorr
\ringcov_\indcov \to \modul \tensorr \ringcov_\indcov$ for all
$\indcovinset$. If $\catopenspax $ is quasicompact, then we may
restrict to finite index sets.
\end{lemma}
\begin{proof}
This follows from Lemma \ref{explicitsite} applied to $(\spax_\topo,
\catopen_\spax, \O^\pre) \to \spax$. Note that in the notation of
that Lemma we have $\opzar = \spax = \obv =\obcov_\indcov$ for all
$\indcovinset$, $\obv_\indcov = \Spec \ringcov_\indcov$ and
$\Gamma(\obv_\indcov, \O^\pre)=\ringcov_\indcov$.
\end{proof}

\begin{remark}
The condition that $\elem \in \Gamma(\spax_\topo, \submod^\topo)$
does only mean that $\elem$ is locally in the image of
$\submod_\topo \to \modul_{\topo}$. The property that there exists a
(global) element $\elet \in \Gamma(\spax_\topo,\submod_\topo )$
mapping to $\elem$ is a stronger condition, because such an element
$\elet$ is represented by compatible elements $\elet_\indcov \in
\submod \tensorr \ringcov_\indcov$ for some covering $\ring \to
\ringcov_\indcov$, $\indcovinset$. Also this deserves interest, and
there are topologies like the Frobenius topology and the surjective
topology where both properties coincide (see Example
\ref{constructibleradicalexample}, Theorem \ref{surjectiveexact} and
Proposition \ref{frobeniusexact}).
\end{remark}

We are going to establish finiteness properties of the closure
operation induced by a Grothendieck topology.

\begin{lemma}
\label{tensorhilfslemma}
\renewcommand{\eles}{t}
\renewcommand{\finsubmod}{\modul'}
Let $\ring$ denote a commutative ring, $\modul$ an $\ring$-module,
$\alg$ an $\ring$-algebra. Let $\elem \in \modul$ and suppose that
$\elem \tensor 1 =0$ in $\modul \tensorr \alg$. Then there exists a
finitely generated $\ring$-submodule $\finsubmod \subseteq \modul$,
$\elem \in \finsubmod$, and an $\ring$-subalgebra $\alg' \subseteq
\alg$ of finite type such that $\elem \tensor 1 =0$ in $\finsubmod
\tensorr \alg'$.
\end{lemma}
\begin{proof}
For this we have to go back to the definition of the tensor product,
which is constructed by finite sums of formal tensor expressions
$\elev \tensor \elea$ modulo certain relations. That an element in
the tensor product is zero means that it is in the submodule
generated by the relations, and for that only finitely many data are
needed.
\end{proof}

\

\ren{\indgen}{{j}} 

\ren{\indgenset}{{J}}

\begin{corollary}
\label{quasicompactfiniteproperty}
Suppose that $\catopenspax $
is a quasicompact Grothendieck topology on a scheme $\spax$. Then
the closure operation induced by the topology fulfills both
finiteness properties from Definition \ref{finitepropdef}.
\end{corollary}
\begin{proof}
Suppose that $\elem \in \submod^\tocl$, where $\submodul $ are
$\ring$-modules. By Lemma \ref{explicit} there exists a covering
$\coveringring$, $\indcovinset$ finite, such that $\elem \tensor 1
\in \modul \tensorr \ringcovindic$ is the image of $\elen_\indcov
\in \submod \tensorr \ringcovindic$.

For the first finiteness property we may assume that $\submod =0$.
Since $\elem \in 0^\tocl$ we have that $\elem \tensor 1 =0$ in
$\modul \tensorr \ringcovindic$ for every $\indcovinset$. By Lemma
\ref{tensorhilfslemma} there exist finitely generated
$\ring$-submodules $\modul_\indcov \subseteq \modul$ such that
$\elem \in \modul_\indcov$ and such that $\elem_\indcov=0$ in
$\modul_\indcov \tensorr \ringcovindic$. Then we can take the
finitely generated $\ring$-submodule $\ring \langle \modul_\indcov,
\indcovinset \rangle \subseteq \modul$.

For the second finiteness property there exist finitely many
elements $\eley_\indgen \in \submod$ and $\eler_\indgen \in
\ringcovindic$, $\indgen \in \indgenset_\indcov$, such that
$\elen_\indcov= \sum_{\indgen \in \indgenset_\indcov} \eley_\indgen
\tensor \eler_\indgen$. This means then that $\elem$ is already in
the closure of the finitely generated submodule $\finsubmod = \ring
\langle \eley_\indgen, \indgen \in \indgenset_\indcov, \indcovinset
\rangle \subseteq \submod$.
\end{proof}

\ren{\ringsec}{{S}}

\begin{theorem}
\label{topclosure}
Let $\grocont$ denote a subcategory of rings, let
$\catopen_\ring$ denote an affine Grothen\-dieck topology of finite
type on $\specring$, $\ring \in \grocont$, and suppose that every
ring homomorphism $ \ring \to \ringsec$ in $\grocont$ induces a site
morphism $(\Spec \ringsec)_\topo \to (\specring)_\topo$. Then the
closure operation $\tocl$ induced by the Grothendieck topology is a
$\grocont$-admissible closure operation for all submodules.
\end{theorem}
\begin{proof}
That $\tocl$ is a closure operation was mentioned in Proposition
\ref{toclprop}(ii).

Admissibility.
Let $\elem \in \submod^\tocl$ and let $\algforc$
denote a forcing algebra (of finite type) for
$(\modul,\submod,\elem)$, where $\submodul$ is a finitely presented
$\ring$-module and where $\elem \in \submod^\tocl$. By Lemma
\ref{explicit} there exists a covering $\coveringring$,
$\indcovinset$, such that $\elem \tensor 1 \in \modul \tensorr
\ringcovindic$ lies in the image of $\submod \tensorr
\ringcovindic$. Then by the universal property (Proposition
\ref{universalforcing}) of the forcing algebra we have
factorizations $\ring \to \algforc \to \ringcovindic$ for all
$\indcovinset$. We have to show that $\algforc$ is $\tocl$-pure. So
let $\submodulzweit$ and $\elemzweit \in \modulzweit$ be given and
suppose that $\elemzweit \tensor 1$ is in the image of $\submodzweit
\tensorr \algforc \to \modulzweit \tensorr \algforc$. Then it is
also in the image of $(\submodzweit \tensorr \algforc)
\tensor_\algforc \ringcovindic \cong \submodzweit \tensorr
\ringcovindic$ in $\modulzweit \tensorr \ringcovindic$ for every
$\indcovinset$, and so $\elemzweit \in \submodzweit^\tocl$ by Lemma
\ref{explicit}.

$\grocont$-persistence.
Again, let $\submodul$ and $\elem \in\modul$.
Let $\ring \to \ringsec$ denote a ring homomorphism in
$\grocont$ with corresponding site morphism $(\Spec \ringsec)_\topo
\to (\specring)_\topo$. The condition $\elem \in \submod^\tocl$
means by Lemma \ref{explicit} that there exists a covering
$\coveringring$, $\indcovinset$, such that there exist
$\elen_\indcov \in \submod \tensorr \ringcovindic$ mapping to $\elem
\tensor 1$ in $\modul \tensorr \ringcovindic$. The pull-back
$\ringsec \to \ringsec_\indcov := \ringcovindic \tensorr \ringsec$,
$\indcovinset$, defines a covering of $(\Spec \ringsec)_\topo$. For
every $\indcovinset$ we have the commutative diagram
$$
\begin{CD}
\submod \tensorr \ringcovindic  @>>>  \modul \tensorr  \ringcovindic  \\
@V \tensorr \ringsec VV  @VV \tensorr \ringsec V     \\
\submod \tensorr  \ringcovindic \tensorr  \ringsec  @>>> \modul
\tensorr \ringcovindic \tensorr \ringsec \, .
\end{CD}
$$
Note that $\submod \tensorr  \ringcovindic \tensorr \ringsec \cong
\submod \tensorr \ringsec_\indcov \cong (\submod \tensorr \ringsec)
\tensor_\ringsec \ringsec_\indcov$ (and also for $\modul$), so the
diagram shows that also $\elem \tensor 1 \in \modul \tensorr
\ringsec $ has the property that $(\elem \tensor 1) \tensor 1  \in
(\modul \tensorr \ringsec) \tensor_\ringsec \ringsec_\indcov$ is in
the image of $(\submod \tensorr \ringsec) \tensor_\ringsec
\ringsec_\indcov$. This means by Lemma \ref{explicit} that $\elem
\tensor 1 \in (\submod_\ringsec)^\tocl$, where $\submod_\ringsec$
denotes the image module under $\submod \tensorr \ringsec \to
\modul\tensorr \ringsec$.
\end{proof}

\ifthenelse{\boolean{book}}{
\begin{bookexample}
If the empty scheme $\emptyset= \Spec 0$ is a covering in a
Gro\-then\-dieck topology, then the closure of a submodul
$\submodul$ is of course $\submod^\tocl = \modul$.
\end{bookexample}}{}

As mentioned in Lemma \ref{rightexactpullback} the sheafification
with respect to a Grothendieck topology is right exact. In general
it is not left exact (Example
\ref{radicalnotleftexact}). It is however exact for flat
topologies and for hereditary torsion theories, see Section
\ref{torsiontheorysubsection}.

\ren{\shem}{{\modm}}

\ren{\sher}{{\modr}}

\ren{\shel}{{\modl}}

\begin{lemma}
\label{leftexact}
If the sheafification in a Grothendieck topology
$\spax_\topo$ is left-exact, then also the closure operation for the
$0$-submodule is left exact.
\end{lemma}
\begin{proof}
Let $\olmrosheaf$ be an exact sequence of $\O_\spax$-modules.
Suppose that $\elem \in 0^\tocl \subseteq \Gamma(\opzar,\shem)$ maps
to $0$ in $\Gamma(\opzar,\sher)$. Then $\elem$ is the image of an
element $\elel \in \Gamma(\opzar,\shel)$. Since by assumption
$\Gamma(\opzar,\shel_\topo) \to \Gamma(\opzar,\shem_\topo)$ is
injective and $\shef(\elem)=0$, it follows that $\shef(\elel)=0$ in
$\Gamma(\opzar, \modul_\topo)$, hence $\elel \in 0^\tocl$ in
$\shel$.
\end{proof}

\ren{\shem}{{\sheaf }}

\ren{\sher}{{\shE}}

\ren{\shel}{{\shG}}

\ren{\elea}{{h}}

An absolute stalk in the affine single-handed case contains the full
information about the closure operation. It is however in general
not so easy to find an absolute stalk explicitly.

\ren{\ringsec}{{\ring'}}

\begin{corollary}
\label{topclosureabsolute}
Let $\spaxeqspecring$ and let an affine
single-handed Grothen\-dieck topology $\catopenspax $ be given on
$\spax$ with induced closure operation $\tocl$ on submodules.
Suppose that there exists an absolute filter
$\filtcatindtocatopenspax$ with absolute stalk
$\ring \to \stalkabs
= \colim_{\indmincat} \Gamma(\ring_\indm, \O_\topo)$.
Let $\submodul$ denote $\ring$-modules, $\elem \in \modul$. Then $\elem \in
\submod^\tocl$ if and only if
$\elem\tensor 1 \in \im(\submod
\tensorr  \stalkabs \to \modul \tensorr \stalkabs)$.
\end{corollary}
\begin{proof}
Suppose that $\elem \in \submod^\tocl$.
This means by Lemma
\ref{explicit} that there exists a cover $\ring \to \ringsec$ such
that $\elem \tensor 1 \in \im(\submod \tensorr \ringsec \to \modul
\tensorr \ringsec)$. By the definition of an absolute filter, there
exists $\indmincat$ and $\ringsec \to \ring_\indm \to
\colim_{\indmincat} \ring_\indm = \colim_{\indmincat}
\Gamma(\ring_\indm, \O_\topo) = \stalkabs$ (Corollary
\ref{stalkpresheaf}). Hence $\elem \tensor 1 \in \im(\submod
\tensorr \stalkabs \to \modul \tensorr \stalkabs)$. Suppose now on
the other hand that this holds. This means that there exist elements
$\elen_\indgen \in \submod$ and $\elea_\indgen \in \stalkabs$,
$\indgeninset$ finite, such that $\elem \tensor 1 =
\sum_\indgeninset \elea_\indgen \elen_\indgen$ in $\modul \tensorr
\stalkabs$. The $\elea_\indgen$ are represented in some
$\ring_{\indm_\indgen}$ and hence also in one $\ring_\indm $,
$\indmincat$. Since we have a filter, the equality $\elem \tensor 1
= \sum _\indgeninset \elea_\indgen \elen_\indgen$ in the colimit
holds already in some covering $\ring_\indm \to \ring_\indl$, hence
$\elem \in \submod^\tocl$ by Lemma \ref{explicit}.
\end{proof}


\ifthenelse{\boolean{book}}{

\begin{bookremark}
Suppose that $c_1$ and $c_2$ are two closure operations (defined on
the same category of submodules). Then also their intersection is a
closure operation. In case that they are induced by Grothendieck
topologies $\catopen_1$ and $\catopen_2$ their intersection is
induced by the Grothendieck topology where $W_j \to U$, $j \in J$,
is a covering if there exist subsets $J_1, J_2 \subseteq J$ such
that $W_j \to U$, $j \in J_1$, is a covering in $\catopen_1$ and
$W_j \to U$, $j \in J_2$, is a covering in $\catopen_2$.
\end{bookremark}}{}

\subsection{Pure Grothendieck topologies}
\label{puregrothendiecksubsection}

\

\medskip
If a Grothendieck topology is given on a subcategory $\grocont$ of
rings, then there will be certain rings $\ring \in \grocont$ such
that the induced closure operation is trivial, like $\frobf $-pure
rings in the context of Frobenius closure/Frobenius topology or
$\frobf$-regular rings in the context of tight closure. Here we
characterize the Grothendieck topologies where the induced closure
operation is trivial.

\begin{definition}
\label{puretopdef}
We call a family of scheme-morphisms
$\morcov_\indcov: \obcov_\indcov \to \spax$, $\indcovinset$,
\emph{pure} if for every quasicoherent $\O_\spax $-module $\modul$
the natural module homomorphism $\Gamma(\spax, \modul) \to
\bigoplus_{\indcovinset} \Gamma(\obcov_\indcov,
\morcov_\indcov^*(\modul))$ is injective.

We call a Grothendieck topology on a scheme $\spax$ \emph{pure} if
every covering of $\spax$ in the topology is pure.
\end{definition}

\begin{remark}
An affine single-handed Grothendieck topology is pure if and only if
every cover $\ring \to \ringsec$ is a pure homomorphism.
\end{remark}

\begin{proposition}
\label{pureequivalent} Let $\spax_\topo$ be an affine scheme endowed
with a Gro\-then\-dieck topology. Then the following are equivalent.

\numiii

\begin{enumerate}

\item
The induced closure operation on submodules is trivial.

\item
The sheafification homomorphism $\shef : \modul \to
\Gamma(\spax_\topo, \modul_\topo)$ is injective for every
$\ring$-module $\modul$.

\item
The Grothendieck topology is pure.

If the Grothendieck topology is quasicompact, then this is also
equivalent to

\item
The induced closure operation on submodules of finitely generated
modules is trivial.

Moreover, if $\spax$ is noetherian, then this is also equivalent to

\item
The induced closure operation on submodules in artinian modules is
trivial.
\end{enumerate}
\end{proposition}
\begin{proof}
(i) $\Leftrightarrow $ (ii) is clear from Proposition
\ref{toclprop}(v). (ii) $\Rightarrow$ (iii). Suppose that $\elem \in
\modul$ becomes $0$ in every
$\Gamma(\obcov_\indcov,\morcov_\indcov^*\modul)$, where
$\morcov_\indcov: \obcov_\indcov \to \spax$, $\indcovinset$, is a
covering. Then this is also true for the element $\elem$ considered
in $\Gamma(\spax_\topo, \modul_\topo)$. Hence $\elem=0$ by (ii).

(iii) $\Rightarrow$ (ii). If $\elem = 0$ in $\Gamma(\spax_\topo,
\modul_\topo)$, then it must be $0$ on the presheaf level on some
covering (by the definition of sheafification).

The equivalence (i) $\Leftrightarrow $ (iv) follows from Corollary
\ref{quasicompactfiniteproperty}. Suppose now that $\spax$ is also
noetherian, and that (v) holds. Let $\elem \in \modul$, a finitely
generated module, $\elem \neq 0$. Then $\elem \not\in \fom^\expor
\modul$ for some power of some maximal ideal $\fom$. Therefore
$\overline{\elem} \neq 0$ in the artinian module $\modul/\fom^\expor
\modul$ and therefore $\overline{\elem} \not\in 0^\tocl$ by (v).
Hence $\elem \not\in (\fom^\expor \modul)^\tocl$ in $ \modul$ by
Proposition \ref{admissibleproperties}(iii) and in particular
$\elem \not\in 0^\tocl$ in $\modul$ by Proposition
\ref{admissibleproperties}(i) (or Proposition \ref{toclprop}).
\end{proof}

\begin{remark}
In a pure topology the identity
$\Gamma(\spax_\topo, \modul_\topo)=\modul$ does not hold in general. This is only true
if for every object in $\catopenspax $ the induced topology is pure.
For example, for the $\frobf$-pure rings, i.e. the rings of positive
characteristic which are pure in the Frobenius topology (see
Proposition \ref{frobeniuspureequivalent} and Theorem
\ref{frobeniusglobal}) we have $\Gamma(\ring_\frob, \O_\frob)=
\ring^\perf$.
\end{remark}

\ren{\ringsec}{{S}}

\begin{example}
\label{fullpuretopology} In the pure affine single-handed
Grothendieck topology only the pure ring homomorphisms $\ring \to
\ringsec$ are covers. Then $\modul \to \modul \tensorr \ringsec$ is
by definition injective for all $\ring$-modules $\modul$. Hence
$(\modul^\pre)_0=0$ and the associated presheaf $\modul^\pre
=(\modul^\pre)_1$ is already separated. Therefore the closure
operation induced by the pure topology is trivial. As soon as an
affine single-handed Grothendieck topology encompasses one non-pure
covering $\ring \to \ringsec$ the closure operation is not trivial
anymore.
\end{example}

\subsection{Grothendieck topologies which induce given closure operations}

\label{grothendieckconstructionsubsection}

\

\medskip In this section we show that every admissible closure operation
for ideals (or submodules) comes from the closure operation induced
by a suitable (not uniquely determined) affine single-handed
Grothen\-dieck topology.

\begin{theorem}
\label{correspondence} Let $\grocont$ denote a subcategory of
commutative rings which is closed under the formation of forcing
algebras. Suppose that we have an admissible $\grocont$-persistent
closure operation $\submod \mapsto \submod^\clop$ defined for
submodules in a certain category $\modcat$ of finitely presented
modules. Then there exists an affine single-handed Grothendieck
topology on every $\spaxeqspecring$, $\ring \in \grocont$, such that
the induced closure operation is $\submod^\tocl = \submod^\clop$ and
such that the homomorphisms $\ring \to \ringsec$ in $\grocont$ yield
site morphisms $(\Spec \ringsec)_\topo \to ( \specring)_\topo$.
\end{theorem}
\begin{proof} We look at the following topology on $\spax = \Spec \ring$, $\ring \in \grocont$.
Let the objects in $\catopenspax$ be given by (morphisms which are
isomorphic to) compositions
$$ \spax_\indk \ltodots \spax_1 \lto \spax\, ,$$
where each $\spax_{\indi+1} \to \spax_\indi$ is an admissible
forcing algebra for the closure operation. We declare such
compositions to be the covers in $\catopen_\spax$ and the only
morphisms.

We have to check that this gives indeed a Grothendieck topology. The
identity $\spax \to \spax$ is the forcing algebra of the zero-data.
The coverings are of course closed under composition. If $\spay \to
\spax$ is a $\grocont$-morphism (in particular if $\spay \in
\catopen_\spax$), then we set $\spay_\indi = \spay \timesx
\spax_\indi $ and because of $ \spay_{\indi-1}
\times_{\spax_{\indi-1}} \spax_\indi =(\spay \timesx
\spax_{\indi-1}) \times_{\spax_{\indi-1}} \spax_\indi = \spay
\timesx \spax_\indi = \spay_\indi$ the chain
$$ \spay_\indk \ltodots \spay_1 \lto \spay \,  $$
is again a composition of forcing algebras over $\spay$. Note that
these forcing algebras are again $\clop$-admissible due to the
$\grocont$-persistence of the closure operation. Therefore
$\catopen_\spax$ is closed under products, the covers respect base
changes and every $\grocont$-morphism yields a site morphism.

We have to show that the closure operation $\tocl$ given by this
topology is the same as the given closure operation $\clop$. Let
$\submodul $, $\elem \in \modul$, $\quotmod$ finitely presented.
Suppose first that $\elem \in \submod^\clop$. Then the forcing
algebra $\ring \to \algforc$ for these data is admissible and yields
cover (by definition of the topology) such that $\elem \tensor 1$ is
in the image of $\submod \tensorr  \algforc$. Therefore $\elem \in
\submod^\tocl$ by Lemma \ref{explicit}.

On the other hand, if $\elem \in \submod^\tocl$, then by Lemma
\ref{explicit} there exists a covering, that is, a composition
$$\ring \lto \algforc_1 \ltodots \algforc_\indk$$
of admissible forcing algebras, such that $\elem \tensor 1$ is in
the image of $ \submod \tensorr  \algforc_\indk  \to \modul \tensorr
\algforc_\indk$. Since by assumption the admissible forcing algebras
are $\clop$-pure and since this is true for a composition
(Proposition \ref{admissiblecomposition}(ii)), it follows that
$\elem \in \submod^\clop$.
\end{proof}

\begin{remark}
Different topologies may yield the same closure operation. The
(somewhat artificial) topology constructed in the proof of Theorem
\ref{correspondence} is the coarsest topology such that the
admissible forcing algebras are covers. In general we will deal with
other more natural topologies which yield the same closure
operation.
\end{remark}

\ifthenelse{\boolean{book}}{For the identical closure operation
$\submod^\clop=\submod$ the construction yields as coverings all
compositions of forcing algebras which admit a ring section. But of
course also the trivial topology consisting only of the isomorphic
coverings $\spax \to \spax$ yield this closure operation. Also the
\'{e}tale and the flat Grothendieck topology yield the trivial
closure operation, as well as the pure Grothendieck topology
(Example \ref{fullpuretopology}).

\begin{remark}
Instead of compositions of $\clop$-admissible forcing algebras one
can take all ring homomorphisms as covers which are universally
(after base change in $\grocont$) $\clop$-pure to get a Grothendieck
topology with induced closure operation $\tocl =\clop$. This would
yield for the identical closure operation all pure ring
homomorphism.
\end{remark}}{}

\ren{\ringcov}{{S'}}

\ren{\homtest}{{\psi}}

\begin{proposition}
\label{testringcatclcl} Let $\subcatring$ be a subcategory of test
rings in the sense of Definition \ref{zielringcategory} and let
$\ring$ be a commutative ring. For $\spax =\Spec \ring$ let
$\catopen_\spax$ be the full subcategory of finitely generated
$\ring$-algebras. Let $(\spax, \catopen_\spax)$ be endowed with an
affine single-handed Grothendieck topology, where the covers
$\morcov : \ringsec \to \ringcov $ are given by one of the following
conditions.

\numiii

\begin{enumerate}

\item
For every ring homomorphism $\homtest  : \ringsec \to \ringtest$,
$\ringtestincat$, the induced ring homomorphism
$\morcov_\ringtest:\ringtest \to \ringcov \tensor_\ringsec\ringtest$
has a ring section.

\item
For every ring homomorphism $\homtest  : \ringsec \to \ringtest$,
$\ringtestincat$, the induced ring homomorphism
$\morcov_\ringtest:\ringtest \to \ringcov \tensor_\ringsec\ringtest$
is an direct summand.

\item
For every ring homomorphism $\homtest  : \ringsec \to \ringtest$,
$\ringtestincat$, the induced ring homomorphism
$\morcov_\ringtest:\ringtest \to \ringcov \tensor_\ringsec\ringtest$
is pure {\rm(}so $\morcov$ is a $\subcatring$-pure
homomorphisms{\rm)}.
\end{enumerate}
Then for finitely presented $\ring$-modules $\submod \subseteq
\modul$ we have $\submod^\topocl = \submod^\subcatringcl$. If the
closure induced by $\subcatring$ fulfills the first finiteness
condition then this holds for arbitrary $\ring$-modules.
\end{proposition}
\begin{proof}
Note that the allowed covers are different, but the allowed forcing
algebras in the topology are the same under (i), (ii) and (iii)
because of Lemma \ref{forcingequivalent}. Hence the three
Grothendieck topologies induce the same closure operation, and we
restrict to the pure case (iii).

\ren{\ringcov}{{A}}

Suppose that $\elem \in \submod^\topocl$. This means by Lemma
\ref{explicit} that there is a cover $\ring \to \ringcov$ such that
$\elem \tensor 1 \in
\im ( \submod \tensorr \ringcov \to  \modul \tensorr \ringcov)$. This means for
$\homtest: \ring \to \testring$,
$\ringtest \in \catringtest$, that
$\homtest(\elem) \in
\im (\submod \tensorr \ringcov \tensorr \ringtest
 \to \modul \tensorr\ringcov \tensorr \ringtest)$.
Since
$\ringtest \to \ringcov \tensorr \ringtest$ is pure it follows that already
$ \homtest (\elem)
\in \im ( \submod \tensorr \ringtest  \to \modul \tensorr \ringtest)$.
Hence
$\elem \in \submod^\subcatringcl$.

If, on the other hand,
$\elem \in \submod^\subcatringcl$, then we
consider the (finitely generated) forcing algebra $\algforc$. It
follows that $\ring \to \algforc$ is an allowed covering, since
after base change to
$\homtest: \ring \to \ringtest$ it has even a
ring section. Therefore $\elem \in \submod^\topocl$ again by Lemma
\ref{explicit}.

The last statement is clear since the topological closure fulfills
the finiteness conditions by Corollary
\ref{quasicompactfiniteproperty} and so the identity reduces to the
finite case.
\end{proof}

\ifthenelse{\boolean{book}}{

\begin{remark}
 This is already relevant for $\subcatring$
being the category of fields. For $\subcatring$ being the discrete
valuation domains, the `pure'- viewpoint gives the universal
submersive morphisms, which we adopt in
\cite{blicklebrennersubmersion}, whereas the
`ring-section'-viewpoint excludes certain ramification phenomena.
\end{remark}

\ren{\ringcov}{{A}}

We call a ring homomorphisms $\ring \to \ringcov$
$\subcatring$-\emph{pure} if for every ring homomorphism $\homtest:
\ring \to \testring$, $\testring \in \subcatring$, the homomorphism
$\testring \to \testring \tensorr \ringcov$ is pure. }{}

\subsection{Exactness and closure operations}
\label{exactnessclosuresubsection}

\

\medskip
We relate now for a complex $\lmr$ of $\ring$-modules several
notions of exactness in a given Grothendieck topology with the
property that the kernel is contained in the closure of the image.

\ren{\ringcov}{{\ringsec  }}

\begin{proposition}
\label{closurexexact} Let $\catopenspax    $ denote an affine
single-handed Grothendieck topology on $\spaxeqspecring$. Let $\lmr$
denote a complex of $\ring$-modules. Consider the following
properties.

\numiii
\begin{enumerate}

\item
The complex of sheaves $\lmrr{_\topo}$ on $\spax_\topo$ is exact.

\item
The complex of sheaves $\lmrr{_\topo}$ on $\spax_\topo$ is \xexact.

\item
The complex of presheaves $\lmrr{^\pre}$ on $\spax_\topo$ is
\xexact.

\item
$\ker \cob \subseteq (\im \coa)^\tocl$, where $\tocl$ denotes the
closure operation of the topology {\rm(}inside $\modm ${\rm)}.
\end{enumerate}
Then {\rm(i)} $\Rightarrow$ {\rm(ii)}, {\rm(ii)} $\Rightarrow$
{\rm(iii)}, {\rm(iii)} and {\rm(iv)} are equivalent. Moreover, if
$\modm \to \Gamma(\spax_\topo, \modm_\topo)$ is surjective and
$\modr \to \Gamma(\spax_\topo, \modr_\topo)$ is injective, then also
{\rm(ii)} is equivalent to {\rm(iii)} and {\rm(iv)}.
\end{proposition}
\begin{proof}
(i) $\Rightarrow$ (ii) and  (ii) $\Rightarrow$ (iii) were given in
Proposition \ref{exactexact}.

(iii) $\Rightarrow$ (iv). Let $\elem \in \modm$, $\cob(\elem)=0$. By
(iii) there exists a cover $\ring \to \ringcov$ and $\elel \in \modl
\tensorr \ringcov$ mapping to $\elem \tensor 1$. We may write $\elel
= \sum_{\indgeninset} \elel_\indgen \tensor \elerp_\indgen$,
$\elel_\indgen \in \modl$, $\elerp_\indgen \in \ringcov$. Then of
course $\coa(\elel_\indgen) \in \im (\coa)$ and so $\elem \tensor 1
\in \im( (\im \coa) \tensorr \ringcov \to \modm \tensorr \ringcov)$,
hence $\elem \in (\im \coa)^\tocl$.

(iv) $\Rightarrow$ (iii). Again let $\elem \in \modm$,
$\cob(\elem)=0$, so that $\elem \in (\im \coa)^\tocl$ by (iv). This
means by Lemma \ref{explicit} that there exists a cover $\ring \to
\ringcov$ such that $\elem \tensor 1 \in \im( (\im \coa) \tensorr
\ringcov \to \modm \tensorr \ringcov)$. So in particular there
exists $\elemp = \sum_\indgeninset \elem_\indgen \tensor
\elerp_\indgen$ mapping to $\elem \tensor 1$, where $\elem_\indgen =
\coa(\elel_\indgen)$, $\elel_\indgen \in \modl$ and $\elerp_\indgen
\in \ringcov$. But then also $\elem \tensor 1 =
\coa(\sum_\indgeninset \elel_\indgen \tensor \elerp_\indgen)$, so
locally $\elem$ is in the image of the presheaf given by $\modl$.

Now suppose that (iii) holds and that $\modm \to \Gamma(\spax_\topo,
\modm_\topo)$ is surjective and $\modr \to \Gamma(\spax_\topo,
\modr_\topo)$ is injective. Then $\elem \in \Gamma(\spax_\topo,
\modm_\topo)$ mapping to $0$ comes from $\elem \in \modm$, which
also maps to $0$ in $\modr$. So there exists a cover $\ring \to
\ringcov$ and $\elelp \in \modl \tensorr \ringcov$ mapping to $\elem
\tensor 1$ in $\modm \tensorr \ringcov$. Then also $\elelp \in
\Gamma( \ringcov, \modl_\topo)$ maps to the restriction of $\elem$
in $\Gamma(\ringcov, \modm_\topo)$.
\end{proof}

\begin{corollary}
\label{closurexexactcor}
Let $\catopenspax $ denote an affine
single-handed Grothendieck topology on $\spaxeqspecring$.
Let $\lmr$
denote a complex of $\ring$-modules. Then the following are
equivalent. \numiii
\begin{enumerate}
\item
The complex of sheaves $\lmrr{_\topo}$ on $\spax_\topo$ is exact.


\item
The complex of presheaves $\lmrr{^\pre}$ on $\spax_\topo$ is
$\opentop$-exact for every $\opentop \to \spax$ in $\catopenspax $.
\item
For every $\ring \to \ringtop  $ in $\catopenspax    $ the complex
of $\ringtop$-modules $\lmrtensor {\ringtop  }$ has the property
that $\ker \cob \subseteq (\im \coa)^\tocl$, where $\tocl$ denotes
the closure operation induced by the topology {\rm(}inside $\modm
${\rm)}.
\end{enumerate}
\end{corollary}
\begin{proof}
Everything follows from Proposition \ref{closurexexact} except (iii)
$\Rightarrow $ (i). We have to show that the complex of sheaves is
exact for every $\opentop \to \spax$ in $\catopenspax    $, and we
may assume that $\opentop = \spax$. Suppose that $\elem \in
\Gamma(\spax, \modm_\topo)$ maps to $0$ in $\Gamma(\spax,
\modr_\topo)$, and assume that it is represented by
$\elem \in \modm \tensorr \ringtop$, $\ring \to \ringtop$
a cover, mapping to $0$ in
$\modr \tensorr \ringtop$. Then by (iii) we have $\elem \in ({\im
\coa})^\tocl$. This means that there exists by Lemma \ref{explicit}
a cover $\ringtop \to \ringsecsec$ such that
$\elem \tensor 1 \in \modul \tensorr \ringsecsec$
is in the image of
$\modl \tensorr \ringsecsec$. This means that
$\elem \in \Gamma(\spax,\modul_\topo)$
comes locally from the left.
\end{proof}

\ren{\ringsec}{{\ring'}}

\begin{corollary}
\label{catringtestclosureexact}
Let $\catringtest$ be a category of
test rings and let $\ring$ be a commutative ring. Let $\spax= \Spec
\ring$ be endowed with the affine single-handed Grothendieck
topology where the covers are given by $\catringtest$-pure
homomorphisms of finite type. Suppose that $\catopen_\spax$ contains
with every ring $\ringsec$ also the full category of finitely
generated $\ringsec$-algebras. Let $\lmr$ denote a complex of
finitely generated $\ring$-modules. Then the following are
equivalent.

\numiii
\begin{enumerate}
\item
The complex of sheaves $\lmrr{_\topo}$ on $\spax_\topo$ is exact.

\item
For every $\ring \to \ringtop $ in $\catopen_\spax$ the complex of
$\ringtop$-modules $\lmrtensor {\ringtop}$ has the property that
$\ker \cob \subseteq (\im \coa)^\tocl$.

\item
For every ring homomorphism
$\homtest:\ring \to \ringtest$,
$\ringtest \in \catringtest$, we have that the complex $ \lmrr{
\tensorr \ringtest}$ is exact.
\end{enumerate}
If the closure induced by $\subcatring$ fulfills the first
finiteness condition then this holds for arbitrary $\ring$-modules.
\end{corollary}
\begin{proof}
The equivalence of (i) and (ii) was given in Corollary
\ref{closurexexactcor}. Suppose that (iii) holds. We may assume that
$\ring = \ringtop$. Suppose that $\elem \in \modm$ maps to $0$ in
$\modr$ and consider $\im (\coa) \subseteq \modm$. Let $\homtest:
\ring \to \ringtest$ be given, $\ringtestincat$. Then by assumption
$\elem \tensor 1 \in \im (  \modl \tensorr \ringtest \to \modm
\tensorr \ringtest )$ and hence also
$\elem \tensor 1
\in \im ( \im(\coa) \tensorr \ringtest
\to \modm \tensorr \ringtest )$.
But this implies by Lemma \ref{testringcatclcl} that $\elem \in (\im
(\coa))^\topocl$.

Suppose now that (ii) holds and let
$\homtest: \ring \to \ringtest$
be given, $\ringtestincat$. Suppose that
$\elema \in  \modm \tensorr \ringtest$ is mapped to $0$ in $\modr \tensorr \ringtest$. There
exists by Lemma \ref{tensorhilfslemma} a finitely generated
$\ring$-algebra $\ringtop$ with a factorization
$\ring \to \ringtop \to \ringtest$ and an element
$\elemaa \in \modm \tensorr \ringtop$
mapping to $\elema$ and also mapping to $0$ in $\modr \tensorr
\ringtop$. Hence by (ii) $\elemaa \in (\im (\coa))^\topocl$.
Therefore by Lemma \ref{testringcatclcl} we know that
$ \elema
=\elemaa \tensor 1
\in \im ( \im (\coa) \tensor_\ringtop \ringtest
\to (\modm \tensorr \ringtop) \tensor_\ringtop \ringtest)$.
Since
$\modl \tensorr \ringtop \to \im (\coa)$ is surjective, it follows
also that
$\elema \in \im ( (\modl \tensorr \ringtop)
\tensor_\ringtop \ringtest
= \modl \tensorr \ringtest
\to (\modm \tensorr \ringtop) \tensor_\ringtop \ringtest)
= \modm \tensorr \ringtest) $.
\end{proof}

\begin{example}
Let $\spaxeqspecring$ be an affine scheme endowed with a
Grothendieck topology $\catopenspax $. Let $\ideal=(\runfuf)
\subseteq \ring$ denote an ideal with ideal closure $\ideal^\tocl
\neq \ideal$. Then the complex $\ring^\numgen
\stackrel{\runfuf}{\lto} \ring \lto \ring/\ideal^\tocl$ is not exact
as a complex of $\ring$-modules, but the corresponding complex of
presheaves is \xexact, because every element in $\ideal^\tocl$ comes
from the left under some covering.
\end{example}

\begin{remark}
Complexes of $\ring$-modules with the property that the kernel is
inside the closure of the image (as in Proposition
\ref{closurexexact}(iv)) were studied by several people with
respect to different closure operations. Katz studies complexes
acyclic up to integral closure \cite{katzdcomplex} and Hochster and
Huneke study phantom homology within the setting of tight closure
\cite[\S 10]{hunekeapplication}.
\end{remark}

\ifthenelse{\boolean{book}}{

\subsection{Some exact sequences}
\label{someexactsequencesubsection}

\

\medskip}{}

An exact sequence of $\ring$-modules leads in general not to an
exact sequence of sheaves. \ifthenelse{\boolean{book}}{In this
section we describe several consequences from this.}{There are
however some cases when they do.}

\ren{\ringsec}{{\ring'}}

\begin{proposition}
\label{presheafexact} Let $\ring$ denote a commutative ring,
$\olmro$ an exact sequence of $\ring$-modules and let an affine
single-handed Grothendieck topology be given on $\spaxeqspecring$.
Then the complex of sheaves $\olmroo {_\topo}$ is exact in the
following cases.

\numiii
\begin{enumerate}

\item
The sequence splits.

\item
$\modr$ is flat.

\item
$\modl$ is pure in $\modm$.
\end{enumerate}
\end{proposition}
\begin{proof}
We only have to show that the corresponding complex of presheaves is
exact. (i) is clear, since the presheafifcation respects direct
sums. (ii),(iii). For $\ring \to \ringsec$ we have the exact
sequence of presheaves
$$ \to \Tor_1^\ring(\modr, \ringsec) \to \lmroo {\tensorr \ringsec} \, .$$
If $\modr$ is flat, then $\Tor_1^\ring(\modr, \ringsec)=0$, and the
presheaf complex is exact. If $\modl \to \modm$ is pure, then again
$\modl \tensorr \ringsec \to \modm \tensorr \ringsec$ is injective.
\end{proof}

\ifthenelse{\boolean{book}}{
\begin{proposition}
\label{toptopprop} Let $\ring$ denote a commutative ring, let
$\spaxeqspecring$ be endowed with an affine Grothendieck topology
and let $\olmro$ be a short exact sequence of $\ring$-modules. Then
the following holds. \numiii
\begin{enumerate}
\item
We have a short exact sequence of $\O_\topo$-modules
$$ 0 \lto \modl^\topo \lto \modm_\topo \lto \modr_\topo \lto 0\, .$$
If $\Gamma(\spax_\topo, \modm_\topo)=\modm$ and $H^1(\spax_\topo,
\modm_\topo)=0$, then we get the exact sequence
$$0 \to \modl^\tocl \to \modm \to \Gamma(\spax_\topo, \modr) \to H^1(\spax_\topo,\modl^\topo)
\to 0 \, .$$

\item
We have a short exact sequence of $\O_\topo$-modules
$$ 0 \lto \shker \lto \modl_\topo \lto \modl^\topo \lto 0\, .$$

\item
If $\modm$ is flat, then $\shker$ is the sheaf associated to the
presheaf $\ringsec   \mapsto \Tor^\ring_1( L, \ringsec  )$, where
$\ring
\to \ringsec  $ is in $\catopenspax    $.

\item
If $H^1(\spax_\topo, \modl_\topo)=0$, then we get
$$ \Gamma(\spax_\topo,\modl^\topo)/ (\im(\Gamma(\spax_\topo, \modl_\topo) \to \Gamma(\spax_\topo,\modl^\topo))
\cong H^1(\spax_\topo, \shT) \, .$$ If also $\Gamma(\spax_\topo,
\modm_\topo)=\modm$ and $\Gamma(\spax_\topo, \modl_\topo)=\modl$,
then $\modl^\tocl/\modl =H^1(\spax_\topo, \shker)$.
\end{enumerate}
\end{proposition}
\begin{proof}
By right exactness of the tensor product we get $ \lmrr {_\topo} \to
0$ and $\modl^\topo$ is by definition the image sheaf of
$\modl_\topo$ inside $\modm_\topo$. (i) and (ii) follow. (iii). In
general we have the long exact sequence
$$
\longto \Tor^\ring_1(\modl, \ringsec  ) \longto \Tor^\ring_1(\modm, \ringsec
)\longto \Tor^\ring_1(\modr, \ringsec  )
$$
$$ \longto \modl \tensorr \ringsec   \longto \modm \tensorr \ringsec    \longto \modr \tensorr \ringsec
\longto 0 \, ,
$$
so $\shker$ is also the image sheaf of the sheafification
$(\Tor^\ring_1(L,-))^\sheafify \to \modl_\topo$.
If $\modm$ is flat,
then $\Tor^\ring_1(\modm,\ringsec  )=0$ and we get the short exact
sequence
$$ 0 \longto \shT = (\Tor_1(L,-))^\sheafify \longto \modl_\topo \longto \modl^\topo \longto 0 \,.$$
(iv). The short exact sequence of (ii) yields the long exact
cohomology sequence (the first line)
$$
\xymatrix{\! 0\! \ar[r]\! &\! \Gamma\!(\spax,\!\shker)\! \ar[r]
\!&\! \Gamma\!(\spax,\!\modl_\topo)\! \ar[r] \!&\! \Gamma\!(\spax,\!
\modl^\topo) \! \ar[r] \!&\!
H^1\!(\!\spax,\!\shker) \! \ar[r] \!&\! H^1\!(\!\spax,\!\modl_\topo )\\
& 0 \ar[u]^{=}\ar[r] & \modl \ar[r]\ar[u]^{=}&\modl^\tocl \ar[r]
\ar[u]^{=}& H^1(\spax,\shker) \ar[r]\ar[u]^{=} &  0 \ar[u]^{=} \,
.}$$ Under the assumptions we have the identities to the bottom row,
and the statements follow.
\end{proof}

\begin{remark}
\label{idealtoptop} Let $\spax_\topo$ denote a Grothendieck topology
on an affine scheme $\spaxeqspecring$, let $\idealsubring$ denote an
ideal. By Proposition \ref{toptopprop}(i) we have the short exact
sequence
$$0 \lto \ideal^\topo \lto \O_\topo \lto (\O/I)_\topo \lto 0\, $$
of sheaves on $\spax_\topo$. Note that on the right we have the
sheafification of the $\ring$-module $\ring/\ideal$ on
$\spax_\topo$. This is in general not the same as the sheafification
of the structure sheaf $\O_Z$ on $Z_\topo$, $Z=V(I)=\Spec
\ring/\ideal$. The point is that a covering $\Spec T \to Z$ (like
the reduction) is in general not extendible to a covering $\Spec
\rings \to \spax$ (it is true for the surjective topology). If
$\ring \to \ring/\ideal$ induces a site morphism $Z_\topo \to
\spax_\topo$, then we get a ring homomorphism $\Gamma(\spax_\topo,
\O) \to \Gamma(Z_\topo, \O)$.

If $\Gamma(\spax_\topo, \O_\topo)=\ring$ and $H^1(\spax_\topo,
\O_\topo)=0$, then we get a short exact sequence $0 \to
\ring/\ideal^\tocl \to \Gamma(\spax_\topo, (\O/I)_\topo) \to
H^1(\spax_\topo,\ideal^\topo) \to 0$.
\end{remark}}{}

\subsection{Sheaf of syzygies and cohomological interpretation}

\label{syzygycohomologysubsection}

\

\medskip
We want to relate the quotient $\submod^\tocl/ \submod$ with the
first cohomology of the syzygy sheaf in the topology. We consider
finitely generated $\ring$-modules $\submodul $ and a homomorphism
$\freisurn: \ring^\numgen \to \modul$ such that the image of
$\freisurn$ is $\submod$. Let a Grothendieck topology $\catopenspax
$ on $\spaxeqspecring$ be given. This yields a surjective
homomorphism of sheaves
$$ \freisurn: \O_\topo^\numgen \lto \submod^\topo \subseteq \modul_\topo \, .$$ The kernel of this sheaf
homomorphism is a sheaf on $\spax_\topo$ which we denote by
$\kersyz$. Note that this kernel sheaf is not the sheafification of
the kernel of the given $\ring$-module homomorphisms.

\begin{proposition}
\label{syzygycohomology} Let $\ring$ denote a commutative ring with
a Grothendieck topology $\catopenspax    $ on $\spaxeqspecring$. Let
$\submodul$ denote finitely generated $\ring$-modules and let
$\ring^\numgen \to \modul$ denote an $\ring$-module homomorphism
with image $\submod$. Let $\kersyz$ denote the kernel sheaf of the
sheaf morphism $\O^\numgen_\topo \to \submod^\topo \subseteq
\modul_\topo$. Then the following holds.

\numiii
\begin{enumerate}
\item
We have a homomorphism $\submod^\tocl\!/\!\submod \!\to \! H^1(\!
\spax_\topo,\!\kersyz)$. If $\Gamma(\spax_\topo,\!
\O_\topo)\!=\!\ring$, $\Gamma(\spax_\topo, \modul_\topo)= \modul$
and $H^1(\spax_\topo, \O_\topo)=0$, then this is an isomorphism.

\item
Suppose that $\openzar \subseteq \spax$ is Zariski open, that
$\catopen_\spax$ is a refinement of the Zariski topology and that
the support of $\quotmod$ is outside of $\openzar$. Set $\closedzar=
\spax- \openzar$. Then we have a commutative diagram
$$
\xymatrix{
\submod^\tocl/\submod  \ar[r] \ar[d] & \quotmod \ar[r] \ar[d] & \modul/\submod^\tocl \ar[d] \\
\Gamma\!(\spax_\topo\! ,\! \submod^\topo)\!/\!\submod\!
\Gamma(\spax_\topo,\O_\topo)\! \ar[r]\! \ar[d]^{\subseteq}\! &\!
\Gamma(\openzar_\topo ,\submod^\topo)\!/\!\submod \!
\Gamma(\openzar_\topo , \O_\topo) \ar[r] \!\ar[d]^{\subseteq} \!
& \!
H^1_\closedzar\!(\!\submod^\topo)\! /\! \!\im\!(\!H^1_\closedzar\!(\!\O_\topo^\numgen)\!)\! \ar[d]^{\subseteq} \! \\
H^1(\spax_\topo, \kersyz) \ar[r] &  H^1(\openzar_\topo, \kersyz)
\ar[r] & H^2_{\closedzar_\topo}(\kersyz) \, .}
$$
\end{enumerate}
\end{proposition}
\begin{proof}
We apply to the short exact sequence of sheaves on $\spax_\topo$,
$$ 0 \lto \kersyz \lto \O_{\topo}^\numgen \lto\submod^{\topo} \lto 0 $$
several cohomological functors, to wit,
$\Gamma(\spax_\topo,-)$, $\Gamma(\openzar_\topo,-)$ and
$\Gamma_{\closedzar_\topo} (-)$, and get the commutative diagram
$$\xymatrix{
\Gamma(\spax_\topo, \O_\topo)^\numgen \ar[r] \ar[d] &
\Gamma(\openzar_\topo, \O_\topo)^\numgen \ar[r] \ar[d] & H^1_\closedzar (\O_\topo)^\numgen  \ar[d] \\
\Gamma(\spax_\topo, \submod^\topo) \ar[r] \ar[d] &
\Gamma(\openzar_\topo, \submod^\topo) \ar[r] \ar[d]&
 H^1_\closedzar(\submod^\topo) \ar[d] \\
H^1(\spax_\topo, \kersyz) \ar[r] \ar[d] & H^1(\openzar_\topo, \kersyz) \ar[r] \ar[d] &
H^2_\closedzar (\kersyz)  \ar[d] \\
H^1(\spax_\topo, \O_\topo)^\numgen \ar[r] & H^1(\openzar_\topo,
\O_\topo)^\numgen \ar[r] & H^2_\closedzar(\O_\topo)^\numgen \, .}$$
Note that the image of
$\Gamma(\spax_\topo, \O_\topo)^\numgen$ in
$\Gamma(\spax_\topo, \submod^\topo)$ is just
$\submod \Gamma(\spax_\topo, \O_\topo)$
(same for the second column).
From this everything in the lower half of the diagram  in (ii) follows.

The first down arrow in the diagram in (ii) on the left comes from
the sheafifying homomorphism $\shef: \modul \to \Gamma(\spax_\topo,
\modul_\topo)$, where by definition $\submod^\tocl$ is the submodule
mapping to
$\Gamma(\spax, \submod^\topo) \subseteq \Gamma(\spax_\topo,\modul_\topo)$.
This gives the first statement in (i).
The other statement in (i) follows from the long exact sequence
(in the first column), since
$\Gamma(\spax_\topo,\modul_\topo)=\modul $ implies that
$\Gamma(\spax_\topo,\submod^\topo)=\submod^\tocl$.

Suppose now that the support of $\modul/\submod$ is outside
$\openzar$, or equivalently that
$\ring^ \numgen \to \modul$ is
surjective on $\openzar$. Then we have
$\submod^\topo= \modul_\topo$
on $\openzar_\topo$. The mapping
$\modul \to\Gamma(\spax_\topo,\modul_\topo)
\stackrel{\rest}{\to} \Gamma(\openzar_\topo, \modul_\topo)
=\Gamma(\openzar_\topo,\submod^\topo)$
induces the first mapping in
the second column and the upper left square commutes.

The induced mapping
$\quotmod \to H^1_\closedzar(\submod^\topo)  /
\im(H^1_\closedzar(\O_\topo)^\numgen)$
sends
$\submod^\tocl \mapsto 0$, so it factors through
$\modul/\submod^\tocl$.
\end{proof}

{\ren{\spau}{{U}}
\begin{example}
\label{onegeneratorexplicit}
We describe explicitly how a class
$[\elem] \in \submod^\tocl/ \submod$ gives a cohomology class
$\classcoho= \homcon ([\elem]) \in H^1(\spax_\topo, \Syz)$
in the easiest case, that of a principal ideal
$\ideal=(\gen)$ in an affine single-handed Grothendieck topology on $\spax =\Spec \ring$. For
$\spau \in\catopen_\spax$
we have the sequence
$0 \to \Gamma(\spau, \Syz) \to \Gamma(\spau, \O_\topo)
\stackrel{\gen}{\to} \Gamma(\spau, \ideal^\topo)$,
so
$\Gamma(\spau, \Syz)= \{ \fuh \in \Gamma(\spau,\O_\topo):\, \fuh \fug=0 \}=\Ann (\fug)$,
which is zero if $\fug $ is a non-zero divisor in
$\Gamma(\spau, \O_\topo)$. If $\fuf \in \ring$ and $\fuf \in
\ideal^\topocl$, then there exists a cover $\spay \to \spax$ and an
element $\eleq \in \Gamma(\spay, \O_\topo)$ such that $\fuf = \eleq
\fug$. Therefore over $\spay \timesx \spay$ the element
$\project_2^*( \eleq) - \project_1^*( \eleq)$ maps to $0$ in
$\Gamma(\spay \timesx \spay, \ideal^\topo)$
and comes from an element in
$\Gamma(\spay \timesx \spay, \Syz)$. This element is a
\v{C}ech representative of $\homcon(\fuf) \in H^1 (\spax, \Syz)$.

The \v{C}ech cohomology of the syzygy sheaf (which is in our case an
annihilator sheaf of one element) with respect to a cover
$\chch^1(\spay \to \spax, \Syz)$
is computed as the homology of the \v{C}ech complex
$\Gamma(\spax, \Syz) \to \Gamma(\spay, \Syz) \to
\Gamma(\spay \timesx \spay, \Syz) \to \Gamma(\spay \timesx \spay
\times \spay, \Syz)$.
We compute an example in the finite topology in Example
\ref{onegeneratorfinite}.
\end{example}
}

\ifthenelse{\boolean{book}}{
\begin{remark}
Suppose in the situation of Proposition \ref{syzygycohomology} that
$\modul/\submod^\tocl \to H^2_{\spay_\topo}(\kersyz)$ is injective.
Then an element $\elem \in \modul$ belongs to $\submod^\tocl$ if and
only if its cohomology class $\delta(\elem) \in H^1(U_\topo,
\kersyz)$ maps to $0$ on the right or equivalently comes from the
left. In particular, the containment $\elem \in \submod^\tocl$ is a
property of the cohomology class $\delta(\elem)$. Note however that
it does depend on $\kersyz$, not only on the restriction
$\kersyz|_U$.
\end{remark}

\begin{remark}
We want to relate the several notions of syzygies for a given
homomorphism $\ring^n \to \modul$, the image being $\submod$. The
short exact sequence of $\ring$-modules $0 \to \ker \to \ring^n \to
\submod \to 0$ yields by Proposition \ref{toptopprop}(i) the exact
sheaf sequence $0 \to \ker^\topo \to \O^n_\topo \to \submod_\topo
\to 0$. The surjection $\submod_\topo \to \submod^\topo$ induces
also a surjection $ \O^n_\topo \to \submod^\topo$, and (since
$\submod^\topo\subseteq \modul_\topo $) the kernel sheaf of this
surjection is the same as $\kersyz$. Thus we get the commutative
diagram
$$ \xymatrix{ & 0 \ar[d] & & \shker \ar[d]^{\subseteq} & \\
0 \ar[r]   & \ker^{\topo} \ar[r] \ar[d] & \O^n_\topo \ar[r] \ar[d]^{=}
& \submod_\topo \ar[r] \ar[d]   & 0   \\
0 \ar[r]   & \kersyz \ar[r] \ar[d] & \O^n_\topo \ar[r]   & \submod^ \topo \ar[r] \ar[d]   & 0  \\
 & \shekok  & & 0  & }$$
Here $\shekok$ is defined as the cokernel of the first column. By
the snake lemma \cite[Exercise A 3.10]{eisenbud} it follows that
$\shker \cong \shekok$. Suppose that
$\Gamma(\spax_\topo,\O_\topo)=\ring$, $\Gamma(\spax_\topo,
\modul_\topo)=\modul$, $\Gamma(\spax_\topo, \submod_\topo)=\submod$,
$H^1(\spax_\topo,\O_\topo)=0$ and moreover $H^1(\spax_\topo,
\submod_\topo)=0$. Then $H^1(\spax_\topo, \Syz)=
\submod^\tocl/\submod$ by Proposition \ref{syzygycohomology}(i), but
also $H^1(\spax_\topo,\shker)= \submod^\tocl/\submod$ by the
sequence in the right column (Proposition \ref{toptopprop}(iv)). The
sequence in the first row shows that $H^1(\spax_\topo,
\ker^\topo)=0$, and the first column yields then the isomorphism $$0
\to H^1(\spax_\topo, \kersyz)= \submod^\tocl/\submod \to
H^1(\spax_\topo,\shekok)= \submod^\tocl/\submod \to $$ followed by
the inclusion $H^2(\spax_\topo,\ker^\topo) \to H^2(\spax_\topo,
\kersyz)$ (we do not know what these groups are).
\end{remark}

\begin{remark}
If $\ideal=(\fgen) \subseteq \ring$ is an ideal, then we denote the
syzygy sheaf by $\Syz (\fgen)$. A natural choice for $U$ is then
$U=D(\fgen) \subseteq \specring$. The assumptions in Proposition
\ref{syzygycohomology} hold in the submersive topology for normal
rings in characteristic zero, see
\cite{blicklebrennersubmersion}
Some other identities depend on the height and on the depth of the
ideal $\ideal$.
\end{remark}}{}

\ren{\varx}{{x}}
\ren{\vary}{{y}}

\begin{remark}
\label{notcohomologyclass}
It is possible that two syzygy modules
are isomorphic as $\ring$-modules and even isomorphic to the
Zariski structure sheaf, but that they define different sheaves of
syzygies in a Grothendieck topology. This is for example the case in
the submersive topology \cite{blicklebrennersubmersion} for
$\Syz(\varx,\vary)$ and
$\Syz (\varx^2,\vary^2)$ on $\spax=\Spec\field[\varx, \vary]$,
which are both Zariski-isomorphic to
$\O_\spax$. This is related to the fact that the containment
$\fuf \in \overline{(\fgen)}$ in the integral closure is not a property of
the Zariski cohomology class
$\delta(\fuf)\! \in \! H^1(D(\runfuf),$ $\Syz(\runfuf))$.
For example, both data sets $(\varx,\vary);1$ and
$(\varx^2, \vary^2);\varx \vary$ yield the same Zariski cohomology
class
$1/\varx \vary \in H^1(D(\varx,\vary) ,\O_\spax)$, but
$1\not\in \overline{(\varx, \vary)}$ and
$\varx \vary \in \overline{(\varx^2,\vary^2)}$.
This is also related to the fact that
there exists no theory of test ideals for the integral closure.
\end{remark}

\ifthenelse{\boolean{book}}{ The closure of $0$ in a cohomology
module.

\begin{definition}
Let $\ring$ be a commutative ring, endowed with a affine
single-handed Grothendieck topology. Let $\modul$ be an
$\ring$-module and let $H^\numcoho_\idealsup(\modul)$ be a local
cohomology module. We say that $c$ belongs \emph{cohomologically} to
the closure of $0$, written $c \in 0^{\ctopc}$, if there exists a
cover $\ring \to \rings   $ such that the image cohomology class
$c'=0$ in $H^\numcoho_{\idealsup S}(\modul \tensorr S)$.
\end{definition}

\begin{remark}
One can also say that $c$ belongs cohomologically to the closure of
$0$ if $c$ lies in the kernel of the natural homomorphism
$H^\numcoho _\idealsup(\modul) \to
H^\numcoho_\idealsup(F_\absol,\modul^\pre) = \colim_\absol
H^\numcoho_\idealsup(\modul \tensorr S)$, where $F_\absol$ is an
absolute filter.
\end{remark}

We compare this new notion with the closure of the submodule $0$ in
$H^\numcoho_\idealsup(\modul)$, considered as an $\ring$-module.

\begin{lemma}
\label{closurecohoclosure}
Let $\idealsup$ be an ideal in a
noetherian ring $\ring$, endowed with an affine single-handed
Grothendieck topology. Let $\modul$ be an $\ring$-module and let
$H^{\numcoho}_\idealsup (\modul)$ be a local cohomology module. Then
$c \in 0^\topc$ implies $c \in 0^\ctopc$. If $\numcoho$ is such that
$\rad(\idealsup) = \rad (f_1 \komdots f_\numcoho)$, then the
converse is also true.
\end{lemma}
\begin{proof}
The containment $c \in 0^\topc$ means that there exists a cover
$\ring \to \rings   $ such that $c \tensor 1 = 0$ in
$H^{\numcoho}_\idealsup (\modul) \tensorr S$. The image cohomology
class $c'$ in $H^\numcoho_{\idealsup S}(\modul \tensorr S)$ is the
image of $c \tensor 1$ under the natural $S$-module homomorphism
$H^\numcoho_\idealsup (\modul) \tensorr S \to H^\numcoho_{\idealsup
S}(\modul \tensorr S)$ \ifthenelse{\boolean{book}}{(see also Remark
\ref{cohomologyfunctoriality})}{}, so it is also $0$.

Let now $f_1 \comdots f_\numcoho$ be elements in $\idealsup$ which
generate $J$ up to radical. Then the exact sequence $ \bigoplus_{j=1
\komdots \numcoho} \modul_{f_1 \cdots \check{f_j} \cdots f_\numcoho}
\to \modul_{f_1 \cdots f_\numcoho} \to H^\numcoho_\idealsup( \modul)
\to 0$ is tensored exactly to the corresponding complex over $S$, so
$H^\numcoho_\idealsup (\modul ) \tensorr S= H^\numcoho_{\idealsup S}
(\modul \tensorr S)$.
\end{proof}

\begin{example}
If $\fom$ is the maximal ideal in a local ring of dimension $d$,
then for the local cohomology modules $H^\numcoho_\fom (\modul)$
with $\numcoho<d$ the two notions do fall apart. For example, let
$\ring$ be a generalized local Cohen-Macaulay ring, which means that
$H^\numcoho_\fom(\ring)$ is finitely generated for $\numcoho < d=
\dim (\ring)$. If $\ring$ is not itself Cohen-Macaulay, then also
$H^\numcoho_\fom (\ring) \neq 0$ for some $0 < \numcoho <d$. In the
surjective topology (but also in the submersive or in the tight
topology) the image of $H^\numcoho_{\fom}(\ring) $ vanishes in
$H^\numcoho_{\fom S}(S)$ for a suitable cover $\ring \to \rings   $,
so in particular $0^\ctopc $ is the whole module. On the other hand,
by Corollary \ref{radicalnotall} we have $ \rad (0) \neq
H^\numcoho_\fom(\ring)$. For $d= \height (\fom) \geq 2$ we also have
$0^\ctopc = H^{d}_\fom (\ring)$ in the surjective and in the
submersive topology (\cite{blicklebrennersubmersion}
), but not in the plus topology or in the tight topology.
\end{example}

\begin{corollary}
Let $\spaxeqspecring$, where $\ring$ is a ring with a maximal ideal
$\fom$ of height $d$. Let $\spax_\topo$ be endowed with a
single-handed Grothendieck topology and with an absolute filter
$F_\absol$ {\rm(}consisting of ring homomorphisms $\ring \to \rings
${\rm)} and absolute stalk $\O_\absol$. Then we have
$$ H^d_\fom(\modul) \tensorr \O_{\absol}= \colim_\absol  H^d_{\fom
S}( \modul \tensorr S) = H^d_\fom (F_\absol, \modul^\pre) \to
H^d_\fom (F_\absol,\modul_\topo) \, .$$
\end{corollary}
\begin{proof}
The equalities on the left are clear from the proof of Lemma
\ref{closurecohoclosure}. The mapping to the right is also clear.
\end{proof}

\begin{proposition}
Let $\spax$ denote a scheme, endowed with a Grothendieck topology
$\spax_\topo$, let $\modul$ be a quasicoherent $\O_\spax$-module
with sheafification $\modul_\topo$. Let $F: \indexcat \to
\catopenspax $ be a filter where every $\obu_\indec \to \obu$ is a
covering. Then we have homomorphisms
$$ \colim_{\indec \in \indexcat} \chch^{\numcoho}( \obu_{\indec, \zar},
\strutosheaf_\indec (\modul)) \longto \colim_{\indec \in \indexcat}
\chch^{\numcoho}( \obu_{\indec, \topo }, \modul_\topo ) =
\chch^\numcoho(F, \modul_\topo) \, .
$$
\end{proposition}
\begin{proof}
Clear.
\end{proof}

\begin{remark}
This homomorphism needs neither be injective nor surjective. We are
in particular interested in the kernel of $H^\numcoho(\spax, \sheaf
)
\to H^\numcoho (F_\abso , \sheaf _\topo)$, where $F_\abso$ is an
absolute filter in a single-handed Grothendieck topology. E.g. for
the top-dimensional local cohomology module, which seems to give an
important invariant of the topology (like $0^*$).

If $B$ is an absolute stalk, then probably $ H^\numcoho(F_\abso,
\sheaf )= H^\numcoho(\spax, \sheaf ) \tensor B$, with a non-trivial
kernel.
\end{remark}

\begin{remark}
We also have an injection
$$H^{i}_\fom (\ring)/0^\ctopc \to \colim_{\ring \to \rings    \mbox{ covering }}
H^{i}_{\fom S}(S) \, ,$$ where the right hand side is the stalk in
an absolute filter of Zariski local cohomology.
\end{remark}

How related is this to $H^{i}_\fom (\O_\topo)$ (the local cohomology
in the topology). We have $H^{i}(V_\zar, \sheaf ) \to H^{i}(V_\topo,
\sheaf _\topo)$ for every $V \in F$. This yields a homomorphism in the
colimit, though $\sheaf $ is only a presheaf.

For an absolute stalk $F$ the colimit of a presheaf equals the
colimit of the associated presheaf. But there will be no relation
like that between $\colim_{U \in F} H^{i}(U_\zar,\struto^*_U (\sheaf
))
\to H^{i}(F,\sheaf _\topo)= \colim_{U \in F} H^i(U_\topo, \sheaf _\topo)$.


\begin{lemma}
Let $F: \indexcat \to \catopenspax    $ be an absolute filter in a
single-handed Grothendieck topology on $\spaxeqspecring$. Let
$\submodul$ be $\ring$-modules and let $\elem \in \modul$. Let $\foa
=\Ann (\modul/\submod)$ and set $U=D(\foa)$. Suppose that $U$ is
such that for every cover $\ring \to \rings   $ we have
$\Gamma(D(\foa S), \O_S)=S$ (at least after normalization -
extreme?). Let $0 \to \Syz \to \O_U^\submodra \to \O_U^\modra \to 0$
be exact and let $\tilde{\elem} \in \Gamma(U, \O_\spax^\modra)$ be
representing $\elem \in \modul$. Let $\delta (\tilde{\elem}) \in H^1
(U, \Syz)=H^2_\foa (\Syz)$ ($\Syz$ is the Zariski module). Then
$\delta(\elem) = 0$ in $\colim_{\ring \to \rings   } H^2_{\foa S} (
\Syz)$ if and only if $\elem \in \submod^\topc$.
\end{lemma}
\begin{proof}
If $\elem \in \submod^\tocl$, then there exists $\ring \to \rings $
such that $\elem \tensor 1 \in \im( \submod \tensorr S \to \modul
\tensor S)$. This means on $D(\foa S)$ that $\tilde{\elem} \in
\Gamma(D(\foa S), \O^\modra)$ comes from the left and so the
cohomology class is $0$. On the other hand, if this cohomology class
is $0$ in the colimit, then it is $0$ for a cover $\ring \to \rings
$. This again means that on $D(\foa S)$ the element $\tilde{\elem}$
comes from the left.
\end{proof}}{}

\ifthenelse{\boolean{book}}{

\subsection{Principal fiber bundles for the syzygy sheaf}
\label{principalsubsection}

\

\medskip
We interpret some of our results in terms of group schemes and
principal fiber bundle. Suppose that $\mata: \ring^\submodra  \to
\ring^\modra \to \quotmod \to 0$ is a presentation of $\submodul$.
We may consider $\mata : \AA_\ring^\submodra \to \AA_\ring^\modra$
as a group homomorphisms between group schemes. The kernel group
scheme $\ker (\mata)$ is then a group scheme, which is given by
$$\syzgroup=\ker (\mata) = \Spec R[\vart_1 \komdots
\vart_\submodra]/(\mata T) \, .$$ We have denoted its sheaf of
sections by $\Syz$ in Section \ref{syzygycohomologysubsection}. An
element $\elem \in \modul$ defines the forcing algebra, on which the
syzygy group scheme acts additively,
$$ \syzgroup \times_{\specring}  \Spec R[\vart_1
\comdots \vart_\submodra]/ (\mata T - \elem) \longto  \Spec
R[\vart_1 \comdots \vart_\submodra]/ (\mata T - \elem)\, ,$$ sending
$(t_1 \komdots t_\submodra ; \tilde{t}_1 \komdots
\tilde{t}_\submodra) \mapsto (t_1 +\tilde{t}_1 \komdots t_\submodra
+ \tilde{t}_\submodra)$. This action is locally trivial in the
Zariski topology on the open set where $\elem \in \submod$. More
generally, if $\ring \to \ringsec  $ is such that $\elem \tensor 1 =
0$ in $\modul/\submod \tensorr \ringsec  $, then the induced action
of $\syzgroup_\ringsec$ is trivial in the sense that the action of
the forcing algebra is equivalent to the action of
$\syzgroup_\ringsec$ on itself.
\begin{corollary}
Let $\ring$ denote a commutative ring, $\submodul$ finitely
presented $\ring$-modules, and let $\mata: \ring^\submodra  \to
\ring^\modra $ denote a presentation of $\quotmod$. Let a
Grothendieck topology be given on $\spaxeqspecring$ with induced
closure operation $\tocl$. Then the spectrum of the admissible
forcing algebra
$$\spay_\elem = \Spec R[\vart_1 \comdots \vart_\submodra]/ (\mata T - \elem), \,
\elem \in \submod^\tocl\, ,$$together with the natural action of the
group scheme $\Spec R[T]/(\mata T)$, is a model for the principal
fiber bundle on $\spax_\topo $ classified by $[\elem] \in
\submod^\tocl/\submod \subseteq H^1(\spax_\topo, \Syz)$.
\end{corollary}
\begin{proof}
Since $\elem \in \submod^\tocl$, there exists a covering $\{R \to
\ring_i\}, i \in I$, such that $ s \tensor 1 \in \im (\submod \tensorr
\ring_i \to \modul \tensorr \ring_i)$. This means that $\spay_\elem$ is
locally trivial with respect to this covering. Hence $\spay_\elem$
is a principal fiber bundle for $\syzgroup$ in the topology.
\end{proof}

\begin{remark}
In general, a $G$-torsor $G \times \spay \to \spay$ on $\spay \to
\spax$ gets trivial when pulled-back to $\spay$. If the spectrum of
the forcing algebra $\spay_\elem$ is for $\elem \in \submod^\tocl$ a
covering in the topology, then it gets trivial on this covering.
\end{remark}}{}

\section{The surjective topology and the radical}

\label{surjectivesection}

In this part we describe our first main example, the radical of an
ideal and its relationship to the surjective Grothendieck topology,
where the coverings are given by affine spec-surjective morphisms.

\subsection{The radical of an ideal and of a submodule}
\label{radicalsubsection}

\

\medskip
The radical of an ideal $\idealsubring$ is by definition
$$ \sqrt{\ideal} = \rad (\ideal)
 = \{ \fuf \in \ring:\, \exists \expon \in \NN
 \mbox{ such that }\,\, \fuf^\expon \in \ideal\}\, .$$
The radical is also the intersection of all prime ideals which
contain $\ideal$ \cite[Corollary 2.12]{eisenbud}. The containment
$\fuf \in \rad (\ideal)$ is thus equivalent with
$\homtest(\fuf) \in \ideal \field$
for all ring homomorphisms
$\homtest: \ring \to \field$ to fields.

\begin{proposition}
\label{radicalequivalent}
Let $\submodul$ denote $\ring$-modules
such that $\quotmod$ is finitely presented, let $\elem \in \modul$.
Then the following are equivalent.

\numiii

\begin{enumerate}

\item
There exists a spec-surjective ring homomorphism
$\ring \to \alg$
{\rm(}of finite type{\rm)} such that
$\elem \tensor 1 \in
\im(\submod \tensorr \alg \to \modul \tensorr \alg)$.

\item
The forcing algebra
$\algforc$ for $(\modul,\submod,\elem)$
induces a surjective morphism
$\Spec \algforc \to \specring$.

\item
For all homomorphisms
$\ring \to \field$ to a field $\field$ we have
$\elem \tensor 1 \in
\im(\submod \tensorr \field \to \modul \tensorr \field)$.
\end{enumerate}
\end{proposition}
\begin{proof}
(i) $\Ra$ (ii).
By the universal property of the forcing algebra
(Proposition \ref{universalforcing})
we have a factorization
$\Spec\alg \to \Spec \algforc \to \specring$,
so also
$\Spec \algforc \to\specring$ is surjective
(and $\algforc$ is of finite type).
(ii) $\Ra$ (i) is clear.

(ii) $\Leftrightarrow$ (iii).
By the functoriality of forcing algebras
(Remark \ref{forcingfunctorremark})
$\algforc \tensorr \field$
is the forcing algebra for
$(\modul\tensorr \field,\submod',\elem\tensor1)$,
$\submod'$ being the image of
$\submod \tensorr \field$ in $\modul\tensorr \field$.
(ii) says that every fiber is non-empty, and this is equivalent to
the property that
$\algforc \tensorr \field \neq 0$
for every $\ring \to \field   $.
This is equivalent to
$\field\to\algforc \tensorr \field $
being a direct summand, hence this is by
Lemma \ref{forcingequivalent}
equivalent to (iii).
\end{proof}

\begin{definition}
\label{radicalmoduledefinition}
For $\ring$-modules $\submodul$ we
define the radical of $\submod$ in $\modul$ to be the submodule
\begin{eqnarray*}
\rad (\submod)\!\!\!\!
 &=&\!\!\!\! \{ \elem \in \modul:
 \elem \tensor 1 \in \im(\submod\tensorr \field \to \modul \tensorr \field)
\mbox{ for all fields } \ring \to \field \} .
\end{eqnarray*}
\end{definition}

\begin{remark}
This definition gives clearly a submodule. If $\submod$ and $\modul$
are such that $\quotmod$ is finitely presented, then the containment
$\elem \in \rad(\submod)$ can be expressed by the equivalent
properties of Proposition \ref{radicalequivalent}. The radical is by
Proposition \ref{subcatringadmissible} an admissible closure
operation with respect to every ring homomorphism. This follows also
from the topological interpretation in Proposition
\ref{surjectiveinduceradical} below.
\end{remark}

\begin{remark}
\label{radicalzero}
Let $\modul$ denote an $\ring$-module and
consider the $0$-submodule. Then the radical of $0$ consists of all
elements $\elem \in \modul$ such that $\elem \tensor 1 =0$ in
$\modul \tensorr \field $ for all fields $\ring \to \field   $.
\end{remark}

\ren{\ind}{{i}}

\begin{proposition}
\label{radicalfiniteness}
Suppose that $\ring$ is a noetherian ring of finite dimension. Then
the radical fulfills both finiteness properties from Definition
\ref{finitepropdef}.
\end{proposition}
\begin{proof}
Let $\submodul$ be $\ring$-modules, $\elem \in \modul$ and $\elem
\in \rad( \submod)$. For the first finiteness property we may assume
that $\submod=0$. Let $\fop_\ind$, $ \indinset$, denote the finitely
many minimal prime ideals of $\ring$. For all $\ind$ we have $\elem
\tensor 1 = 0$ in
$\modul \tensorr \quotfield( \ring/ \fop_\ind)$.
Then there exist by Lemma \ref{tensorhilfslemma} finitely generated
$\ring$-submodules $\modul_\ind \subseteq \modul$ such that $\elem
\in \modul_\ind$ and $\elem \tensor 1 = 0$ in $\modul_\ind \tensorr
\quotfield(\ring / \fop_\ind)$. The same holds for the finitely
generated submodule $\tilde{\modul}$ generated by the $\modul_\ind$.
This means that there exist $\fuf_\ind \not\in \fop_\ind$ such that
$ \elem \tensor 1 = 0 $ in $\tilde{\modul} \tensorr
(\ring/\fop_\ind)_{\fuf_\ind}$. Now if $\ring \to \field   $ is a
homomorphism to a field $\field$, then either $\fuf_\ind$ maps to a
unit for some $\ind$ or all $\fuf_\ind$ map to $0$. In the first
case the factorization $\ring \to \ring_{\fuf_\ind} \to \field $
shows that $\elem \tensor 1 =0$ in $\tilde{\modul} \tensorr \field$.
In the second case we look at the minimal primes of
$\ring/(\fuf_\ind, \indinset)$ as before, which has smaller
dimension, so we can (increasing the submodule $\tilde{\modul}$)
argue inductively.

The proof of the second finiteness condition is similar.
\end{proof}

\begin{definition}
For an $\ring$-module $\modul$ we denote the factor module by
$\modul_1=\modul/\rad (0)$ and call it the $\ring$-\emph{reduction}
of $\modul$.
\end{definition}

The $\ring$-reduction of an $\ring$-module $\modul$ yields indeed
the separated presheaf in the surjective topology, as we will see in
a minute.

\begin{example}
\label{pidradical}
Let $\ring$ denote a principal ideal domain. Then
every finitely generated $\ring$-module $\modul$ is isomorphic to
$\modul \cong \ring^\rak \times \modul(\prim_1) \timesdots
\modul(\prim_\num)$, where $\prim_1 \komdots \prim_\num$ are non
associated prime elements in $\ring$ and where
$\modul(\prim_\ind)
=(\ring/\prim_\ind)^{\ulm_{\ind 1}} \times
(\ring/\prim_\ind^2)^{\ulm_{\ind 2}} \timesdots
(\ring/\prim_\ind^{\expon_\ind})^{\ulm_{\ind \expon_\ind}}$
is the
$\prim_\ind$-primary component of $\modul$, see \cite[Satz 61.8 and
Satz 61.10]{schejastorch2}. The radical of $0$ is then given by the
direct product (Proposition \ref{admissibleproperties}(iv)) of the
radicals of $0$ in $\ring^\rak$ (which is $0$) and in
$\modul(\prim_\ind)$. The radical of $0$ in
$(\ring/\prim)^{\ulm_{1}}
\times (\ring/\prim^2)^{\ulm_{2}}
\timesdots (\ring/\prim^{\expon})^{\ulm_{\expon}}$
is
$0 \times (\prim \ring/\prim^2)^{\ulm_2}
\timesdots (\prim \ring/\prim^\expon)^{\ulm_\expon}$
and the $\ring$-reduction is
$(\ring/\prim)^{\ulm_1+\ulm_2 \plusdots \ulm_\expon}$.
\end{example}

\begin{example}
Let $\ring = \field[\varx,\vary]$ ($\field$ a field) and lets
consider the $\ring$-module $\modul =
\ring^2/((\varx,\vary),(\varx^2\vary,\varx))$. We claim that the
residue class $\elem$ of $(\varx(1-\vary^2),0)$ is not $0$ in
$\modul$, but it is in the radical of $0$. The determinant of the
generators is $\varx^2(1- \vary^2)$. If $\ring \to \fieldl$ is a
homomorphism to a field, and if the determinant (in $\fieldl$) is
$\neq 0$, then $\modul_\fieldl=0$. If the determinant is $0$, then
also $\elem =0$, so $\elem \in \rad(0)$. Assume now that
$(\varx(1-\vary^2),0) = \fuf(\varx,\vary)+ \fug(\varx^2 \vary,\varx)$
with
$\fuf,\fug \in \ring$.
Then $\fuf=\fuh \varx $ and $\fug=- \fuh
\vary$ by the second column and the first gives
$\varx(1-\vary^2)=\fuh\varx^2-\fuh \varx^2 \vary^2$,
which has no solution.
\end{example}

\ifthenelse{\boolean{book}}{
\begin{bookexample}
Let $\ring$ be a domain with quotient field $\quotfield(\ring)$. Then $\rad
(0)=0$ and $\rad (\submod)=\quotfield(\ring)$ for every
$\ring$-submodule
$\submod \neq 0$. In the first case the tensoration with $\quotfield(\ring)$
yields $0 = 0 \tensorr \quotfield(\ring) \to \quotfield(\ring)
\tensorr \quotfield(\ring) =\quotfield(\ring)$, so an element $0
\neq \elem \in \quotfield(\ring)$ does not belong to $\rad (0)$. In
the second case the tensoration with $\quotfield(\ring)$ yields the
identity $\submod \tensorr \quotfield(\ring) =\quotfield(\ring) \to
\quotfield(\ring)$. The same holds for fields $\quotfield(\ring)
\subseteq \field$, and for other fields we have $\quotfield(\ring)
\tensorr \field =0 $.
\end{bookexample}}{}

\subsection{Surjective topology}
\label{surjectivetopologysubsection}

\

\ren{\obcov}{{V}}

\begin{definition}
\label{surjectivetopology}
Let $\spax$ be a scheme and let $\catopenspax$
be the full subcategory of schemes (or schemes of finite type)
which are affine over $\spax$. The coverings in the (\emph{affine}-)
\emph{surjective topology} on $\catopenspax    $ are the affine
surjections $\obcov \to \obcovp$ of finite type (over $\spax$). We
denote $\spax$ endowed with this topology by $\spax_\surj$.
\end{definition}

\begin{proposition}
\label{surjectivebasicprop}
The surjective topology is an affine single-handed Gro\-then\-dieck
topology, and it is a refinement of the single-handed
Zariski-topo\-lo\-gy. Every scheme-morphism
$\spay\to \spax$ induces a site morphism $\spay_\surj \to \spax_\surj$.
\end{proposition}
\begin{proof}
This follows immediately from the definitions. The last statement
follows since surjections are universally surjective.
\end{proof}

\begin{proposition}
\label{surjectiveinduceradical}
Let $\ring$ denote a ring and let
$\submodul$ denote $\ring$-modules such that $\quotmod$ is finitely
presented. Then the radical of $\submod$ is the closure operation
induced by the surjective topology. If $\ring$ is noetherian of
finite dimension, then this is true for arbitrary submodules
$\submodul$.
\end{proposition}
\begin{proof}
The first statement follows from Proposition \ref{radicalequivalent}
and Lemma \ref{explicit}. The second statement follows from this,
from Corollary \ref{quasicompactfiniteproperty} and from Proposition
\ref{radicalfiniteness}.
\end{proof}

We are going to describe the sheafification in the surjective
topology, in particular for the structure sheaf. We start with the
separated presheaf.

\begin{lemma}
\label{surjseparatedreduction}
Let $\ring$ denote a commutative
ring, $\modul$ a finitely generated $\ring$-module yielding the
presheaf
$\ringsec   \mapsto \modul\tensorr \ringsec $
for $\ring\to \ringsec  $ of finite type. Then
$(\modul \tensorr\ringsec )_0=\rad (0)$
and the $\ringsec  $-reduction
$(\modul \tensorr\ringsec)_1
= \modul \tensorr \ringsec  /(\modul \tensorr \ringsec)_0$
is the associated separated presheaf.
\end{lemma}
\begin{proof}
This follows from Proposition \ref{radicalequivalent} and Section
\ref{grothendieckschemesubsection}.
\end{proof}

\begin{corollary}
\label{surjseppresheaf}
Let $\ring$ denote a commutative ring
endowed with the surjective topology. Then the separated presheaf
$\O_1$ associated to the Zariski structure sheaf is the reduction,
i.e.
$\Gamma(\ringsec  , \O_1) = \ringsec_{\red}$.
\end{corollary}
\begin{proof}
This follows from Lemma \ref{surjseparatedreduction}.
\end{proof}

\begin{remark}
It follows that an element $\fuglob \in \rings$, $\ring \to \rings $
spec-surjective, yields a global section in $\Gamma(\ring,
\O_\surj)$ if and only if $\differ (\fuglob) = \fuglob \tensor 1 - 1
\tensor \fuglob \in \rings \tensorr \rings$ is nilpotent. We will
compute the ring of global sections in Section
\ref{surjectiveglobalsubsection}.
\end{remark}

\subsection{The Jacobson radical and the Jacobson topology}
\label{jacobsonsubsection}

\

\medskip
In this section we describe briefly a variant of the radical and of
the surjective topology, the Jacobson radical and the Jacobson
topology. The \emph{Jacobson radical} of a commutative ring $\ring$
is by definition the intersection of all maximal ideals
\cite[Chapter 4]{eisenbud}.
Accordingly, we make the following definition.

\begin{definition}
The \emph{Jacobson radical} of a submodule $\submodul$ is
$$\jac (\submod)\!
=\!\{ \elem \in \modul:\, \elem \tensor 1 \in
\im( \submod \tensorr \fieldres(\fom)
\to \modul \tensorr \fieldres(\fom)) \mbox{ for all } \fom \mbox{ maximal}\}  .$$
\end{definition}

\begin{remark}
\label{jacobsonremark}
Note that this closure operation is not
persistent with respect to all ring homomorphisms. It is however
persistent with respect to homomorphisms between $\field$-algebras
of finite type over a field $\field$, in which case it is the same
as the radical, see Proposition \ref{radicaljacobson}.
For
$0 \subseteq \modul$ we have
$\jac (0)= \bigcap_{\fom \mbox{ maximal }}\fom \modul$
(note that for the radical it is not true that
$\rad (0)= \bigcap \fop \modul$) and for $\submod$ we have
$\jac (\submod) = \bigcap_{\fom \mbox{ maximal }} (\submod +\fom \modul)$.
\end{remark}

\begin{definition}
Let $\spax$ be a scheme and let $\catopenspax$ be the
subcategory of schemes of finite type which are affine over $\spax$.
The \emph{Jacobson topology} on $\spax$ is the Grothendieck topology
on $\catopenspax $,
where $\morcov: \obcov \to \obcovp$ is a covering
if $\morcov$ maps onto the closed points of $\obcovp$. We denote
$\spax$ endowed with this topology by $\spax_\jac$.
\end{definition}

\begin{remark}
\label{jacobsontopologyremark}
The Jacobson topology is a refinement of the surjective topology.
For $\spaxeqspecring$,
a morphism $\Spec \rings  \to \specring$ is a covering if and only
if every maximal ideal of $\ring$  lies in the image.
For a local ring $\ring$, the morphism
$\Spec \fieldres(\fom) \to \specring$ is a covering,
and more generally, if $\ideal \subseteq \jac $, then
$\Spec \ring/\ideal \to \specring$
is a covering in the Jacobson topology.

A scheme morphisms $\moryx$ which sends closed
points to closed points induces a site morphism between the Jacobson
sites. This holds in particular for morphisms between schemes of
finite type over a field \cite[Corollaire 10.4.7]{EGAIV}.
The closure operation induced by the
Jacobson topology is the
Jacobson radical. Also, $\rad (\submod) \subseteq \jac (\submod)$.
\end{remark}

Recall that a commutative ring is called a \emph{Jacobson ring} if
every prime ideal is the intersection of maximal ideals
\cite[D\'{e}finition 10.4.1]{EGAIV}.

\begin{proposition}
\label{radicaljacobson}
For a noetherian Jacobson ring $\ring$ and
for finitely generated $\ring$-modules $\submod \subseteq \modul$
the identity
$\jac (\submod) =\rad(\submod)$ holds.
\end{proposition}
\begin{proof}
Let $\ring \to \rings   $ be of finite type and suppose that all
maximal ideals of $\ring$
 are in the image of $\Spec \rings \to \specring$. We want to show
that the morphism is surjective. For this we may assume that $\ring$
is a domain and that $\fop =0$. The closure of the image is
$\specring$ due to the Jacobson property. Hence by the Theorem of
Chevalley (\cite[Th\'{e}or\`{e}me 1.8.4]{EGAIV}
or \cite[Corollary 14.7]{eisenbud} the image contains a
non-empty set and in particular it contains $0$. This shows that a
Jacobson cover is spec-surjective. The stated result follows then
from Remark \ref{jacobsontopologyremark} and Proposition
\ref{surjectiveinduceradical}.
\end{proof}

\begin{definition}
\label{faithfullnakayama}
We say that a Grothendieck topology
$\mortop: \spax_\topo \to \spax$ on a scheme $\spax$
is \emph{faithful} (or has the \emph{Nakayama property}) if for
every coherent
$\O_\spax$-module $\sheaf$ the following holds: if
$\mortopmod(\sheaf)=\sheaf_\topo= 0$, then $\sheaf=0$.
\end{definition}

The following Lemma is a version of the Lemma of Nakayama
\cite[Corollary 4.8]{eisenbud}.

\begin{lemma}
\label{nakayamalemma}
The Jacobson topology
{\rm(}and hence the surjective topology{\rm)}
on a noetherian scheme $\spax$ is faithful.
\end{lemma}
\begin{proof}
We may assume that $\spaxeqspecring$ is affine.
Let $\modul$ be a finitely generated $\ring$-module and suppose that
$\modul_\jac =0$.
This means that for every
$\elem \in \modul$ there exists a Jacobson cover
$\Spec \rings \to \specring$ such that $\elem \tensor 1 = 0$ in
$\modul \tensorr \rings$.
Since $\modul$ is finitely generated, there exists also such an
$\rings$ with
$\modul \tensorr \rings = 0$. Since $\ring \to \rings $
is a Jacobson cover, for every maximal ideal $\fom \subset \ring$
the homomorphism
$\ring/\fom \to \rings/\fom \rings$ is faithful, and so
$\modul \tensorr \ring/\fom = \modul/ \fom \modul= 0$
for every maximal ideal.
But then
$\modul_\fom =0$ by the Lemma of Nakayama and so $\modul =0$.
\end{proof}

\begin{corollary}
\label{radicalnotall}
Let $\ring$ denote a noetherian commutative
ring and $\submod \subset \modul$ a submodule such that
$\modul/\submod \neq 0$ is finitely generated. Then $\jac (\submod)
\neq \modul$ and $\rad( \submod) \neq \modul$.
\end{corollary}
\begin{proof}
This follows from Lemma \ref{nakayamalemma}.
\end{proof}

\ifthenelse{\boolean{extra}}{ Let $ \ring^\submodra
\stackrel{\mata}{\lto} \ring^\modra$ be a presentation. $\mata$ does
not define a surjective homomorphism, hence there exists also $\ring
\to \field   $ a field such that the matrix $\mata_K$ over $\field$
is not surjective. This means that the rank of $\mata_K$ is smaller
than $\modra$, and this means that there exists a linear dependence
between the rows. Say that the $i$-th row has non-zero coefficient
in a linear relation. We consider the forcing algebra for
$\tilde{s}=$ the $i$-th unit vector. Then the linear relation over
the field shows that over $\field$ the ideal defining the forcing
algebra is the unit ideal, and so the forcing algebra does not
define a surjective map. }{}

\begin{example}
\label{quotientmodulodvd}
Corollary \ref{radicalnotall} is not true
for modules which are not finitely generated. If $\dvd$ is a
discrete valuation domain, then the quotient module
$\quotfield(\dvd)/\dvd$ has the property that its tensoration both
with the residue class field $\dvd/ \fom_\dvd$ and with
$\quotfield(\dvd)$ is zero, and hence the tensoration with every
field. Hence the radical of $0$ in this
module is the whole module.
\end{example}

\begin{remark}
Let $\ring$ be a principal ideal domain with infinitely many maximal
ideals (like $\ZZ$). Then by Proposition \ref{radicaljacobson} for a
finitely generated $\ring$-module the radical is the same as the
Jacobson radical. This is not true for not finitely generated
modules, and the Jacobson radical does not fulfill the finiteness
property \ref{finitepropdef}(i). For example, for the quotient field
$\quotfield= \quotfield(\ring)$ the radical of $0$ is $0$, since
$\quotfield$ itself is a test ring, but the Jacobson radical of $0$
is $\quotfield$, since
$\quotfield \tensorr \ring/\fom =0$ for every
maximal ideal $\fom$. Every finitely generated submodule
$\modul'\subseteq \quotfield$ is isomorphic to
$\ring^\rak$, and there the Jacobson radical is $0$.
\end{remark}

\begin{lemma}
\label{surjectivefinitetype} If $\spax$ is a scheme of finite type
over an algebraically closed field $\field$, and $\submodul$ are
coherent $\O_\spax$-modules, then $\elem \in \rad(\submod)$ if and
only if this is true for all $\field$-points of $\spax$.
\end{lemma}
\begin{proof}
The surjectivity of the map from the forcing scheme to $\spax$ can
be tested on $\field$-points alone.
\end{proof}

\begin{example}
\label{radicalarbitraryfunctions}
Let $\spax$ be a scheme of finite
type over $\CC$ and consider the ringed site $\spax_\arb $ given by
the point set $\spax(\CC)$ (with the complex topology) and the sheaf
of all (not necessarily continuous) complex-valued functions on
$\spax(\CC)$.
The natural mapping $\spax_\arb \to \spax_\zar$ is a morphism of
ringed sites (in fact, of ringed spaces).
We claim that the closure operation induced by this morphism
(it is not a Grothendieck topology on $\spax$)
equals the radical.
This follows from the simple fact that
$\spay \to \spax$ is surjective if and only if
$\spay(\CC) \to \spax(\CC)$ is surjective,
and the later is equivalent to the existence of a
(not necessarily continuous)
section $\spax(\CC) \to \spay (\CC)$.
Applying this fact to the forcing schemes $\moryx$ gives the claim.
\end{example}

We relate briefly our notion of the radical of a submodule with
another one found in the literature, see
\cite{mccaslandmoore}, \cite{marcelomarcelorodriguez},
\cite{pusatyilmazsmith},
\cite{maanishirazisharif}
(there are also more abstract categorial definitions of a radical).
A submodule $\modprim \subseteq \modul$
of an $\ring$-module is called a \emph{prime submodule} if whenever
$\fuglob \elem \in \modprim$
(where $\fuglob \in \ring$, $\elem \in \modul$),
then $\elem \in \modprim$ or
$\fuglob \modul \subseteq \modprim$.
Maximal submodules
$\modprim \subseteq \modul$,
which are characterized by
$\modul/\modprim \cong \ring/\fom$ for
some maximal ideal $\fom$, are prime. For an arbitrary submodule
$\submod
\subseteq \modul$ one defines the (prime)-radical by
$\varrad(\submod) = \bigcap_{\submod \subseteq \modprim \mbox{ prime }}
\modprim$.

\begin{proposition}
Let $\ring$ denote a noetherian Jacobson ring and let $\modul$
denote a finitely generated $\ring$-module, let $\submod \subseteq
\modul$ be a submodule. Then
$$\rad (\submod) = \jac (\submod)
= \bigcap_{\fom \mbox{ maximal ideal }} (\submod + \fom \modul)
= \bigcap  _{\submod \subseteq \modprim \mbox{ maximal }}
\modprim \supseteq \varrad(\submod) \, .$$
In particular, $\varrad(\submod) \subseteq \rad(\submod)$.
\end{proposition}
\begin{proof}
The first equation is Proposition \ref{radicaljacobson},
the second is the definition
(Remark \ref{jacobsonremark})
and the forth inclusion is trivial,
so we only have to look at the third equation.
We may assume that $\submod =0$.
Assume first that $\elem \in \jac (\submod)$ and let
$\modprim \subseteq \modul$ be a maximal submodule.
Then $\modul/ \modprim \cong \ring/\fom$ as
$\ring$-modules for some maximal ideal $\fom \subset \ring$ and so
$\modul/ \modprim$ is annihilated by $\fom$.
Then $\fom \modul \subseteq \modprim$ and so $\elem \in \modprim$.

Suppose now that
$\elem \in \modprim$ for all maximal submodules
$\modprim \subseteq \modul$, and let $\fom$ be a maximal ideal of
$\ring$.
Then $\modul/\fom \modul$ is a finitely generated
$\ring/\fom = \field$-module,
hence isomorphic to $\field^\expon$
(we may assume $\expon \geq 1$).
The maximal submodules of
$\field^\expon$ are the subspaces of dimension $\expon-1$, and their
intersection is $0$.
Pulling this back to $\modul$ we see that $\fom
\modul$ is the intersection of maximal submodules,
hence $\elem \in \fom \modul$.
\end{proof}

\subsection{The surjective topology over a field}
\label{surjectivefieldsubsection}

\

\medskip
We consider the case of a field $\field$. Although an algebra
$\field \to \alg$ is spec-surjective if and only if it is faithfully
flat (if and only if it is $\neq 0$), the surjective topology is not
the same as the faithfully flat topology. This is because the
compatibility condition in the surjective topology takes into
account the reduction of
$\alg \tensor_\field \alg$.
Recall that the
\emph{perfect closure} of a field $\field$ in characteristic zero is
$\field$ itself and $\field^\perf = \colim_\expoe \frobe \field $ in
positive characteristic, see also Part \ref{frobeniussection}.

\begin{proposition}
\label{globalsurjectivefield}
Let $\field$ denote a field endowed
with the surjective topology. Then
$\Gamma(\field_\surj,\O_\surj) = \field^{\perf} $,
the perfect closure of $\field$.
If $\field$ is a perfect field, then
$\Gamma(\field_\surj,\O_\surj)=\field$.
\end{proposition}
\begin{proof}
Theorem \ref{frobeniusglobal} below shows that an element in the
perfect closure
$\field^\perf$ gives an element in $\Gamma(\field_\surj,\O_\surj)$.
Therefore we may assume that $\field$ is perfect and we have to
prove the second statement.

Let $\field \to \alg$ denote an $\field$-algebra $\alg \neq 0$ of
finite type.
We may assume that $\alg$ is a domain, since any
integral component yields a refined covering.
Suppose that an element
$\fuglob \in \alg$ is compatible,
that is,
$\differ(\fuglob) = \fuglob \tensor 1 - 1 \tensor \fuglob$ is nilpotent
in $\alg \tensor_\field \alg$.
By
\cite[Proposition 4.6.1]{EGAIV} this algebra is reduced,
hence $\differ(\fuglob) =0 $ in $\alg \tensor_\field
\alg$.
Therefore $\fuglob$ is also compatible in the faithfully flat topology,
hence $\fuglob \in \field $ by \cite[Proposition
I.2.18]{milne}.
\end{proof}

\ren{\point}{{x}}

\subsection{Absolute stalks in the surjective topology}
\label{surjectiveabsolutesubsection}

\

\medskip
Let $\ring$ be a ring and set
$\rad(\ring)= \prod_{\point \in \specring}\overline{\fieldres(\point)}$
(we call it the \emph{radicalization} of $\ring$).
The homomorphisms
$\ring \to \ring_{\fop_\point}
\to \fieldres(\point) \to\overline{\fieldres(\point)}$
induce a homomorphism
$\ring \to \rad(\ring)$.

{
\ren{\alg}{{S}}
\ren{\rings}{{\ringt}}
\ren{\vart}{{U}}

\begin{proposition}
\label{surjectiveabsolutestalk}
Let $\ring$ denote a commutative ring,
considered in the surjective topology. Then the fixation given
by $\ring \to \rad (\ring)$ defines an absolute filter and the
radicalization is an absolute stalk.
\end{proposition}
\begin{proof}
Let $\ring \to \alg$ in $\catopen_\spax$ be given with a fixation
$\alg \to \rad (\ring)$. Then $\Spec \alg \to \Spec \ring$ is
surjective, hence a cover. Let $\alg \to \rings$ be a cover in the
surjective topology. We have to show that there exists an
$\alg$-algebra homomorphism $\rings \to \rad(\ring)$, and this can
be checked for every component $\cloalg{\fieldres (\point)}$
separately. For $\point \in \specring$ the geometric fiber ring
$\rings_\point:= \rings \tensor_\alg \overline{\fieldres(\point)}\neq 0$
is of finite type over $\overline{\fieldres(\point)}$,
hence there exists a homomorphism
$\rings_\point \to \overline{\fieldres(\point)}$.
So this quasifilter is absolute. For $\alg \to \rad(\ring)$ the
image $\overline{\alg} = \im(\alg) \subseteq \rad(\ring)$ is also
(spec-surjective and) of finite type over $\ring$, so $\alg \to
\overline{\alg}$ belongs to $\catopenspax $. Hence by Proposition
\ref{filtertest} this is a filter. Since $\ring \to \ring[\vart]$ is
a spec-surjective covering, the corresponding homomorphism $\O_\filt
=\colim_{\indmincat} \alg_\indm \to \rad (\ring)$ is an isomorphism.
\end{proof} }

\begin{example}
Let $\dvd$ denote a discrete valuation domain.
Then by Proposition \ref{globalsurjectivefield} and Corollary
\ref{constructibleglobalsectioncorollary} we have
$\Gamma(\dvd,\O_\surj)
= \quotfield(\dvd)^\perf \oplus (\dvd/\fom)^\perf
=: \field \oplus \fieldmax$.
An ideal $\ideal= (\pipo)$ extends in
$\Gamma(\dvd, \O_\surj)$ to the ideal $ \field \oplus \fieldmax$ for
$\popi =0$ and to $\field \oplus 0$ for $\popi \geq 1$. The
zero-ideal in $\dvd$ defines the zero ideal. The ideal $0 \oplus
\fieldmax \subseteq \Gamma(\dvd,\O_\surj)$ defines an ideal sheaf in
the radical topology, but it is not the sheafification of an ideal
in $\dvd$. \ifthenelse{\boolean{book}}{Look at the exact
$\dvd$-sequence ($\popi \geq 1$)
$$0 \lto \dvd \stackrel{\pipo}{\lto} \dvd \lto \dvd/(\pipo)
\lto 0 \,.$$
The global evaluation in the surjective topology gives
the mapping
$$ \field \oplus \fieldmax \stackrel{\pipo}{\lto} \field \oplus \fieldmax \lto 0 \oplus \fieldmax
 = \Gamma(\dvd,(\dvd/ (\pipo))) \lto 0 \, ,$$
where the multiplication is bijective in the first component and the
zero map in the second component.}{}
\end{example}

\begin{example}
Let $\dvd$ denote a discrete valuation domain with maximal ideal
$\fom=(\pi)$, consider $\ideal = \dvd \pi^2 \subset \idealsec =\dvd
\pi \subset \dvd$. Then the radical of $\dvd\pi^2$ inside $\dvd$ is
$\dvd \pi$. The inclusion
$\dvd \pi^2 \subset \dvd \pi$ is as
$\dvd$-modules isomorphic to $\dvd \pi \subset \dvd$. Hence the
radical of $\dvd \pi^2 $ inside $ \dvd \pi$ is just $\dvd \pi^2$. So
the radical of a submodule depends on the surrounding module.

This also shows that ideals $\ideal \subseteq \idealsec$ with the
same radical do not induce isomorphisms $\ideal_\surj \to
\idealsec_\surj$. This mapping is in our example just $\O_{\surj}
\to \O_{\surj}$, $1 \mapsto \pi$. However the image ideal sheaves of
$\ideal_\surj $ and $ \idealsec_\surj$ in $\O_\surj$, that is,
$\ideal^\surj$ and $ \idealsec^\surj$, are the same. Also, the
surjective homomorphism
$\ring/\ideal \to \ring/\idealsec$ induces an
isomorphism
$(\O/\ideal)_\surj \to (\O/\idealsec)_\surj$.
\ifthenelse{\boolean{book}}{ The exact complex $0 \to \ideal \to
\idealsec$ is not exact in the surjective topology, since the
tensoration with $\fieldmax = \dvd/\fom$ gives the non-exact complex
$0 \to \field \stackrel{0}{\to} \field$.}{}
\end{example}

\subsection{Constructible partitions}
\label{constructiblesubsection}

\

\ren{\fieldres}{{k}}

\medskip
For a discrete valuation domain $\dvd$ with quotient field $\field$
and residue class field $\fieldres$ the ring homomorphism
$\dvd \to\field \oplus \fieldres$
defines a surjection
$$\Spec \field \uplus \Spec \fieldres \to \Spec \dvd \, ,$$
which is neither flat nor submersive. Since
$(\field \oplus \fieldres )\tensordvd (\field \oplus \fieldres )
= \field \oplus \fieldres $,
every element in
$\field \oplus \fieldres $ is compatible and defines a global
section on $\Spec V$ in the surjective topology. We treat here
similar constructible partitions of a scheme and the global sections
in the surjective topology which arise from them.

Recall that a subset in a scheme $\spax$ is called
\emph{constructible} if it is the finite union of intersections of a
closed and of an open subset
(see \cite[\S 0.9]{EGAIII}, \cite[Exercise II. 3.18]{haralg}
or \cite[Section 14.3]{eisenbud}).

\begin{definition}
\label{constructpartdef}
We call a finite disjoint union
$\spax=\biguplus_{\indpartinset} \obpart_\indpart$,
where each
$\obpartindic$ factors as
$\obpartindic \to \clopart_\indpart \to
\spax$ with an open embedding
$\obpartindic \to \clopart_\indpart$ and
a closed embedding
$\clopart_\indpart \to \spax$
 a \emph{constructible partition} of $\spax$.
 If all $\obpartindic$ are affine, then we
call this an \emph{affine constructible partition} of $\spax$.
\end{definition}

\begin{remark}
We can refine every constructible partition to an affine
constructible  partition. We can further refine so that each patch
$\obpartindic$ is integral, i.e. irreducible and reduced. If $\spax$
is a scheme of finite type over a field,
then we can get also a regular constructible partition \cite[Corollaire 6.12.5]{EGAIV}.
Note that every constructible partition
$\biguplus_{\indpartinset} \obpartindic \to \spax$
is a bijection,
and that for every $\point \in \spax$ the induced
homomorphism of residue class fields is an isomorphism.
\end{remark}

We introduce the constructible topology, where essentially the
constructible partitions are the coverings.

\begin{definition}
The \emph{constructible topology} on a scheme $\spax$ is given by
the full subcategory $\catopenspax$ of schemes of finite type which
are affine over $\spax$.
A morphism $\obtopa \to \obtopo$ in $\catopenspax$
is a cover if there exists a constructible partition
$\biguplus _{\indpartinset} \obtopo_\indpart \to \obtopo$
and a factorization
$\biguplus _{\indpartinset} \obtopo_\indpart \to \obtopa \to \obtopo$.
We denote a scheme $\spax$ endowed with the constructible
topology by $\spax_\cons$.
\end{definition}

\begin{remark}
The constructible topology is an affine single-handed Grothendieck
topology. It is a refinement of the affine single-handed
Zariski topology. Every scheme morphism
$\spay\to\spax$ induces a site-morphism
$\spay_\cons \to \spax_\cons$.
The surjective topology is a refinement of the
constructible topology, so there exists a site morphism $\spax_\surj
\to
\spax_\cons$. We will see below (Theorem
\ref{surjectiveglobalsection}, Corollary
\ref{surjectivecoherentmoduleglobal}) that in characteristic $0$
these two topologies yield the same modules of global sections.
\end{remark}

\begin{lemma}
\label{constructibleproductlemma}
Let $\spax$ denote a scheme and let
$\obpart_\indpart=\Spec \ring_\indpart$, $ \indpartinset$,
denote an affine constructible partition of $\spax$.
Then
$
(\biguplus_\indpartinset \obpart_\indpart) \timesx
(\biguplus_\indpartinset \obpart_\indpart) \cong
(\biguplus_\indpartinset \obpart_\indpart)$,
and the projections are isomorphisms.
\end{lemma}
\begin{proof}
Note that $\spaz \timesx \spaz \cong \spaz$ for a closed immersion
$\spaz \to \spax$ and
$\openzar \timesx \openzar \cong \openzar$ for
an open subscheme $\openzar \subseteq \spax$.
Hence this is also true for a composition
$\openzar \to \spaz \to \spax$.
Moreover, for disjoint immersions
$\obpart_\indpart,\obpart_\indpartsec \to \spax$
we have
$\obpart_\indpart \timesx \obpart_\indpartsec = \emptyset$.
So
$$(\biguplus_\indpartinset \obpart_\indpart) \timesx (\biguplus_\indpartinset \obpart_\indpart)
\cong \biguplus_{(\indpart,\indpartsec) \in \setindpart \times
\setindpart} \obpart_\indpart \timesx \obpart_\indpartsec
\cong \biguplus_\indpartinset \obpart_\indpart \timesx \obpart_\indpart
\cong \biguplus_\indpartinset \obpart_\indpart \, .$$
\end{proof}

\begin{corollary}
\label{constructibleglobalsectioncorollary}
Let $\spax$ denote a scheme endowed with the constructible or the surjective topology.
Let $\obpartindic$, $ \indpartinset$,
denote a constructible partition of $\spax$, and let $\sheaf$ denote
a sheaf of abelian groups on $\spax_\topo$,
where $\topo =\surj$ or $= \cons$.
Then the restriction homomorphism
$$\Gamma(\spax_\topo,\sheaf )
\to \Gamma(\biguplus_\indpartinset \obpart_{\indpart\, \topo}, \sheaf )$$
is an isomorphism.
\end{corollary}
\begin{proof}
This isomorphism follows from Lemma \ref{constructibleproductlemma}
and Lemma \ref{antistrictlemma}.
\end{proof}

\begin{corollary}
\label{constructibleglobalsectiondisjoint}
Let $\spax$ denote a scheme and let $\modul$ be an $\O_\spax$-module on $\spax$.
Let
$\morpart_\indpart: \obpartindic \to \spax$, $ \indpartinset$,
denote a constructible partition of $\spax$. Then in the
constructible or in the surjective topology we have an isomorphism
$$\Gamma(\spax_\topo,\modul_\topo)
\to \bigoplus_\indpartinset \Gamma(\obpart_{\indpart\, \topo}, (\morpart_\indpart ^*(\modul))_\topo) \, .$$
In particular, any $(\eles_\indpart)_\indpartinset \in
\bigoplus_\indpartinset \Gamma(\obpart_{\indpart\, \topo},
\morpart_\indpart ^*( \modul))$
yields a global section in
$\Gamma(\spax_\topo , \modul_\topo)$.
\end{corollary}
\begin{proof}
Due to Corollary \ref{constructibleglobalsectioncorollary} we may
assume that $\spax= \biguplus_\indpartinset \obpartindic$ is a
disjoint union of open subschemes. Then this follows from
Lemma
\ref{pullbacksheafdisjoint}
\end{proof}

\ifthenelse{\boolean{book}}{
\begin{bookexample}
\label{disjointcover} Let $g \in \ring$. Then the ring homomorphism
$\ring \to \rings   =\ring_g \oplus \ring/(g)$ is a constructible
partition and hence every element $(r/g^k, \bar{s}) \in $ defines a
compatible element and hence an element in $\Gamma(\ring_\surj,
\O_\surj)$.
\end{bookexample}
}{}

\begin{proposition}
\label{constructibleradical}
Let $\spaxeqspecring$ be noetherian of finite dimension.
Then the closure operation for finitely generated
modules induced by the constructible topology is the radical.
\end{proposition}
\begin{proof}
Since the surjective topology is a refinement of the constructible
topology and because of Proposition \ref{surjectiveinduceradical} we
only have to show for
$\eles \in \rad(\submod)$, $\submod \subseteq \modul$
finitely generated modules,
that $\eles $ belongs locally to
$ \submod$ in a constructible partition.
We may assume that $\submod=0$ and we may pass to a constructible partition where the
patches are affine and integral.
So we may assume that $\ring$ is a noetherian domain of finite dimension.
Let $\eles \in \rad(0)$, so
that $\eles \tensor 1=0$ in $\modul \tensorr \field $ for every
field $\ring \to \field   $.
Then in particular $\eles \tensor 1=0$ in
$\modul \tensorr \quotfield(\ring)$ and so
$\eles \tensor 1 =0$ in $\modul \tensorr \ring_\fug$ for some $0 \neq \fug \in \ring$.
Considering the constructible partition $D(\fug) \uplus V(\fug)$ we
can go on inductively.
\end{proof}

\begin{example}
\label{constructibleradicalexample}
In the case of an ideal we can give a constructible partition as predicted by
Proposition \ref{constructibleradical} immediately.
Let $\fuf \in \rad (\ideal)$, $\ideal = (\runfuf)$ an ideal in a commutative ring
$\ring$.
An affine constructible partition is given by the
reductions of
$$ D(\fuf_1) \uplus (D(\fuf_2) \cap V(\fuf_1)) \uplus \ldots
\uplus ( D(\fuf_\numfuf)\cap V(\fuf_1 \komdots \fuf_{\numfuf-1}) )
\uplus V(\runfuf)\, .$$
The corresponding rings are
$\ring_\numfufbreak
=(\ring/(\fuf_1 \komdots \fuf_\numfufbreak)_{\fuf_{\numfufbreak+1}})_\red$,
$\numfufbreak = 0 \komdots \numfuf$
(setting $\fuf_{\numfuf+1}=1$).
In each of the
$\ring_\numfufbreak$ we clearly have
$\fuf \in  \ideal \ring_\numfufbreak$,
since this extended ideal is the unit ideal for
$\numfufbreak \leq \numfuf -1$ and for $\numfufbreak = \numfuf$ it
is the zero-ideal,
but also $\fuf$ becomes $0$ in
$\ring_\numfuf=(\ring/(\runfuf))_\red$.

This constructible partition is perhaps even more canonical that the
forcing algebra $\algforc$ of these data. By the universal property
of the forcing algebra (Proposition \ref{universalforcing}) there
exists a ring homomorphism
$\algforc \to \bigoplus_{\numfufbreak=0}^\numfuf \ring_\numfufbreak$.
Since the global sections of this constructible partition
are global sections in the surjective structure sheaf by Corollary
\ref{constructibleglobalsectiondisjoint}, this also shows that
$\fuf \in \ideal \Gamma(\ring, \O_\surj)$,
and not only that $\fuf$ belongs locally to this ideal sheaf.
\end{example}

\subsection{Surjections and constructible partitions}
\label{surjectiveconstructiblesubsection}

\

\medskip
In this section we gather together several statements which relate
surjections to constructible partitions and their finite extensions
which will later on help us to compute modules of global sections
and cohomology in the surjective topology.

{
\ren{\closedzar}{{Z}}

\begin{lemma}
\label{surjectiveconstructible}
Let $\spax$ denote a noetherian
scheme of finite dimension and let $\morcov: \spay \to \spax$ be a
spec-surjective morphism of finite type. Then there exists an affine
constructible partition
$\biguplus_\indpartinset \obpartindic \to \spax$
such that there is an {\rm(}$\spax$ -covering{\rm)}
refinement
$\biguplus_\indpartinset \obpartsec_\indpart \to \biguplus_\indpartinset \obpartindic' \to \spay$,
where $\obpartindic' = \obpartindic \times_\spax \spay$, $\obpartindic=\Spec \ring_\indpart$ and
$\obpartsec_\indpart=\Spec \rings_\indpart$, such that $\ring_\indpart \to \rings_\indpart$ is a finitely
generated faithfully flat extension of domains. If $\spax$ is a
scheme over a field of characteristic $0$, then one can also achieve
that $\rings_\indpart \tensor_{\ring_\indpart} \rings_\indpart$ is
reduced.
\end{lemma}
\begin{proof}
We do induction on the dimension of $\spax$ and we construct step by
step an affine integral constructible partition of $\spax$ with the
required properties. We pass first to an affine constructible
partition of $\spax$ such that the patches are affine and integral.
We may also assume that $\spay$ is affine, by replacing $\spay$ by a
disjoint union of an open affine cover.
So we may assume that $\ring$ is a domain,
and that
$\morcov:\spay=\Spec \rings  \to \spaxeqspecring$ is
surjective. By the Theorem of Chevalley
\cite[Corollary 14.7]{eisenbud}
the images of the irreducible components of $\Spec \rings $ are
constructible. In particular there exists at least one component
such that its image contains an open affine subset $\openzar$ of
$\specring$. We pass to the constructible affine partition
$\openzar \uplus \closedzar$ ($\closedzar$ the closed component).
Then, on the one hand, the morphism
$\morcov^{-1}(\openzar) \to \openzar$ has a refinement by an affine integral
scheme, and, on the other hand,
the image of the surjection
$\morcov^{-1}(\closedzar) \to \closedzar$
has lower dimension and can be treated as before.

So we may assume that $\ring \to \rings $ are both domains.
By generic flatness
(\cite[Th\'{e}or\`{e}me 6.9.1]{EGAIV} or \cite[Theorem
14.4]{eisenbud})
there exists
$0 \neq \fug \in \ring$ such that
$\ring_\fug \to \rings_\fug$ is faithfully flat.
We pass to the corresponding constructible partition,
the complement being again of lower dimension.
So we get a partition and a refinement such that the mappings
$\obtopa_\indpart \to \obtopo_\indpart$ are faithfully flat morphisms between domains.

Suppose now that $\ring$ is a domain which contains a field of
characteristic zero and let $\ring \subseteq \rings $ be a
faithfully flat extension of finite type, $\rings$ a domain. We
consider the commutative diagram
\[ \xymatrix{\ring \ar[r] \ar[d] & \rings \ar[d] \ar[r]^{i_1\,\,\,\,\,\,\,}
 & \rings \tensorr  \rings  \ar[d]& \\
\quotfield(\ring) \ar[r] & \rings_\nzdring
\ar[r]^{i_1\,\,\,\,\,\,\,\,\,\,} & (\rings \tensorr \rings)_\nzdring
&\!\!\!\!= \rings_\nzdring \tensor_{\quotfield(\ring)}
\rings_{\nzdring} .}\]
Since the characteristic is zero,
$\quotfield(\ring)$ is a perfect field and so
$\rings_{\nzdring}\tensor_{\quotfield(\ring)} \rings_{\nzdring} $
is reduced by
\cite[Proposition 4.6.1]{EGAIV}.
If $\foa \subseteq \rings \tensorr \rings$ is the
nilradical, then it maps to $0$ in
$(\rings \tensorr \rings)_{\nzdring}$, and
so there exists $\fug \in \nzdring$ with $\fug \foa =0$. Then
$(\rings \tensorr\rings)_\fug
= \rings_\fug \tensor_{\ring_\fug} \rings_\fug$
is reduced and we look at the partition
$\spax=D(\fug) \uplus V(\fug)$ and use the induction hypothesis on $V(\fug)$.
\end{proof}
}

\ifthenelse{\boolean{extra}}
{If $g \in A$ is such that all nilpotent
elements of $\alg$ are $0$ in $A_g$, the $A_g$ is reduced: if
$a/g^s$ is nilpotent, say $a^n/g^{ns} = 0 $ in $A_g$, then $g^ua^n =
0 $ in $\alg$. Then also $(g^va)^n=0$ in $\alg$ for suitable $v$ and
so $g^va$ is nilpotent. But then $g^{v+m}a=0$ and $a=0$ in $A_g$.}{}

\begin{lemma}
\label{surjectionfinitelemma}
Let $\spay \to \spax$ be a surjection
of finite type between noetherian schemes of finite dimension.
Then there exists an affine integral constructible partition
$\obpartindic \to \spax$ and finite surjective morphisms
$\obpart'_\indpart \to \obpartindic$, where
$\obpart'_\indpart$ are also affine and integral,
and sections
$\obpart'_\indpart \to \spay$.
\end{lemma}
\begin{proof}
The statement is stable under refining by constructible partitions,
so we may assume that $\spaxeqspecring$, where $\ring$ is a
noetherian domain with quotient field $\quotfield(\ring)$.
For a closed point
$\point \in \spay_{\quotfield(\ring)}$
the extension
$\quotfield(\ring) \subseteq \fieldres(\point)$
is of finite type and hence finite. Let
$\ring \subseteq \ringsec \subset \fieldres(\point)$
be finite such that
$\quotfield(\ringsec) = \fieldres(\point)$.
Then there exists an open subset
$D(\fuh) \subseteq \Spec \ringsec$
such that
$\embpoint: \point = \Spec (\fieldres(\point)) \to \spay$
extends to
$D(\fuh) \to \spay$.
For this reduce to
$\spay=\Spec \alg$ affine and find a common denominator $\fuh$ for the images
$\embpoint^*(\vart_\indj)= \fuq_\indj \in \fieldres(\point)$,
where the $\vart_\indj$ are $\ringsec $-algebra generators of $\alg$.
There exists also $\emptyset \neq D(\fug) \subseteq \spax$ such that the preimage of
$D(\fug)$ lies inside $D(\fuh) \subseteq \Spec \ringsec $. Therefore for
$D(\fug) \subseteq \spax$ we have a finite extension $\spax'= \Spec
\ringsec |_{D(\fug)}$ and a section $\spax' \to \spay$. We go on on
$V(\fug) \subset \spax$ inductively.
\end{proof}

{
\ren{\closedzar}{{V}}

\ren{\indcomp}{{j}}

\ren{\setindcomp}{{J}}

\begin{lemma}
\label{finiteconstructibledisjoint}
Let $\submor :\spay \to \spax$
be a finite surjective morphism between noetherian schemes,
where $\spax$ is integral of finite dimension.
Then there exists an affine integral constructible partition
$\obpartindic \to \spax$ such that
$\spay_\indpart = \obpartindic \times _\spax \spay$ is a disjoint union of
irreducible schemes.
\end{lemma}
\begin{proof}
Let $\obcompindic$, $\indcompinset$, be the irreducible components of $\spay$.
The union of the intersections
$\obcompindic \cap \obcomp_\indcompsec$
($\indcomp \neq \indcompsec$) has smaller dimension than $ \dim (\spax)$,
so its image in $\spax$ is a proper closed subset
$\closedzar \subset \spax$.
For an open non-empty affine subset
$\openzar \subseteq \spax-\closedzar$
in its complement the components
$\openzar_\indcomp = \submor^{-1}(\openzar) \cap \obcompindic$ are disjoint.
So we get inductively a constructible partition of $\spax$ with the required properties.
\end{proof}

}

\subsection{Global sections in the  surjective topology}
\label{surjectiveglobalsubsection}

\

\medskip
We characterize now the ring of global sections
in the constructible topology and in the surjective topology. We
first describe
$\Gamma(\spax_\cons, \O_\cons)$ as a subring of the
radicalization of $\ring$.

\ren{\fieldres}{{\kappa}}

\begin{proposition}
\label{constructiveglobalsection}
Let $\ring$ denote a noetherian ring,
and let $\spaxeqspecring$ be endowed with the constructive
topology. Then
$$\Gamma(\spax_\cons,\O_\cons)\!
=\!\! \{\fuglob \in \prod_{\point \in \spax} \overline{\fieldres(\point)} \!: \,
\mbox{there exists an affine constructible partition}$$
$$\,\,\,\,\,\,\,\,\,\,\,\,\,\,\,\,\,\,\,\,\,\,\,\,\,\,\,\,\,\,\,\,\,\,\,\, \,\,\,\,\,\,\,
\ring \to \bigoplus_\indpartinset \ring_\indpart
\mbox{ and } \fuglob_\indpart \in \ring_\indpart \mbox{ such that }
\fuglob_\indpart(\point)= \fuglob(\point)
\mbox{ in } \overline{\fieldres(\point)} \} \, .$$
\end{proposition}
\begin{proof}
We have natural homomorphisms
$$\Gamma(\spax_\cons, \O_\cons)
\lto \Gamma(\spax_\surj,\O_\surj)
\lto \prod_{\point \in \spax} \overline{\fieldres(\point)}$$
by Proposition \ref{surjectiveabsolutestalk}.
Since the radicalization is the absolute stalk, the second homomorphism is
injective by Lemma \ref{globalstalkinjective}.
The first homomorphism is induced by the site morphism
$\spax_\surj \to \spax_\cons$.
This homomorphism is also injective: if
$\fuglob \in \Gamma(\spax_\cons, \O_\cons)$ is
represented by
$\fuglob_\indpart \in \Gamma(\obpartindic,\O_\indpart)$,
$\spax= \biguplus_\indpartinset \obpartindic$ an affine integral constructible partition and maps to
$0$ in $\Gamma(\spax_\surj,\O_\surj)$, then there exists a
surjection $\spay \to \spax$ such that $\fuglob$ becomes $0$. But
then also the pull-backs of $\fuglob_\indpart$ to $\obpartindic \times_\spax \spay$
are $0$ for all $\indpartinset$. This implies, since the
$\obpartindic$ are integral, that $\fuglob_\indpart=0$ and so $\fuglob=0$
in
$\Gamma(\spax_\cons, \O_\cons)$.

The image of $\Gamma(\spax_\cons, \O_\cons)$ inside the absolute
stalk is the described subring: an element
$\fuglob \in \Gamma(\spax_\cons, \O_\cons)$
is represented by
$\fuglob_\indpart \in \ring_\indpart$,
where
$\obpartindic = \Spec \ring_\indpart$ is an affine constructible partition of
$\spax$,
and every such tuple yields a global element.
\end{proof}

\begin{theorem}
\label{surjectiveglobalsection}
Let $\ring$ denote a noetherian ring of finite dimension endowed with the surjective topology,
$\spaxeqspecring$. Suppose that $\ring$ contains a field of
characteristic zero. Then
$$\Gamma(\spax_\surj,\O_\surj) = \Gamma(\spax_\cons, \O_\cons)
 \, .$$
\end{theorem}
\begin{proof}
We have proven already in Proposition
\ref{constructiveglobalsection} that the natural homomorphism
$\Gamma(\spax_\cons, \O_\cons) \to \Gamma(\spax_\surj,\O_\surj)$ is
injective. For the surjectivity suppose that
$\fuglob \in \Gamma(\spax,\O_\surj)$ is represented by a compatible element
$\fuglob \in \rings$,
$\ring \to \rings   $ being spec-surjective. We have
to construct a constructible partition of $\spax$ such that
$\fuglob$ comes from there.
By Lemma \ref{surjectiveconstructible}
there exists an affine constructible partition
$\spax= \biguplus_\indpartinset \obpartindic$ and affine schemes
$\obpartsec_\indpart \to \obpartindic \times_\spax \spay$ such that
$\biguplus \obpartsec_\indpart \to \spay \to \spax$ is a surjection,
that $\obpartsec_\indpart \to \obpartindic$ are faithfully flat
morphisms of finite type between integral affine schemes and that
$\obpartsec_\indpart \times_\obpartindic \obpartsec_\indpart$
is reduced.
We can pass to the refinement
$\biguplus \obpartsec_\indpart$ and may assume that
$\fuglob$ is given as
$\fuglob_\indpart \in \rings_\indpart$,
$\obpartsec_\indpart = \Spec \rings_\indpart$,
which are compatible over $\obpartindic$.
We can treat the $\fuglob_\indpart$ separately.

The compatibility of
$\fuglob_\indpart \in \rings_\indpart$ means that
$\fuglob_\indpart \tensor 1 - 1 \tensor \fuglob_\indpart$
is nilpotent by Corollary \ref{surjseppresheaf}.
Hence this difference is $0$, since
$\rings_\indpart \tensor_{\ring_\indpart} \rings_\indpart$ is reduced.
But this means that $\fuglob_\indpart$
is also compatible in the faithfully flat topology; hence
$\fuglob_\indpart
\in \ring_\indpart$ \cite[Proposition I.2.18]{milne}, where
$\obpart_\indpart= \Spec \ring_\indpart$.
\end{proof}

\begin{remark}
Theorem \ref{surjectiveglobalsection} and Proposition
\ref{constructiveglobalsection} show that every element $\fuglob \in
\Gamma(\spax_\surj, \O_\surj)$ has at every point
$\point \in \spax$
a well-defined value
$\fuglob(\point) \in \fieldres(\point)$,
not only in
$\overline{\fieldres(\point)}$. This is also clear considering
$\Spec \fieldres(\point) \to \spax$.
This scheme morphism defines a site morphism
$(\Spec \fieldres(\point))_\surj \to \spax_\surj$ and
therefore a homomorphism of the rings of global sections,
$\Gamma(\spax,\O_\surj) \to \Gamma(\Spec \fieldres(\point),
\O_\surj) = \fieldres(\point)$ (in characteristic zero). This
reasoning shows also that in positive characteristic every global
function over $\spax_\surj$ has for every point $\point$ a
well-defined value in $\fieldres(\point)^\perf$.
\end{remark}

We extend some of the results in this section to coherent modules.

\begin{lemma}
\label{coherentpartition}
Let $\spax$ be a noetherian scheme of
finite dimension and let $\modul$ be a coherent sheaf on $\spax$.
Then there exists an affine constructible partition
$\biguplus_\indpartinset \obpartindic \to \spax$ such that
$\morpartindic^*(\modul)$ is free, where
$\morpartindictospax$ is the natural embedding.
\end{lemma}
\begin{proof}
We may assume that $\spaxeqspecring$ is affine and integral.
Then we have an isomorphism
$\modul \tensorr \quotfield(\ring) \cong
\quotfield(\ring)^\ran$. Therefore
$\modul_\fug \cong \ring_\fug^\ran$ for some
$0 \neq \fug \in \ring$ and the statement is true on $\openzar=D(\fug)$.
Then we may go on on $V(\fug)$ inductively.
\end{proof}

\begin{corollary}
\label{surjectivecoherentmoduleglobal}
Let $\spax$ be a noetherian scheme of finite dimension over a field of characteristic $0$ and
let $\modul$ be a coherent sheaf on $\spax$.
Then the global sections in
$\Gamma(\spax_\surj, \modul_\surj)$ are constructible sections,
i.e. they are of the form $\elem=(\elem_\indpart)_\indpartinset$,
$\elem_\indpart \in \Gamma( \obpartindic, \morpart^*_\indpart(\modul))$,
where
$\spax= \biguplus_\indpartinset \obpartindic$
is an affine constructible partition.
\end{corollary}
\begin{proof}
We may assume by Lemma \ref{coherentpartition} that
$\spax= \biguplus \obpartindic$ is such that all
$\morpart_\indpart^*(\modul)$ are free. Then
$\Gamma(\spax_\surj, \modul_\surj)= \bigoplus_\indpartinset
\Gamma(\obpart_{\indpart\,\surj}, \modul_{\indpart\, \surj}) = \bigoplus_\indpartinset
\Gamma(\obpart_{\indpart\, \surj}, \O_\surj^{\ran_\indpart})$  (by Corollary
\ref{constructibleglobalsectiondisjoint}) and the result follows
from Theorem \ref{surjectiveglobalsection}.
\end{proof}

\subsection{Exactness in the surjective topology}
\label{surjectiveexactsubsection}

\

\medskip
We characterize when a complex of modules is exact in the surjective
topology. First examples were already given in Proposition
\ref{presheafexact}.

\begin{theorem}
\label{surjectiveexact}
Let $\ring$ denote a commutative ring and
let $ \lmr $ denote a complex of finitely generated $\ring$-modules.
Then the following are equivalent {\rm(}suppose in the last
statement that $\ring$ is noetherian of finite dimension and
contains a field of characteristic $0${\rm)}.
\numiii
\begin{enumerate}

\item
The complex of sheaves
$ \lmrr{_\surj}$ on $(\specring)_\surj$ is
exact {\rm(}where $\catopen_\ring$ contains all $\ring$-algebras or
all $\ring$-algebras of finite type{\rm)}.

\item
For every ring homomorphism
$\ring \to \ringsec  $ {\rm(}of finite
type{\rm)}, the complex
$ \lmrtensor {\ringsec  }$ is such that the
kernel lies inside the radical of the image.

\item
The complex
$\lmrtensor {\rad(\ring)}$ of $\rad(\ring)$-modules is
exact.

\item
For every homomorphism $\ring \to \field $ to a field the complex
$\lmrtensor{\field}$ of $\field$-vector spaces is exact.

\item
The complex
$\lmrglobal {\ringsec}{_\topo}$
of global sections is exact for all ring homomorphism
$\ring \to \ringsec$ of finite type.
\end{enumerate}
\end{theorem}
\begin{proof}
The conditions (i), (ii) and (iv) are equivalent by Corollary
\ref{catringtestclosureexact}.
(i) $\Rightarrow $ (iii) follows from
Corollary \ref{irreduciblefiltermodulestalk}, Lemma
\ref{exactirreducible} and Proposition
\ref{surjectiveabsolutestalk}.
(iii) $\Rightarrow $ (iv).
A homomorphism
$\ring \to \field   $ factors as $\ring \to
\fieldres(\point) \to \field$, where $\point \in \specring$. By
restricting the exactness of (iii) to the component given by
$\point$, it follows that
$\lmrtensor {\overline{\fieldres(\point)}}$
is exact. This is then also true if we replace
$\overline{\fieldres(\point)}$ by $\fieldres(\point)$ and
then by $\field$.

(iv) $\Ra$ (v).
The condition in (iv) holds also for $\ringsec  $,
and
$(\modul \tensorr \ringsec  )_\topo =(\modul_\topo)|_{\spax'}$ on
$\spax'_{\surj}$ ($\spax' = \Spec \ringsec $) by Corollary
\ref{pullbacksheafify}, hence we may assume that $\ring =\ringsec $.
Let $\elem \in \Gamma(\spax_\surj, \modm_\topo)$ mapping to $0$ in
$\Gamma(\spax_\surj, \modr_\topo)$.

By Corollary \ref{surjectivecoherentmoduleglobal} we may assume that
$\elem$ is represented by $\elem_\indpart \in \modul \tensorr \ring_\indpart$, where
$\morpart_\indpart: \obpartindic =\Spec \ring_\indpart \to \spax$ is an affine integral
constructible partition. We have to find a global section $\elel \in
\Gamma(\spax_\surj, \modl_\surj)$ which maps to $\elem$. Since
$\Gamma(\spax_\surj, \modl_\surj)
 = \bigoplus_\indpartinset
\Gamma( \obpart_{\indpart\,\surj}, \morpart_\indpart^*(\modl)_\surj)$
(Corollary \ref{constructibleglobalsectiondisjoint})
we may assume that
$\spax=\biguplus_\indpartinset \obpartindic$ is the disjoint union of open subschemes
and we can treat the components separately; hence we skip the index
$\indpart$. Since $\elem$ maps to $0$, there exist spec-surjective
homomorphisms $\ring \to \rings   $ such that
$\cob(\elem) \tensor 1 = 0$ in $\modr \tensorr \rings$.
Then by Proposition
\ref{radicalequivalent} for every field $\ring \to \field   $ we
have
$\cob(\elem) \tensor 1 = 0$ in $\modr \tensorr \field$, so in
particular for the quotient field $\quotfield(\ring)$.

By condition (iv) there exists an element $\elel \in \modl \tensorr
\quotfield(\ring)$ mapping to $\elem$, and this is then also true
for $\elel \in \modl \tensorr \ring_\fug$, $0 \neq \fug \in \ring$.
Therefore
$\elem$ comes from the left on the open set $D(\fug) \subseteq \obpart$ and
we pass to the constructible partition
$D(\fug) \uplus V(\fug)$. The
$\ring/(\fug)$-complex $\lmrtensor {(\ring/(\fug))}$ fulfills also
property (iv), so we can construct inductively a constructible
partition
$\obpart_\indpartsec \to \obpart$ such that $\elem$ comes from an
element
$\elel_\indpartsec \in \modl \tensorr \ring_\indpartsec$. Such a tuple is compatible and
defines a global section in $\Gamma(\spax_\surj,\modl_\surj)$ by
Corollary \ref{constructibleglobalsectiondisjoint}.

(v) $\Ra$ (i) is trivial,
since a global solution is also a local solution.
\end{proof}

\begin{corollary}
\label{surjectiveglobalsurjective}
Let $\ring$ denote a commutative ring
of finite dimension containing a field of characteristic $0$,
$\spaxeqspecring$. Let $\modm \to \modr \to 0$ be an exact sequence
of finitely generated $\ring$-modules. Then also
$\Gamma(\spax_\surj, \modm_\surj) \to \Gamma(\spax_\surj,
\modr_\surj)$ is surjective.
\end{corollary}
\begin{proof}
The complex of sheaves $\modm_\surj \to \modr_\surj \to 0$ is exact,
hence by Theorem \ref{surjectiveexact}(v) also $\Gamma(\spax_\surj,
\modm_\surj) \to \Gamma(\spax_\surj, \modr_\surj) \to 0$ is exact.
\end{proof}

\begin{example}
Let $\ring$ be a commutative ring with reduction $\ring_\red$. Then
the complex of $\ring$-modules $0 \to \ring \to \ring_\red$ is exact
only when $\ring$ is reduced. But tensoring with $\fieldres(\fop)$
(and so with every field) yields always an exact complex, since
$\ring_\red \tensorr \ring/\fop = \ring/ \fon_\ring \tensorr
\ring/\foq= \ring/\foq $, where $\fon_\ring$ is the nilradical. So
this complex is exact in the surjective topology.
\end{example}

\begin{remark}
\label{surjectiveexactclosedpoint}
Using the characterization in
Theorem \ref{surjectiveexact}(ii) one can show that for a ring of
finite type over a field $\field$ the exactness in the surjective
topology is equivalent to the exactness after every ring
homomorphism $\ring \to \cloalg{\field}$.
\end{remark}

\begin{proposition}
\label{surjuntensurjoben}
Let $\ring$ be a noetherian ring of finite dimension over a field of characteristic zero
and let
$\submod \subseteq \modul$ be finitely generated
$\ring$-modules. Let
$\submod_\surj \to \submod^\surj \subseteq \modul_\surj$
be the surjection of sheaves in the surjective
topology on
$\spax = \Spec \ring$. Then also
$\Gamma(\spax_\surj, \submod_\surj) \to \Gamma(\spax_\surj, \submod^\surj)$
is surjective.
\end{proposition}
\begin{proof}
Let
$\elem \in \Gamma(\spax_\surj, \submod^\surj)
\subseteq \Gamma(\spax_\surj, \modul_\surj)$.
Recall that $\submod^\surj$ is the image sheaf of
$\submod_\surj \to \modul_\surj$,
so $\elem$ comes locally from $\submod_\surj$.
Due to Corollary \ref{surjectivecoherentmoduleglobal}
we can represent $\elem$ by
$\elem_\indpart \in \Gamma(\obpartindic, \morpart_\indpart^*( \modul))
= \modul \tensorr \ring_\indpart$,
($\morpart_\indpart^*$ is the Zariski pull-back), where
$\partitionofspax$ is an affine partition of $\spax$.
As the global sections of $\submod_\surj$ also respect this partition
we may assume that $\elem \in \modul$. Since
$\elem \in \Gamma(\spax_\surj,\submod^\surj)$, we know by
Proposition \ref{constructibleradical} that there exists a
constructible partition $\partitionofspax$ and sections
$\elel_\indpart \in \Gamma(\obpartindic, \submod_\surj)$ mapping to the
restrictions of $\elem$. Again by Corollary
\ref{surjectivecoherentmoduleglobal} this represents a global element
$\elel \in \Gamma(\spax, \submod_\surj)$
mapping to $\elem$.
\end{proof}

\ifthenelse{\boolean{book}}{
\begin{bookexample}
The complex $0 \to I=(x,y) \to R=K[x,y]$ is exact, but the
tensoration with $\tensorr \field $ yields $0 \to K^2 \to K$, which
is not exact, hence also the corresponding complex of sheaves in the
surjective topology is not exact.
\end{bookexample}}{}

\ifthenelse{\boolean{book}}{
\begin{bookexample}
Let $\modul$ denote an $\ring$-module. Then the short exact sequence
(see Remark \ref{radicalzero} and Lemma
\ref{surjseparatedreduction})
$$0 \lto \rad (0) = \modul_0 \lto
\modul \lto \modul_1 \lto 0
$$ yields the exact complex of sheaves on $(\specring)_\surj$
$$ (\modul_0)_\surj \stackrel{0}{\lto} \modul_\surj  \lto (\modul_1)_\surj  \lto 0 \,.$$
Note that $(\modul_0)_\surj \neq 0$ in general.
\end{bookexample}

\begin{bookexample}
\label{quotientmodulodvdtwo} Let $\dvd$ denote a discrete valuation
domain. The complex of $\dvd$-modules, $0 \to Q(\dvd)/\dvd \to 0$
yields an exact complex in the surjective topology, because
$(Q(\dvd)/\dvd)_\surj = 0$. This is due to the fact that it is $0$
on the disjoint covering $\Spec \dvd/ \fom_\dvd \uplus \Spec
Q(\dvd)$, see Example \ref{quotientmodulodvd}. Also, the other
conditions in Theorem \ref{surjectiveexact} are fulfilled, though
$Q(\dvd)/\dvd$ is not finitely generated.

The short exact sequence $0 \to \dvd \to Q(\dvd) \to Q(\dvd)/\dvd
\to 0$ yields the exact complex of sheaves $\dvd_\surj \to
Q(\dvd)_\surj \to 0$. Its complex of global sections is the
projection $\Gamma((\Spec V)_\surj, \dvd_\surj)= Q(\dvd)^\perf
\oplus (\dvd/ \fom_\dvd)^\perf \to \Gamma((\Spec V)_\surj,
Q(\dvd)_\surj)= Q(\dvd)^\perf \to 0$.
\end{bookexample}}{}

For free modules the exactness in the surjective topology is quite a
strong condition.

\begin{lemma}
\label{freecomplexradicalexactlemma}
Let
$\ring^\ral \stackrel{\coa}{\to} \ring^\ram \stackrel{\cob}{\to} \ring^\rar$
be complex of free $\ring$-modules, so that $\coa$ and $\cob$ are
given by matrices. Let
$\idealminor_\indseize(\coa)$ denote the
$\indseize$th minor ideal, i.e. the ideal generated by all
determinants of the $\indseize \times \indseize$-submatrices of
$\coa$. Then the following are equivalent.

\numiii
\begin{enumerate}

\item We have
$ \sum_{\indseize+\indseizesec =\ram} \idealminor_\indseize(\coa) \idealminor_\indseizesec(\cob) =(1)$.

\item
For every ring homomorphism $\ring \to  \ringsec  $ the tensored complex is exact.

\item
The complex is exact in the surjective topology.

\end{enumerate}
\end{lemma}
\begin{proof}
(i) $\Rightarrow $ (ii).
Since the minor condition is stable under
arbitrary base change, we may assume that $\ringsec  = \ring$ and
that $\ring$ is local. For a local ring the condition means that
there exists an invertible submatrix of $\coa$ of seize $\indseize$
and an invertible submatrix of $\cob$ of seize $\indseizesec$ with
$\indseize+\indseizesec=\ram$. After a base change in the first two
modules we may assume that the first $\indseizesec$ columns in the
matrix $\coa$ are the first $\indseizesec$ unit vectors. It follows
that the first $\indseizesec$ columns of $\cob$ must be zero.
Therefore the remaining
$\ram -\indseizesec = \indseize$ columns of
$\cob$ contain an invertible $\indseize \times \indseize$-submatrix. 
It follows that the first $\indseizesec$ unit vectors generate
exactly the kernel of $\cob$, and this is also the image of $\coa$.

The implication (ii) $\Rightarrow $ (iii) is trivial.
(iii) $\Rightarrow $ (i).
If the minor condition is not fulfilled, then
$\sum_{\indseize+\indseizesec =\ram} \idealminor_\indseize(\coa)
\idealminor_\indseizesec(\cob)$
is contained in a maximal ideal and so there exists a field
$\ring \to \field $
such that
$\sum_{\indseize+\indseizesec
=\ram} \idealminor_\indseize(\coa) \idealminor_\indseizesec(\cob) =0$.
Let $\indseize = \rk (\coa)$
over $\field$.
Then
$\rk (\cob) < \ram - \indseize$ and
$\dim (\ker \cob)= \ram - \dim(\im \cob)> \ram - (\ram -\indseize) = \indseize$,
so the tensorized complex is not exact.
\end{proof}

\begin{example}
Let $\fuf \in \ring$ be a non-zero divisor and a non unit. Then the
sequence
$\ring \stackrel{0}{\to} \ring \stackrel{\fuf}{\to} \ring$
is exact, but not exact in the surjective topology. For a maximal
ideal $\fom $ containing $\fuf$ the tensoration with
$\field=\ring/\fom$ yields
$0 \to \field \stackrel{0}{\to} \field$.
\end{example}

\begin{remark}
Lemma \ref{freecomplexradicalexactlemma} says in particular that for
given $\coa: \ring^\ras \to \ring^\rat$ it can be checked in the
surjective topology (i.e. after tensoration with fields alone)
whether it is an isomorphism or not. This statement can be found for
a homomorphism between flat modules over a noetherian ring in
\cite{hashimotoapproximation}.
\end{remark}

\ifthenelse{\boolean{book}}{
\begin{bookexample}
Lemma \ref{freecomplexradicalexactlemma} has the following special
cases. A homomorphism $\coa: \ring^\ras \to \ring^\rat$ is
surjective if and only it is surjective in the surjective topology.
It is injective in the submersive topology if and only if
$\ideal_\ras (\coa) = (1)$; in this case the image is a free direct
sum.
\end{bookexample}

\begin{bookexample}
Consider the complex $R \stackrel{(1,0)}{\to}  R \oplus R
\stackrel{(0,f)}{\to} \ring$, where $\fuf$ is not a unit, but a
non-zero divisor. Then the complex is exact and the image is a
direct summand of the module in the middle. But this is not exact in
the surjective topology.
\end{bookexample}}{}

\ifthenelse{\boolean{book}}{
\begin{bookexample}
An exact complex of free modules is not necessarily exact in the
sense of Lemma \ref{freecomplexradicalexactlemma}. For example, if
$x,y$ are a regular sequence in $\ring$, then $R
\stackrel{(y,-x)}{\to} \ring^2 \stackrel{(x,y)}{\to} \ring$ is
exact, but for $\ring \to \ringsec  $, where $x,y$ is not a regular
sequence anymore, the tensored complex is not exact. In particular,
the tensored complex for a field where $x$ and $y$ are $0$ is not
exact. The ideal of the minor products is $1 (x,y) + (x,y) (x,y) +
(x,y) 1 = (x,y) \neq (1)$.
\end{bookexample}}{}

{ \ren{\ind}{{k}}

\ren{\modfree}{{G}}

\ren{\num}{{n}}

For a homomorphism $\coa$ between free modules (or for a matrix) we
call the maximum $\indseize$ such that
$ \idealminor_\indseize (\coa)$
contains a non zero divisor its rank
(cf. \cite[Proposition 1.4.11]{brunsherzog}).

\begin{corollary}
\label{freecomplexradicalexactcor}
Let
$ \modfree_\bullet = \modfree_\num \stackrel{\coa_{\num-1}}{\to} \modfree_{\num-1}
\todots \modfree_2 \stackrel{\coa_1}{\to} \modfree_1 \stackrel
{\coa_0} \to \modfree_0 \to 0$
be a complex of free $\ring$-modules
$\modfree_\ind$ of finite ranks $\raa_\ind$.
Suppose that
$\ramat_\ind = \rk (\coa_\indk) = \raa_\ind - \raa_{\ind-1}$
{\rm(}this is the standard condition on ranks{\rm)}. Then the
following are equivalent.

\numiii \begin{enumerate}

\item The minors are
$\idealminor_{\ramat_\indk}(\coa_\ind) = (1)$
for all $\ind=0 \komdots \num$.

\item
The complex
$\modfree_\bullet \tensorr \ringsec$ is exact for
every $\ring \to \ringsec  $.

\item The complex is exact in the surjective topology.
\end{enumerate}
\end{corollary}
\begin{proof}
(i) $\Leftrightarrow $ (iii).
Note that
$\ramat_\ind +\ramat_{\ind-1} = \raa_\ind $
and that
$\idealminor_\indseize(\coa_\ind) \neq (1)$
for $\indseize > \ramat_\ind$.
Then the minor condition of Lemma
\ref{freecomplexradicalexactlemma} for the $\ind$th place reduces to
$\idealminor_{\ramat_\ind}(\coa_\ind) \idealminor_{\ramat_{\ind -1}}(\coa_{\ind -1})=(1)$.
Hence this equivalence follows by
induction on $\ind$.

(ii) $\Rightarrow$ (iii) is contained in Theorem
\ref{surjectiveexact},
so assume that (i) holds.
Since the minor condition respects base changes, we may assume that $\ring =
\ringsec $ and that $\ring$ is local. Since
$\idealminor_{\ramat_0} (\coa_0)=(1)$ and $\ramat_0=\raa_0$ the homomorphism $\coa_1$ is
surjective. Then the kernel $\ker (\coa_0)$ is free of rank $\raa_1
- \raa_0$ and the homomorphism
$\ring^{\raa_2} \to \ring^{\raa_1-\raa_0} \cong\ker (\coa_0) \subseteq \ring^{\raa_1}$
fulfills also the minor condition. So we can go on inductively.
\end{proof}

\begin{remark}
\label{standardconditionremark}
The Buchsbaum-Eisenbud criterion
asserts that if in the situation of Corollary
\ref{freecomplexradicalexactcor} (with also a $0$ on the left) the
so-called standard condition on depth is fulfilled, namely that
$\depth (\idealminor_{\ramat_\ind}( \coa_\ind)) \geq \ind + 1$
, then the complex itself is exact. A variant of Hochster and Huneke
\cite[Theorem 10.3]{hunekeapplication} says that if the height
condition
$\height (\idealminor_{\ramat_\ind}( \coa_\ind)) \geq \ind+1$
is true, then the complex has phantom homology.
\end{remark}
}

\ifthenelse{\boolean{extra}} {how does for a domain $\ring$ the
sheaf $\quotfield(\ring)_\topo$ look like?}{}

\ifthenelse{\boolean{book}}{
\begin{bookexample}
Let $\dvd$ denote a discrete valuation domain with local parameter
$\pi$. Then multiplication by $\pi^n$, $n \geq 1$, yields the exact
sequence $0 \to \dvd \stackrel{\pi^n}{\to} \dvd \to \dvd/(\pi^n) \to
0$ and therefore the exact sequence $\dvd_\surj
\stackrel{\pi^n}{\to} \dvd_\surj \to (\dvd/(\pi^n))_\surj \to 0$ of
sheaves in the surjective topology. Let $n \geq m \geq 1$. Then we
have the following commutative diagram of $\dvd$-modules
\[ \xymatrix{0 \ar[r] &\dvd \ar[r]^{\pi^n} \ar[d]_{\pi^{n-m}} & \dvd \ar[d]^{=} \ar[r] & \dvd/(\pi^n) \ar[d]^{} \ar[r] & 0\\
0 \ar[r]& \dvd \ar[r]^{\pi^m} & \dvd   \ar[r] & \dvd/(\pi^m) \ar[r]
& 0 \, .}\] The ideal sheaves defined by the ideals $\dvd \pi^n$ and
$\dvd \pi^m$ inside $\O_{\surj}= \dvd_\surj$ are the same and the
downarrow on the right becomes an isomorphism in the surjective
topology, but the injective multiplications are not injective
anymore in the surjective topology. This behavior is mirrored by the
corresponding diagrams under a homomorphism $V \to \field $ to a
field. If $V \to \field $ sends $\pi$ to $0$ (and $n >m$), then the
resulting diagram is
\[ \xymatrix{K \ar[r]^{0} \ar[d]_{0} & \field \ar[d]^{=} \ar[r]^{=} & \field \ar[d]^{=} \ar[r] & 0\\
K \ar[r]^{0} & \field   \ar[r]^{=} & \field \ar[r] & 0 \, }\] and if
$\pi$ is not $0$ in $\field$, it looks like
\[ \xymatrix{K \ar[r]^{=} \ar[d]_{=} & \field \ar[d]^{=} \ar[r]^{0} & 0 \ar[d]^{=} \ar[r] & 0\\
K \ar[r]^{=} & \field   \ar[r]^{0} & 0 \ar[r] & 0 \, . }\]
\end{bookexample}}{}

\ifthenelse{\boolean{book}}{
\begin{bookremark}
\label{generatingsurjectivenotara} The surjective topology exhibits
some strange but trivial properties with respect to generating. If
$\idealsubring$ is an ideal with zero-set $\spay=V(I) \subseteq
\spaxeqspecring$, then the ideal sheaf $\ideal^\surj \subseteq
\O_\surj$ is generated by the global indicator function $e_\spay \in
\Gamma(\spax, \O_\surj)$. For $\spay \neq \emptyset, \spax$, this
function $e_\spay$ is a zero divisor, and an infinite resolution is
given by
$$\to \O_\surj \stackrel{e_\spay}{\to}  \O_\surj \stackrel{1-e_\spay}{\to}
  \O_\surj \stackrel{e_\spay}{\to}  \O_\surj \to \O_\surj/\ideal^\surj \to 0 \, .$$
This shows in particular that generating in the surjective topology
is not at all the same as generating up to radical (the number of
elements in $\ring$ needed to generate the ideal up to radical,
called its \emph{arithmetical rank} or $\ara$, is of course not $1$
in general). Generators $f_1 \komdots f_n \in \ring$ up to radical
yield a surjection $\O^n_\surj \to \ideal^\surj$ in the surjective
topology. Resolution up to radical has hardly any reasonable
uniqueness property.
\end{bookremark}}{}

\subsection{Cohomology in the surjective topology}
\label{surjectivecohomologysubsection}

\

\medskip
In this section we show that the first cohomology of coherent
modules in the constructible and in the surjective topology is zero.

\begin{proposition}
\label{constructiblecohomology}
Let $\spax$ denote a scheme endowed
with the constructible or the surjective topology. Let
$\spax'=\biguplus _\indpartinset \obpartindic \to \spax$
denote a constructible partition.
Let $\presh$ denote a presheaf of abelian groups. Then the
\v{C}ech cohomology groups for this covering are
$\chch^0(\spax' \to\spax, \presh) = \Gamma(\spax', \presh)$ and
$\chch^\numcoho (\spax' \to \spax, \presh)=0$ for $\numcoho \geq 1$.
\end{proposition}
\begin{proof}
This follows from Lemma \ref{constructibleproductlemma} and Lemma
\ref{chchisoproj}.
\end{proof}

The first cohomology for coherent sheaves vanishes also in the
surjective topology, but the proof is more involved. We restrict to
characteristic zero.

\ren{\spau}{{U}}

\begin{theorem}
\label{surjectivecohomologytheorem}
Let $\spax$ be a noetherian
scheme of finite dimension over a field of characteristic $0$. Then
$\chch^1(\spax_\surj, \O_\surj)=0$.
\end{theorem}
\begin{proof}
Let $\spay \to \spax$ be a surjection of finite type. We have to
show that a \v{C}ech cocycle represented by this covering is $0$. To
show this we may use refinements $\spaz \to \spay$, where
$\spaz\to\spax$ is also a covering.
Let $\obpartindic \to \spax$ be a
constructible partition and set
$\spay_\indpart =\obpartindic\times_\spax \spay \to \obpartindic$.
The first \v{C}ech cohomology (with respect to the given cover) is
given as the homology of the first row in the following commutative
diagram
\newcommand{\surjweg}{{}}
$$
 \xymatrix{
\Gamma(\spay ,\O_\surj) \ar[r] \ar[d] &
\Gamma(\spay\times_\spax \spay , \O_\surj) \ar[r] \ar[d] &
\Gamma(\spay\times_\spax \spay \times _\spax \spay,\O_\surj) \ar[d] \cr
\Gamma(\biguplus \spay_\indpart,\O_\surj) \ar[r] \ar[d]
&\Gamma(\biguplus(\spay_\indpart\times_{\obpartindic}\spay_\indpart), \O_\surj) \ar[r] \ar[d]
&\Gamma(\biguplus (\spay_\indpart\times_{\obpartindic}
 \spay_\indpart \times_{\obpartindic}\spay_\indpart), \O_\surj) \ar[d] \cr
\bigoplus_{i}\Gamma(\spay_{i},\O_\surj) \! \ar[r] \!
&\bigoplus_{i}\Gamma(\spay_\indpart\times_{\obpartindic}\spay_\indpart, \O_\surj)
\!\ar[r] \!
& \bigoplus_{i} \Gamma(\spay_\indpart\times_{\obpartindic}
\spay_\indpart \times _{\obpartindic}\spay_\indpart, \O_\surj)
 \, \! ,}
 $$
where the down-arrows are isomorphisms by Corollaries
\ref{constructibleglobalsectioncorollary} and
\ref{constructibleglobalsectiondisjoint}.
An element
$\elem \in \Gamma((\spay\times_\spax \spay)_\surj, \O_\surj)$ is a cocycle if
and only if for every $\indpart$ its component
$\elem_\indpart
\in \Gamma((\spay_\indpart\times_{\obpartindic} \spay_\indpart)_\surj, \O_\surj)$
is a cocycle. We have to show that $\elem$ comes after a refinement
$\spaz \to \spay$ from the left, and for this we can look at the components
separately. Therefore we will repeatedly simplify the situation by
using a constructible partition of $\spax$.

By Lemma \ref{surjectionfinitelemma} and by the preceding
consideration we may assume that $\spay \to \spax$ is a finite
surjective morphism between noetherian affine integral schemes. We
may also assume by Lemma \ref{surjectiveconstructible} that this
morphism is flat. Suppose that
$\elem \in \Gamma(\spay \times_\spax \spay, \O_\surj)$
maps to $0$ on the right.
This $\elem$ is by Theorem \ref{surjectiveglobalsection} represented by
$\elem_\indpartsec \in \Gamma(\spaz_\indpartsec,\O)$,
$\indpartsecinset$, where
$\biguplus_\indpartsecinset \spaz_\indpartsec \to \spay\times_\spax \spay$
is a constructible partition into affine integral schemes.

By Lemma \ref{finiteconstructibledisjoint} we find a constructible
partition of $\spax$ such that we may assume that all irreducible
components $\spav_\indpart$, $\indpartinset$, of
$\spay\times_\spax\spay$ are disjoint.
We may furthermore assume by refining that the partition
$\spaz_\indpartsec$ of $\spay \times_\spax \spay$
(where the cocycle $\elem = \elem_\indpartsec$ lives)
is a subpartition of the disjoint irreducible components.
So we write the cocycle as
$\elem_{\indpart \indpartsec} \in \Gamma(\spaz_{\indpart \indpartsec}, \O)$, where
$\spaz_{\indpart \indpartsec} \to \spav_\indpart $,
$\indpartsec \in \setindpartsec_\indpart$, are constructible partitions.
Let
$\spaw_\indpart = \spaz_{\indpart \indpartsec_\indpart}$
be the open subsets of the components $\spav_\indpart$ where
the generic elements
$\elem_{\indpart \indpartsec_\indpart} \in \Gamma(\spaw_\indpart, \O)$ live.
Then the images of $\spav_\indpart - \spaw_\indpart$ are closed and their union
over $\indpartinset$ is a closed subset $\closed$ of $\spax$ of
smaller dimension.
We consider the constructible partition
$\spax= \spau \uplus \spaa$. The generic tuple
$\elem_{\indpart \indpartsec_\indpart}$, $\indpartinset$,
is represented by the element
$$
(\elem_{\indpart \indpartsec_\indpart})_\indpartinset
\in \Gamma( \biguplus_\indpartinset \spav_\indpart \cap
(\project^{-1}(\spau) \times_\spau \project^{-1}(\spau)),\O)
=\Gamma( \project^{-1}(\spau) \times_\spau
\project^{-1}(\spau), \O)\, .$$
Note that
$ \spav_\indpart \cap (\project^{-1}(\spau)\times_\spau \project^{-1}(\spau))
\subseteq  \spaw_\indpart$,
so we can consider the restriction of
$\elem_{\indpart \indpartsec_\indpart}$.
If we set $\spau= \specring$
(possibly after replacing $\spau$ by a smaller affine non-empty subset)
and $\project^{-1}(\spau)=\Spec \algb$,
then we have
$\Gamma( \project^{-1}(\spau) \times_\spau \project^{-1}(\spau), \O)
= \algb \tensorr\algb$
and we may consider the generic part of the cocycle in this
tensor product.

So we may assume now that $\ring$ is a noetherian domain,
$\ring\subseteq \algb$ is finite and flat and that we have an element
$\elem \in \algb \tensorr \algb$ which maps to $0$ in
$\Gamma((\spay \times_\spax \spay \times_\spax \spay)_\surj, \O_\surj)$.
With a similar argument as in
the proof of Lemma \ref{surjectiveconstructible} we find an affine
constructible partition of
$\biguplus \obpartindic \to \spax$
such that the
$\spay_\indpart \times_{\obpartindic} \spay_\indpart \times_{\obpartindic}\spay_\indpart$
are reduced. But then the image of $\elem$ is not only $0$ in the
surjective topology, but already in
$\algb_\indpart \tensor_{\ring_\indpart} \algb_\indpart \tensor_{\ring_\indpart} \algb_\indpart$.
Then $\elem$ comes from the left by the exactness of
the \v{C}ech complex in the flat topology.
\end{proof}

\begin{corollary}
\label{surjectivecohomologycoherent}
Let $\spax$ be a noetherian
scheme of finite dimension over a field of characteristic $0$
endowed with the surjective topology. Let $\modul$ be a coherent
module on $\spax$ with corresponding module $\modul_\surj$ in the
surjective topology. Then
$\chch^1(\spax, \modul_\surj)=0$ and
$H^1(\spax, \modul_\surj)=0$.
\end{corollary}
\begin{proof}
The vanishing of the first cohomology follows from the vanishing of
the \v{C}ech cohomology. By Lemma \ref{coherentpartition} we may
pass to a constructible partition
$\morpart_\indpart: \obpartindic\to\spax$
such that
$\morpart_\indpart^*(\modul)$ is finitely generated and
free. Then this follows from Theorem
\ref{surjectivecohomologytheorem}.
\end{proof}

\begin{corollary}
\label{surjectivecohomologycorollary}
Let $\spaxeqspecring$ be an
affine noetherian scheme of finite dimension over a field of
characteristic $0$. Let $\olmro$ be a short exact sequence of
finitely generated $\ring$-modules with corresponding short exact
sequence of sheaves
$0 \to \ker = \modl^\surj \to \modm_\surj \to\modr_\surj \to 0$
on $\spax_\surj$. Then
$H^1(\spax_\surj, \ker )=0$.
In the situation of Proposition \ref{syzygycohomology} we get
$H^1(\spax_\surj, \Syz)=0$.
\end{corollary}
\begin{proof}
The exactness of
$\modm_\surj \to \modr_\surj \to 0$ induces by
Corollary \ref{surjectiveglobalsurjective}
a global surjection
$\Gamma(\spax_\surj,\modm_\surj) \to  \Gamma(\spax_\surj,\modr_\surj)$.
Since
$H^1(\spax_\surj, \modm_\surj)=0$
by Corollary \ref{surjectivecohomologycoherent},
it follows
$H^1(\spax_\surj,\ker)=0$
by the long exact sequence of cohomology.

For the second statement let $\ring^\numgen \to \submod \to 0$ be
given. By
Corollary \ref{surjectiveglobalsurjective}
and Proposition \ref{surjuntensurjoben} we get a surjection
$\Gamma(\spax_\surj, \O_\surj^\numgen)
\!\!\to\! \! \Gamma(\spax_\surj,\submod_\surj)
\!\!\to\!\! \Gamma(\spax_\surj, \submod^\surj)$. Hence the result follows
from $H^1(\spax_\surj, \O_\surj^\numgen)=0$.
\end{proof}

\ifthenelse{\boolean{book}}{
\begin{bookremark}
The results of this section hold probably also over positive
characteristic. Also, we expect that in the situation of Theorem
\ref{surjectivecohomologytheorem} the higher \v{C}ech cohomologies
vanish as well.
\end{bookremark}}{}

\section{Frobenius topology and Frobenius closure}
\label{frobeniussection}

\subsection{Frobenius morphism and Frobenius closure}

\label{frobeniussubsection}

\

\medskip
Let $\ring$ denote a ring containing a field of positive
characteristic $\posp$. Then the \emph{Frobenius homomorphism} is
given by
$\Frob: \ring \to \ring , \fuf \mapsto \fuf^\posp$.
This is  a ring homomorphism due to
$(\fuf+\fug)^\posp=\fuf^\posp+\fug^\posp$ and so are its iterations
$\Frob^\expoe: \ring \to \ring, \fuf \mapsto \fuf^\potq$,
where $\potq=\posp^\expoe$. The $\expoe$-th
\emph{Frobenius power} of an ideal $\idealsubring$ is the extended ideal
$$ \ideal^{[\potq]} = \Frob^{\expoe}(\ideal) = \{\fuf^\potq: \fuf \in \ideal \}
 = (\fuf_1^\potq \komdots \fuf_\numgen^\potq)\, ,$$
where $\ideal = (\runfuf)$.

\begin{definition}
\label{frobeniusclosuredef} The \emph{Frobenius closure} of an ideal
$\idealsubring$ in a commutative ring containing a field of positive
characteristic $\posp$ is given as
$$ \ideal^\frobf = \{\fuf \in \ring: \exists \expoe \in \NN \mbox{ such that } \fuf^\potq \in
\ideal^\frobpotq , \, \potq = \posp^\expoe\}\,.$$

For a ring $\ring$ of positive characteristic and an $\ring$-module
$\modul$ we denote $\frobe \modul = \modul \tensorr  \frobe \ring$,
where $\frobe \ring$ is the ring $\ring$, but considered as an
$\ring$-algebra via the $\expoe$-th iteration of the Frobenius
homomorphism. For submodules $\submodul$ we get $\frobe
\ring$-homomorphisms $\frobe \submod \to \frobe \modul$, and the
Frobenius closure of a submodule $\submod \subseteq \modul$ is
defined by
$$\submod^\frobf
= \{ \elem \in \modul:\, \exists \expoe \mbox{ such that } \elem
\tensor 1 \in \im (\frobe \submod \to \frobe \modul) \} \, .$$
\end{definition}

We consider the category of schemes of positive characteristic
$\posp$, that is, the category of schemes $\spax$ with a (unique)
structure morphism $\spax \to \Spec \ZZ /(\posp)$. The
\emph{absolute Frobenius morphism} is the identity on the underlying
set $\spax$ and on every open affine subset $\specring = \openzar
\subseteq \spax$ it is the scheme morphism given by the Frobenius
homomorphism. We denote the Frobenius by ${}^{1} \spax \to \spax$
and the $\expoe$-th iteration by $\frobe \spax \to \spax$.

\begin{definition}
\label{frobeniuscoverdef} A (finite) morphism $\obv \to \obu$ of
finite type of schemes of positive characteristic $\posp$ is called
a (\emph{finite}) \emph{Frobenius cover} if there exists $\expoe \in
\NN$ and a factorization  $\frobe \,\obu \to \obv \to \obu$  such that
$\frobe \,\obu \to \obv$ is surjective and the composition is an
iteration of an absolute Frobenius morphism.
\end{definition}

\begin{lemma}
\label{frobeniuscoverlemma} The composition of two {\rm(}finite{\rm)}
Frobenius covers and an arbitrary base change of a {\rm(}finite{\rm)}
Frobenius cover is again a {\rm(}finite{\rm)} Frobenius cover.
\end{lemma}
\begin{proof}
Let $ \frobe \obu \to \obv \to \obu $ and $\frobesec \obv \to \obw
\to \obv$ be given. Then we also have ${}^{\expoe +\expoesec}\! \obu
\to \frobesec \obv$ and so $ {}^{\expoe +\expoesec}\! \obu \to
\frobesec \obv \to \obw \to \obv \to \obu$. If $\spay \to \spax$ is
a base change, then the pull-back of $\frobe \obu \to \obv \to \obu$
and the universal property of the product give  $\frobe \spay
\times_{\frobe \spax} \frobe \obu \to  \spay \timesx \frobe \obu \to
\spay \timesx \obv \to \spay \timesx \obu $.
\end{proof}

\begin{lemma}
\label{frobeniusdiagonalreduziert} Let $\spax$ denote a scheme in
positive characteristic, $\expoe \in \NN$. Then the projections
$\frobe \spax \timesx \frobe \spax \to \frobe \spax$ induce
{\rm(}identical{\rm)} isomorphisms in the reductions, $(\frobe \spax
\timesx \frobe \spax)_\red \cong (\frobe \spax)_\red$. The same is
true for a finite Frobenius cover $\spay\to \spax$, and the diagonal
$\spay \to \spay \times_\spax \spay$ is also a finite Frobenius cover.
\end{lemma}
\begin{proof}
We may assume that $\spaxeqspecring$ is affine. We have the diagonal
mapping $\frobe \spax \stackrel{\diag}{\to} \frobe \spax \times
_\spax \frobe \spax \stackrel{\project_1}{\to}
 \frobe \spax$ resp. $\frobe \ring \stackrel{\inject_1}{\to} \frobe \ring
\tensorr \frobe \ring \stackrel{\diag^*}{\to} \frobe \ring$, $\fuf
\mapsto \fuf \tensor 1$, $\fua \tensor \fub \mapsto \fua\fub$. The
composition is of course the identity. The other composition sends
$\fua \tensor \fub \mapsto \fua \fub \tensor 1$. We have $(\fua
\tensor \fub)^\potq= \fua^\potq \tensor \fub^\potq =
\fua^\potq\fub^\potq \tensor 1$, because $\fub^\potq \in \ring$.
Hence $ \fua \tensor \fub - \fua \fub \tensor 1$ is nilpotent in
$\frobe \ring \tensorr \frobe \ring$, and so after reduction we have
an isomorphism. The two projections yield the same isomorphism,
because the inverse is $\diag_\red$ in both cases.

If $\rings$ is finite over $\ring$ with a spec-surjective homomorphism $\rings \to
\frobe \ring$, then again we look at $\rings \tensorr \rings
\stackrel{\diag^*}{\to} \rings \stackrel{\inject_1}{\to} \rings
\tensorr \rings$. The difference between $ \fua \tensor \fub$ and
$\fua \fub \tensor 1$ maps to a nilpotent element under the finite
spec-surjective homomorphism $ \rings \tensorr \rings \to \frobe
\ring \tensorr \frobe \ring$. Since this homomorphism is
spec-surjective, the difference itself must be nilpotent.

For the last statement we look at
$ \rings \tensorr \rings
\stackrel{\diag^*}{\to}
\rings \stackrel{\inject_1}{\to} \rings \tensorr \rings
\to  \frobe ( \rings \tensorr \rings )$.
The composition maps
$\fua \tensor \fub \mapsto \fua^\potq \fub^\potq \tensor 1$,
and for $\expoe$ big enough this equals
$(\fua \tensor \fub)^\potq$,
as the difference
$ \fua \tensor \fub -\fua \fub \tensor 1$
is nilpotent.
\end{proof}

\subsection{Frobenius topology}

\label{frobeniustopologysubsection}

\

\begin{definition}
\label{frobeniustopologydef} The \emph{Frobenius topology} on a
scheme $\spax$ of characteristic $\posp$ is given by schemes over
$\spax$ where the morphisms $\obv \to \obu$ (including the
structural morphisms $\obv \to \spax$) are finite Frobenius covers.
We denote $\spax$ endowed with the Frobenius topology by
$\spax_\frob$.
\end{definition}

\begin{proposition}
The Frobenius topology is an affine single-handed covering
Gro\-then\-dieck topology. Every scheme morphism $\spay \to \spax$
of schemes over $\Spec \ZZ/(\posp)$ yields a site morphism
$\spay_\frob
\to \spax_\frob$.
\end{proposition}
\begin{proof}
This follows from Lemma \ref{frobeniuscoverlemma}.
\end{proof}


\begin{proposition}
\label{frobeniusfrobenius}
The Frobenius closure is the closure
operation induced by the Frobenius topology.
\end{proposition}
\begin{proof}
This follows from the Definitions \ref{frobeniusclosuredef} and
\ref{frobeniustopologydef} and from Lemma \ref{explicit}.
Just note that if
$\elem \tensor 1 \in \im (\submod \tensorr \frobe \ring \to \modul
\tensorr \frobe \ring)$,
then there exists also a finite subalgebra
$\ring \to \rings \to \frobe \ring$
such that
$\elem \tensor 1 \in \im (\submod \tensorr \rings \to \modul \tensorr \rings )$.
\end{proof}

\begin{proposition}
\label{frobeniuspureequivalent}
Let $\ring$ denote a ring containing a field of positive characteristic $\posp$.
Then the following are equivalent.

\numiii
\begin{enumerate}

\item
$\ring$ is $\frobf$-pure, that is, the Frobenius homomorphism is a
pure ring homomorphism.

\item
The Frobenius topology on $\specring$ is pure.

\item
The Frobenius closure is trivial on arbitrary submodules.

\item
If $\ring$ is a local Gorenstein ring, then this is also equivalent
with the property that the Frobenius closure is trivial on ideals.
\end{enumerate}
\end{proposition}
\begin{proof}
The equivalence (i) $\Leftrightarrow$ (ii) follows immediately from
Definition \ref{frobeniustopologydef} of the Frobenius topology and
from Definition \ref{puretopdef}. The equivalence of (ii) and (iii)
follows from Proposition \ref{pureequivalent}. (iii) implies (iv).
Under the Gorenstein hypothesis (iv) implies (iii) by
\cite[Corollary 1.5]{fedderfpure}.
\end{proof}

$\frobf$-pure rings have found a great deal of interest, including
their classification in low dimension, see e.g.
\cite{hochsterrobertsfrobenius}, \cite{fedderwatanabe},
\cite{fedderfpure}, \cite{mehtasrinivasfpure},
\cite{haraclassification}, \cite{harawatanabe}.

\subsection{Rings of global sections in the Frobenius topology}
\label{frobeniusglobalsubsection}

\

\medskip
The \emph{prefect closure} of a ring $\ring$ over a field of
positive characteristic $\posp$ is given as the colimit
$\ring^\perf= \colim_{\expoe \in \NN} \frobe \ring$, see
 \cite{sharpnossem}. For an $\ring$-module
$\modul$ we have $\colim \frobe \modul = \colim \modul \tensorr
\frobe \ring = \modul \tensorr \ring^\perf$. We will see that this
is also the module of global sections of $\modul_\frob$ in the
Frobenius topology.

\begin{proposition}
\label{frobeniusstalk} Let $\ring^\perf$ be the perfect closure of
$\ring$. The fixation given by $\ring^\perf$ defines an absolute
filter with absolute stalk $\ring^\perf$.
\end{proposition}
\begin{proof}
If $\ring \to \rings$ is a Frobenius cover and
$\rings \to \ring^\perf$ an $\ring$-algebra homomorphism,
then this factors through some $\frobe \ring$.
For a Frobenius cover $\rings \to \ringt$ over $\ring$ and a given fixation
$\rings \to \frobe \ring \to \ring^\perf$
there exists by definition
$\ringt \to \frobesec \rings \to \frobeplussec \ring$.
Hence there exists also an $\rings$-algebra
homomorphism $\ringt \to \ring^\perf$.
Hence this quasifilter is absolute (Definition
\ref{absolutefilterdef}). The image of $\rings$ under $ \varphi:
\rings \to \frobe \ring \to \ring^\perf$ gives also a Frobenius
cover (over $\Spec \ring$ and $\Spec \rings$), since the kernel of
$\varphi$ is nilpotent as $\varphi$ is spec-surjective.
Hence by Proposition \ref{filtertest} it is a filter and
$\ring^\perf$ is the absolute stalk.
\end{proof}

\begin{remark}
For a scheme $\spax$ of positive characteristic such that the
Frobenius morphism is finite ($\frobf$-finite schemes), the
Frobenius powers
$\frobe \spax \to \spax$, $\expoe \in (\NN,\geq)$
belong to the category and they also form an absolute filter with
absolute stalk $\ring^\perf$.
\end{remark}

\begin{theorem}
\label{frobeniusglobal}
Let $\ring$ denote a ring of positive
characteristic $\posp$ endowed with the Frobenius topology,
$\spax_\frob=(\specring)_\frob$. Let $\modul$ denote an
$\ring$-module with corresponding sheaf $\modul_\frob$ in the
Frobenius topology. Then for every finite Frobenius cover {\rm(}i.e.
for every object in $\catopen_\spax${\rm)} $\spay \to \spax$ we have
$\Gamma(\spay_\frob, \modul_\frob) = \modul \tensorr \ring^\perf $.
In particular, $\Gamma(\spax_\frob, \O_\frob) = \ring^\perf $.
\end{theorem}
\begin{proof}
Let $\spay= \Spec \rings$.
By Proposition \ref{frobeniusstalk} and Corollary
\ref{stalkpresheaf} we have a natural homomorphism $\Gamma(\spay,
\modul_\frob) \to \colim_{\indm \in \indexcat} \Gamma(\openfil_\indm ,
\modul_\frob)  = \modul \tensorr \ring^\perf $, where $\indexcat$ is
the indexing category corresponding to the fixation given by the
perfect closure of $\ring$. Let $\fuf \in \frobe \ring$ be
represented ($\expoe$ sufficiently large) by $\fuf \in \rings'$
finite over $\rings$,
$\ring \to \rings \to \rings' \to \frobe \ring$ finite.
By Lemma \ref{frobeniusdiagonalreduziert} the element
$\fuf \tensor 1-1 \tensor \fuf$ is nilpotent in
$\rings' \tensor_\rings \rings' $.
But nilpotent elements vanish in some Frobenius power, so $\fuf$
defines a global element in $\Gamma(\spay, \O_\frob)$.
It follows that
$\Gamma(\spay_\frob, \O_\frob) \to \ring^\perf$
is surjective and then also that
$\Gamma(\spay_\frob, \modul_\frob)\to \modul\tensorr \ring^\perf$
is surjective.
On the other hand, if
$\elem \in \Gamma(\spay_\frob, \modul_\frob)$
maps to $0$ in the stalk, then it is $0$ in some
$\modul \tensorr\frobe \ring$
and also in a finite Frobenius cover $\rings \to\rings'$, so
$\elem$ itself must be $0$.
\end{proof}

\ifthenelse{\boolean{book}}{ Alter Beweis

 First note that
$(\modul \tensorr S)_0 = \{ \eles \in \modul \tensorr S: \eles
\tensor 1 = 0 \in \modul \tensorr \frobe S \mbox{ for some } e \in
\NN\}$ and in particular $(\modul \tensorre )_0$ consists of all
elements which vanish in some higher Frobenius power $\frobesec
\modul$ or, equivalently, vanish in the limit $\modul \tensorr
\ring^\perf$. In particular, $\Gamma(R,{\O}_1)$ is the reduction of
$\ring$. Let $\spay= \frobe \spax$ and consider the commutative
diagram
$$ \xymatrix{\modul \ar[r] \ar[d] &  \modul \tensorr \frobe R \ar[r]^{p_1^*- p_2^*} \ar[d] &
 \modul \tensorr (\frobe R \tensorr \frobe R) \ar[d]   \\
\modul \ar[r] \ar[d] &  \modul \tensorr (\frobe R)_\red \ar[r]^{p_1^*- p_2^*} \ar[d]
& \modul \tensorr (\frobe R \tensorr \frobe R)_\red \ar[d]   \\
\modul \ar[r]  & \Gamma(\spay,\modul_1)=(\modul \tensorr \frobe R)_1
\ar[r]^{p_1^*-p_2^*\,\,\,\,\,\,\,\,\,\,\,\,\,\,\,\,\,\,\,\,\,\,\,} &
\Gamma(\spay \times_\spax \spay, \modul_1)\! =\! ( \modul \tensorr
(\!\frobe R \tensorr \frobe R)\!)_1 \, .
\\}
$$
Here the last row is the relevant complex to look at for compatible
sections. The arrows from the second to the last row exist, because
if $\eles \in \modul \tensorr S$ (where $S=\frobe \ring$ or $=\frobe
R \tensorr \frobe \ring$) vanishes in $\modul \tensorr \rings_\red$,
then it also vanishes in $\modul \tensorr S^\perf$ and so it
vanishes in some $ \modul \tensorr \frobesec S$. Therefore $\eles
\in ( \modul \tensorr S)_0$ and we get the factorization.

Let $\eles \in \modul \tensorre $. By Lemma
\ref{frobeniusdiagonalreduziert} the two projections induce
(identical) isomorphisms in the reductions, hence the mapping in the
middle row is $0$ and so $p_1^*(\eles)$ and $p_2^*(\eles)$ are
compatible. Therefore $\eles$ gives an element in
$\Gamma(\spax_\frob , \modul_\frob)$.

This gives compatible homomorphisms $\frobe \modul \to
\Gamma(\spax_\frob , \modul_\frob)$ and therefore a homomorphism $
\colim_{e \in \NN} \frobe \modul \to \Gamma(\spax_\frob ,
\modul_\frob)$. Every global section must be represented by an
element in a covering, so this homomorphism is surjective. On the
other hand, if a global section represented by $\eles \in \modul
\tensorre $ is $0$ in $\Gamma(\spax_\frob, \modul_\frob)$, then it
must be $0$ is some covering, that is, in some higher Frobenius
power. But then it is also $0$ in the limit. Hence we altogether get
an isomorphism $\colim_{e \in \NN} \frobe \modul =
\Gamma(\spax_\frob, \modul_\frob)$.
\end{proof}}

\begin{remark}
Note that already for a (non-perfect) field $\field$
the ring of global sections in the Frobenius topology is not the
field itself, as it is in the flat topology, though the Frobenius
homomorphism on a field is flat. This rests on the fact that the
Frobenius homomorphism on
$\frobe \field \tensor_\field \frobe \field$,
which enters the computation of the
global section ring in the Frobenius topology, is not flat.
\end{remark}

\ifthenelse{\boolean{book}}{
\begin{bookremark}
\label{frobeniusflabby} It follows also from Theorem
\ref{frobeniusglobal} that for all coverings $\frobe \spax \to \spax
$ the restrictions $\Gamma(\spax,\modul_\frob) \to \Gamma(\frobe
\spax, \modul_\frob)$ are isomorphisms. So these Frobenius sheaves
are ``flabby'' in some sense (the definition in \cite[Example III.
1.9b,c]{milne} is however different). We will see below (Proposition
\ref{frobeniuscohomology}) that these sheaves have no cohomology.
\end{bookremark}}{}

\subsection{Exactness in the Frobenius topology}
\label{frobeniusexactsubsection}

\begin{proposition}
\label{frobeniusexact} Let $\ring$ denote a ring of positive
characteristic and let $\lmr$ denote a complex of $\ring$-modules.
Let $\lmrr {_\frob}$ denote the corresponding complex of sheaves on
$(\specring)_\frob$. Then the following are equivalent.

\numiii

\begin{enumerate}

\item
The complex of sheaves $\lmrr {_\frob}$ is exact in the Frobenius
topology.

\item
The complex
$ \lmrtensor {\ring^\perf}$
is exact
{\rm(}this is the global complex as well as the
corresponding complex in the absolute stalk given by
$\ring^\perf${\rm)}.

\item
For every $\expoe \in \NN$ the complex $\lmrfrobe$ has the property
that $\ker (\cob) \subseteq (\im \coa)^\frobf$.

\item
For every
$\expoe \in \NN$
and every
$\elem \in \modm \tensorr\frobe \ring$
mapping to
$0$ in $\modr \tensorr \frobe \ring$
there exists
$\expoesec \geq \expoe$
and an element
$\elel \in \modl \tensor \frobesec \ring$
mapping to the image of
$\elem$ in $ \modm \tensorr \frobesec \ring$.
\end{enumerate}
\end{proposition}
\begin{proof}
The equivalence (i) $\Leftrightarrow $ (ii) was given in Proposition
\ref{exactnessabsolute}.
The equivalence (ii) $\Leftrightarrow $
(iii) (for which (iv) is an explicit version) follows directly and
follows also in the $\frobf$-finite case from Corollary
\ref{closurexexactcor}.\end{proof}
%
%

\begin{example}
The sequence $0 \to \ring \to \ring_\red$ is exact in the Frobenius
topology.
\end{example}

\begin{example}
Let $\ring = \field[\varx,\vary,\varz]/(\varx^3+\vary^3+\varz^3)$,
where $\field$ is a field of positive characteristic $\posp \neq
2,3$. It is known that $ \varz^2 \in (\varx,\vary)^\frobf$ if and
only if $\posp=2 \modu 3$, and in this case already $(\varz^2)^\posp
\in (\varx^\posp,\vary^\posp)$ \cite[Example 2.2]{hunekeparameter}.
Look at the complex of $\ring$-modules
$$\ring^2 \stackrel{\varx, \vary}{\lto} \ring \lto \ring/(\varx, \vary, \varz^2) \, .$$
This complex is of course not exact as $\ring$-modules. For $\posp=2
\modu 3$ it is however Frobenius-exact, since $(\varz^2)^\posp$ is
in the image of $\varx^\posp$ and $\vary^\posp$. Therefore
$(\varx,\vary)\ring^\perf=(\varx,\vary,\varz^2)\ring^\perf$ and so
it is exact. It is not Frobenius exact for $\posp=1 \modu 3$.
\end{example}

\begin{example}
Consider the domain $\ring = \field
[\varx,\vary,\varu]/(\vary^\posp-\varu \varx^\posp)$ with
normalization $\field[\vary/ \varx, \varx]$ and the ideal
$\ideal=(\varx,\vary)$. The first Frobenius power of $\ideal$ is
$\ideal^{[\posp]}=(\varx^\posp,\vary^\posp)=(\varx^\posp)$, so it is
the principal ideal generated by $\varx^\posp$. The ideal $\ideal$
is as an $\ring$-module isomorphic to $\ring^2/((\vary,-\varx),
(\varu\varx^{\posp-1},-\vary^{\posp-1})) \cong \ideal$, where
$\vecunit_1 \mapsto \varx$, $\vecunit_2 \mapsto \vary$. Therefore we
get isomorphisms (say $\posp$ is odd, $\potq= \posp^\expoe$)
$$ \ring^2/((\vary^\potq ,-\varx^\potq),
(\varu^\potq \varx^{\potq(\posp -1)},-\vary^{\potq (\posp-1)}))
\cong \ideal \tensorr \frobe \ring \, .$$ For $\expoe=1$
($\potq=\posp$), the element $(\varu,-1)$ on the left corresponds to
$\varx \tensor \varu - \vary \tensor 1$ in $\ideal \tensorr \!\!\!\
^{1}\! \ring$, and the natural mapping $\ideal \tensorr \!\!\! \
^{1}\! \ring \to \ideal^{[\posp]}$ sends this to $\varu
\varx^\posp-\vary^\posp= 0$. The element $\varx \tensor \varu -
\vary \tensor 1$ is however not $0$ in $\ideal \tensorr \!\!\!\
^{1}\! \ring$ nor does it become $0$ in $\ideal \tensorr \frobe
\ring$ for some $\expoe$. This is because $(\varu , -1) \in \ideal
\tensorr \!\!\!\ ^{1}\! \ring$ is mapped to $(\varu^{\potq/\posp},
-1^{\potq/\posp})$ in $ \ideal \tensorr \frobe \ring$, but this does
not belong to $((\vary^\potq ,-\varx^\potq), (\varu^\potq
\varx^{\potq(\posp -1)},-\vary^{\potq (\posp-1)}))$. Therefore  $0
\to \ideal \to \ring$ is exact, but not Frobenius-exact in the sense
of Proposition \ref{frobeniusexact}.
\end{example}

\ifthenelse{\boolean{extra}}{ Example where the complex of $\ring$
-modules is exact, but not Frobenius-exact. Due to right-exactness
this is only possible if the second mapping is not a surjection.
Look at injections $0 \to I \to \ring$ such that $\ideal \tensorr
\frobe R \to \frobe R $ is not injective.

What about $R=K[x,y,z,u,v,s,t]/(zu-x^p,zv-y^p,z-x^ps-y^pt)$ and the
ideal given by $(x,y)$. Its first Frobenius power is by construction
a principal ideal generated by $z$.}{}

\subsection{Cohomology in the Frobenius topology}

\label{frobeniuscohomologysubsection}

\

\begin{proposition}
\label{frobeniuscohomology}
Let $\ring$ be a ring of positive characteristic,
$\spax_\frob= \specring$
endowed with the Frobenius topology,
and let $\modul$ denote an
$\ring$-module with corresponding Frobenius sheaf $\modul_\frob$ on
$\spax_\frob$.
Then
$\chch^\numcoho(\spax_\frob, \modul_\frob)=0$ for $\numcoho \geq 1$
and in particular $H^{1}(\spax_\frob, \modul_\frob)=0$.
\end{proposition}
\begin{proof}
For every finite Frobenius cover $ \spay \to \spax$ we have $\Gamma(\spay,
\modul_\frob) = \modul \tensorr \ring^\perf$ by Theorem
\ref{frobeniusglobal}.
In particular we have
$\Gamma(\spay \timesx \cdots \timesx \spay, \modul_\frob) =
\modul \tensorr \ring^\perf $, and all restriction morphisms are
isomorphisms. Hence the \v{C}ech complex is
$$\Gamma(\spax_\frob,\modul_\frob)
=\ring^\perf \tensorr  \modul \stackrel{\id}{\lto}
\ring^\perf \tensorr  \modul \stackrel{0}{\lto} \ring^\perf \tensorr
\modul \stackrel{\id}{\lto} \ring^\perf \tensorr  \modul \to \cdots
\, .$$
Therefore
$\check{H}^\numcoho (\ring_\frob,\modul_\frob)=0$
for
$\numcoho \geq 1$ and $H^1(\ring_\frob,\modul_\frob)=0$ by Lemma
\ref{chchisoproj}.
\end{proof}

\begin{corollary}
\label{frobeniuscohomologycorollary}
Let $\ring$ be a ring of
positive characteristic.
Let
$\olmro$ be a short exact
sequence of $\ring$-modules with corresponding short exact sequence
of sheaves
$0 \to \ker = \modl^\frob \to \modm_\frob \to \modr_\frob\to 0$
on $\spax_\frob$
. Then $H^1(\spax_\frob, \modl^\frob)=0$. In the situation of
Proposition \ref{syzygycohomology} we get
$H^1(\spax_\frob, \Syz) = 0$.
\end{corollary}
\begin{proof}
The exactness of $\modm_\frob \to \modr_\frob \to 0$ induces by
Proposition \ref{frobeniusexact} the global exactness
$\Gamma(\spax_\frob,\modm_\frob)\!=\! \modm \tensorr \ring^\perf
\!\to\! \Gamma(\spax_\frob, \modr_\frob) \!=\! \modr \tensorr
\ring^\perf \to 0$. Since $H^1(\spax_\frob, \modm_\frob)$ is $0$ by
Proposition \ref{frobeniuscohomology}, it follows $H^1(\spax_\frob,
\modl^\frob)=0$ by the long exact cohomology sequence.

For the second statement let $\submod \subseteq \modul$ be finitely generated
$\ring$-modules and let $\ring^\numgen \to \submod$ be a surjection.
Then $\Gamma(\spax_\frob, \O_\frob^\numgen)
=(\ring^\perf)^\numgen
\to
\Gamma(\spax_\frob, \submod_\frob)= \submod\tensorr \ring^\perf$
is a surjection but also
$\submod \tensorr \ring^\perf
\to \Gamma(\spax_\frob, \submod^\frob)
\subseteq \Gamma(\spax_\frob, \modul_\frob)
= \modul \tensorr\ring^\perf$
is a surjection, as
$\Gamma(\spax_\frob, \submod^\frob)$ consists of
the elements in
$\Gamma(\spax_\frob, \modul_\frob)$ which are in the
Frobenius closure of $\submod$.
Hence
$\Gamma(\spax_\frob, \O_\frob^\numgen) \to \Gamma(\spax_\frob,\submod^\frob)$
is surjective and so again $H^1(\spax_\frob, \Syz)=0$.
\end{proof}

\ifthenelse{\boolean{extra}}{

\begin{remark}
Relation to phantom homology. Here the two things do not fit
together, because the vanishing is on the Frobenius category but the
existence is in the tight category.
\end{remark}

\begin{remark}
If $\ring \to \rings $ is not cyclic pure, then it cannot be
universally injective. For if $x \in \ring$, $x \notin I$ and $x \in
IS$, then the base change $\ring \to \ring/\ideal$ yields
$\ring/\ideal \to (\ring/\ideal) \tensorr S = S/IS$, where $\bar{x}
\neq 0$ is sent to $0$.
\end{remark}}{}

\section{Finite topology and plus closure}
\label{plussection}

\subsection{Plus closure for ideals and for submodules}
\label{plusclosuresubsection}

\

\medskip
The \emph{plus} (or \emph{finite}) \emph{closure} of an ideal
$\idealsubring$ in a domain $\ring$ is given by $$\ideal^+= \{\fuf
\in \ring:\, \exists \ring \subseteq \rings \mbox{ finite domain
such that } \fuf \in \ideal \rings\} .$$ The properties of this
closure operation depend heavily on the characteristic. If $\ring$
contains a field of characteristic zero and is normal, then
$\ideal^+=\ideal$
(this is due to the trace map, see
\cite[Remark 9.2.4]{brunsherzog}).
If however $\ring$ contains a field of positive characteristic, then
we have
$\ideal^\frobf \subseteq \ideal^+ \subseteq \ideal^*$,
where $\ideal^\frobf$ is the Frobenius closure
(Section \ref{frobeniussubsection})
and where $\ideal^*$ is the tight closure
\cite{hunekeapplication}, \cite{hochsterhunekebriancon}.
Moreover, it is conjectured that
$\ideal^+=\ideal^*$ \cite{hochstertightsolid}.
The plus closure may
also be characterized with the help of the so-called \emph{absolute
integral closure} $\ring^+$ of $\ring$. This is the integral closure
of $\ring$ inside the algebraic closure
$\overline{\quotfield(\ring)}$ of its fraction field
(first considered by M. Artin in \cite{artinjoinshensel}).
The plus closure is then just
$\ideal^+=\ring \cap \ideal\ring^+$
\cite[Theorem 1.7]{hunekeapplication}.
$\ring^+$ is a big Cohen-Macaulay algebra for excellent local domains of positive
characteristic
\cite[Theorem 7.0]{hunekeapplication}, \cite{hunekelyubeznik} or \cite{hochsterhunekeinfinitebig}.

We first generalize this notion for rings which are not domains and
for arbitrary submodules. For this we introduce finite covers.

\begin{definition}
A scheme morphism $\spay \to \spax$ is called a \emph{finite cover}
if it is finite and surjective.
\end{definition}

\begin{definition}
\label{finiteclosuremoduledefinition} For $\ring$-submodules
$\submodul $ we define the \emph{finite closure} of $\submod$ in
$\modul$ by
$$\submod^+ \!\!=\!\! \{ \elem \in \modul:\, \exists \, \ring \to \rings    \mbox{ a finite cover such that }
\elem \tensor 1 \in \im(\submod \tensorr \rings \to \modul \tensorr \rings) \}
\, .$$
\end{definition}

This definition coincides for an ideal in a domain with the standard
definition due to Proposition \ref{finitecoverproperties}(v).

\begin{proposition}
\label{finitecoverproperties}
Finite covers have the following properties.

\numiii
\begin{enumerate}

\item
The composition of finite covers is a finite cover.

\item
The base change of a finite cover is a finite cover.

\item
The reduction $\spax_\red \to \spax$ is a finite cover.

\item
If $\spax_\indcomp$ are the integral components of $\spax$, then $
\biguplus_\indcompinset \spax_\indcomp \to \spax$ is a finite cover.

\item
If $\spax$ is integral and $\spay \to \spax$ is a finite cover, then
there exists an integral component $\spay_0$ of $\spay$ such that
$\spay_0 \to \spax$ is a finite cover. This means that every finite
cover of an integral scheme can be refined by an integral finite
cover.
\end{enumerate}
\end{proposition}
\begin{proof}
These are clear. For (v) let $\pointy \in \spay$ map to the generic
point of $\spax$ and let $\pointy \in \spay_0 \subseteq \spay$ be an
integral component. Then $\spay_0 \to \spax$ is finite but also
surjective by going-up \cite[Proposition 4.15]{eisenbud}.
\end{proof}




\begin{definition}
\label{plustopologydefinition}
Let $\spax$ denote a scheme. The
\emph{finite topology} $\catopenspax $ on $\spax$ is given by the
schemes over $\spax$,
where all morphisms are finite and where the finite covers are covers.
The structural morphisms are declared to be covers. We denote
$\spax$ endowed with this topology by $\spax_\fin$.
\end{definition}

\begin{proposition}
\label{finitefinitetopology}
The finite topology is an affine single-handed Grothen\-dieck
topology. Every scheme morphism $\spay \to \spax   $ yields a
corresponding site morphism $\spay_\fin \to \spax_\fin$. The closure
induced by the finite topology is the plus closure.
\end{proposition}
\begin{proof}
This is clear from Proposition \ref{finitecoverproperties} and Lemma
\ref{explicit}.
\end{proof}

\begin{remark}
\label{qfh}
In the $qfh$-\emph{topology} introduced by Voevodsky
(\cite{voevodskyhomology}, \cite{suslinvoevodskysingularhomology})
a cover is given by a topological epimorphism
(or universal submersion)
$\biguplus \spax_\indcov \to \spax$
such that each $\spax_\indcov \to \spax$ is a quasi-finite morphism
\cite[D\'{e}finition 6.2.3]{EGAII}. So this quasifinite topology is
a refinement of the finite topology. On the other hand, over a
complete local domain it follows from
\cite[Corollaire 6.2.6]{EGAII} that a quasifinite cover
contains a connected component which is a finite cover of the
$\spax$ (this is not true for a quasifinite surjection without the submersive property).
\end{remark}

\subsection{The ring of global sections in the finite topology}
\label{finiteglobalsubsection}

\

\medskip

\begin{definition}
\label{compatibleintegralclosure} Let $\ring$ denote a ring and let
$\ring \to \ringl$ be an $\ring$-algebra. Then we set
$$\ring^{\compic/\ringl} = \{ \fuglob \in \ringl: \mbox{ there exists a finite }
\ring\mbox{-algebra } \ringt \subseteq \ringl \mbox{ such that }
$$
$$\hspace{3.1cm}
\fuglob \in \ringt \mbox{ and } \tensortwistfuglob \in \ringt \tensorr
\ringt \mbox{ is nilpotent}\}\, $$
\and call it the \emph{compatible
integral closure} of $\ring$ in $\ringl$.
\end{definition}

Compare \cite[Section I]{manaresiseminormal} for this definition. By
the
\emph{compatible normalization}
of a domain $\ring$ we mean
the compatible integral closure of $\ring$ in
$\ringl=\quotfield(\ring)^\perf$, and denote it by $\ring^\comp $.
Recall that the \emph{seminormalization} of a ring is the biggest
ring inside the normalization which does not separate points and
does not change the residue class fields, see \cite[Definition
7.2.1]{kollarrational}. The following characterization of the
seminormalization is due to Manaresi, see \cite[Theorem
I.6]{manaresiseminormal}. We include a
proof.

\begin{lemma}
\label{finiteglobalcompatible} Let $\ring$ denote a noetherian
domain in characteristic zero. Then the seminormalization is the
compatible normalization.
\end{lemma}
\begin{proof}
The seminormalization $\Spec \ring^\seminorm \to \specring$ is a
homeomorphism and induces isomorphisms $\fieldres(\point) \to
\fieldres(\point')$ for $\point \in \specring$. The same is true for
$\ring \subseteq \rings \subseteq \ring^\seminorm$. Let $\fuglob \in
\rings \subseteq \ring^\seminorm $ be given, $\ring \subseteq
\rings$ finite, and consider $ \tensortwistfuglob \in \rings \tensorr
\rings \to \ring^\seminorm \tensorr \ring^\seminorm $. For every
ring homomorphism $\ring \to \fieldres(\point)$ there exists exactly
one ring homomorphism $\ring^\seminorm \to \fieldres(\point)$.
Therefore by Lemma \ref{nilpotenttest} we know that
$\tensortwistfuglob $ is nilpotent (in $\rings \tensorr \rings$ and
in $ \ring^\seminorm \tensorr \ring^\seminorm $). Hence $\fuglob \in
\ring^\comp$.

Now suppose that an element $\fuglob$ in the normalization is
compatible. We know by Lemma \ref{compatiblevalue} that it has a
uniquely determined value in $\fieldres(\point)$ at each point
$\point \in \specring$. Therefore $\fuglob$ belongs to the
seminormalization.
\end{proof}

\begin{proposition}
\label{finiteglobalsection} Let $\ring$ denote a noetherian domain,
$\spaxeqspecring$. Then the ring of global sections in the finite
topology is the compatible normalization inside the perfect closure
$\quotfield(\ring)^\perf$, i.e. $\Gamma(\spax_\fin, \O_\fin)=
\ring^\comp$.
\end{proposition}
\begin{proof}
The natural morphism $(\Spec \quotfield(\ring))_\fin \to \spax_\fin$
induces the ring homomorphism $\Gamma(\spax_\fin, \O_\fin) \to
\Gamma(\quotfield(\ring)_\fin, \O_\fin)$. First note that
$\Gamma(\quotfield(\ring)_\fin, \O_\fin)= \quotfield(\ring)^\perf $,
by Theorem \ref{frobeniusglobal} and Proposition
\ref{globalsurjectivefield}.

An element $\fuglob \in \Gamma(\spax, \O_\fin)$ is represented by
$\fuglob \in \rings$, where $\ring \to \rings$ is a finite cover. The
compatibility condition for the sheafification means exactly that
$\tensortwistfuglob \in \rings \tensorr \rings$ is nilpotent. Hence
the natural ring homomorphism maps into the compatible
normalization.

Suppose that
$\fuglob \in \ring^\comp \subseteq \quotfield(\ring)^\perf$
and that it is represented by
$\fuglob \in \rings$, $\ring \to \rings $ a finite cover.
Then
$\tensortwistfuglob \in \rings \tensorr \rings$ is nilpotent.
Since reductions are finite covers
(Proposition \ref{finitecoverproperties}(iii)),
this means that $\fuglob \in \rings$ is compatible in the sense of
sheafification and gives a global element in $\Gamma(\spay,
\O_\fin)$. Hence the homomorphism is surjective.

Suppose that
$\fuglob \in \Gamma(\spax_\fin, \O_\fin)$ maps to $0$. We may
assume that $\fuglob$ is represented by
$\fuglob \in \rings$, where $\rings$ is a domain
(Proposition \ref{finitecoverproperties}(v)).
Hence the image of $\fuglob $ in $\Gamma(\quotfield(\ring), \O_\fin)$
is represented by $\fuglob$ considered in $\quotfield(\rings)$. So
$\fuglob=0$.
\end{proof}

\ifthenelse{\boolean{book}}{
\begin{bookexample}
Consider $R=K[x,y]/(x^2-y^3)$ with the normalization
$$R \lto K[t],\, \, \, x \longmapsto t^3, y \longmapsto t^2 \, .$$
We show that $t$ is compatible. We
have to consider the difference of $1 \tensor t$ and $t \tensor 1$
in $K[t] \tensor_{K[t^2,t^3]} K[t]$. We have $(1 \tensor t)^3=x =(t
\tensor 1)^3$ and $(1 \tensor t)^2=y =(t \tensor 1)^2$. Therefore we
have
$$0=(1 \tensor t)^2 - (t \tensor 1)^2
= (1 \tensor t-t \tensor 1)(1 \tensor t+t \tensor 1)$$ and
$$ 0=(1 \tensor t)^3 - (t \tensor 1)^3
= (1 \tensor t-t \tensor 1)
((1 \tensor t)^2+  (1\tensor t)(t \tensor 1) +(t \tensor 1)^2)\, .$$
On the other hand, $(1 \tensor t - t \tensor 1)^2$ lies in the ideal
generated by
$1 \tensor t+t \tensor 1$ and by $(1 \tensor t)^2+
(1\tensor t)(t \tensor 1) +(t \tensor 1)^2$. Therefore $(1\tensor
t-t \tensor 1)^3=0 $ and so $ 1 \tensor t - t \tensor 1$ is
nilpotent, hence $t$ is compatible in the finite topology
(and also in the surjective or the submersive topology).
Therefore the ring of
global section for $\ring$ in the finite topology is just its
normalization, if $\field$ has characteristic zero.
\end{bookexample}}{}

\begin{remark}
The seminormalization is also the ring of global sections in the
submersive topology resp. $h$-topology
(see \cite[Proposition 3.2.10]{voevodskyhomology} and \cite{blicklebrennersubmersion})
The normalization cannot be the ring of global sections for a
Grothendieck topology which is functorial for all scheme morphisms,
since the normalization is not functorial for quotient rings, see
the following example.
\end{remark}

{\ren{\ring}{{R}}
\begin{example}
\label{normalnonfunctorialexample} Let $\ring \subset \CC[[\varx]]$
denote the subring consisting of formal power series with constant
term in $\RR$. Then $\frac{\varx \imu}{\varx}= \imu$ and $\imu^2=-1
\in \ring$, hence $\imu$ is both rational and integral over $\ring$
and so $\CC[[\varx]]$ is the normalization of $\ring$. The residue
class field of $\ring$ with respect to the maximal ideal
$\fom_\ring=(\varx,\imu \varx)$ is the (normal) field $\RR$, but
there does not exist an $\ring$-algebra homomorphism between the
normalizations $\CC[[\varx]] \to \RR$. In particular $\imu$ is not
compatible, and $\Gamma(\ring_\fin,\O_\fin)=\ring$.
\end{example}}


%
%

\subsection{Absolute stalks in the finite topology}
\label{finitestalksubsection}

\

\medskip
There are several relevant absolute filters in the finite topology
to study. Already for an integral scheme it is not clear which one
is the best choice.

\begin{proposition}
\label{plusabsolutestalk}
Let $\ring$ denote a noetherian ring with integral components
$\ring_\indcomp$, $\indcompinset$.
Set $\ring^+:= \prod_\indcompinset (\ring_\indcomp)^+$. Then
$\prod_\indcompinset \overline{\quotfield(\ring_\indcomp)}$ and also
$\ring^+$ define an absolute filter in $\spax_\fin
=(\specring)_\fin$ with absolute stalk $\ring^+$.
\end{proposition}
\begin{proof}
By Proposition \ref{finitecoverproperties}(v) we may assume that $\ring$
is a domain, because we can treat the components separately.
Every ring homomorphism
$\rings \to \cloalg{\quotfield(\ring)}$, $\rings$ a finite cover,
factors uniquely through $\ring^+$, so we restrict to the fixation
given by $\ring^+$.

Note first that for every finite cover $\ring \to \rings$
there exists an $\ring$-algebra homomorphism $\rings \to \ring^+$.
According to Proposition \ref{finitecoverproperties}(v) there exists
a homomorphism $\rings \to \ringcomp$ to a domain which is a finite
extension of $\ring$. Hence we may assume that $\rings$ is a domain.
Then $\quotfield(\ring) \to \quotfield(\rings)$ is a finite field
extension and hence there exists an inclusion $\quotfield(\rings)
\subseteq \cloalg{\quotfield(\ring)}$. This gives an inclusion
$\rings \subseteq \ring^+$. So in the quasifilter given by the
fixation with $\ring^+$ every finite cover occurs. It is a filter by
Proposition \ref{filtertest}, since the image of a finite cover
$\rings$ in $\ring^+$ is itself a finite (integral) cover over
$\ring$.

To show that it is absolute let $\rings \to \ring^+$ be an
$\ring$-algebra homomorphism and let $\rings \to \ringt$ be in
$\catopen_\ring$. Then let $\bar{\rings}$ be the image of $\rings$
in $\ring^+$ and let $\bar{\ringt}$ be a $\ringt$-algebra which is a
finite domain over $\bar{\rings}$. Then there exists an
$\bar{\rings}$-embedding of $\bar{\ringt}$ into $\ring^+$. This
gives $\ringt \to \ring^+$ over $\rings$. The stalk in this absolute
filter is of course $\ring^+$.
\end{proof}

\begin{example}
\label{finabsoluteintegral} It is also useful to consider the filter
where the objects are only the domains $\rings$ with a fixation
$\ring \to \rings \to \ring^+$. This is a quasifilter since for two
finite domains $\ring \subseteq \rings, \rings'$ there exists a
common finite domain in $\ring^+$ containing both. It is again an
absolute filter by similar arguments as in the proof of Proposition
\ref{plusabsolutestalk}.
The only morphisms in this filter are inclusion of domains.
Moreover, it is a covering filter
(Remark \ref{filterabsoluteremark}), since every finite
$\ring$-homomorphism $\rings \to \ringt$ is a cover
when $\rings$ is a domain (because $\ringt \in \catopen_\ring$ is a
finite cover over $\ring$).
\end{example}

{\ren{\var}{{x}}
\begin{example}
\label{finabsoluteextra}
Consider the affine line
$\AA^1_\field=\Spec \ring$, $\ring=\field[\var]$ in the finite topology and
consider the fixation given by
$\ring^+ \oplus \field^+$
(where in the second component $\var$ is sent to $0$).
This is again an absolute filter with absolute stalk
$\ring^+ \oplus \field^+$.
The sequence
$0 \to \ring \stackrel{\var}{\to} \ring$ is not exact
in this absolute filter, since it is not exact in $\field^+$.
Proposition \ref{exactnessabsolute}(i)-(v) do not hold, but (vi)
and (vii) hold, since for given $\spay $ the sequence is not
$\catopen$-exact in general, but its restriction to an integral component of
dimension one is. Note that these inclusions belong to the category,
but are not always indicated by the filter
(e.g. $\AA^1 \subset \AA^1 \uplus \{\point\}$
is not indicated).
\end{example}}

\begin{proposition}
Let $\ring$ denote a noetherian ring with integral components
$\ring_\indcomp$, $\indcompinset$. Let $\submodul $ be
$\ring$-modules, $\elem \in \modul$. Then the following are
equivalent.

\numiii

\begin{enumerate}

\item
$\elem \in \submod^+$.

\item
$\elem \in \submod_\indcomp^+$, where $\submod_\indcomp$ is the
image of $\submod \tensorr \ring_\indcomp$ inside
$\modul_\indcomp=\modul \tensorr \ring_\indcomp$.

\item
There exists a finite and spec-surjective ring homomorphism $\ring
\to \rings $ such that $\elem \tensor 1 \in \im(\submod \tensorr \rings \to
\modul \tensorr \rings)$.

\item
$\elem \in \submod^\tocl$, where $\tocl$ denotes the closure
operation induced by the finite topology on $\specring$.

\item
$\elem \tensor 1 \in \im(\submod \tensorr  \ring^+ \to \modul \tensorr
\ring^+)$, where $\ring^+$ is the absolute stalk.
\end{enumerate}
\end{proposition}
\begin{proof}
(i) and (iii) are equivalent by Definition
\ref{finiteclosuremoduledefinition}. The defining property (iii)
might be checked on the components by Proposition
\ref{finitecoverproperties}, hence (i) and (ii) are equivalent. The
equivalence of (iii) and (iv) was mentioned in Proposition
\ref{finitefinitetopology}. The equivalence of (iv) and (v) follows from
Corollary
\ref{topclosureabsolute} and Proposition \ref{plusabsolutestalk}.
\end{proof}

\begin{remark}
The exactness of a complex in the finite topology is a strong
condition, since $\spax \cup \{\point\} \to \spax$ is a finite
cover, where $\point $ is a closed point of $\spax$. If a sheafified
complex is exact on $\spax_\fin$, then this implies that it is exact
on every closed point! Hence the exactness in the finite topology
for a scheme of finite type over a field implies the exactness in
the surjective topology (Theorem \ref{surjectiveexact} and Remark
\ref{surjectiveexactclosedpoint}). The converse is not true, as
$\ring^\numgen \stackrel{\runfuf}{\to}  \ring \to \ring/ \rad
(\ideal) \to 0$ shows, where $\ideal =(\runfuf) \subseteq
\ideal^+ \subset \rad(\ideal)$. 
%
\end{remark}

We characterize now the exactness in the absolute filter given by
$\ring^+$.

{
\ren{\ringsec}{{\rings}}

\begin{proposition}
\label{plussexact} Let $\ring$ denote a domain endowed with the
finite topology and let $\lmr$ denote a complex of $\ring$-modules.
Let $\lmrr {_\fin}$ denote the corresponding complex of sheaves on
$(\specring)_\fin$. Then the following are equivalent.

\numiii

\begin{enumerate}

\item
The complex of sheaves $\lmrr {_\fin}$ is exact in the absolute
filter given by $\ring^+$.

\item
The complex of sheaves $\lmrr {_\fin}$ is exact in the integral
absolute filter given by $\ring^+$ where only integral rings occur
{\rm(}Example \ref{finabsoluteintegral}{\rm)}.

\item
The complex $ \lmrtensor {\ring^+} $ is exact.

\item
For every finite cover $\ring \to \rings \to \ring^+$ and every $\elem \in \modm
\tensorr \rings$ mapping to $0$ in $\modr \tensorr \rings$ there
exists a finite $\rings$-algebra $\ringt \to
\ring^+$ {\rm(}which might be chosen to be integral{\rm)} such that there exists
$\elel \in \modl \tensorr \ringt$ mapping to the restriction of
$\elem$.


\item
The complex of sheaves $\lmrr {_\fin}$ is $\catopen$-exact in the
finite topology on every integral object.

\item
For every finite cover $\ring \to \rings$ and every $\elem \in \modm
\tensorr \rings$ mapping to $0$ in $\modr \tensorr \rings$ there
exists a morphism $\rings \to \ringt$ in $\catopen_\ring$ and
$\elel \in \modl \tensorr \ringt$ mapping to the restriction of
$\elem$.

\item
For every finite cover $\ring \to \ringsec $ with $\ringsec$
integral the complex $\lmrr {\tensorr \ringsec  }$ has the property
that $\ker (\cob) \subseteq (\im \coa)^+ $.

\end{enumerate}
\end{proposition}
\begin{proof}
This follows from Proposition \ref{exactnessabsolute} applied to the
(integral) absolute filter given by $\ring^+$. The equivalence of
(v) and (vii) was proven in Proposition \ref{closurexexact}.
\end{proof}
}

\subsection{Cohomology in the finite topology}

\label{finitecohomologysubsection}

\

\medskip
The computation of cohomology in the finite topology , e.g. of the
structure sheaf, provides us with several difficulties. So we give
here only a result concerning the normalization of a curve and some
examples.

{\ren{\point}{{P}}

\begin{proposition}
\label{finitenormalcohomology}
Let $\curve = \Spec \ring$ be an affine integral curve over an
algebraically closed field $\field$ and let $\curve^\nor \to
\curve$ its normalization. Then
$\chch^1(\curve^\nor \to \curve,\O_\fin)=0$.
\end{proposition}
\begin{proof}
Let $\sing \subset \curve$ be the set of singular points. For
$\point \in \sing$ let $\point_1 \comdots \point_{\numnor_\point}$
be the points in $\curve^\nor$ mapping to $\point$. Then
$$\curve^\nor \times_\curve \curve^\nor
 = \curve^\nor \uplus \{\point_{\inda \indo}, 1 \leq \inda , \indo \leq \numnor_\point,
 \inda \neq \indo :\point \in \sing\} \, $$
(the points $\point_{\inda \inda}$ are the points
$\point_\inda \in \curve^\nor$),
and similarly
$\curve^\nor \times_\curve \curve^\nor \times_\curve\curve^\nor
 = \curve^\nor \uplus \{ \point_\indaou \}$.
The indices of the points show where to one has to project, e.g.
$\projcech_1(\point_\indaou)= \project_{2,3}(\point_\indaou )=
\point_\indou$. Let now a cocycle $ \fuh \in \ring^\nor \tensorr \ring^\nor$ (this is seminormal)
be given (Section \ref{cohomologysubsection}).
On the diagonal
$\diag=\curve^\nor$ the cocycle condition gives immediately
$\fuh|_{\curve^\nor} =0$.
On the isolated points the cocycle condition gives
$$(\projcech_1^* - \projcech_2^* + \projcech_3^*) (\fuh) (\point_\indaou)
= \fuh (\point_\indou)- \fuh(\point_\indau) +
\fuh(\point_\indao)=0 \, .$$
For every $\point \in \sing$ this is exactly the cocycle condition for the
complex
$\field^\setind \to \field^{\setind^2} \to\field^{\setind^3}$,
where
$\setind=\{\point_1 \comdots\point_{\numnor_\point}\}$
(the \v{C}ech complex of the cover $\setind \to \{\point\})$.
Hence there exist values $\fuh(\point_1)
\comdots \fuh(\point_{\numnor_\point})$
such that the differences of the pull-backs under the two
projections are
$\fuh(\point_\indao)$.
Since $\curve^\nor$ is affine there exists a global function
$\fug \in \Gamma(\curve^\nor, \O_{\curve^\nor})$
such that
$\fug( \point_\inda) =\fuh ( \point_\inda) $ for all
$\point_\inda$, $\point \in \sing$,
$\inda =1 \comdots \numnor_\point$.
Then the derivation maps $\fug$ to $\fuh$, and so $\fuh$
is a coboundary.
\end{proof}

}

{\ren{\varz}{{z}}
\begin{example}
\label{selfproductnotseminormal}
One of the difficulties in computing cohomology in the finite topology
is that the ring of global sections is the seminormalization. This
might bring in new cocycles which are not seen on the presheaf
level. Let $\field$ be an algebraically closed field and consider
the (flat) extension
$\ring =\field[\varx]
\to \field[\vary]=\rings$,
$\varx \mapsto \vary^\expon$.
Then
$\rings \tensorr \rings
= \field [\vary,\varz]/(\vary^\expon-\varz^\expon)
=\field [\vary,\varz]/(\prod_{\zeta_\expon} (\vary-\zeta_\expon \varz))$.
This ring consists of $\expon$ lines in the
plane meeting in the origin
(say the characteristic does not divide $\expon$).
So $\Spec \rings \tensorr \rings$ has embedding dimension $2$ and is for
$\expon \geq 3$ not seminormal,
its seminormalization is the scheme of $\expon$ axes
with embedding dimension $\expon$
\cite[Section 3]{brennercontinuous}.
Hence
$\Gamma(\rings\tensorr\rings,\O_\fin) \neq \rings\tensorr\rings$.
\end{example}

}

{\ren{\varu}{{v}} \ren{\point}{{P}}

\begin{example}
\label{onegeneratorfinite}
Let
$\ring = \field[\vara,\varu]/(\vara\varu)$ with normalization
$\ringnorm=\field[\varasec] \oplus \field[\varusec]$,
$\vara=(\varasec , 0)$ and $\varu =(0,\varusec)$.
Consider the principal
ideal $\ideal=(\fug)=(\vara + \varu)$
(cf. Example \ref{onegeneratorexplicit}).
Then in $\ringnorm$ we have
$\fuf = \vara = (\varasec,0) = (1,0)(\varasec, \varusec) =(1,0) (\vara+\varu)$.
So $\fuf \in \ideal^+$ (but $\fuf \not\in \ideal$).
Set $\spax =\Spec \ring$ and
$\spay=\Spec \ringnorm$. We want to compute
$\chch^1(\spaxnorm \to \spax, \Syz)$, where $\Syz$ is given as the kernel of
$\fug: \O_\fin\to \ideal^\fin \subseteq \O_\fin$
(so $\Syz$ is the annihilator of $\fug$),
and
$\homcon(\fuf) \in \chch^1(\spaxnorm \to \spax, \Syz)$
(the cohomology of the structure sheaf is $0$ as in Proposition
\ref{finitenormalcohomology}).

$ \spaxnorm \timesx \spaxnorm$ consists of two lines $\linel_1$,
$\linel_2$ and two points $\point_{1,2}$ and $\point_{2,1}$.
$ \spaxnorm \timesx \spaxnorm \timesx \spaxnorm$
consists again of two lines
$\linel_1$, $\linel_2$ and
$8$ points $\point_{\indi,\indj,\indk}$, where $\point_{1,1,1} \in
\linel_1$ and $\point_{2,2,2} \in \linel_2$ and where the other six
points are isolated components.

The pull-back of the function $\gen= \vara + \varu$ to
$ \spaxnorm \timesx \spaxnorm$ and to
$\spaxnorm \timesx \spaxnorm \timesx \spaxnorm$ is the function which is
$\varasec$ on $\linel_1$ and $\varusec$ on $\linel_2$
and has value $0$ on the isolated points. So the annihilator of
$\fug$ consists of the functions which are zero on the two lines and
have arbitrary values on the isolated points. Let
$\fuh \in \Gamma(\spay \timesx \spay,\Syz)$
have values
$\fuh(\point_{1,2})=\elea$ and $\fuh(\point_{2,1})=\eleb$
(and $=0$ on the lines).
The cocycle condition is
$(\project_{2,3}^* - \project_{1,3}^* +\project_{1,2}^*) (\fuh)=0$.
This means evaluated at $\point_{1,2,1}$ that
$0 = \fuh(\point_{2,1})- \fuh(\point_{1,1}) + \fuh(\point_{1,2})
= \elea -0 + \eleb$
($\point_{2,1,2}$ yields the same, the other points yield no
condition).
So the cocycle condition reduces to $\eleb =- \elea$, and
since there are no coboundaries at all, we infer that $\chch^1(\spay
\to \spax, \Syz)$ has dimension one.

We come to the computation of $\homcon (\fuf)$. In $\ringnorm$ we
have $\fuf=\eleq \gen$, $\eleq = (1,0)$, so we have to compute the
difference between the pull-backs of $\eleq$ under the two
projections.
The pull back of $\eleq$ under both projections gives on $\linel_1$
the function $1$ and on $\linel_2$ the function $0$. The values of
the pull backs under the projections at the points are
$\project_1^*(\eleq)(\point_{1,2})= 1$,
$\project_1^*(\eleq)(\point_{2,1})= 0$,
$\project_2^*(\eleq)(\point_{1,2})= 0$,
$\project_2^*(\eleq)(\point_{2,1})= 1$. Hence the difference is $0$
on $\linel_1$ and $\linel_2$ and has value $1$ at $\point_{1,2}$ and
$-1$ at $\point_{2,1} $. So this difference gives an element in
$\Gamma(\spaxnorm \timesx \spaxnorm, \Syz)$
which represents the cohomology class
$\homcon (\fuf) \in H^1(\spax,\Syz)$.
\end{example}

}

\section{Topologies defined by Zariski filters and torsion theories}
\label{zariskifiltersection}

Some properties of a variety, or objects associated to it, depend
only on the points of codimension one, i.e. they are unchanged if we
remove a subset of codimension at least two, like e.g. the divisor
class group. Other objects get simpler after removing small subsets,
e.g. torsion-free sheaves are locally free in codimension one.
A (hereditary) torsion theory for $\ring$-modules studies which part
of a module can be neglected in a certain sense.
The notions of a pure resp. parafactorial pair are defined by requiring that the restriction
of an \'{e}tale covering resp. an invertible sheaf to an open subset
gives an equivalence of categories
(\cite[Expos\'{e}s X.3 and XI.3]{SGA2}, \cite[D\'{e}finition
21.13.1]{EGAIV}).
In birational geometry one may restrict to any non-empty open
subset. So smaller open subsets may have the full or the essential
information, hence they are covers in a certain sense.
We formalize these observations in this part.

\ren{\obcov}{{U}}

\subsection{Zariski filters}
\label{zariskifiltersubsection}

\

\begin{definition}
\label{zariskifilterdef}
A \emph{Zariski filter} on a scheme $\spax$ is a \emph{topological
filter} $\filtzar $ in the Zariski topology, i.e. a non-empty
collection of open subsets $\obfil \subseteq \spax$ such that
$\obfil   \in \filtzar$ and $\obfil  \subseteq \obfilu$ implies $\obfilu \in \filtzar$ and such that
$\obfil , \obfilu \in \filtzar$ implies $\obfil  \cap \obfilu \in \filtzar$.
\end{definition}

A Zariski filter
on $\spax$ defines for every morphism
$\morbc :\spay \to \spax$ a
Zariski filter
$\morbc^*(\filtzar)$ on $\spay$ by setting
$\morbc^*(\filtzar)
= \{\obfilw \subseteq \spay :\, \obfilw \supseteq
\morbc^{-1}(\obfil) \mbox{ for some } \obfil \in \filtzar\}$.
In particular, for an open subset
$\openzar \subseteq \spax$
this is
$\filtzar_\openzar =\filtzar|_\openzar
=\{\openzar \cap \obfil:\obfil \in\filtzar\}$ (we will also write $\obfilw \in \filtzar_\openzar$
if we mean $\openzar \cap \obfilw \in \filtzar_\openzar$).

\begin{definition}
Let $\spax$ denote a scheme and let $\filtzar $ denote a Zariski
filter on $\spax$. The Grothendieck topology $\spax_\zftop$
\emph{induced} by $\filtzar $ is given by all (Zariski-)open subsets of
$\spax$, and
$\obcovindic \to \openzar$, $\indcovinset$, is a covering if and only if
$\bigcup_\indcovinset \obcovindic \in \filtzar_\openzar$, that is, $\openzar \cap
\obfil \subseteq \bigcup_\indcovinset \obcovindic$ for some $\obfil \in \filtzar$.
\end{definition}

\begin{remark}
Here we design the Grothendieck topology given by a Zariski filter
to be a non-affine topology which is a refinement of the Zariski
topology. The open sets $\obfil \in \filtzar$ define the single
covers of $\spax$. This topology is flat, but not faithfully flat.
If $\spax$ and $\spay$ are schemes with Zariski filters $\filtzar
_\spax$ and $\filtzar _\spay $, then a morphism $\moryx$ induces a
site morphism
$\spay_\zftop \to \spax_\zftop$
if and only if
$\morbc ^*(\filtzar_\spax) \subseteq \filtzar_\spay$.
\end{remark}

\begin{example}
Consider a pair $(\spax_0,\spax)$ consisting of a scheme $\spax$ and
a subset $\spax_0 \subseteq \spax$.
All open neighborhoods
$\obfil \supseteq \spax_0$ form a filter denoted by $\filtzar (\spax_0)$.
A morphism
$\morbc: \spay \to \spax$ which maps
$\spay_0\subseteq\spay$ into $\spax_0$ induces a site morphism.
The disjoint union
$\biguplus_{\point \in \spax_0} \Spec \O_\point$
(but also
$\biguplus_{\point \in \spax_0} \Spec \fieldres(\point)$)
gives an absolute cover
(Definition \ref{absolutecoverdef})
for this topology.
\end{example}

\begin{example}
\label{zariskifiltercodimension}
For a natural number $\numk \in
\NN$ and a scheme $\spax$ denote by $\spax^{(\numk)}$ the points of
height at most $\numk$, that is,
$$\spax^{(\numk)} := \{ \point \in \spax: \, \dim \O_{\spax,\point} \leq \numk \} \, .$$
Let $\filtzar ^\numk_\spax$ denote the filter consisting of all open
subsets which contain $\spax^{(\numk)}$.
A morphism
$\morbc:\spay\to\spax$ with the property that
$\morbc(\spay^{(\numk)})\subseteq \spax^{(\numk)}$
(so that there are no contractions of points of height $\numk$ to points of bigger height)
induces a site morphism.
For $\numk=0$ see Example \ref{densetopology}, for
$\numk=1$ see Section \ref{divisorialsubsection}.
\end{example}

\begin{proposition}
\label{zariskifilterproperties}
Let a Zariski filter $\filtzar $ be given on a topologically
noetherian scheme $\spax$ given rise to the Grothendieck topology
$\spax_\zftop$. Let $\modul$ be a quasicoherent $\O_\spax$-module
and let $\submodul$ be a submodule. Then the following holds.

\numiii
\begin{enumerate}

\item
The module of global sections is
$\Gamma(\spax_\zftop, \modul_\zftop)
= \colim_{\obfil \in \filtzar} \Gamma(\obfil, \modul)$.

\item
The assignment $\openzar \mapsto \Gamma(\openzar_\zftop, \modul_\zftop)$ for
$\openzar \subseteq \spax$ open is a quasicoherent $\O_\spax$-module.

\item
The sheafification $\modul \mapsto \modul_\zftop$ is exact.

\item
The Zariski filter $\filtzar $ is an absolute filter in the
Grothendieck topology, and the absolute stalk of a module is the
same as its module of global sections.

\item
The closure operation induced by $\spax_\zftop$ is {\rm(}on an open
affine subset{\rm)}
$$\submod^{\zftop - \clop}
= \{ \elem \in \modul:\rest^\spax_\obfil (\elem)
\in \Gamma(\obfil, \submod) \mbox{ for some }
\obfil \in \filtzar \} \, .$$
If $\filtzar = \filtzar(\spax_0)$
for a subset $\spax_0 \subseteq \spax$, then
$$\submod^{\zftop - \clop}
= \{ \elem \in \modul:
\elem_\point \in \submod_\point \mbox{ for all } \point \in \spax_0\} \, .$$

\item
If $\modul$ is coherent, then $\modul_\zftop = 0$ if and only if
$\modul |_\obfil=0$ for some $\obfil \in \filtzar$, and this is true if
and only if the support of $\modul$ does not meet $\obfil$.

\item
We have
$H^\numcoho (\spax_\zftop, \modul_\zftop)
= \colim_{\obfil \in \filtzar} H^\numcoho (\obfil, \modul)$.
\end{enumerate}
\end{proposition}
\begin{proof}
(i). The (Zariski-)module
$\openzar \mapsto \Gamma(\openzar, \modul)$ is the
associated $\O_\zftop$-premodule in the sense of
Proposition \ref{sheafzariskicompare}(iv), which we have to sheafify.
We have
$$\Gamma(\openzar, \modul_0)
= \{\elem \in \Gamma(\openzar,\modul):\, \rest^\openzar_{\openzar \cap \obfil}(\elem)
=0 \mbox{ for some } \obfil \in \filtzar \}$$
and the associated separated presheaf is
$\Gamma(\openzar, \modul_1)
= \Gamma(\openzar,\modul)/\Gamma(\openzar, \modul_0)$
(see Section \ref{sheafifysubsection}).
We claim that there are natural homomorphisms
$$ \colim_{\obfil \in \filtzar} \Gamma(\obfil,\modul)
\lto \colim_{\obfil \in \filtzar} \Gamma(\obfil,\modul_1)
\lto \Gamma(\spax_\zftop, \modul_\zftop) \, .  $$
It is clear that
the first mapping exists and that it is a bijection.
For the existence of the second map, let
$\elem \in \Gamma(\obfil,\modul_1)$, $\obfil\in \filtzar$.
Since
$\obfil \times_\spax \obfil = \obfil \cap \obfil=\obfil$, the section
$\elem$ is compatible by Lemma \ref{antistrictlemma}, so we get a
mapping $\Gamma(\obfil, \modul_1) \to \Gamma(\spax_\zftop,
\modul_\zftop)$ and also from the colimit.

Every global section
$\elem \in \Gamma(\spax_\zftop, \modul_\zftop)$
is given by a compatible family $\elem_\indcov \in
\Gamma(\obcovindic, \modul)$,
$\indcovinset$,
in a covering $\obcovindic \to \openzar$,
so that the differences are
$0$ in
$\Gamma( \obcov_\indcov \cap \obcov_\indcovsec ,\modul_1)$.
We may assume that $\ideal$
is finite, since $\openzar$ is quasicompact, and suppose that
$\elem_\indcov - \elem_\indcovsec
\in \Gamma(\obcov_\indcov \cap \obcov_\indcovsec, \modul)$
restricts to $0$ in
$\obcov_\indcov \cap \obcov_\indcovsec \cap \obfil_{\indcov  \indcovsec}$,
$\obfil_{\indcov  \indcovsec} \in \filtzar$.
Then
$\spaw:=(\bigcup_\indcovinset \obcovindic)
\cap \bigcap_{(\indcov, \indcovsec) \in \setindcov \times \setindcov}
\obfil_{\indcov \indcovsec} \in \filtzar_\openzar$
and $\elem_\indcov$
is represented by an element $\elem \in \Gamma(\spaw, \modul)$.
Hence the second mapping is surjective. It is injective as it is the
sheafification homomorphism of a separated presheaf.

(ii). By (i) we have
$\Gamma(\openzar_\zftop, \modul_\zftop)
=\colim_\obfilinfil \Gamma(\openzar \cap \obfil, \modul) $.
For fixed $\obfil$, the assignment
$\openzar \mapsto \Gamma(\openzar \cap \obfil, \modul)$ is a sheaf,
as it is the push forward of $\modul|_\obfil$ under
$\obfil\subseteq \spax$.
Therefore also the colimit
$\openzar \mapsto \colim_\obfilinfil \Gamma(\openzar \cap \obfil,
\modul)$ is a sheaf \cite[Exercise II.1.11]{haralg}.
For quasicoherence we may assume hat $\spax$ is affine and that
$\fuf \in \Gamma(\spax, \O_\spax)$.
Then we have
$\colim_\obfilinfil \Gamma(D(\fuf) \cap \obfil, \modul)
= \colim_\obfilinfil \Gamma(\obfil,\modul)_\fuf
= (\colim _\obfilinfil \Gamma(\obfil, \modul))_\fuf $
(see \cite[Lemma II.5.14]{haralg} for the first equality).

(iii).
The sheafification is
$\shG \mapsto \colim_\obfilinfil (\incl_*(\shG|_\obfil)$, ($\incl:
\obfil \to \spax$)
so the left exactness of the sheafification follows
from the left exactness of push forward and from the exactness of
taking the colimit. The right exactness follows from Lemma
\ref{rightexactpullback}.

(iv).
The filter $\filtzar $ is a small cofiltered category
(Definition \ref{cofilteredcatdef}) with an arrow $\obfilsec \to \obfil$
corresponding to an inclusion $\obfilsec \subseteq \obfil$. The
inclusions
$\obfilsec \to\obfil$ and $\obfiltri \to \obfil$ come together in the intersection $\obfilsec \cap \obfiltri \to
\obfilsec,\obfiltri$. The identical mapping gives then a filter in the sense of
Definition \ref{filterdefinition}.
It is clearly an absolute filter
in the sense of Definition \ref{absolutefilterdef}
(but not irreducible, as the topology is not single-handed).
Then we have
$\colim_\obfilinfil \Gamma(\obfil_\zftop,\modul_\zftop)
= \colim_\obfilinfil(\colim_{\obfilw \in \filtzar}
\Gamma(\obfil \cap \obfilw ,\modul))
=\colim_\obfilinfil \Gamma(\obfil ,\modul)
=\Gamma(\spax_\zftop, \modul_\zftop)$.

(v) follows from Lemma \ref{explicitsite} or directly.
For the second part suppose first that
$\rest^\spax_\obfil(\elem) \in \Gamma(\obfil,\submod)$
for some $\obfil \in \filtzar$.
Then $\spax_0 \subseteq \obfil$ and so
$\elem_\point = \rest_\point(\rest^\spax_\obfil(\elem) )
\in\submod_\point$
for all $\point \in \spax_0$.
On the other hand, if
$\elem_\point \in \submod_\point$, then there exists an open
neighborhood of $\point \in \obfil_\point$ such that
$\rest^\spax_{\obfil_\point} (\elem) \in \Gamma(\obfil_\point, \submod)$.
Then
$\obfil= \bigcup_{\point \in \spax_0} \obfil_\point$ is an open neighborhood of
$\spax_0$ and $\rest^\spax_\obfil(\elem) \in \Gamma(\obfil, \submod)$, since
it is locally true.

(vi).
We may reduce to $\spax$ affine by (ii).
Let $\elem_1 \komdots \elem_\numgen$ be generators of $\modul$. Since $\modul_\zftop= 0$, there
exists by (v) for every $\indob$ an $\obfil_\indob \in \filtzar$
such that
$\rest^\spax_{\obfil_\indob}(\elem_\indob)=0$
in $\Gamma(\obfil_\indob, \modul)$, so all
generators become $0$ in $\Gamma(\obfil_1 \capdots \obfil_\numgen,
\modul)$.
On the other hand, if $\modul|_\obfil=0$, $\obfilinfil$, then the module is $0$ in a
covering,
hence it is $0$.

(vii). The right hand side has all the properties of a universal
$\delta$-functor for $\spax_\zftop$ \cite[Corollary III.1.4]{haralg}.
\end{proof}

\begin{remark}
\label{zariskiabelianremark}
The statements (i), (iii) and (vii) in
Proposition \ref{zariskifilterproperties} hold for arbitrary sheaves
of abelian groups on $\spax$ in the sense that
$\Gamma(\spax_\zftop,\sitemorinv (\shG)^\sheafify)
= \colim_\obfilinfil \Gamma(\obfil, \shG)$.
\end{remark}

\ren{\algfilt}{{T}}

\begin{example}
Let $\spaxeqspecring$ denote an affine scheme and let $\algfilt$
denote a multiplicatively closed subset in $\ring$ defining the
topological filter $\filtzar _\algfilt$ consisting of open subsets
$V \supseteq
D(\fuf)$ for some $\fuf \in \algfilt$. Then for an $\ring$-module
$\modul$ the module of global sections is
$$ \Gamma(\spax_\zftop, \modul_\zftop)
 = \colim_{V \in \filtzar_\algfilt} \Gamma(V,\modul)
 = \colim_{\fuf \in \algfilt} \modul_\fuf
 = \modul_\algfilt \, .$$
The closure of $\submod \subseteq \modul$ is
$\submod^\clop
=\{\elem \in\modul: \elem \in \submod_\algfilt \mbox{ in } \modul_\algfilt \}
= \{ \elem \in \modul:\, \exists \fuf \in \algfilt \mbox{ such that }
\fuf \elem \in \submod \}$.
In the case of such an algebraically
defined filter one can also define an affine Grothendieck topology
with the same closure operation by declaring $\ring \to \ring_\fuf$
to be a cover if $\fuf$ belongs to the saturated multiplicative
system given by $\algfilt$ (that is, $\fuf$ divides some element in
$\algfilt$).

The affine scheme $\Spec \ring_\algfilt$ gives an absolute cover for
this topology. If $\spax_0= \{\fop \in \Spec \ring_\algfilt\}=
\{\fop \in \specring,\, \fop \cap \algfilt = \emptyset\}$, then also
$\filtzar (\spax_0)= \filtzar_\algfilt$. On the other hand, if
$\spax_0$ is a finite collection of prime ideals $\fop_1 \komdots
\fop_\numn$, then the filter of open neighborhoods $\filtzar (\spax_0)$
is given algebraically by $\algfilt = \ring-
\bigcup_{\ind=1}^\numn
\fop_\ind$.
\end{example}

\begin{example}
\label{nzdtopology}
Let $\spaxeqspecring$ be affine and let $\algfilt \subseteq \ring$
be the multiplicative system of all non-zero divisors. The open
subsets $\{\obfil \subseteq \spax: \obfil \supseteq D(\fuf) \mbox{
for some }
\fuf \in \algfilt \}$ form a Zariski filter, and we denote the
corresponding Grothendieck topology by $\spax_\nzd$. The induced
closure operation on submodules is the (non-zero divisor-)
\emph{torsion},
$$\submod^{ \nzd - \clop} =\tors(\submod)
= \{ \elem \in \modul: \overline{\elem}
 = 0 \mbox{ in} (\modul/ \submod)_\algfilt  \}\, .$$
The quotient module $\modul/ \tors (0)$ is torsion-free.
The module of global sections is then
$\Gamma(\ring_{\nzd}, \modul_{\nzd})= \modul _\algfilt$,
and
$\modul/\tors (0) \subseteq \modul_\algfilt$
is a torsion-free submodule.
$\modul_\nzd$ is trivial if and only if $\modul$ is a torsion
module. If $\ring$ is noetherian, then $\algfilt =\ring -
\bigcup_{\ind=1}^\numn \fop_\ind$,
where $\fop_\ind$ runs through the associated primes of $\ring$
\cite[Theorem 3.1]{eisenbud}.
Therefore
$\submod^{\nzd-\clop} = \{\elem \in \modul: \elem \in \submod_\fop
\mbox{ for all associated primes } \fop \}$ (Proposition \ref{zariskifilterproperties}(v)).
Here one can also restrict to the maximal members among
$\Ass (\ring)$. For an ideal $\ideal \neq 0$ in a domain $\ring$ always
$\ideal^{\nzd -\clop}=\ring$.

For a scheme $\spax$ the same definitions apply. For an
$\O_\spax$-submodule $\shG \subseteq \sheaf $ the closure
$\shG^{\nzd -\clop}$ defined by the property that
$\sheaf / \shG^{\nzd -\clop}$
is torsion-free is sometimes called the saturation of $\shG$.
\end{example}

\begin{example}
\label{densetopology}
For a scheme $\spax$ let $\filtzar _{\dense}$
be the topological filter of all open dense subsets $\obfil$ in
$\spax$.
This means for $\spax$ irreducible that $\obfil$ is non-empty, and in
general that $\obfil$ meets every irreducible component. This is the
special case $\numk=0$ of Example \ref{zariskifiltercodimension},
and we call this topology the \emph{dense topology} and denote it by
$\spax_\dense$.

Let $\spaxeqspecring$ be affine, and let $\submodul$ be
$\ring$-modules. The closure operation induced by the dense topology
is the ``dense'' torsion of a submodule, in the sense that
$$\submod^{\dense -\clop} = \{ \elem \in \modul: \overline{\elem} = 0 \mbox{
in } (\modul/ \submod)_\fop \mbox{ for all minimal prime ideals }
\fop \}\, .$$ For $\ring$ noetherian, the module of global sections
is $\Gamma(\ring_{\dense}, \modul_{\dense}) =
\bigoplus_{\height{\fop}
= 0} \modul_{\fop}$.
\end{example}

\ren{\idealsup}{{J}}

\ren{\expoid}{{k}}
\ren{\closedzar}{{Z}}

\subsection{Submodules with support condition}
\label{submodulesupportsubsection}

\

\medskip
The easiest topological filters are the ones given by just one open
subset $\openzar \subseteq \spax$, so that $\filtzar (\openzar)=
\{\obfil:\,\obfil \supseteq \openzar\}$.
The closure operations coming from these filters are directly
related to support conditions and to local cohomology. Recall the
following definition.

\begin{definition}
Let $\ring$ be a ring, $\idealsup$ an ideal and  $\modul$ an
$\ring$-module. Then
$$\Gamma_\idealsup(\modul)= \{ \elem \in \modul: \idealsup^\expoid \elem = 0 \mbox{ for some } \expoid \in \NN\}\, $$
is called the \emph{module of sections in $\modul$ with support in}
$\idealsup$.

If $\spax$ is a scheme with an open subset $\openzar \subseteq \spax$ and
closed complement $\closedzar = \spax-\openzar$, and if $\sheaf $ is
a sheaf on
$\spax$,
then
$$\Gamma_\closedzar(\spax,\sheaf ) = \{\elem \in \Gamma(\spax,\sheaf ):\, \support(\elem) \subseteq \closedzar
\}$$ is called the \emph{module of sections with support in}
$\closedzar$
\cite[Expos\'{e} I]{SGA2}.
\end{definition}

\begin{lemma}
Let $\spax$ be a scheme and let $\openzar \subseteq \spax$ be an open
subset. Let $\sheaf $ be an $\O_\spax$-module. Then
$\Gamma_{\spax-\openzar}(\spax, \sheaf )$
is the closure of $0$ in
$\Gamma(\spax,\sheaf )$ in the Grothendieck topology induced by $\openzar$
{\rm(}that is, induced by the filter $\filt(\openzar)${\rm)}. In
particular, if $\ring$ is a noetherian ring with an ideal
$\idealsup$, then $\Gamma_\idealsup(\modul)$ is the closure of $0$
in the topology induced by $D(\idealsup)$.
\end{lemma}
\begin{proof}
In general,
$\support(\elem)
= \{\point \in \spax: \elem_\point \neq 0\} \subseteq \closedzar$
if and only if $\elem = 0$ in $\Gamma(\openzar, \sheaf )$. So the
statement follows from Proposition \ref{zariskifilterproperties}(v).
In the case of a noetherian ring, if
$\elem \in \modul$ restricts to
$0$ on
$\openzar = D(\idealsup)$, then, letting
$\idealsup =(\fug_1 \comdots \fug_\numm)$,
we also have
$\fug_\ind^{\expoid_\ind} \elem =0$ and so
$\idealsup^\expoid\elem =0$ for
$\expoid=\expoid_1 \plusdots \expoid_\numm$.
On the other hand, if
$\idealsup^\expoid \elem = 0$, then
$\elem |_\openzar =0$.
\end{proof}

Recall that $\Gamma_\idealsup(-)$ is left exact
(this follows also from Proposition \ref{zariskifilterproperties}(iii)
and Lemma \ref{leftexact})
and induces a right derived functor,
called local cohomology
(see \cite{SGA2}, \cite{brodmannsharp}).
This is true more generally for a Zariski filter
$\filtzar $ (but not for an arbitrary closure operation,
see Example \ref{radicalnotleftexact}).
In this case we also write
$$\Gamma_\filtzar(\modul)
= \{ \elem \in \Gamma(\spax,\modul):\, \elem|_\obfil= 0 \mbox{ for some} \obfil \in \filtzar\}
= \colim_{\obfil \in \filtzar} \Gamma_{\spax-\obfil}(\modul)
\, .$$
The long exact cohomology sequences \cite[Corollaire
I.2.9]{SGA2} for
$\obfil\in \filtzar$ induce a long exact cohomology sequence in the colimit, namely
$$0 \lto \Gamma_\filtzar(\modul)
\lto \Gamma(\spax,\modul) \lto \colim_{\obfil \in \filtzar} H^0 (\obfil, \modul)
= H^0(\spax_\zftop, \modul_\zftop) \lto $$
$$ H^1_\filtzar (\modul) \lto H^1(\spax, \modul )
\lto \colim_{\obfil \in \filtzar} H^1(\obfil ,\modul)
= H^1(\spax_\zftop, \modul_\zftop) \lto \,.$$

{
\ren{\numvar}{{d}}

\begin{example}
For a commutative ring $\ring$ with an ideal $\idealsup$ the closure of
an $\ring$-submodule $\submodul$ is given as
$\submod^\clop
= \{\elem \in \modul:\, \idealsup^\expoid \elem \subseteq \submod \mbox{ for some } \expoid \}$.
If $\ring = \field[\varx_1 \komdots \varx_\numvar]$
is a polynomial ring, and
$\idealsup= \fom =(\varx_1 \komdots \varx_\numvar)$,
then for an ideal
$\idealsubring$ this closure yields
$\ideal^\clop
 = \{ \fuf \in \ring:\, \fom^\expoid \fuf
\subseteq \ideal \mbox{ for some } \expon \}$, which
is also called the
\emph{saturation} of the ideal. A homogeneous ideal $\ideal$ and its
saturation $\ideal^\sat$ define the same projective variety
$\Proj \ring/\ideal= \Proj \ring/\ideal^\sat$.
\end{example}
}

\begin{example}
\label{radicalnotleftexact}
The closure of $0$ is not left exact in
general. Let $\field$ be a field and let
$\ring = \field[\epsilon]/(\epsilon^2)$,
and suppose we are in the surjective Grothendieck topology
(Part \ref{surjectivesection}).
Then the short exact sequence of $\ring$-modules
$0 \to  \field \cong (\epsilon) \to \ring \to \ring/ (\epsilon)= \field \to 0$
yields the non-exact sequence of closures of $0$ (radicals),
$0 \to 0 \to (\epsilon) \to 0 \to 0$.
So the left-exactness of $\Gamma_\idealsup( - )$ is a
particularity of the functor ``sections with support in''.
\end{example}

\subsection{Torsion theories, Zariski filters and Grothendieck topologies}

\label{torsiontheorysubsection}

\

\medskip
In this section we relate Grothendieck topologies induced by a
Zariski filter to hereditary torsion theories and show that they are
essentially the same. Torsion theories for abelian categories have
gained a lot of attention, in particular for the category of
left-modules for non-commutative rings. We recall briefly the basics
of torsion theory until we are able to translate them into our
setting. Introductions to torsion theories can be found in
\cite{blandtorsiontheory}, \cite{dicksontorsiontheory},
\cite{golantorsiontheory} or
\cite[Chapter 3]{vasthesis}.

A \emph{torsion theory} on an abelian category $\catc$
(think of the category of $\ring$-modules)
consists in two classes $\torclass$ and
$\freeclass$ of objects (the ``torsion class'' and the
``torsion-free class'') of $\catc$ such that
$\torclass \cap \freeclass = \{ 0 \}$, $\torclass$ is closed under quotients,
$\freeclass$ is closed under subobjects, and for every object
$\modul \in \catc$ there exists a short exact sequence
$0 \to \modtor \to \modul \to \modfree \to 0$ with
$\modtor \in \torclass $ and $\filtzar  \in
\freeclass$,
see \cite{dicksontorsiontheory}. There exists then also such a short
exact sequence where $\modtor$ is maximal
(\cite[Proposition 2.4]{dicksontorsiontheory} or \cite[Proposition
3.1]{vasthesis}), called the \emph{torsion submodule} of $\modul$.

We restrict to the case of $\ring$-modules and we denote the torsion
submodule of $\modul$ by $\modtor_\modul$. The class $\torclass$ is
characterized by the property
$\modtor_\modul = \modul$
\cite[Proposition 3.2]{vasthesis}. A torsion theory is immediately a
closure operation for the $0$-submodule and it becomes a closure
operation for submodules $\submodul$ by pulling back the torsion
submodule from
$\modul/\submod$ \cite[Chapter 3.4]{vasthesis}.
A torsion theory is called \emph{hereditary} if $\torclass$ is also
closed under subobjects.

\begin{theorem}
\label{grothendieckleftexacttorsiontheory}
Let $\ring$ be a noetherian commutative ring. Then the following
structures are equivalent.

\numiii
\begin{enumerate}

\item
A closure operation on $\ring$-Mod coming from a Grothendieck
topology on $\Spec \ring$
{\rm(}i.e. an admissible closure operation{\rm)}
such that the functor
$\modul \mapsto 0^\clop_\modul$ is left exact.

\item
A hereditary torsion theory on $\ring$-Mod.

\item
A Gabriel filter $\filtgab$ of ideals.
This is a set of ideals fulling {\rm(}a{\rm)} if
$\idealsup \in \filtgab$ and $\fug \in \ring$,
then
$(\idealsup:\fug) \in \filtgab$ and
{\rm(}b{\rm)} if $\idealsup \in \filtgab$ and if
$(\idealsupsec :\fug) \in \filtgab$ for all
$\fug \in \idealsup$, then also $\idealsupsec \in \filtgab$ \cite[Definition 3.4]{vasthesis}.

\item
A Zariski filter of open subsets in $\Spec \ring$.
\end{enumerate}
\end{theorem}
\begin{proof}
We describe the constructions leading from one structure to another.
From (i) to (ii).
One declares
$\torclass$ to be the class of $\ring$-modules $\modul$ with
$0^\clop = \modul$ and
$\freeclass$ to be the $\ring$-modules with $0^\clop =0$.
We go through the axioms of a torsion theory.
$\torclass \cap \freeclass = 0$ is clear.
Let $\homsurj: \modtor \to \modtorsec \to 0$ be a surjection and
suppose that $\modtor$ is $\catopen $-torsion. Then
$0^\clop_\modtorsec \supseteq \homsurj (0^\clop_\modtor)
= \homsurj (\modtor) = \modtorsec$.
If $0 \to \modfreesec \to \modfree$ is an injection and
$0^\clop_\modfree = 0$, then by the module persistence
(Proposition \ref{admissibleproperties}(i)) of the closure
operation also $0^\clop_\modfreesec = 0 $, so
$\modfreesec \in \freeclass$.

So far we have not used the left-exactness. For $\modul$ we consider
the short exact sequence
$0 \to 0^\clop \to \modul\stackrel{\homsurj}{\to} \modul/ 0^\clop = \modfree \to 0$.
By left exactness we have also the intersection property
$\modl \cap 0^\clop_\modm = 0^\clop_\modl$ for $\modl \subseteq \modm$. For if
$\elel \in \modl$ maps to $0^\clop_\modm$, then it maps to $0$ in
$0^\clop_{\modm/\modl}$ and by left exactness it must belong to
$0^\clop_\modl$.
From this intersection property it follows that
$0^\clop_{0^\clop_\modul}= 0^\clop_\modul$, so this module belongs
to $\torclass$. By the independence of presentation (Proposition
\ref{admissibleproperties}(ii)) we have
$\homsurj^{-1} (0^\clop_\modfree)
=(\homsurj^{-1}(0))^\clop = (0^\clop)^\clop
=0^\clop$,
hence $0^\clop_\modfree =0$ and
$\filtzar  \in \freeclass$.
By the intersection property it also follows for
$0 \to \modtorsec \to \modtor$ and
$\modtor \in \torclass $ that $\modtorsec \in \torclass$, so $(\torclass,\freeclass)$ is a
hereditary torsion theory.

The equivalence of (ii) and (iii) is a well known fact from torsion
theory (see \cite[Proposition 3.6]{vasthesis}). In a hereditary
torsion theory for
$\ring$-modules one has
$\modtor_\modul = \{\elem \in \modul: \ring/ \Ann (\elem) \in \torclass \}$,
hence a hereditary torsion theory is characterized by the set of
ideals $\idealsup \subseteq \ring$ such that
$\ring/\idealsup$ belongs to $\torclass$.
These ideals form the \emph{Gabriel filter}
$\filtgab =\filtgab(\torclass,\freeclass)$
of the torsion theory.
A Gabriel filter $\filtgab$ defines a hereditary torsion theory by
setting
$\modtor_\modul = \{\elem \in \modul : \Ann (\elem) \in \filtgab\}$.

From (iii) to (iv).
A Gabriel filter is closed under bigger ideals
and contains with $\idealsup$ and $\idealsupsec$ also $\idealsup
\idealsupsec$.
In particular, it contains arbitrary powers of $\ideal$.
Therefore, in the case of a commutative noetherian ring, a Gabriel
filter contains with an ideal also all ideals with the same radical.
Hence it corresponds to a Zariski filter.

From (iv) to (i) follows from Proposition
\ref{zariskifilterproperties}(iii) and Lemma \ref{leftexact}.

We have to show that the constructions are inverse to each other. A
left exact closure operation
$\clop$ yields the Zariski filter
$D(\idealsup)$, $\idealsup^\clop = \ring$ and the corresponding
Grothendieck topology, and we have to show that $\clop =
\zftop - \clop$. We can restrict to the closure of $0$. If $\elem \in
0^\clop$, then
$\elem \in \modtor_\modul$, so
$\Ann(\elem) \in\filtgab$ and so
$D(\Ann(\elem)) \in \filtzar$. Since $\elem|_{D(\Ann(\elem))}=0$,
we have $\elem \in 0^{\zftop - \clop}$. If
$\elem \in 0^{\zftop -\clop}$,
then $\elem|_{D(\idealsup)}$ for some
ideal $\idealsup$ with $\idealsup^\clop = \ring$. Then
$\support (\elem) \cap D(\idealsup)= \emptyset$
and so
$D(\idealsup) \subseteq D(\Ann(\elem))$,
hence $\Ann(\elem) \in \filtgab$. Therefore
$\Ann(\elem)^\clop =\ring$ and $\ring/\Ann(\elem)$ is a torsion module.
We have
$\ring/ \Ann(\elem) \to \modul$, $1\mapsto \elem$, and
$1 \in 0^\clop$ implies by persistence that $\elem \in 0^\clop$.
The other directions are similar.
\end{proof}

\begin{example}
We give an example of a non-hereditary torsion theory \cite[Example
3.1(8)]{vasthesis}. Consider the torsion theory for $\ring$ -modules
which is cogenerated by $\ring$. This means that an $\ring$-module
$\modul$ is torsion if and only if
$\Hom_\ring(\modul,\ring)=0$. If $\ring$
is a domain, but not a field, then
$\Hom_\ring (\quotfield(\ring),\ring)=0$,
so the quotient field is torsion, but the submodule $\ring \subset
\quotfield(\ring)$ is not-torsion, so the hereditary property does
not hold. This torsion theory does not fulfill the first finiteness
condition \ref{finitepropdef}, so it can not be induced by a
quasicompact Grothendieck topology by Corollary
\ref{quasicompactfiniteproperty}.
\end{example}

\subsection{The divisorial topology,  reflexive hull and symbolic powers}
\label{divisorialsubsection}

\

\medskip
In this section we deal with the following special case of Example
\ref{zariskifiltercodimension}.

\begin{definition}
Let $\spax$ be a scheme and let $\spax^{(1)}$ be its points of
height at most one. The Grothendieck topology on $\spax$, where the
open sets are the Zariski open sets and where
$\obcovindic \to \openzar$, $\indcovinset$,
is a covering if
$\openzar^{(1)} \subseteq \bigcup_\indcovinset \obcovindic$, is
called the \emph{divisorial topology}, denoted by $\spax_\divtop $.
\end{definition}

We will relate the divisorial topology with symbolic powers of
ideals and with the reflexive hull of a module. Recall that the dual
of an $\ring$-module is
$\modul^\dual=\Hom_\ring(\modul,\ring)$ and
that the
\emph{bidual} of a module $\modul$ is
$\modul^{\dual \dual}$.
We denote the natural homomorphism
$\modul \to \modul^{\dual \dual}$ given by
$\elem \mapsto (\varphi \mapsto \varphi(\elem))$ by $\shef$.
$\modul$ is called \emph{torsionless}
if $\shef$ is injective and \emph{reflexive} if $\shef$ is
bijective; see \cite[1.4]{brunsherzog} for these definitions. For a
domain $\ring$ a finite torsion-free
(with respect to non-zero divisors, see Example \ref{nzdtopology})
module $\modul$ is also torsionless.

\begin{lemma}
\label{divisorialglobalreflexive}
Let $\ring$ be a noetherian normal
domain and let $\modul$ be a finite torsion-free $\ring$-module.
Then
$\Gamma(\ring_\divtop, \modul_\divtop) = \modul^{\dual \dual}$,
the reflexive hull of $\modul$.
\end{lemma}
\begin{proof}
We have a natural $\ring$-module homomorphism
$$ \modul \lto \colim_{\spax^{(1)} \subseteq \obfil} \Gamma(\obfil, \modul)
=\Gamma(\ring_\divtop,\modul_\divtop) \, .$$
We first claim that the colimit on the right is reflexive,
using \cite[Proposition 1.4.1(b)]{brunsherzog}.
For prime ideals $\fop$ with depth $\! \ring_\fop \leq 1$ also
the height is $\leq 1$
(as $\ring$ is normal)
and so $\ring_\fop$ is a field or a discrete
valuation domain. Hence $\modul_\fop$ is free, since it is
torsion-free, and so it is reflexive.

Let now $\fop$ be a prime ideal with depth $\! \ring_\fop \geq 2$
and let $\elet_1,\elet_2 \in \ring$ be a regular sequence in
$\ring_\fop$.
We have to show that they are a regular sequence for
$\modul$ as well. Let $\elem_1,
\elem_2 $ be two sections of the colimit and suppose that $\elet_1
\elem_1 + \elet_2 \elem_2=0$.
We may assume
that $\elem_1, \elem_2 \in \Gamma(\obfil, \modul)$, that
$\obfil \subseteq D(\elet_1,\elet_2)$
and that the equation holds over $\obfil$.
Then $\frac{\elem_1}{\elet_2} = - \frac{\elem_2}{\elet_1} =: \elem$
is a global element of
$\Gamma(\obfil,\modul)$ and hence $\elem_1= \elet_2 \elem$.

The reflexivity implies that we get a homomorphism
$\modul^{\dual \dual} \to \Gamma(\ring_\divtop, \modul_\divtop)$. There
exists an open subset $\spaw \supseteq \spax^{(1)}$ where
$\modul|_\spaw$ is locally free, since $\modul_\fop$ is free for every
prime ideal of height one. Hence the colimit equals
$\Gamma(\spaw,\modul)$ and it is a finite $\ring$-module.
A homomorphism between two finite reflexive
$\ring$-modules is
by \cite[Lemma 24.2]{schejastorchverzweigung} a bijection if this is
true after localization at every prime ideal of height one.
For such a prime ideal we have
$(\colim_{\spax^{(1)}\subseteq \obfil}\Gamma(\obfil, \modul))_\fop
= \colim_{\fop \in \obfil} \Gamma(\obfil,\modul)
= \modul_\fop =(\modul^{\dual\dual})_\fop$.
\end{proof}

We come to the closure operation induced by the divisorial topology.

\begin{proposition}
\label{divisorialnormalclosure}
Let $\ring$ be a normal noetherian domain.
The closure operation for finite module $\submod \subseteq \modul$
induced by the divisorial topology is the {\rm(}\!`internal'{\rm)}
reflexive hull, namely
$$\submod^{\divtop -\clop}
= \{ \elem \in \modul : \shef(\elem) \in \submod^{\dual \dual} \,\,\,(\mbox{ inside } \modul^{\dual \dual})\}
\, .$$
For submodules of a finite reflexive module, the closure is the
reflexive hull $\submod^{\dual \dual}$.
For ideals $\idealsubring$ we have
\begin{eqnarray*}
\ideal^{\divtop - \clop}
&=& \ring \cap \bigcap_{\ideal \subseteq \fop \mbox{ prime of height one}} \ideal \ring_{\fop} \cr
&=& \{\fuf \in \ring: \nu_\fop( \fuf) \geq \min\{ \nu_\fop(\fug) : \fug \in \ideal\} \mbox{ for all }
\fop \mbox{ of height one} \} \, ,
\end{eqnarray*}
where $\nu_\fop$ denotes the corresponding discrete valuation. In
particular, for a prime ideal $\fop$ of height one we have
$(\fop^\numn)^{\divtop - \clop} = \fop^{(\numn)}$,
the $\numn$th symbolic power of $\fop$.
For ideals $\idealsubring$ of height $\geq 2$ we have
$\ideal^{\divtop -\clop}= \ring$.
\end{proposition}
\begin{proof}
The inclusion $\submod \subseteq \modul$ yields by left exactness
(Proposition \ref{zariskifilterproperties}(iii))
an inclusion
$\submod_\divtop \subseteq \modul_\divtop$ of sheaves in
the divisorial topology. The modules of global sections are by Lemma
\ref{divisorialglobalreflexive}
$\Gamma(\ring_\divtop,\submod_\divtop)
= \submod^{\dual \dual}
\subseteq \Gamma(\ring_\divtop, \modul_\divtop)
= \modul^{\dual \dual}$.
For $\modul$ reflexive we have the inclusions
$\submod \subseteq \submod^{\dual \dual }
\subseteq \modul^{\dual \dual} = \modul$.

By Proposition \ref{zariskifilterproperties}(v) we have
$\ideal^{\divtop-\clop}
= \ring \cap \bigcap_{\fop \mbox{ prime of height one}} \ideal \ring_{\fop}$,
and we can skip the prime ideals which do not contain $\ideal$.
The extended ideal $\ideal \ring_\fop$ for a prime ideal $\fop$ of height one equals
$(\fop \ring_\fop)^\nu$ for
$\nu = \min \{\nu_\fop( \fug): \fug \in \ideal\}$,
which gives the second equality.
Also, $\fop^\numn \subseteq \fop$ is the
only minimal prime ideal, and
$\fop^{(\numn)}= \ring \cap \fop^\numn \ring_\fop$
is by definition the $\numn$-th symbolic power of $\fop$
\cite[Section 3.9]{eisenbud}.
Finally, if
$\ideal$ has height $\geq 2$, meaning that every minimal prime ideal
of $\ideal$ has height $\geq 2$, then the index set is empty.
\end{proof}

\ifthenelse{\boolean{book}}{

\begin{remark}
For a finite reflexive $\ring$-module $\modul$ and a submodule
$\submodul$ one can also characterize the (internal) reflexive hull
as
$$\{\elem \in \modul: \mbox{ for every } \ring-\mbox{linear }
\varphi: \submod \to \ring \mbox{ there exists a unique extension
to} (\submod, \elem) \} .$$
If
$\elem \in \submod^{\dual\dual}$,
then
$\submod \subseteq (\submod,\elem) \subseteq
\submod^{\dual\dual} \subseteq \modul$
and since
$\submod^\dual \to \submod^{\dual\dual\dual}$
is a bijection, we have a unique extension to
$\submod^{\dual \dual}$.
This gives an extension to
$(\submod, \elem)$, which is also unique (replacing $\submod$ by
$(\submod, \elem)$). On the other hand, suppose that every
$\varphi \in \submod^\dual $ is uniquely extendible to $(\submod,\elem)$.
Then defining
$\elem(\varphi):= \tilde{\varphi}(\elem)$, where
$\tilde{\varphi}$ is the extension, shows that
$\elem \in \submod^{\dual \dual}$.
\end{remark}
}{}

\begin{proposition}
\label{divisorialdivisor}
Let $\spax$ be a normal noetherian
integral scheme endowed with the divisorial topology, let
$\quotfield = \quotfield (\spax)$ be its function field. Then the
following holds.

\numiii
\begin{enumerate}

\item
The coherent submodules of $\quotfield$ in the divisorial topology
are all invertible and are in one to one correspondence with Weil
divisors. The coherent ideal sheaves correspond to effective Weil
divisors.

\item
We have $H^1(\spax_\divtop, \O^\times_\divtop) = \Cl \spax$, the
divisor class group of $\spax$
{\rm(}compare \cite[Proposition XI.3.7.1]{SGA2}{\rm)}.
\end{enumerate}
\end{proposition}
\begin{proof}
(i). By coherent we mean here just coherent as a module in the
Zariski topology. Let $\shinv \subseteq \quotfield$ be a coherent
submodule. For $\fop \in \spax$ of height one the localization
$\shinv_\fop \subseteq \quotfield$ is free of rank one, and by
coherence there exists an open affine neighborhood
$\fop \in \openzar$ such that $\shinv|_\openzar$ is free
(in the Zariski topology and in the divisorial topology).
These neighborhoods do cover
$\spax$ in the divisorial topology (not necessarily in the Zariski
topology), hence
$\shinv$ is invertible. In general, a Weil divisor corresponds to an invertible submodule of
$\quotfield$ on some open subset $\obfil \supseteq \spax^{(1)}$. So the
result follows.

(ii).
The exact sequence
$0 \to \O^\times \to \quotfield^\times \to \quotfield^\times/\O^\times \to 0 $
of sheaves of abelian groups
yields the corresponding exact sequences on $\spax_\divtop$
(Proposition \ref{zariskifilterproperties}(iii) and Remark
\ref{zariskiabelianremark})
and the global sequence
$$0\! \to\! \Gamma(\spax, \O^\times)
\!\to \!\quotfield^\times \!\to\!\! \!\colim_{\spax^{(1)} \subseteq \obfil}\!\!
\Gamma(\obfil,\quotfield^\times/ \O^\times)\! \to\!\! \!
\colim_{\spax^{(1)}\subseteq \obfil} \!\!H^1(\obfil, \O^\times)
\! =\! H^1(\spax_\divtop,
\O_\divtop^\times)\! \to\! 0 \, .$$
A section
$\elem \in \Gamma(\obfil, \quotfield^\times/\O^\times)$
is given by a covering
$\obfil = \bigcup _\indinset \obfil_\ind$ and
$\fuq_\ind \in \quotfield^\times$ with
$\fuq_\ind \fuq_\indsec^{-1}
\in \Gamma(\obfil_\ind \cap\obfil_\indsec, \O^\times)$.
It defines clearly a submodule which is invertible
in codimension one. On the other hand, a Weil divisor corresponds to
such a submodule
$\shinv \subseteq \quotfield$, and $\shinv$ is for
every $\point \in \spax^{(1)}$ generated by
$\fuq_\point \in \quotfield$, and this is true in a neighborhood of $\point$.
The union of these neighborhoods together is then an open set
$\obfil \supseteq \spax^{(1)}$.
This gives an identification
$Div \spax =
\colim_{\spax^{(1)}
\subseteq \obfil} \Gamma(\obfil,\quotfield^\times/\O^\times)$.
So the written down cohomology sequences is just the
defining sequence of the divisor class group
\cite[Section II.6]{haralg}.
\end{proof}

\section{Further examples}
{
\ren{\point}{{P}}
\ren{\homring}{{\psi}}

\subsection{Completion}
\label{completionsubsection}

\

\medskip
Also the completion of a ring with respect to an ideal may be
considered as given a Grothendieck topology, but quite in a
different flavor. In particular the coverings there are not
coverings in the set-theoretical sense, and it is not a
single-handed Grothendieck topology nor a quasicompact one.

Let $\ring$ denote a commutative ring and $\foa \subseteq \ring$ an
ideal. We consider the category over
$\spaxeqspecring$ consisting of
$\spax$ and of all
$\spax_\numk := \Spec \ring/ \foa^\numk$, $\numk \geq 0$, and
all $\spax$-morphisms. The coverings in this topology are given by
families which contain the identity
($\spax_\numk \to \spax_\numk$ or $\spax \to \spax$)
and by families which contain
$ \spax_\numk \to \spax$, $\numk \in \setind \subseteq \NN$,
where $\setind$ is infinite.
We denote this Grothendieck topology by
$\spax_{\foa-\topo}$
(similar definitions are possible for ideal families
$\foa_\numk$, $\numk \in \NN$, with $\foa_0=\ring$ and
$\foa_\numk \cdot \foa_\numsec \subseteq \foa_{\numk+\numsec}$;
one can also formulate scheme versions).
An $\ring$-module $\modul$
defines the presheaf $\modul^\pre$ by
$$ \Gamma(\spax_\numk, \modul^\pre)
:= \modul \tensorr  \ring/\foa^\numk
= \modul/ \foa^\numk \modul \, .$$

\begin{proposition}
\label{completeproperties}
Let $\spaxeqspecring$ be an affine scheme
endowed with the topology
$\spax_{\foa-\topo}$ for an ideal $\foa
\subseteq \ring$.
Let $\modul$ denote an $\ring$-module.
Then the following holds.

\numiii
\begin{enumerate}

\item
The presheaf $\modul^\pre$ is a sheaf with module of global sections
$$ \Gamma(\spax_{\foa-\topo}, \modul^\pre)
= \hat{\modul}
= \liminv_{\numk \in \NN}  \modul/ \foa^\numk \modul \, ,$$
the $\foa$-adic completion of $\modul$. In particular, the ring of
global sections in this topology is just the completion
$\hat{\ring}$ of
$\ring$ with respect to the ideal $\foa$.

\item
The global sheafification homomorphism
$\shef: \ring \to \Gamma(\ring_{\foa-\topo}, \O)
= \hat{\ring}$
has
$\fob = \bigcap_{\numk \in \NN}\foa^\numk$
as its kernel, which is zero for $\foa \neq \ring$ in a
noetherian ring.

\item
The mapping
$\overline{\NN} = {\NN} \cup \{\infty \} \to \catopenspax$,
$\numk \mapsto \spax_\numk$ {\rm(}$\spax_\infty = \spax${\rm)},
is a non-absolute irreducible filter.
The stalk of a module is $0$.

\item
The mapping
$\{\point\} \to \catopenspax$,
$ \point \mapsto \spax$,
is an absolute filter which is not irreducible.
Its stalk is the same as the module of global sections.

\item
The closure operation induced by this topology is just
$\ideal \to \ideal^\tocl = \ideal + \fob$,
so it is the identical operation in the
noetherian case.

\item
A ring homomorphism
$\homring : \ring \to \ringsec$ where
$\homring (\foa) \subseteq \foa'$ induces a site morphism
$\spax'_{\foa'-\topo} \to \spax_{\foa-\topo}$
and the corresponding homomorphism of rings of global sections is
just
$\hat{\ring} \to \hat{\ringsec }$.
\end{enumerate}
\end{proposition}
\begin{proof}
(i)
An element in $\hat{\modul}$ is given by a sequence
$\elem_\numk \in \modul/\foa^\numk \modul$ with the compatibility condition
$\varphi_{\numk\numsec}(\elem_\numk)= \elem_\numsec$
for all $\numk \geq \numsec$.
This is exactly the compatibility
condition in the topology.
(ii) is clear.

(iii)
For $\numk \leq \numsec$ we have
$\foa^\numk \supseteq \foa^\numsec$,
hence
$\ring/ \foa^\numsec \to \ring/ \foa^\numk$ and so
$\spax_\numk \to \spax_\numsec$.
As $\spax_\numk \to \spax$
are not covers, this filter is not absolute.
Since
$\emptyset =\Spec \ring/ \foa^0$ is the initial
object of the filter, all stalks are $0$.

(iv)
The only single cover of $\spax$ is the identity, so the filter
is clearly absolute.
The covering
$\spax_\numk \to \spax$, $\numk\in\NN$,
is such that the identity $\spax \to \spax$ factors through no
$\spax_\numk$,
hence this filter is not irreducible
(Definition \ref{irreduciblefilterdef}).
(v) and (vi) are clear.
\end{proof}
}

\subsection{Proper topology}
\label{propersubsection}

\

\medskip
Proper topologies seem to have first been considered by S.-I. Kimura
\cite[Definition 3.1]{kimuraalexander} in order to give a
cohomological characterization of Alexander schemes. Here we treat
the proper topology in its relation to the integral closure of
ideals.

\begin{definition}
Let $\spax$ be a scheme. Let $\catopenspax$ denote the category
of all schemes over $\spax$ and define a covering to be a proper
surjective morphism $\obv \to \obu$ (we do not insist that the
structure morphisms are proper).
We denote $\spax$ endowed with this
Grothendieck topology by $\spax_\prop$.
\end{definition}

\begin{remark}
Recall that a proper morphism is separated, of finite type and
universally closed \cite[Section II.4]{haralg}. Proper morphisms are
stable under base change and under composition, hence we get indeed
a Grothendieck topology. Note that the proper topology is not an
affine topology and that it is not a refinement of the
(single-handed) Zariski topology. Every scheme morphism induces a
site morphism in the proper topology. The proper topology is a
refinement of the finite topology. A proper and surjective
homomorphism is universally submersive.
\end{remark}

Important examples of proper morphisms are given by blow-ups of
ideals, see \cite[D\'{e}fintion 8.1.3]{EGAII} or
\cite[Propositions II.7.10, II.7.13]{haralg}.
Recall that the blow-up of an ideal
$\idbla \subseteq \ring$ is given by
$\spax_\idbla= \Proj \algb_\idbla$,
$\algb_\idbla= \oplus_{\expok \in \NN} \idbla^\expok$.
We characterize first which blow-ups yield surjective morphisms.

\ren{\indcompsec}{{j}}
\ren{\indr}{{r}}

\begin{lemma}
\label{blowupsurjective}
A blow-up $\spax_\idbla \to \spaxeqspecring$ of
an ideal $\idbla \subseteq \ring$ in a noetherian ring is surjective
if and only if $\idbla$ is not contained in any minimal prime ideal of
$\ring$.
\end{lemma}
\begin{proof}
Let
$\compindic = \Spec \ring_\indcomp$ denote the integral components of
$\specring$. If $\idbla$ is not contained in the minimal prime $\foq$,
then the extended ideal in $\ring/\foq$ is not the $0$-ideal.
For a domain the blow-up of a non-zero ideal defines a surjective
morphism, since the image is closed and the morphism is an
isomorphism on a non-empty open set
\cite[Proposition II.7.13]{haralg}.
Since the blow-up of the component factors through
$\spax_\idbla \times_\spax\obpartindic$ \cite[Corollary II.7.15]{haralg},
it follows that the morphism is surjective on every component, hence
surjective.

On the other hand, if
$\idbla\ring_\indcomp=0$, then over the open set
$\obu_\indcomp= \spax_\indcomp- \bigcup_{\indcompsec \neq \indcomp}\spax_\indcompsec$
we have
$\spax_\idbla \times _\spax \obu_\indcomp = (\obu_\indcomp)_{\idbla|\obu_\indcomp} = \emptyset$.
\end{proof}

The \emph{extended ideal} of $\ideal \subseteq \ring$ under a scheme
morphism $\spay \to \Spec \ring$ is on every affine subset $\obv
=\Spec \rings \subseteq \spay$ given by $\ideal \rings$. This gives
an ideal sheaf on $\spay$. The following lemma is motivated by
Ratliff's $\catideal$-closure, see Section \ref{ratliffsubsection}.

\begin{lemma}
\label{ratliffextended} Let $\ring$ denote a noetherian ring and
suppose that $\fuf \idbla \subseteq \ideal \idbla$, where $\fuf \in
\ring$ and
where $\ideal$ and $\idbla$ are ideals.
Then $\fuf$ belongs to the extended
ideal of
$\ideal$ under the blow-up
$\spax_\idbla \to \spax$.
\end{lemma}
\begin{proof}
Set $\ideal = (\runfuf)$ and $\idbla=(\fuh_1 \komdots \fuh_\num)$. The condition
$\fuf\idbla \subseteq \ideal \idbla$ means explicitly that for every $\fuh_\indi$ we may write
$\fuf \fuh_\indi
= \sum_{ 1 \leq \indr \leq \numgen, 1 \leq \inds \leq \num} \fua_{\indr\inds \indi}\fuf_\indr\fuh_\inds$,
$\fua_{\indr\inds\indi} \in \ring$.
This translates to
$$\fuf = \sum_{\indr=1}^\numgen
(\sum_{\inds=1}^\num \frac{\fua_{\indr\inds\indi}\fuh_\inds}{\fuh_\indi})
\fuf_\indr \,\,\,\,\, \mbox { in } \,\,\,\,\, (\algb_\idbla)_{\fuh_\indi} \, ,$$
so $\fuf \in \ideal ((\algb_\idbla)_{\fuh_\indi})_0$
for $\indi =1 \komdots \num$. This means that $\fuf$ belongs to the
extended ideal of $\ideal$ on
$D_+(\fuh_i) \subseteq \spax_\idbla$. These open
sets cover $\spax_\idbla$,
therefore $\fuf$ belongs also globally to the extended ideal on $\spax_\idbla$.
\end{proof}

\begin{proposition}
\label{properintegral}
Let $\spaxeqspecring$ denote an affine scheme endowed with the proper topology.
Then the closure  induced by this topology is the integral closure.
\end{proposition}
\begin{proof}
Assume first that $\fuf$ belongs to the extended ideal $\ideal
\O_\spay $ for a proper surjective morphism $\spay \to \spax$. Such
a morphism is universally
submersive \cite{blicklebrennersubmersion}
. Under the base change $\patha :\Spec \dvd \to \spax$, where
$\dvd$ is a discrete valuation domain, the image $\patha^*(\fuf)$ is
also in the extended ideal of $ \ideal \O_\dvd$ in $\spay_\dvd$. Thus we may
assume that $\ring =\dvd$ is a discrete valuation domain with local
parameter $\pi$, that $\fuf=\pi^\expon$ and $\ideal=(\pi^\expom)$.
Then by
\cite{blicklebrennersubmersion} 
it follows that there exists a discrete valuation domain $\dwd$ such
that $\Spec \dwd \to \spay \to\dvd$ is surjective, which shows that
$\fuf \in \intclo{\ideal}=\ideal$.

\ren{\indcomp}{{j}}
\ren{\setindcomp}{{J}}
\ren{\idealsec}{{L}}

Suppose now that $\fuf \in \intclo{\ideal}$.
Let $\spax_\indcomp=\Spec \ring_\indcomp$,
$\indcompinset$, denote the integral components of $\spax$. The morphism
$\biguplus \spax_\indcomp \to \spax$ is proper and surjective. We will
construct $\spay \to \spax$ proper and surjective on each component
$\spax_\indcomp$ separately according to whether the extended ideal $\ideal
\O_{\spax_\indcomp}$ is $0$ or not. In the first case we also have $\fuf_\indcomp=0$,
and we do not need to modify $\spax_\indcomp$.

So suppose that $\ideal$ extends to a non-zero ideal in $\ring_\indcomp$, or
that $\ring$ is a domain and $\ideal \neq 0$. We use the theory of
reductions, see \cite[Chapter 5]{hunekeapplication}.
If $\fuf \in
\intclo{\ideal}$, then $\ideal$ is a reduction of $\idealsec=(\ideal,\fuf)$, that is,
$\idealsec^{\expor+1}= \ideal \idealsec^\expor$ for some $\expor$.
So in particular
$\fuf \idealsec^\expor \subseteq \ideal \idealsec^\expor$.
Hence $\fuf$ belongs to the extended ideal of the blow-up
along the ideal $\idealsec^\expor$ by Lemma \ref{ratliffextended}.
Since
$\idealsec \neq 0$, this blow-up is a surjection.
\end{proof}

\begin{remark}
This result is related to the known fact that $\fuf$ belongs to the
integral closure of an ideal $\ideal$ in a domain if and only if it
belongs to the extended ideal in the normalization of the blow-up of
$\ideal$, see \cite[Proposition 1.4]{faridiblowup}.
\end{remark}

\ifthenelse{\boolean{book}}{

\begin{bookremark}
In terms of forcing algebras, the containment $\fuf \in \ideal^\prop$,
$\ideal = (\runfuf)$, holds true if and only if $\Spec \algforc \to
\specring = \spax$, $\algforc = \ring[\vart_1 \komdots
\vart_\numgen]/(\fuf_1\vart_1 \plusdots \fuf_\numgen\vart_\numgen+
\fuf)$, admits locally (that is, for a covering of affine Zariski
open sets) a section for some proper and surjective morphism $\spay
\to \spax$, see \ref{forcingglobalremark}.
\end{bookremark}}{}

\subsection{The $\catideal$-closure of Ratliff}
\label{ratliffsubsection}

\

\medskip
In \cite{ratliffdeltaclosure}, Ratliff introduced the concept of the
$\catideal$-\emph{closure} of an ideal
(see also \cite{ratliffrushdeltareduction} and \cite{ratliffrushasymptotic}).
Here $\catideal $ is a fixed multiplicatively closed set of ideals
in a commutative ring $\ring$ and the closure operation for an ideal
$\idealsubring$ is given by
$$\ideal_\catideal
= \bigcup_{\idbla \in \catideal}  (\ideal \idbla:\idbla) \, .$$
Thus an element $\fuf \in \ring $ belongs to $\ideal_\catideal$ if
and only if there exists an ideal $ \idbla \in \catideal$ such that
$\fuf \idbla \subseteq \ideal \idbla$.
Under mild conditions the
$\catideal$-closure lies inside the integral closure (see
\cite[Theorem 3.2]{ratliffdeltaclosure} or Proposition
\ref{ratliffintegral} below).

We will translate this closure operation into our language and we
will show how some of Ratliff's basic results follow immediately.
It will turn out that the most natural Grothendieck topology to induce
the $\catideal$-closure is given by certain non-affine proper
mappings.
The basic relationship is given by the following
proposition
(compare with \cite[Fact 2.1]{heinzerlantzshahratliffrush} for the case of the Ratliff-Rush
closure).
\ren{\expot}{{t}}

\begin{proposition}
\label{ratliffblowup}
Let $\ring$ denote a commutative noetherian ring,
let $\catideal$ be a multiplicatively closed collection of ideals in
$\ring$. Let
$\spax_\idbla = \Proj \algb_\idbla$,$\algb_\idbla
=\bigoplus_{\expok \in \NN} \idbla^\expok$,
$\idbla \in \catideal$,
denote the corresponding blow-ups.
Then an element $\fuf \in \ring$ belongs to Ratliff's $\catideal$-closure of an ideal
$\idealsubring$,
i.e. $\fuf\in \ideal_\catideal$,
if and only if $\fuf$ belongs for some
$\spax_\idbla \to \spaxeqspecring$, $\idbla\in \catideal$,
to the extended ideal of $\ideal$.
\end{proposition}
\begin{proof}
Set $\ideal = (\runfuf)$.
Suppose first that $\fuf \in \ideal_\catideal$,
so that there exists an ideal
$\idbla \in \catideal$
with $\fuf \idbla \subseteq \ideal \idbla$.
Then $\fuf $ belongs to the extended ideal under the blow-up $\spax_\idbla \to \spax$ by
Lemma \ref{ratliffextended}.

For the other direction suppose that $\fuf$ belongs to the extended
ideal of $\ideal$ under the blow-up $\spax_\idbla \to \spax$,
$\idbla =(\fuh_1\komdots \fuh_\numm)$.
This means that we have locally $\fuf \in \ideal $
on the affine covering sets $D_+(\fuh_\ind)$, that is,
$\fuf \in\ideal((\algb_\idbla)_{\fuh_\ind})_0$.
This means explicitly that
$\fuf = \sum_{\indsec=1}^\numgen
\frac{\eleb_{\ind \indsec}}{\fuh_\ind^{\expot_{\ind \indsec}}} \fuf_\indsec$,
where
$\eleb_{\ind \indsec} \in \idbla^{\expot_{\ind \indsec}}$.
Hence there exists also a $\expot$ such that
$\fuf \fuh_\ind^\expot \in \setind \idbla^\expot$ for all
$\ind$. It follows that
$\fuf \idbla^{\numm \expot} \subseteq \setind \idbla^{\numm \expot}$,
since for a monomial
$\fuh_1^{\nu_1} \cdots \fuh_\numm^{\nu_\numm}
\in \idbla^{\numm \expot}$ we have
$\nu_\ind \geq \expot$ for some $\ind$ and hence
$$\fuf \fuh_1^{\nu_1} \cdots \fuh_\numm^{\nu_\numm}
= \fuf \fuh_\ind^\expot \fuh_1^{\nu_1} \cdots \fuh_\ind^{\nu_\ind- \expot} \cdots \fuh_\expom^{\nu_\expom}
\in \setind \idbla^\expot \cdot \idbla^{\numm \expot - \expot} \subseteq \setind \idbla^{\numm \expot}\, .$$
Since $\catideal$ is multiplicatively closed we have
$\idbla^{\numm \expot} \in \catideal$ and we get $\fuf \in \ideal_\catideal$.
\end{proof}

\begin{definition}
\label{deltatopologydef}
Let $\spax$ denote a scheme and let
$\catideal$ denote a multiplicatively closed set of ideal sheaves.
For $\moryx$ denote the set of extended ideal
sheaves $\idbla \O_\spay $, $\idbla \in \catideal$, by
$\catideal_\spay$. Then the $\catideal$-\emph{topology} on $\spax$
is given by the full category of schemes of finite type over
$\spax$, and we declare $\spaz\to\spay$ to be a cover if there
exists a factorization $\spay_\idbla \to \spaz\to\spay$, where
$\idbla \in \catideal_\spay $ and where $\spay_\idbla \to \spay$ is
the blow-up of $\idbla \O_\spay $.
\end{definition}

\begin{remark}
If $\spay_\idbla \to \spaz \to \spay$ is a covering and $\spau \to \spay$ an
$\spax$-morphism, then we get a factorization
$\spau_\idbla \to \spay_\idbla
\times_\spay \spau \to \spaz \times_\spay  \spau \to \spau$ by the universal
property of blow-up \cite[Proposition 7.14]{haralg}. Hence the
pull-back of a covering is a covering.
If $\spaz \to \spay$ and
$\spaw  \to\spaz$ are coverings with factorizations
$\spay_\idbla \to \spaz \to \spay$ and
$\spaz_\idblasec \to \spaw \to \spaz$, then we get
$\spay_{\idbla \idblasec} \cong (\spay_\idbla)_\idblasec \to
\spaz_\idblasec$, which shows that the composition of coverings is a covering.
So we have indeed a Grothendieck topology.
\end{remark}

\begin{corollary}
Let $\ring$ denote a noetherian commutative ring and let $\catideal$
denote a multiplicatively closed set of ideals. Then the
$\catideal$-closure equals the closure induced by the
$\catideal$-topology on $\spaxeqspecring$.
\end{corollary}
\begin{proof}
This follows immediately from Lemma \ref{explicit}, Proposition
\ref{ratliffblowup} and Definition \ref{deltatopologydef}.
\end{proof}

\begin{proposition}
Let $\varphi: \spay_{\catideal'} \to \spax_\catideal$ be a scheme
morphism such that $\idbla \O_\spay  \in \catideal'$ for all ideal
sheaves $\idbla \in \catideal$. Then $\varphi$ induces a site
morphism in the $\catideal$-topology. In particular, the persistence
property holds for the $\catideal$-closure for these morphisms
{\rm(}as stated in \cite[Theorem 5.1]{ratliffdeltaclosure}{\rm)}.
\end{proposition}
\begin{proof}
If $\spav \to \spax$ is a covering with factorization
$\spax_\idbla \to \spav \to \spax$, then the factorization
$\spay_\idbla \to \spax_\idbla \timesx \spay
\to \spau \timesx \spay \to \spay$
shows that the pull-back is also a covering. The second statement
follows from Theorem
\ref{topclosure}.
\end{proof}

\begin{proposition}
\label{ratliffintegral}
Let $\ring$ denote a noetherian ring and let
$\catideal$ denote the set of all ideals which are not contained in
any minimal prime ideal. Then for every ideal
$\ideal \in \catideal$ we have
$\ideal_\catideal = \intclo{\ideal}$,
the integral closure of $\ideal$.
\end{proposition}
\begin{proof}
Let $\fuf \in \ideal_\catideal$, and suppose that $\fuf$ belongs to the
extended ideal of $\ideal$ under
$\spax_\idbla \to \spax$, $\idbla \in \catideal$. This blow-up is surjective by Lemma
\ref{blowupsurjective}, hence $\fuf$ belongs to the proper closure
of $\ideal$ and hence to the integral closure by Proposition
\ref{properintegral}.

On the other hand, if $\fuf \in \intclo{\ideal}$, then $\fuf$ belongs to the
extended ideal of $\ideal$ under the blow-up of a power of
$(\ideal,\fuf)$ as in the proof of Proposition \ref{properintegral}.
\end{proof}

\ren{\funum}{{c}}
\ren{\fudeno}{{b}}

Ratliff also associates to a commutative ring endowed with
$\catideal$ another ring $\ring^\catideal$, the $\catideal$-closure
of the ring \cite[Definition 6.2]{ratliffdeltaclosure}. In the case
where $\catideal$ does not contain the zero-ideal and $\ring$ is
noetherian this is (see \cite[Theorem 6.2]{ratliffdeltaclosure})
$$\ring^\catideal
= \{\funum/\fudeno: \fudeno \mbox{ non-zero divisor and }
\funum \in (\fudeno)_\catideal \}
\subseteq \quotfield(\ring) \, ,$$
where
$\quotfield(\ring)$ denotes the total quotient ring. We relate this
construction to absolute stalks and rings of global sections in the
$\catideal$-topology.

We endow a given multiplicatively closed set $\catideal$ of ideals
with the structure of a cofiltered set by declaring that an arrow
$\idbla \to \idblasec$ exists if $\idblasec$
extends to an invertible ideal in $\spax_\idbla$. In this case there
is a unique morphism
$\spax_\idbla \to \spax_\idblasec$
by
\cite[Proposition II.7.14]{haralg}.

\begin{lemma}
\label{ratlifffilterstalk}
Let $\ring$ denote a commutative ring,
and let $\catideal$ denote a multiplicatively closed set of ideals.
Then $\catideal$ is a cofiltered set and
$\catideal \to\catopen_\spax$, $\idbla \mapsto \spax_\idbla$, is an absolute
filter. If
$\ring \to \rings$
is a ring homomorphism such that for every $\idbla \in
\catideal $ the extended ideal $\idbla \rings$ is invertible, then
there exists a ring homomorphism
$\colim_{\idbla \in \catideal}\Gamma(\spax_\idbla, \O) \to \rings$.
\end{lemma}
\begin{proof}
Two arrows
$\idbla \to \idblatri$, $\idblasec \to \idblatri$ have the common refinement $\idbla \idblasec$,
because in
$\spax_{\idbla \idblasec}$ both $\idbla$ and $\idblasec$ become invertible.
Since there is at most one arrow, $\catideal$ is
cofiltered in the sense of Definition \ref{cofilteredcatdef}.

By the definition of a cover in the $\catideal$-topology every
morphism $\spax_\idbla \to \spax$ indexed by this filter is a cover.
If
$\spay \to \spax_\idbla$ is a cover, then there exists an ideal
$\idblasec \in \catideal$ and a factorization
$(\spax_\idbla)_\idblasec \to \spay \to \spax_\idbla$,
and so this filter is absolute
(Definition \ref{absolutefilterdef}).
The corresponding absolute stalk of the
structure sheaf is hence
$\colim_{\idbla \in \catideal} \Gamma(\spax_\idbla,\O)$.

If $\ring \to \rings$ is given and the extended ideal $\idbla
\rings$ is invertible, then there exists a unique morphism
$\Spec \rings \to \spax_\idbla$ and hence
$\Gamma (\spax_\idbla ,\O)\to \rings$,
which extends to a  ring homomorphism
$\colim_{\idbla\in \catideal} \Gamma(\spax_\idbla, \O) \to \rings$.
\end{proof}

\ren{\expot}{{t}}

\begin{proposition}
\label{ratliffglobalstalk}
Let $\spaxeqspecring$ denote an affine scheme and let $\catideal$
denote a multiplicatively closed set of ideals in $\ring$ such that
every
$\idbla \in \catideal$ contains a non-zero divisor.
Then we have identities
$\Gamma(\spax_\catideal, \O_\catideal)
= \colim_{\idbla\in \catideal} \Gamma(\spax_\idbla, \O)
= \ring^\catideal$.
\end{proposition}
\begin{proof}
We have natural homomorphisms
$ \Gamma(\spax_\catideal, \O_\catideal)
\to \colim_{\idbla\in \catideal} \Gamma(\spax_\idbla, \O)
\to \quotfield (\ring)$, where the second homomorphism exists due to
Lemma \ref{ratlifffilterstalk} applied to
$\rings= \quotfield(\ring)$.
The first homomorphism exists by Lemma \ref{stalkpresheafsheaf}
and is injective by Lemma \ref{globalstalkinjective}.
So let
$\fuglob \in \Gamma(\spax_\idbla, \O)$, $\idbla \in
\catideal$, be given. The diagonal
$\spax_\idbla \to \spax_\idbla \times_\spax \spax_\idbla$
is a cover in the $\catideal$-topology,
since we have a factorization
$(\spax_\idbla \timesx \spax_\idbla)_\idbla
\to \spax_\idbla \to  \spax_\idbla \timesx \spax_\idbla$.
Hence $\fuglob$ is compatible (Lemma \ref{antistrictlemma}) and
$\fuglob \in \Gamma(\spax_\catideal,\O_\catideal)$.

Let now
$\funum/\fudeno \in \quotfield(\ring)$
be such that
$\funum \in (\fudeno)_\catideal$, say
$\funum\idbla \subseteq \fudeno \idbla$, $\idbla=(\fuh_1 \komdots \fuh_\numm)$.
This means that for every $\fuh_\ind$ there exists $\fug_\ind
\in \idbla$ such that
$\funum \fuh_\ind =\fudeno \fug_\ind$.
For $\ind,\indsec$ it follows that
$\fudeno \fug_\ind\fuh_\indsec
=\funum \fuh_\ind \fuh_\indsec
=\fudeno \fug_\indsec\fuh_\ind$
and so
$\fug_\ind\fuh_\indsec = \fug_\indsec\fuh_\ind$,
since $\fudeno$ is a non zero divisor.
Therefore
$\fug_\ind/\fuh_\ind= \fug_\indsec/\fuh_\indsec$
in
$\Gamma(D_+(\fuh_\ind \fuh_\indsec), \O)$
and so we get a global element in $\Gamma(\spax_\idbla,\O)$.
This element maps to
$\funum/\fudeno$ under the natural homomorphism.

On the other hand, suppose that
$\funum/\fudeno \in \quotfield(\ring)$ is the image of
$\fuq \in \Gamma(\spax_\idbla,\O)$, $\idbla \in \catideal$,
$\idbla =(\fuh_1\komdots \fuh_\numm)$.
Let $\fuh \in \idbla$ be a non zero
divisor, and let
$\fuq =\fug/\fuh^\expor \in \Gamma(D_+(\fuh), \O)$.
This maps to
$\fug/\fuh^\expor \in \quotfield(\ring)$
and so
$ \funum/\fudeno = \fug/\fuh^\expor$.
So we may assume that
$\fudeno\in \idbla$. We have
$\funum/\fudeno = \fug_\ind/\fuh_\ind$ in $\Gamma(D_+(\fudeno \fuh_\ind), \O)$,
and so
$\funum \fuh_\ind^\expot \in \fudeno\idbla^\expot$
for some $\expot$ and all $\ind$.
Then the ideal
$\idbla^{\expot \numm}$
shows that
$\funum \in (\fudeno)_\catideal$
(as in the proof of Proposition \ref{ratliffblowup}).
\end{proof}

\bibliographystyle{plain}

\end{document}